\documentclass[12pt,a4paper]{amsart}

\usepackage{amssymb,amsmath,amsfonts,amsthm,amscd}
\usepackage{leftidx}
\usepackage{tikz}

\usetikzlibrary{matrix,positioning,decorations.pathreplacing}

\usepackage[utf8]{inputenc}

\usepackage{mathrsfs}

\usepackage{hyperref}
\hypersetup{
    colorlinks=true,       % false: boxed links; true: colored links
    linkcolor=blue,          % color of internal links
    citecolor=blue,        % color of links to bibliography
    filecolor=magenta,      % color of file links
    urlcolor=cyan           % color of external links
}
\usepackage[all]{hypcap}

\usepackage{subfigure,graphicx}
\usepackage{xcolor}
\usepackage{xspace,accents}

\makeatletter
\renewcommand\normalsize{%
    \@setfontsize\normalsize{11.7}{14pt plus .3pt minus .3pt}%
    \abovedisplayskip 10\p@ \@plus4\p@ \@minus4\p@
    \abovedisplayshortskip 6\p@ \@plus2\p@
    \belowdisplayshortskip 6\p@ \@plus2\p@
    \belowdisplayskip \abovedisplayskip}
\renewcommand\small{%
    \@setfontsize\small{9.5}{12\p@ plus .2\p@ minus .2\p@}%
    \abovedisplayskip 8.5\p@ \@plus4\p@ \@minus1\p@
    \belowdisplayskip \abovedisplayskip
    \abovedisplayshortskip \abovedisplayskip
    \belowdisplayshortskip \abovedisplayskip}
\renewcommand\footnotesize{%
    \@setfontsize\footnotesize{8.5}{9.25\p@ plus .1pt minus .1pt}%%
    \abovedisplayskip 6\p@ \@plus4\p@ \@minus1\p@
    \belowdisplayskip \abovedisplayskip
    \abovedisplayshortskip \abovedisplayskip
    \belowdisplayshortskip \abovedisplayskip}
\ifdefined\pdfpagewidth
\else
\fi
\calclayout
\makeatother

\let \Re \relax
\DeclareMathOperator{\Re}{Re}
\let \Im \relax
\DeclareMathOperator{\Im}{Im}
\DeclareMathOperator*{\esssup}{ess\,sup}
\newcommand{\mb}[1]{\ensuremath{\mathbb{#1}}}
\newcommand{\N}{{\mb{N}}}
\newcommand{\ZZ}{{\mb{Z}}}
\newcommand{\R}{{\mb{R}}}
\newcommand{\C}{{\mb{C}}}

\newcommand{\MM}{{\mb{M}}}

\DeclareMathOperator{\trace}{tr}
\newcommand{\eps}{\varepsilon}

\newcommand{\unitfunction}[1]{\bld{1}_{#1}}

\newcommand{\D}{\ensuremath{\mathscr D}}
\newcommand{\E}{\ensuremath{\mathcal E}}
\newcommand{\G}{\ensuremath{\mathcal G}}
\renewcommand{\H}{\ensuremath{\mathcal H}}
\newcommand{\I}{\ensuremath{\mathcal I}}
\renewcommand{\L}{\ensuremath{\mathcal L}}
\newcommand{\M}{\ensuremath{\mathcal M}}
\renewcommand{\O}{\ensuremath{\mathcal O}}

\newcommand{\A}{\ensuremath{\mathsf A}}

\newcommand{\ESC}{\mathcal B_{esc}}

\renewcommand{\k}{\ensuremath{\kappa}}
\newcommand{\y}{\ensuremath{\varrho}}

\newcommand{\para}[1]{\leftidx{^\parallel}{#1}{}}

\newcommand{\ppi}{\pi_{\parallel}}

\newcommand{\py}{\para{\y}}

\newcommand{\symb}{\Sigma}
\newcommand{\symbo}{\symb_0}
\newcommand{\symbc}{\symb_c}
\newcommand{\symbH}{\symb^{\mathcal H}_0}
\newcommand{\Symbol}[3]{\symb^{#1}(\jp{\xi}^{#2};\R^{#3})}
\newcommand{\Symbolo}[3]{\symb_0^{#1}(\jp{\xi}^{#2};\R^{#3})}

\newcommand{\T}{\mbox{\tiny ${\mathsf T}$}}
\newcommand{\symbt}{\symb_{\T}}
\newcommand{\symbto}{\symb_{\T,0}}
\newcommand{\Symbolt}[4]{\symbt^{#1}(\jp{\eta'}^{#2};\R^{#3}\times
  \R^{#4})}
\newcommand{\Symbolto}[4]{\symbto^{#1}(\jp{\eta'}^{#2};\R^{#3}\times
  \R^{#4})}

\DeclareMathOperator{\Span}{span}

\newcommand{\sd}{{\mathrm{d}}}
\newcommand{\gl}{{\mathrm{3}}}
\newcommand{\sg}{{\mathrm{g}}}

\newcommand{\sdGb}{\G^{\sd}_\d}
\newcommand{\glGb}{\G^{\gl}_\d}
\newcommand{\sgGb}{\G^{\sg}_\d}

\newcommand{\Hb}{\H_\d}
\newcommand{\Gb}{\G_\d}

\newcommand{\Equiv}{\Leftrightarrow}
\newcommand{\imp}{\Rightarrow}
\DeclareMathOperator{\Hamiltonian}{H}
\newcommand{\Hp}{\Hamiltonian_p}
\newcommand{\HpG}{\Hamiltonian_p^{\G}}
\newcommand{\Hpz}{\Hp z}
\newcommand{\Hppz}{\Hp^2 z}
\newcommand{\Hz}{\Hamiltonian_z}

\newcommand{\Hzzp}{\Hz^2 p}

\newcommand{\n}{\mathsf n}
\newcommand{\nx}{\n_x}

\newcommand{\Con}{\ensuremath{\mathscr C}}
\newcommand{\Conc}{\ensuremath{\mathscr C}_{c}}
\newcommand{\Cinf}{\ensuremath{\mathscr C^\infty}}
\newcommand{\Cinfc}{\ensuremath{\mathscr C^\infty_{c}}}

\newcommand{\X}{\mathcal X} 
\newcommand{\Y}{\mathcal Y} 

\renewcommand{\d}{\ensuremath{\partial}}

\newcommand{\transp}[1]{\ensuremath{\leftidx{^t}{\!#1}{}}}

\newcommand{\wrt}{w.r.t.\@\xspace}
\newcommand{\rhs}{r.h.s.\@\xspace}
\newcommand{\lhs}{l.h.s.\@\xspace}

\newcommand{\resp}{resp.\@\xspace}

\newcommand{\pp}{a.e.\@\xspace}
\newcommand{\scm}{s.c.m.\@\xspace}
\newcommand{\nhd}{neighborhood\xspace}

\newcommand{\suff}{sufficiently\xspace}

\newcommand{\cst}{{\mathrm{Cst}}}

\DeclareMathOperator{\sgn}{sgn}

\DeclareMathOperator{\Char}{Char}

\DeclareMathOperator{\obs}{obs}
\newcommand{\Cobs}{C_{\obs}}

\DeclareMathOperator{\Obs}{\mathsf L}

\newcommand{\bichar}{bicharacteristic\xspace}
\newcommand{\bichars}{bicharacteristics\xspace}
\newcommand{\gbichar}{generalized bicharacteristic\xspace}
\newcommand{\gbichars}{generalized bicharacteristics\xspace}

\newcommand{\gammaG}{\ensuremath{\leftidx{^{\mathsf G}}{\gamma}{}}}
\newcommand{\GammaG}{\ensuremath{\leftidx{^{\mathsf G}}{\bar{\gamma}}{}}}

\newcommand{\XG}{\ensuremath{\leftidx{^{\mathsf G}}{X}{}}}

\newcommand{\bld}[1]{\mbox{\boldmath $#1$}}

\newcommand{\tb}{\tilde{b}}

\newcommand{\tchi}{\tilde{\chi}}

\newcommand{\hw}{\hat{w}}

\newcommand{\ty}{\tilde{y}}
\newcommand{\tz}{\tilde{z}}
\newcommand{\hy}{\hat{y}}

\newcommand{\tk}{\tilde{\kappa}}
\newcommand{\tzeta}{\tilde{\zeta}}
\newcommand{\tpsi}{\tilde{\psi}}
\newcommand{\tvarphi}{\tilde{\varphi}}

\newcommand{\hpsi}{\hat{\psi}}
\newcommand{\hvarphi}{\hat{\varphi}}

\newcommand{\ttheta}{\tilde{\theta}}
\newcommand{\htheta}{\hat{\theta}}

\newcommand{\bb}{\bar{b}}
\newcommand{\bQ}{\bar{Q}}
\newcommand{\bq}{\bar{q}}
\newcommand{\bG}{\bar{G}}

\let \div \relax
\DeclareMathOperator{\div}{div}
\newcommand{\nablag}{\nabla_{\!\! g}}
\newcommand{\divg}{\div_{\! g}}
\newcommand{\mug}{\mu_g}
\newcommand{\mugb}{\mu_{g_\d}}

\newcommand{\kk}{\k_k}
\newcommand{\tkk}{\tk_k}
\newcommand{\hk}{h_k}
\newcommand{\gk}{g_k}
\newcommand{\mugk}{\mu_{\gk}}
\newcommand{\mugbk}{\mu_{{\gk}_{\hspace*{-0.07em}\mbox{\tiny$\d$}}}}

\newcommand{\kn}{\k_n}

\newcommand{\hn}{h_n}
\newcommand{\gn}{g_n}
\newcommand{\mugn}{\mu_{\gn}}
\newcommand{\mugbn}{\mu_{{\gn}_{\hspace*{-0.07em}\mbox{\tiny$\d$}}}}

\newcommand{\chart}{\mathcal C}
\newcommand{\hchart}{\mathcal C}
\newcommand{\hO}{O}
\newcommand{\chdiff}{\phi}
\newcommand{\cdiff}{\phi}
\newcommand{\cdiffL}{\phi_\L}

\newcommand{\tM}{\tilde{\M}}
\newcommand{\hM}{\hat{\M}}
\newcommand{\ThM}{T^* \hM}

\newcommand{\hL}{\hat{\L}}
\newcommand{\tL}{\tilde{\L}}
\newcommand{\ThL}{T^* \hL}

\newcommand{\TM}{T^* \M}
\newcommand{\TL}{T^* \L}
\newcommand{\dTM}{\d(\TM)}
\newcommand{\dTL}{\d(\TL)}
\newcommand{\pdTL}{\para{\dTL}}

\newcommand{\pTL}{\para{\TL}}

\newcommand{\pEb}{\para{\mathcal E}_\d}
\newcommand{\pGb}{\para{\mathcal G}_\d}
\newcommand{\pHb}{\para{\mathcal H}_\d}

\newcommand{\interval}{J}
\DeclareMathOperator{\dist}{dist}

\DeclareMathOperator{\Op}{Op}
\newcommand{\OpH}{\Op^h}

\DeclareMathOperator{\supp}{supp}

\newcommand{\et}{\ensuremath{\text{and}}}
\newcommand{\ou}{\ensuremath{\text{or}}}
\newcommand{\avec}{\ensuremath{\text{with}}}
\newcommand{\si}{\ensuremath{\text{if}}}
\newcommand{\pour}{\ensuremath{\text{for}}}
\newcommand{\dans}{\ensuremath{\text{in}}}

\newcommand{\where}{\ensuremath{\text{where}}}

\newcommand{\inp}[2]{(#1, #2)} 
\newcommand{\biginp}[2]{\big(#1, #2 \big)}

\newcommand{\scp}[2]{#1 \cdot  #2}

\newcommand{\bigdup}[2]{\big\langle #1, #2 \big\rangle}
\newcommand{\dup}[2]{\langle #1, #2 \rangle}

\newcommand{\ovl}[1]{\overline{#1}}
\newcommand{\udl}[1]{\underline{#1}}

\newcommand{\Rdp}{\R^d_+}

\newtheoremstyle{note}{} {}{\itshape}{-6pt}{\bf}{. --}{ }{}

\newtheorem{theorem}{Theorem}[section]
\newtheorem{proposition}[theorem]{Proposition}
\newtheorem{lemma}[theorem]{Lemma}
\newtheorem{corollary}[theorem]{Corollary}
\newtheorem{assumption}[theorem]{Assumption}
\newtheorem*{theo*}{Theorem}

\newtheorem{lemmabis}{Lemma}
\newcounter{theorembiss}

\newcounter{lemmabiss}

\newtheorem{definition}[theorem]{Definition}
\newtheorem{remark}[theorem]{Remark}

\newcommand{\jp}[1]{\langle #1 \rangle}

\DeclareMathOperator{\loc}{loc}

\renewcommand{\S}{\ensuremath{\mathscr S}}
\newcommand{\K}{\ensuremath{\mathscr K}}

\newcommand{\Ldct}{dominated-convergence theorem\xspace}

\newcommand{\Norm}[2]{{\| #1 \|}_{#2}}
\newcommand{\bigNorm}[2]{{\big\| #1 \big\|}_{#2}}
\newcommand{\BigNorm}[2]{{\Big\| #1 \Big\|}_{#2}}

\newcommand{\norm}[2]{{|#1|}_{#2}}
\newcommand{\bignorm}[2]{{\big|#1\big|}_{#2}}

\newcommand{\bt}{_{|t=0}}

\newcommand{\ubt}{\underaccent{\bar}{t}}
\numberwithin{equation}{section}


\title[Measure and continuous vector fields] 
{Measure and continuous vector field at a boundary~I:
  propagation equations and wave observability} 

\author{Nicolas Burq}
\address{Nicolas Burq. Laboratoire de Math\'ematiques d'Orsay, Universit\'e
  Paris-Sud, Universit\'e Paris-Saclay, B\^atiment~307, 91405
  Orsay Cedex \& CNRS UMR 8628 \& Institut Universitaire de France}
 \email{nicolas.burq@math.u-psud.fr}
\author{Belhassen Dehman}
  \address{Belhassen Dehman. Universit\'e de Tunis El Manar, Facult\'e
  des Sciences de Tunis, 2092 El  Manar \& Ecole Nationale d'Ing\'enieurs de Tunis, ENIT-LAMSIN, B.P. 37, 1002 Tunis, Tunisia. }
\email{Belhassen.Dehman@fst.utm.tn}
 \author{J\'er\^ome Le Rousseau}
 \address{J\'er\^ome Le Rousseau.
Universit\'e Sorbonne Paris Nord, Laboratoire Analyse,
G\'eom\'etrie et Applications, LAGA, CNRS, UMR 7539, F-93430,
Villetaneuse, France.}  \email{jlr@math.univ-paris13.fr}
\date{\today}

\begin{document}
\begin{abstract}
  The celebrated geometric control condition of Bardos, Lebeau, and
  Rauch is necessary and sufficient for wave
  observability~\cite{BLR:92,BG:1997} and exact controllability. It
  requires that any point in phase-space be transported by the
  generalized geodesic flow to the region of observation in some
  finite time. The initial smoothness ($\Cinf$) required on the
  coefficients of the metric to prove that exact control and geometric
  control are essentially equivalent was subsequently relaxed to
  $\Con^2$-metrics/coefficients and $\Con^3$-domains~\cite{Burq:1997},
  which is close to the optimal smoothness required to preserve a
  generalized geodesic flow. In this article, we investigate a
  natural generalization of the geometric control condition that
  makes sense for $\Con^1$-metrics and we prove that wave
  observability holds under this condition.  Moreover, we establish that
  the observability property is stable under rougher (Lipschitz)
  perturbation of the metric. We also provide a geometric necessary
  condition for wave observability to hold. Transport equations that describe the
  propagation of semi-classical measures are at the heart of the
  arguments. They are natural extensions to geometries with boundaries
  of usual transport equations. This article is mainly dedicated to
  the proof of such propagation equations in this very rough context.
\end{abstract}

\maketitle
\tableofcontents
%%%%%%%%%%%%%%%%%
% section
%%%%%%%%%%%%%%%%%
\section{Introduction}
An observability inequality for the wave equation is  an estimate of
the energy of a solution  by a ``recording'' of this
solution in a restricted domain for some time $T>0$. This restricted domain can
be in the interior of the region $\Omega$ where waves propagate or in
a part of its boundary.
One interest in such an inequality lies in its consequences in
terms of exact controllability and stabilization.
The observability property has been intensively
studied during the last decades. Until the end of the 80's, most of
the results were proven under
a (global) geometrical assumption called
$\Gamma$-condition and introduced by J.-L.~Lions~\cite{Lions:88},
essentially based on a multiplier method. Later, following Rauch and
Taylor~\cite{RT:74}, Bardos, Lebeau, and Rauch proved 
observability inequalities from part of the boundary in their seminal article
\cite{BLR:92}, and as a
consequence,  boundary stabilization,  under a microlocal condition,
that is, a property in the cotangent bundle
$T^{\ast}(\R \times \Omega)$, the so-called geometric control
condition (GCC in short), exhibiting a connection between the set on
which observation is performed
and the generalized geodesics of the wave operator. In
addition, taking into account the work of \cite{BG:1997}, it is now
classical that observability (with stability with respect to
the observation set) is equivalent to the GCC. In terms of geodesics,
the GCC reads as follows:
\begin{center}
  {\em for any point $x$ and any tangent vector $v$,
  the generalized geodesic initiated at $(x,v)$
  enters the observation region in some time $T>0$.}
\end{center}
Generalized geodesics follow the laws of geometrical optics at
boundary points: reflection if the boundary is hit transversally and
possible gliding if hit tangentially.

The proofs of the results in \cite{BLR:92,BG:1997} are based on microlocal tools, namely,
the propagation of wavefront sets or that of microlocal defect measures. Let us notice here that despite their high efficiency and
robustness, these methods present the great disadvantage of requiring
a lot of regularity for the domain and for the metric/coefficients. Starting from the
original result developed in the framework of the Melrose-Sj\"ostrand
$\Cinf$-singularity propagation results, thus requiring
$\Cinf$-smoothness, the theory has been subsequently developed in
the framework of microlocal defect measures allowing one to relax the
assumptions down to a $\Con^2$-metric \cite{Burq:1997}, which barely misses
the natural smoothness ($W^{2, \infty}$) required to define a geodesic
flow (away from any boundary). An important remark is that below this
smoothness threshold, for instance for $\Con^1$-metrics that we will consider, generalized
geodesics may still exist as integral curves of the
$\Con^0$-Hamiltonian vector field
in the interior of the domain but {\em uniqueness is lost} in
general. A natural question lies in the understanding of the
relationship between those nonunique integral curves and the
observability property for such a rough metric. 

Many attempts were made in the last years (see, for instance, the works
\cite{Fan-Zua} in dimension 1, and
\cite{Dehman-Ervedoza}). In the present article, we reach the lowest
possible regularity level for the mere existence of geodesics (with a
gain of a full derivative with respect to all previous geometric
results) and lowest possible regularity level for observability to
hold as exhibited by the counter-example in \cite{CastroZuazua02, CZ07}.

\medskip
We prove the following result: for a $\Con^1$-metric, observability holds, and
    consequently exact controllability, if a generalized geometric control
    condition is satisfied, that is,
    \begin{center}
     {\em  all generalized geodesics enter the
      observation region in some time $T>0$.}
  \end{center}
  In particular, if considering a point $x$ and some tangent vector $v$,
  {\em all} generalized geodesics initiated at this point with $v$ as initial direction
  fulfill this property.  When
   uniqueness of generalized geodesics holds
   the above  condition coincides with the usual GCC.  We
   thus keep the ``GCC'' denomination. 

   \medskip
   Moreover, our proof allows us to go beyond the $\Con^1$-threshold
   in a perturbative regime
   and consider cases where the notion of geodesics is lost.
  %  We prove
%    that under GCC for a reference $\Con^1$-metric,  the
% observability property for the wave equation remains valid up to small
% Lipschitz ($W^{1, \infty}$)  perturbations. More
% precisely, we show
   We prove 
that if a reference $\Con^1$-metric $g$
satisfies the GCC for some time $T>0$, for  any other metric $\tilde{g}$ close enough to
$g$ in the Lipschitz topology the observability property holds also,
moreover in the same time $T$.  We insist once more on the fact that  Lipschitz metrics
are too rough to even define geodesics since the associated Hamiltonian vector field
is only $L^\infty$ and hence integral
curves do not make sense in general.

The GCC stated above stands as a sufficient condition for
observability to hold. We also provide a  necessary condition that is
natural in the sense that it coincides with the usual necessary
condition if uniquess of
the flow holds \cite{BG:1997}.

\medskip
Our proof of the observability property relies on two key results on
semi-classical measures that appear naturally if 
the energy of waves concentrates asymptotically:
\begin{itemize}
\item The derivation of a  transport equation that describes the
  propagation of such semi-classical measures; this is
 one of the main contribution of the present article; see
  Theorem~\ref{thm: equationmesure}. Its proof relies on commutator analysis and finely
  tuned properties of semi-classical opertors to address the low
  regularity level of metric/coefficients and the boundary/manifold.
\item The description  of the support of measures
  that are solutions to the above propagation equation in terms of
  \gbichars (whose projection on the base manifold  are
  geodesics) even if uniqueness of such curves fails to hold; this is the main
  contribution of the companion article~\cite{BDLR2} stated 
  in Theorem~\ref{theo-propagation} here.
\end{itemize}
In section~\ref{sec: Method to prove observability and outline} below, we describe how these two results are used in
the structure of the proof of observability inequalities.

\medskip
An important difficulty in the present article is the presence of a
boundary. In the case of a compact manifold without boundary, we
refer to our much less technically demanding article~\cite{BDLR0}
where both parts of this program are achieved in that simpler
setting.

Our proof of a necessary condition for wave observability to hold  is quite similar and based on:
\begin{itemize}
\item the derivation of a  transport equation for semi-classical
  measures across an isochrone $\{ t=\cst\}$; see Theorem~\ref{theorem: measure equation at t=0};
\item the use of this measure equation to ensure that a maximal 
  \gbichar lies entirely in the support of the measure.
\end{itemize}

%%%%%%%%%%%%%%%%%
% subsection
\subsection{Metrics, elliptic operators and wave equations}
\label{sec: Wave equation, regularity levels}

Consider a compact connected Riemannian manifold $\M$ of dimension $d$
with boundary, endowed with a metric $g = (g_{ij})$. At first $\M$ and its
boundary are assumed $W^{2,\infty}$ and the metric is assumed
Lipschitz. Denote by $\mug$ the canonical positive Riemannian density on
$\M$, that is, the density measure associated with the density
function $(\det g)^{1/2}$. We also consider a positive Lipschitz function $\k$
and we define the density $ \k \mug$. 

The $L^2$-inner product and norm are considered with respect to this
density $\k \mug$, that is, 
\begin{align}
  \label{eq: L2-inner product and norm}
  \inp{u}{v}_{L^2(\M)} = \int_\M u \bar{v} \, \k \mug, 
  \qquad 
  \Norm{u}{L^2(\M)}^2 = \int_\M |u|^2\, \k \mug.
\end{align}
We denote by $L^2 V(\M)$ the space of $L^2$-vector fields on $\M$,
equipped with the norm
\begin{align*}
  \Norm{v}{L^2 V(\M)}^2
  = \int_\M  g(v, \bar{v}) \, \k \mug, \qquad v \in L^2 V(\M).
\end{align*}
We write
\begin{align*}
  \inp{.}{.}_{L^2(\M,\k \mug)},\ \ 
  \Norm{.}{L^2(\M,\k \mug)}, \ \ 
  \et  \ \
  \Norm{.}{L^2 V(\M, \k \mug)}^2,
\end{align*}
if needed for clarity in particular if different metrics and functions
$\k$ are considered simultaneously.

Recall that the Riemannian gradient and divergence are given by 
\begin{align*}
  g (\nablag f, v) = v(f)
  \ \ \et \ \ 
  \int_\M f \divg v \mug = - \int_\M v(f) \, \mug,
\end{align*}
for $f$ a  function and $v$ a vector field with supports away from the
boundary, yielding in local
coordinates
\begin{align*}
  %\label{eq: riemann gradient-divergence}
  (\nablag f)^i= \sum_{1\leq j \leq d} g^{i j}\d_{x_j} f, 
  \qquad 
  \divg v = (\det g)^{-1/2} 
  \sum_{1\leq i\leq d} \d_{x_i}\big( (\det g)^{1/2} v^i \big),
\end{align*}
with $(g_x^{ij}) = (g_{x,i j})^{-1}$.
With the Poincar\'e inequality a norm on $H^1_0(\M)$ is 
\begin{align*}
  %\label{eq: norm H1-0}
  \Norm{u}{H^1_0(\M)} = \Norm{\nablag u}{L^2 V(\M)}.
\end{align*}

We introduce the elliptic operator $A = A_{\k,g}= \k^{-1} \divg (\k
\nablag)$, that is, in local coordinates
\begin{align}
  \label{eq: elliptic operator}
  A f = \k^{-1} (\det g)^{-1/2} 
  \sum_{1\leq i, j \leq d}  \d_{x_i}\big( 
  \k  (\det g)^{1/2} g^{i j}(x)\d_{x_j} f
  \big).
\end{align}
The operator $A$  is unbounded on $L^2(\M)$. With the domain $D(A) =H^2(\M) \cap H_0^1(\M)$ one
finds that $A$ is selfadjoint, with respect to the $L^2$-inner product
given in \eqref{eq: L2-inner product and norm}, and
negative.

With the elliptic operator $A = A_{\k,g}$ one also defines
 the wave operator 
 \begin{equation*}
   %\label{wave}
   P = P_{\k,g}=\d_t^2 - A_{\k,g}.
\end{equation*}

Consider the wave equation
\begin{align}
  \label{eq: wave equation-intro}
  \begin{cases}
     P_{\k, g}  \, u = f
     & \dans\ \R\times\M,\\
     u=0 &  \dans\ \R \times\d\M,\\
      u_{|t=0}  = \udl{u}^0, \ \d_t u_{|t=0}  = \udl{u}^1 
      & \dans\ \M.
  \end{cases}
\end{align}
Solutions are given by the following result. 
%%%%%%%%%%%%%%%%%%%%%%%%
% proposition          %
%%%%%%%%%%%%%%%%%%%%%%%
\begin{proposition}
  \label{prop: well-posedness nonhomogeneous wave equation}
  Consider $\k$ and $g$ both Lipschitz. 
  Let $(\udl{u}^0, \udl{u}^1) \in H_0^1(\M) \times L^2(\M)$ and $f \in L^1_{\loc}(\R, L^2(\M))$.  There exists a unique  
  \begin{align*}
    &u \in \Con^0 \big(\R; H_0^1(\M)  \big)
    \cap \Con^1 \big(\R; L^2(\M)\big),
  \end{align*}
  that is, a weak solution of \eqref{eq: wave
    equation-intro}, meaning 
  $ u\bt = \udl{u}^0$ and $\d_t u\bt  = \udl{u}^1$ and
  $P_{\k, g}  u =f$ holds in $\D' \big(\R \times \M\big)$.
  The map
  \begin{align}
  \label{eq: continuity wave initial conditions}
  H^1_0(\M) \oplus L^2(\M) \oplus L^1_{\loc}(\R, L^2(\M))
  &\to \Con^0\big(\R; H_0^1(\M)\big)  \cap
                             \Con^1\big(\R; L^2(\M)\big)\\
  (\udl{u}^0, \udl{u}^1,f) & \mapsto u \notag,
  \end{align}
  is continuous.
\end{proposition}
One denotes by 
\begin{align*}
  \E_{\k, g} (u) (t)
  &= \frac12 \big(  \Norm{\nablag u(t)}{L^2 V (\M)}^2
    + \Norm{\d_tu(t)}{L^2(\M)}^2 \big)
\end{align*}
the energy of $u$ at time $t$. For any $T >0$ there exists $C_T>0$ such that
\begin{align*}
  %\label{eq: evol energy}
  \sup_{|t|\leq T} \E_{\k, g} (u) (t)^{1/2}
  \leq C_T\big( \E_{\k, g} (u) (0)^{1/2}
  + \Norm{f}{L^1(-T,T; L^2(\M))}\big).
\end{align*}
If $f=0$, then equation \eqref{eq:
  wave equation-intro} is homogeneous.  This is the case we will most
often consider for the issue of observability. Then, for the weak
solution $u$, the energy is independent of time $t$ that is,
\begin{align*}
  \E_{\k,
  g} (u) (t) = \E_{\k, g} (u) (0) = \frac12 \big( \Norm{\nablag
    \udl{u}^0}{L^2 V(\M)}^2 + \Norm{\udl{u}^1}{L^2(\M)}^2 \big).
\end{align*}
In such case, we simply write $\E_{\k, g} (u)$.

%%%%%%%%%%%%%%%%%
% subsection
\subsection{Regularity levels for manifolds, metrics and coefficients}

Two classes of regularity levels will be of importance in what
 follows. A first one for  which microlocal methods apply (the spaces
 $\X^1$ and $\X^2$  below) and a second one for which basic results
 (uniqueness, traces, etc)  remain true (the space $\Y^1$ below). More
 precisely denote by $X^k$ (\resp $Y^k$) the sets of manifolds $\M$ as
 above of class  $\Con^k$ (resp. $W^{k, \infty}$). This regularity of
 the manifold includes that of its boundary. Set
 \begin{align*}
  &\X^k= \{ (\M, \k,g);\ \M \in X^{1+k}, \ \k\in  \Con^k(\M) \ \et \  \ g \ \text{is a}\
    \Con^k\text{-metric on}\ \M \},
   \\
  &\Y^k= \{ (\M, \k,g);\ \M \in Y^{1+k}, \ \k\in  W^{k, \infty} (\M)  \ \et \ \ g \ \text{is a}\
    W^{k, \infty}\text{-metric on}\ \M \}, 
 \end{align*}
  $k\in \N$.
 The levels of regularity we use in what follows are $\X^1$ and
 $\X^2$ on the one hand, and $\Y^1$ on the other hand.

Observe that having the manifold $\M$ with one order of regularity higher
than that for $g$ is consistent with the transformation
rules of a 2-covariant tensor on $\M$. Next, the regularity of the
function $\k$ is set equal to that of $g$ because of the definition of the
elliptic operator $A_{\k, g}$ in \eqref{eq: elliptic operator}.

With $\M$ of class $\Con^{1+k}$ (\resp $W^{1+k,\infty}$) the same
holds for $\d\M$. Once an atlas is given on $\M$ as in
Section~\ref{sec: Local coordinates}, this is quite clear.
 
 \begin{remark}
   \label{rem: diff X1 Y1}
 Note that $\Y^1$ exhibits a `tiny' loss of regularity if compared to
 $\X^1$. Yet, this loss is more like an abyss as far as the geometry underlying
 wave propagation is concerned. In fact, if considering a
 $W^{1,\infty}$-metric $g$ the Hamiltonian vector field that defines the
\bichars at higher levels of regularity is only $L^\infty$
here. Hence, the existence of \bichars is not guaranted. As a result a 
$W^{1,\infty}$-metric is too rough to state the usual GCC and consequently also to implement  standard microlocal
tools. 
\end{remark}

Based on the previous remark we will exploit the geometry of wave
propagation available for some $(\M, \k, g) \in \X^1$ yet
consider some  $(\tilde{\M}, \tilde{\k}, \tilde{ g}) \in\Y^1$ \suff close to   $(\M, \k,
g)$. Such closedness will be understood as follows. 
%%%%%%%%%%%%%%%%%%%%%%%%
% definition           %
%%%%%%%%%%%%%%%%%%%%%%%%
\begin{definition}
  \label{def: eps-close}
 Consider on the one hand $(\M, \k, g) \in \X^1$ and $\omega$
 an open subset of $\M$ (\resp $\Gamma$ an open subset of
 $\d\M$) and, on the other hand, $(\tilde{\M}, \tilde{\k}, \tilde{g}) \in \Y^1$ and $\tilde{\omega}$
 an open subset of $\tilde{\M}$ (\resp $\tilde{\Gamma}$ an open subset of $\d \tilde{\M}$).
Let $\eps>0$. One says that $(\tilde{\M}, \tilde{\k}, \tilde{ g}, \tilde{\omega})$ (\resp $(\tilde{\M}, \tilde{\k}, \tilde{ g},
\tilde{\Gamma})$) is $\eps$-close to $(\M, \k, g, \omega)$ (\resp
$(\M, \k, g, \Gamma)$)  in the $\Y^1$-topology if the
following holds 
\begin{enumerate}
 \item There exists a $W^{2, \infty}$-diffeomorphism $\psi: \M
   \rightarrow \tilde{\M}$ such that $\psi (\omega) = \tilde{\omega}$ (\resp $\psi (
   \Gamma) = \tilde{\Gamma})$).
 \item One has $\Norm{\psi^* \tilde{\k} - \k}{W^{1,\infty}(\M)} 
     +  \Norm{\psi^* \tilde{ g}  - g}{W^{1,\infty} \mathcal T^0_2 (\M)}
     \leq \eps$,
where $\Norm{.}{W^{1,\infty} \mathcal T^0_2 (\M)}$ denotes the
$W^{1,\infty}$-norm for 2-covariant tensors on $\M$. 
 \end{enumerate}
\end{definition}
%%%%%%%%%%%%%%%%%
% subsection
 \subsection{Interior and boundary observability}
 \label{Interior and boundary observability}

 Consider the following homogeneous version of the wave equation:
 \begin{align}
  \label{eq: wave equation-observation}
  \begin{cases}
     P_{\k, g}  \, u = 0
     & \dans\ \R\times\M,\\
     u=0 &  \dans\ \R \times\d\M,\\
      u_{|t=0}  = \udl{u}^0, \ \d_t u_{|t=0}  = \udl{u}^1 
      & \dans\ \M.
  \end{cases}
\end{align}
Let $\omega$ be a nonempty open subset of $\M$ and $T>0$.
Observability of the wave equation from $\omega$
in time $T$ is the following notion.
%%%%%%%%%%%%%%%%%%%%%%%%
% definition           %
%%%%%%%%%%%%%%%%%%%%%%%%
\begin{definition}[interior observability]
  \label{def: observability interior-intro}
  Let $\omega$ be a nonempty open subset of $\M$. 
  One says that the homogeneous wave equation is observable from
  $\omega$ in time $T>0$ if there exists $\Cobs>0$ such that 
  for any
  $(\udl{u}^0, \udl{u}^1) \in H_0^1(\M) \times L^2(\M)$ one has 
  \begin{equation}
    \label{eq: interior observability-intro}
    \E_{\k, g} (u) 
    \leq \Cobs
    \Norm{\bld{1}_{]0,T[ \times \omega}\, \d_t u}{L^2(\R\times\M)}^2,
  \end{equation}
  for the weak solution $u$ to \eqref{eq: wave equation-observation}.
  \end{definition}

\medskip Let $\Gamma$ be a nonempty open subset of $\d\M$. Observability of the wave equation from $\Gamma$
in time $T$ is the following notion.
%%%%%%%%%%%%%%%%%%%%%%%%
% definition           %
%%%%%%%%%%%%%%%%%%%%%%%%
\begin{definition}[boundary observability]
  \label{def: observability boundary-intro}
   Let $\Gamma$ be a nonempty open subset of $\d\M$. 
  One says that the homogeneous wave equation is observable from
  $\Gamma$ in time $T>0$ if there exists $\Cobs>0$ such that 
  for any
  $(\udl{u}^0, \udl{u}^1) \in H_0^1(\M) \times L^2(\M)$ one has 
  \begin{equation}
    \label{eq: boundary observability-intro}
    \E_{\k, g} (u) 
    \leq \Cobs
    \Norm{\bld{1}_{]0,T[ \times \Gamma}\,
      \d_{\n} u_{|\R\times \d\M}}{L^2(\R\times\d\M)}^2,
  \end{equation}
  for the weak solution $u$ to \eqref{eq: wave equation-observation}.
  \end{definition}

%%%%%%%%%%%%%%%%%
% subsection
\subsection{Main results and open questions}
\label{results}

The following observability  results were proven in~\cite{Burq:1997}.
\begin{theorem}[Burq, 97]
  \label{thm: Burq97}
  Let $(\M, \k, g) \in \X^2$.
  
  {\em \bfseries Interior observability.}
  Let $\omega$ be an open subset of $\M$ that satisfies the interior geometric control condition (see
  Definition~\ref{control-geo-interior} below for a precise
  description) associated with the infimum time $T_{GCC}(\omega)$.
  Let $T>T_{GCC}(\omega)$. Then, the wave equation is observable
  from $\omega$ in time $T$.
 
  {\em \bfseries  Boundary observability.} Let $\Gamma$ be an open subset of $\d \M$ such that
   $\Gamma$ satisfies the boundary geometric control condition
   (see Definition~\ref{control-geo-boundary} below for a precise
   description) associated with the infimum time $T_{GCC}(\Gamma)$.
Let $T>T_{GCC}(\Gamma)$. 
   Then, the
   wave equation is observable from $\Gamma$ in time $T$.
\end{theorem}
The proof of these results essentially relies on pseudo-differential calculus and microlocal tools, namely the microlocal defect measures and their localization and propagation properties. 
 
One of our main contributions  in the present article is to improve upon  the
regularity assumptions on the metric $g$, the function $\k$ and the
manifold $\M$.
\begin{theorem}
  \label{theoprinc}
  Let $(\M, \k, g) \in \X^1$. The two conclusions of Theorem~\ref{thm: Burq97} hold.
\end{theorem}
In fact, we prove a stronger result, namely that these observability results are stable
by small perturbations in $\mathcal{Y}^1$, that is, perturbations
that are slightly less smooth. 
\begin{theorem}\label{theo-perturbation}
  Let $(\M, \k, g) \in \X^1$.
  
  {\em \bfseries Interior observability.}
  Let $\omega$ be an open subset of $\M$ that satisfies the interior geometric control condition associated with the infimum time $T_{GCC}(\omega)$.
  Let $T>T_{GCC}(\omega)$. Then, there exists $\eps>0$ such that if  $(\tilde{\M}, \tilde{\k}, \tilde{ g},
  \tilde{\omega})$ is $\eps$-close to $(\M, \k, g,\omega)$ in the $\Y^1$-topology in the
  sense of Defintion~\ref{def: eps-close} for $(\tilde{\M}, \tilde{\k}, \tilde{ g}) \in \Y^1$
  and $\tilde{\omega}$ an open subset of $\tilde{\M}$, then the wave equation
  associated with $P_{\tilde{\k}, \tilde{ g}}$ on $\tilde{\M}$  is interior observable
  from $\tilde{\omega}$ in time $T$.

  {\em \bfseries  Boundary observability.}
  Let $\Gamma$ be an open subset of $\d \M$ such that
   $\Gamma$ satisfies the boundary geometric control condition
  associated with the infimum time $T_{GCC}(\Gamma)$.
Let $T>T_{GCC}(\Gamma)$. 
   Then, there exists $\eps>0$ such that if  $(\tilde{\M}, \tilde{\k}, \tilde{ g},
  \tilde{\Gamma})$ is $\eps$-close to $(\M, \k, g,\Gamma)$ in the $\Y^1$-topology in the
  sense of Defintion~\ref{def: eps-close} for $(\tilde{\M}, \tilde{\k}, \tilde{ g}) \in \Y^1$
  and $\tilde{\Gamma}$ an open subset of $\d \tilde{\M}$, then the wave equation
  associated with $P_{\tilde{\k}, \tilde{ g}}$ on $\tilde{\M}$  is boundary observable
  from $\tilde{\Gamma}$ in time $T$.
\end{theorem}
\begin{remark}
First, as pointed out in Remark~\ref{rem: diff X1 Y1}
above, note that standard microlocal
tools cannot be used at the $W^{1,\infty}$ level of regularity of $\tilde{\k}$ and
the metric $\tilde{ g}$. 

Second, Theorem~\ref{theoprinc} shows that the observability property
is stable by small Lipschitz perturbations around rough ($\Con^1$)
metrics satisfying GCC. We exhibited in~\cite[Remark 1.13]{BDLR0} an
example showing that the observation property is {\em not} stable by
small (even smooth) perturbations of the geometry/metric around geometries satisfying only the obervation property.  This counter example is actually connected to an example due to
G.~Lebeau \cite{Lebeau:96}. In particular, this  shows that our perturbation argument
will have to be performed on the proof of the fact that the geometric
control condition implies observability and not on the final property itself.
Since Theorem~\ref{theoprinc} is a straightforward consequence of
Theorem~\ref{theo-perturbation} we will hence focus on the proof of
Theorem~\ref{theo-perturbation} in what follows.
\end{remark}

\medskip
In the two  kind of observability-inequality results stated above, the GCC of
Definitions~\ref{control-geo-interior} and \ref{control-geo-boundary}
appear as sufficient conditions. We also formulate the following weak
GCC condition
\begin{center}
     {\em  For all point in the tangent bundle, at least one
       generalized geodesics intitiated at this point enters any
       region larger than the
      observation region in some time $T>0$.}
  \end{center}
  A more precise definition is given in Section~\ref{sec: Necessary
    geometric control condition}. Observe that this condition reduces to  the classical necessary condition in the case of uniqueness of
  \gbichars.
  
We prove the following result proven in Section~\ref{sec: Measure equation at isochrones and necessary geometric control condition}.
\begin{theorem}[Necessary geometric control conditions]
  \label{theorem: necessary control conditions-intro}
  Let $(\M, \k, g) \in \X^1$.
  
  {\em \bfseries Interior observability.}
  Let $\omega$ be an open subset of $\M$ and $T>0$ such that interior
  observability holds from $\omega$ in time $T$. Then the weak
  geometric control condition holds.

    {\em \bfseries  Boundary observability.}
  Let $\Gamma$ be an open subset of $\d \M$ and $T>0$ such that
  boundary  observability holds from $\Gamma$ in time $T$. Then the weak
  geometric control condition holds.
\end{theorem}

\medskip
In the framework of the present article, sharpness of
  GCC as a sufficient condition for observability is an open question.
  Sharpness of the weak GCC as a necessary condition for observability
  is also open.  In the case where generalized geodesics are not unique,
  there is quite a gap between the GCC and the weak GCC. Note that
  this gap closes as soon as uniqueness holds. Note also that a lack of
  uniqueness of generalized geodesics can be connected to a low
  regularity of the coefficients, as in the present setting. However,
  even in the case of smooth coefficients and a smooth manifold, an infinite 
  contact order of a generalized geodesics with the boundary can be a
  source of nonuniqueness. We refer to the Taylor
  example~\cite{Taylor:75} (see also \cite[Example
    24.3.11]{Hoermander:V3}). No thorough study of nonuniqueness
  issues at boundary has been carried out for $\Con^k$ coefficients, $k\geq 2$, up
  to our knowledge.

%%%%%%%%%%%%%%%%
\subsection{Method to prove observability and outline}
\label{sec: Method to prove observability and outline}
We present here the scheme of the proof of
the observability inequalities even
though some material needs to be introduced. Yet, this will provide the
reader with a road map.

The strategy of the proof of observability follows the following
steps:
\begin{itemize}
\item
First, we reduce the observation estimate for general data to a
high-frequency observation estimate for semi-classically localized
initial data. This step is quite classical \cite{BDLR-abs-arg}.
\item Second, we proceed by contradiction and obtain sequences
of $L^2$-normalized initial data that are spectrally localized and
vanishing asymptotically in the observation region.
\item Third, associated with
these sequences is a semi-classical measure $\mu$ that characterises in phase space 
points where mass concentrates asymptotically. The main result of the
present article is a  propagation equation provided  in  
Theorem~\ref{thm: equationmesure} and fulfilled by the measure $\mu$.

\item Fourth, we exploit the result of the companion article~\cite{BDLR2} stated here in Theorem~\ref{theo-propagation} and we deduce that the support of $\mu$
is a union of maximal \gbichars. This leads to a
contradiction with $\mu$ vanishing in the observation region.
\end{itemize}

The remainder of the article is organized as follows.

In Section~\ref{geometry} we introduce the necessary geometrical
notions to precisely state the geometric control condition (GCC) in our low
regularity framework and to state the main result of the companion article~\cite{BDLR2}.
In Section~\ref{sec.4} we
perform some classical {\em a priori} estimates for the normal derivatives
of solutions to wave equations and we recall that an observability
inequality is equivalent to an exact controllability result for the
wave equation.  Section~\ref{Semi-classical reduction} is devoted to a
semi-classical reduction of observability estimates and the definition
of semi-classical measures. In Section~\ref{sec: Semi-classical operators and
  measures}, we recall and introduce some aspects of semi-classical
pseudo-differential operators with minimal regularity properties of
the symbols. We also recall the notions of semi-classical measures and
some of their properties.  In Section~\ref{sec: Proof of the observability
  results} we write the contradiction argument that leads to the proof
of a semi-classical observability inequality. This generates a
semi-classical measure associated with a sequence of solutions to the
wave equation and a semi-classical measure associated with their
normal derivatives on the boundary. The measure-propagation equation
that links these two measures is stated in Theorem~\ref{thm:
  equationmesure} of Section~\ref{sec: The
  measure propagation equation} allowing one to conclude the proof of
the semi-classical observability inequality. Sections~\ref{sec: Proof
  of the propagation equation1} to \ref{sec: Proof of the propagation
  equation2} are dedicated to the proof of the measure-propagation
equation of Theorem~\ref{thm:
  equationmesure}. Section~\ref{sec: Proof of the propagation equation1} exposes
the commutator argument that is the foundation of the
measure-propagation equation.  In Section~\ref {sec: More on
  semi-classical symbols and operators} we present a Weierstrass
division argument to be applied to the test functions to prove the
measure-propagation equation. This leads to symbols with low regularity
and low decay in the conormal direction. Further analysis for such
symbols and associated operators is provided.  Finally, in
Section~\ref{sec: Proof of the propagation equation2} the different
symbols obtained in the Weierstrass division are quantized leading to
the proof of the measure-propagation equation.

In Section~\ref{sec: Measure equation at isochrones and necessary
  geometric control condition} we prove that observability implies a
weak GCC (Theorem~\ref{theorem: necessary control
  conditions-intro}). The proof is based on the measure equation that
is stated in Theorem~\ref{theorem: measure equation at t=0} and proven
in Section~\ref{sec: Proof the measure equation at isochrones}.

%%%%%%%%%%%%%%%%%
% subsection
\subsection{Acknowledgements}
This research was partially supported by Agence Nationale de la
Recherche through project ISDEEC ANR-16-CE40-0013 (NB), by the
European research Council (ERC) under the European Union’s Horizon
2020 research and innovation programme, Grant agreement 101097172 -
GEOEDP (NB), and by the Tunisian Ministry for Higher Education and
Scientific Research within the LR-99-ES20 program (BD). The authors
acknowledge GE2MI (Groupement Euro-Maghr\'ebin de Math\'ematiques et
de leurs Interactions) for its support.

%%%%%%%%%%%%%%%%%
% Section
%%%%%%%%%%%%%%%%%
\section{Geometry}
\label{geometry}
In this section we define the basic notions required to understand the
statements of our results.  We will work in special
{\em quasi-normal geodesic} coordinates near the boundary. We refer
to~\cite[Section 5]{BDLR2} for a more thorough and intrinsic
presentation of the geometry.

\medskip
Here, $\M$ is a $\Con^2$-compact connected Riemannian manifold
of dimension $d$ with boundary, endowed with a $\Con^1$-metric $g$. An example would
be a bounded open subset $\Omega$ of $\R^d$ with a $\Con^{2}$-boundary,
that is, with the boundary given locally by $\varphi(x)=0$ with
$\varphi \in \Con^{2}(\R^d)$ and $d \varphi \neq 0$. Then
$\M = \Omega \cup \d\Omega$ and one can simply consider the Euclidean
metric.  In the spirit of this simple example, consider an open
$d$-dimensional manifold\footnote{The manifold $\tM$ can be
  constructed by embedding $\M$ in $\R^{2d}$ thanks to the
  Whitney theorem \cite{Whitney:36}.} $\tM$ such that $\M \subset \tM$ and 
extend the metric $g$ to a \nhd of $\M$ is a $\Con^1$-manner.

%%%%%%%%%%%%%%%%%
% subsection
\subsection{Local coordinates}
\label{sec: Local coordinates}
Equip a compact \nhd $\hM$ of $\M$ in $\tM$ with a finite
$\Con^2$-atlas. A local chart is denoted $(\hO, \chdiff)$ with $\hO$
an open subset of $\hM$ and $\chdiff$ a one-to-one map from $\hO$ onto an open
subset of $\R^d$. Charts can be chosen so that 
\begin{align*}
  %\label{eq: local chart boundary}
  &\chdiff(\hO \cap \M) = \chdiff(\hO) \cap \{x_d\geq 0\}
    \ \text{is an open subset of}\  \ovl{\Rdp}, \\
  &\chdiff(\hO \cap \d\M) = \chdiff(\hO) \cap \{x_d=0\},
    \  \et \  
    \chdiff(\hO \setminus \M) = \chdiff(\hO) \cap \{x_d<0\},
    \notag
\end{align*}
if $\hO \cap \d\M \neq \emptyset$. Denote the local coordinates by $x=(x', x_d)$ with $x'
\in \R^{d-1}$. Note that $\M$ being compact it
contains its boundary $\d\M$.

In a  local chart, the metric $g$ is
given by $g_x = g_{ij} (x) dx^i \otimes dx^j$, where  $g_{ij}\in
\Con^1(\chdiff(\hO))$. We use the classical notation $(g^{ij}(x) )_{i,j}$ for  the
inverse of  $(g_{ij}(x) )_{i,j}$. 
The metric $g_x = (g_{ij}(x) )_{i,j}$ provides an inner  product
on $T_x \M$. The metric $g_x^* = g^{ij}(x) d \xi_i \otimes d \xi_j$ provides an
inner product
on $T^*_x \M$, denoted $g^*_x(\xi,\tilde{\xi})$, for $\xi, \tilde{\xi} \in T^*_x
\M$. Define the associated norm
\begin{align*}
  \norm{\xi}{x} = g^*_x(\xi,\xi)^{1/2}.
\end{align*}
Near a
boundary point, local coordinates  are chosen according to the following proposition. They simplify
the exposition of some geometrical notions and are key in arguments developed in what follows.
%%%%%%%%%%%%%%%%%%%%%%%%
% proposition          %
%%%%%%%%%%%%%%%%%%%%%%%%
\begin{proposition}[quasi-normal geodesic coordinates]
  \label{prop: quasi-normal coordinates}
  Suppose $m^0 \in \d\M$. There exists a
  $\Con^2$-local chart $(\hO, \chdiff)$ such that $m^0 \in \hO$,  $\chdiff(m) = (x',z)$, with $x' \in \R^{d-1}$ and $z \in \R$, and 
  \begin{enumerate}
   \item $\chdiff(\hO\cap \M) =  \{ z\geq 0\} \cap \chdiff(\hO)$, $\chdiff(\hO\cap
      \d\M) =  \{ z= 0\} \cap \chdiff(\hO)$,  and $\chdiff(\hO\setminus \M) =  \{ z<0\} \cap \chdiff(\hO)$ ;
    \item at the boundary, the representative of the metric has the form 
      \begin{align*}
        g(x',z=0) = \sum_{1\leq i, j \leq d-1} g_{i j}(x',z=0)  dx^i \otimes dx^j
        + |dz|^2.
      \end{align*}
  \end{enumerate}
\end{proposition}
In other words the matrix of $g=(g_{ij})$ has the block-diagonal form
   {\em at the boundary}
   \begin{equation}
     \label{eq: structure metric quasi-geodesic coordinates}
      g (x',z=0) = \begin{tikzpicture}[
  every left delimiter/.style={xshift=.75em,yshift=0.4em},
  every right delimiter/.style={xshift=-.75em,yshift=0.4em},
  style1/.style={
  matrix of math nodes,
  every node/.append style={text width=#1,align=center,minimum height=1.8ex},
  nodes in empty cells,
  left delimiter=(,
  right delimiter=)
  },
baseline=(current  bounding  box.center)
]
\matrix[style1=0.27cm] (1mat)
{
  & &  & \\
  &  & &\\
  &  & &\\
  & &  &  \\
  & &  &  \\
};
\node[font=\normalsize] 
  at (1mat-1-4.north) {$0$};
  \node[font=\normalsize] 
  at (1mat-2-4.north) {$\vdots$};
\node[font=\normalsize] 
  at (1mat-3-4) {$0$};
\node[font=\normalsize] 
  at (1mat-4-4.south) {$1$};
 \node[font=\normalsize] 
  at (1mat-4-1.south) {$0$};
\node[font=\normalsize] 
  at (1mat-4-2.south) {$\ \cdots$};
  \node[font=\normalsize] 
  at (1mat-4-3.south) {$\ 0$};
  \node[font=\Huge] 
  at (1mat-2-2) {$*$};
\end{tikzpicture}.
    \end{equation}
  Naturally, the same form holds for $g_x^* =
  (g^{ij}(x))$ at the boundary. 
One deduces that 
  \begin{align*}
    g_{j d}(x',z)  = z h_{j d}(x',z)  \ \ \et \ \ 
    g_{d d}(x',z)  = 1 +  z h_{d d}(x',z),
  \end{align*}
  for some continuous functions $h_{j d}$, $j=1, \dots, d$.

  Proposition~\ref{prop: quasi-normal coordinates} can be found in
  \cite{BM:21} with a different regularity level. A proof of
  Proposition~\ref{prop: quasi-normal coordinates}  at the regularity
  level we consider is written in 
  Appendix~B of \cite{BDLR2} with a generalization to other levels of
  regularity.
%%%%%%%%%%%%%%%%%%%%%%%%
% remark               %
%%%%%%%%%%%%%%%%%%%%%%%%
\begin{remark}
  \label{remark: no normal geodesic coodinates}
  Because of the low regularity of $g$ and $\M$ one {\em
    cannot choose} normal geodesic coordinates, that is, local
  coordinates for which $g_{j d} = g_{d j} = 0$ for $j \neq d$ and
  $g_{d d}=1$ near a point $m^0$ of the boundary. The
  coordinates that Proposition~\ref{prop: quasi-normal coordinates}
  provides only have this property in a \nhd of $m^0$ {\em within} the
  boundary $\d\M$.
\end{remark}

\medskip One sets $\L = \R \times \M$ and
$\hL = \R \times \hM$ .  From a local chart $(\hO, \chdiff)$
in the atlas for $\hM$ one defines a map
$\cdiffL: (t,m) \mapsto (t, \chdiff(m))$ from $\O = \R\times \hO$ onto
$\R\times \chdiff(\hO)$, yielding a local chart $(\O, \cdiffL)$  for $\hL$ and thus a
finite atlas.

For $x = \chdiff(m)$, $m \in \hO \cap \M$, denote by $v = (v',v^d)$
and $\xi =(\xi',\xi_d)$ the associated coordinates in $T_m \M$ and
$T_m^* \M$, with $v', \xi'\in \R^{d-1}$ and $v^d, \xi_d \in \R$.  We
write $T_x\M$ and $T_x^*\M$ by abuse of notation.  In what follows, it
will be convenient to write $z$ in place of $x_d$, in particular for
the local coordinates given by Proposition~\ref{prop: quasi-normal
  coordinates}. Accordingly we denote the associated cotangent
variable $\xi_d$ by the letter $\zeta$, that is, $\xi= (\xi', \zeta)$.
We however do not change the notation for the associated tangent
variable $v^d$.  With local charts at the boundary given by
Proposition~\ref{prop: quasi-normal coordinates}, if $x \in \d\M$ and
$v \in T_x\d\M$ then $v = (v',0)$ and we use the bijective map
$(\xi',0) \mapsto \xi'$ to parameterize $T_x^* \d\M$.

Also classically set
\begin{align*}
  T \M = \bigcup_{x \in \M} \{x\} \times T_x \M, 
  \quad \TM = \bigcup_{x \in \M} \{x\} \times T_x^* \M\\
  \big(\text{\resp}\ T \hM = \bigcup_{x \in \M} \{x\} \times T_x \hM, 
  \quad \ThM = \bigcup_{x \in \M} \{x\} \times T_x^* \hM\big).
\end{align*}
With $\M$ containing its boundary $\d\M$, one sees that $ T \M$ (\resp
$\TM$) contains  $\{x\} \times T_x \M$ (\resp $\{x\} \times T_x^* \M$)
for $x \in \d\M$. We denote by $\dTM$ the boundary of $\TM$ that is
the set of $(x,\xi)$ with $x \in \d\M$. 
In the local coordinates, $\dTM$ is
given by $\{ z=0\}$ and $\TM$ by $\{ z\geq 0\}$.

In the associated local chart on $\L$, the representative of $(t,m)
\in \L$ is $(t,x) = (t,x',z)$.  We use the letter $\y$ to denote
an element of $\TL$, that is, $\y = (t,x; \tau, \xi)$ with $(t,x) \in
\L$, $\tau \in \R$ and $\xi \in T^*_x\M$.  Classically, we write $\TL
\setminus 0$ for the set of points $\y = (t,x; \tau, \xi)$ with
$(\tau, \xi)\neq0$.  The boundary $\dTL$ is the set of points $\y=
(t,x;\tau,\xi)$ such that $x \in \d\M$. Note that $\dTL$ is locally
given by $\{ z=0\}$ and $\TL$ is locally given by $\{ z\geq 0\}$.

%%%%%%%%%%%%%%%%%
% subsection
\subsection{Wave operators and bicharacteristics}
\label{sec: wave operator}

On the manifold $\M$ consider the elliptic operator $A = A_{\k,g}= \k^{-1} \divg (\k
\nablag)$, that is, in local coordinates
\begin{align*}
  A f = \k^{-1} (\det g)^{-1/2} 
  \sum_{1\leq i, j \leq d}  \d_{x_i}\big( 
  \k  (\det g)^{1/2} g^{i j}(x)\d_{x_j} f
  \big).
\end{align*}
Its principal symbol is simply
$a(x,\xi) = - g^*_x(\xi,\xi)  = - g_x^{i j} \xi_i \xi_j = - \norm{\xi}{x}^2$. Note
that for $\k=1$, one has $A = \Delta_g$, the Laplace-Beltrami
operator associated with $g$ on $\M$. Together with $A$ consider
the wave operator $ P_{\k,g} = \d_t^2 - A_{\k,g}$. Its
principal symbol in a local chart is given by 
\begin{align*}
  p(\y) = -\tau^2 + \norm{\xi}{x}^2.
\end{align*}
Note that $p(\y)$ is smooth in the variables $(\tau,\xi)$ and $\Con^1$
in $x$. 

For a function $f$ of the variable $\y$, the Hamiltonian vector field
$\Hamiltonian_f$ is defined by $\Hamiltonian_f(h)  = \{ f, h\}$, where
$\{.,.\}$ is the Poisson bracket. 
In local coordinates one has 
\begin{align*}
  %\label{eq: Hp}
  \Hp (\y) &= \d_\tau p (\y) \d_t 
             + \nabla_\xi p(\y) \cdot \nabla_x
             -\nabla_x p (\y) \cdot \nabla_\xi \\
           &=  - 2 \tau \d_t + 2 g^{ij}(x) \xi_i \d_{x_j} 
             - \d_{x_k} g^{ij}(x) \xi_i \xi_j  \d_{\xi_k}.\notag
\end{align*}
Recall the following definition. 
%%%%%%%%%%%%%%%%%%%%%%%%
% definition           %
%%%%%%%%%%%%%%%%%%%%%%%%
\begin{definition}
  \label{def: bichar}
  Suppose $V$ is an open subset of $\TL \setminus \dTL$ and $J\subset \R$ is an interval. A $\Con^1$-map $\gamma: J \to V \cap \Char p$ is called a
  \bichar in $V$ if 
  \begin{align*}
    \frac{d}{d s} \gamma(s) = \Hp \big( \gamma(s) \big), 
    \qquad  s \in J. 
  \end{align*}
  It is called {\em maximal} in $V$ if it cannot be extended by another
  \bichar also valued in~$V$. 
\end{definition}
Note that $\transp{\Hp} f (\y) = 2 \tau \d_t f (\y)  - 2 \d_{x_j}
  \big(g^{ij}(x) \xi_if (\y)  \big) 
  + \d_{\xi_k} \big(\d_{x_k} g^{ij}(x) \xi_i \xi_j f (\y) \big)$ and 
deduce 
\begin{align*}
  %\label{eq: transp Hp}
  \transp{\Hp}  = - \Hp.
\end{align*} 
Recall also that 
\begin{align}
  \label{eq: derivative along bichar}
  \Hp f (\gamma(s)) = \frac{d}{d s} f (\gamma(s)),
  \ \ \si \ \gamma \ \text{is a \bichar}. 
\end{align}

%%%%%%%%%%%%%%%%%
% subsection
\subsection{A partition of the cotangent bundle at the boundary}
\label{sec: A partition of the cotangent bundle at the boundary-intro}

Denote  by $\pdTL \subset \dTL$ the bundle of points $\y = (\y',0)= 
(t,x',z=0,\tau,\xi',0) \in \TL$ for $\y' =
(t,x',z=0,\tau,\xi') \in T^* \d\L$. 
 Identifying $\y'$ and $(\y',0)$ as presented
above  thanks to the chosen  local coordinates  allows one to
indentify $\pdTL$ and $T^* \d\L$.

Denote by $\ppi$ the map from $\dTL$ into $\pdTL$ given by
\begin{align*}
  \ppi (t,x',z=0,\tau,\xi',\zeta) = (t,x',z=0,\tau,\xi',0).
\end{align*}
%%%%%%%%%%%%%%%%%%%%%%%%
% definition           %
%%%%%%%%%%%%%%%%%%%%%%%%
\begin{definition}[elliptic, glancing, and hyperbolic regions]
  \label{def: E', H', G'-intro}
  One partitions $\pdTL$  into three homogeneous regions.
  \begin{enumerate}
  \item The elliptic region $\pEb = \pdTL \cap \{ p >0\}$; 
    if
    $\y \in \pEb$ it is called an elliptic point.          	
  \item The glancing region $\pGb = \pdTL \cap \{ p =0\}$;   if $\y
    \in \pGb$ it is called a glancing point.
  \item The hyperbolic region $\pHb =\pdTL \cap \{ p < 0\}$; if
    $\y \in \pHb$ it is called a   hyperbolic point.
  \end{enumerate}
\end{definition}
Since $p (\y) = -\tau^2 + \zeta^2 + g_x(\xi',\xi')_x$ by \eqref{eq: structure metric quasi-geodesic coordinates} if $\y \in
\dTL$, one has the following properties:
\begin{enumerate}
  \item If $\y \in \pEb$ then $\ppi^{-1}\big( \{ \y\}\big) \cap \Char p
    = \emptyset$.
  \item  If $\y \in \pGb$ then $\ppi^{-1}\big( \{ \y\}\big) \cap \Char p
    = \{ \y\}$.
  \item If $\y = (t,x',z=0,\tau,\xi',0) \in \pHb$ then $\ppi^{-1}\big( \{ \y\}\big) \cap \Char p
    = \{ \y^-, \y^+\}$, where
    \begin{align}
      \label{eq: relevement H-G-intro}
      \y^\pm = (t,x',z=0,\tau,\xi^\pm), \ \ \where \ \xi^\pm=(\xi',\zeta^\pm)
       \   \avec \  \zeta^\pm = \pm \sqrt{- p(\y)}.
    \end{align}
\end{enumerate}
Associated with the previous partition of $\pdTL$  is a partition of
$\Char p \cap \dTL$. Indeed, if $\y \in \Char p \cap \dTL$ then
$\ppi(\y) \in \pdTL$ and 
$p \big( \ppi(\y)\big ) \leq 0$. 
Note that having $\y \in \Char p \cap \dTL$ and $p \big( \ppi(\y)\big
) = 0$ is equivalent to having $\y \in \pGb$. 
%%%%%%%%%%%%%%%%%%%%%%%%
% definition           %
%%%%%%%%%%%%%%%%%%%%%%%%
\begin{definition}[partition of $\Char p$ at the boundary]
  \label{def: H, G-intro}
  One partitions $\Char p \cap \dTL$  into two homogeneous regions
  $\Gb$ and $\Hb$:
  \begin{enumerate}
  \item $\Gb = \pGb$; $\y \in \Gb \ \Equiv \ \y \in \Char p$ and $\ppi(\y) = \y$.
  \item $\y \in \Hb$ if $\y \in \Char p$ and $\ppi(\y) \in \pHb$. It
    is also called a hyperbolic point. If $\y =
    (t,x',z=0,\tau,\xi',\zeta)$ one says that $\y  \in \Hb^+$ if $\zeta>0$ and
    $\y  \in \Hb^-$ if $\zeta<0$.
    \end{enumerate}
\end{definition}
Thus, if $\y \in \pHb$ then $\ppi^{-1}\big( \{ \y\}\big) \cap \Char p
= \{ \y^-, \y^+\}$ with $\y^+ \in \Hb^+$ and $\y^- \in \Hb^-$, with
$\y^\pm$ as given in \eqref{eq: relevement H-G-intro}.

Introducing the following involution on $\dTL$
\begin{align*}
  %\label{eq: involution-intro}
  \Sigma(t,x',z=0,\tau,\xi',\zeta) = (t,x',z=0,\tau,\xi',-\zeta),
\end{align*}
one finds that $\Sigma (\y^-) = \y^+$ if $\y \in \pHb$. Thus, $\Sigma$
is a one-to-on map from  $\Hb^-$ onto $\Hb^+$.  

%%%%%%%%%%%%%%%%%
% subsection
\subsection{Glancing region, gliding vector field, and generalized bicharacteristics}
\label{sec: glancing region, generalized bichar-intro}

One computes
\begin{align*}
  \Hpz (\y) = \Hpz (x,\xi)= 2 g^{dj}(x) \xi_j\qquad  \text{(recall $x_d = z$ \ \et \ $\xi_d = \zeta$)}. 
\end{align*}
Observe that $\Hpz$ is a $\Con^1$-function and that $\Hpz_{|z=0} = 2
\zeta$ in the present local coordinates. Hence, locally one has
\begin{align*}
  \pGb= \Gb=\{ z=\Hpz=p=0 \}
  \ \ \et \ \
  \Hb^\pm =\{ z=p=0, \ \Hpz \gtrless 0\}.
\end{align*}
With \eqref{eq: derivative along bichar} this means that a \bichar
going through a point $\y\in \Hb$ has a contact of order exactly one
with the boundary: it is transverse to $\dTL$. A \bichar going through
a point $\y\in \Gb$ has a contact of order greater than or equal to
two: it is tangent to $\dTL$.

One can further compute $\Hppz$.  It is a continuous and gives the
following partition of $\Gb$.
%%%%%%%%%%%%%%%%%%%%%%%%
% definition           %
%%%%%%%%%%%%%%%%%%%%%%%%
\begin{definition}[partition of $\Gb$]
  \label{def: diffrative}
  Introduce
  \begin{align*}
    &\sdGb = \{ \y \in \Gb; \Hppz (\y) >0\},\\
    &\glGb =\{ \y \in \Gb; \Hppz (\y) =0\}, \\ 
    &  \sgGb =\{ \y \in \Gb; \Hppz (\y) <0\}.
  \end{align*}
    One calls $\sdGb$ the diffractive
    set, $\sgGb$ the gliding set.
    One calls $\glGb$ the glancing set of order three: if $\y^0 \in \glGb$ a \bichar that goes through $\y^0$ has
    a contact with the boundary of order greater than or equal to three.
\end{definition}

On $\pdTL$ one defines 
\begin{align*} 
    \HpG (\y)  = \Big(\Hp + \frac{\Hppz}{\Hzzp}\Hz\Big)  (\y) ,
  \end{align*}
referred to as the {\em gliding vector field}. In the present coordinates
one has $\Hzzp=2$. Define the following  vector field on $\TL$ 
\begin{align*}
  \XG(\y)=
  \begin{cases}
    \Hp (\y) & \si \ \y \in \TL\setminus \sgGb,\\
    \HpG (\y) & \si \ \y \in \sgGb,
  \end{cases}
\end{align*} 
that is, $\XG = \Hp + \unitfunction{\sgGb} (\HpG - \Hp) $. More
explainations on the vector field $\HpG$ are given in
Section~5 in the companion article~\cite{BDLR2}.  
%%%%%%%%%%%%%%%%%%%%%%%%
% definition           %
%%%%%%%%%%%%%%%%%%%%%%%%
\begin{definition}[\gbichar]
	\label{def: generalized bichar-intro}
  Let $J \subset \R$ be an interval, $B$ a discrete subset of $J$, and
  \begin{align*}
    \gammaG: J\setminus B 
    \to \Char p \cap \TL.
  \end{align*}
  One says that $\gammaG$ is a \gbichar if  the
  following properties hold:
  \begin{enumerate}
  \item For $s \in J \setminus B$, $\gammaG(s) \notin \Hb$ and the map $\gammaG$ is differentiable
    at $s$ with
	\begin{equation*}
		\frac{d}{d s} \gammaG(s) = \XG\big( \gammaG(s)\big).
	\end{equation*}
  \item If $S\in B$, then $\gammaG (s) \in \TL \setminus \dTL$ for $s
    \in J\setminus B$ \suff
    close to $S$ and moreover
    \begin{enumerate}
      \item  if $[S-\eps,S] \subset J$ for some $\eps>0$, then $\gammaG(S^-) = \lim_{s \to S^-} \gammaG(s) \in
      \Hb^-$;
       \item if $[S,S+\eps] \subset J$ for some $\eps>0$, then $\gammaG(S^+) = \lim_{s \to S^+} \gammaG(s) \in
      \Hb^+$;
    \item and if $[S-\eps,S+\eps]\subset J$ for some $\eps>0$, then
    $\gammaG(S^+) = \Sigma \big(\gammaG(S^-)\big)$.
    \end{enumerate}   
  \end{enumerate}
\end{definition}
Recall that $\TL$ contains its boundary $\dTL$; as a result a
\gbichar $\gammaG (s)$ may lie in the boundary for $s$ in
some interval. Details on \gbichars can be found in
Section~5 of the companion article~\cite{BDLR2}.

\medskip
When one refers to a (generalized) \bichar one often means the points visited in
$\TL$ by
$s\mapsto \gammaG(s)$ as $s$ varies, that is, 
\begin{align*}
	\{ \gammaG(s); \ s \in J \setminus B\}.
\end{align*}
Observe however that this set may not be a closed set if $B\neq
\emptyset$ as its intersection with $\Hb$ is empty. 
Consequently, we rather use its closure to describe the set of reached points.
%%%%%%%%%%%%%%%%%%%%%%%%
% definition           %
%%%%%%%%%%%%%%%%%%%%%%%%
\begin{definition}[\gbichar]
  \label{def: generalized bichar 2-intro}
	By \gbichar one also refers to 
	\begin{align*}
	\GammaG = \ovl{\{ \gammaG(s); \ s \in J \setminus B\}} 
	= \{ \gammaG(s); \ s \in J \setminus B\} 
	\cup \bigcup_{s\in B} \{ \gammaG(s^-), \gammaG(s^+)\}.
\end{align*}
\end{definition}

\bigskip
The following theorem states that for every point of  $\TL$ one can find a maximal \gbichar that goes through this point. 
%%%%%%%%%%%%%%%%%%%%%%%%
% theorem              %
%%%%%%%%%%%%%%%%%%%%%%%%
\begin{theorem}
  \label{theorem: existence bichar}
  Suppose $J\setminus B \ni s \mapsto \gammaG(s) = (t(s), x(s),
  \tau(s), \xi(s))$ is a \gbichar. If $\gammaG$ is maximal
  then $J = \R$. Moreover, $t(\R) = \R$ if $\tau(s) = \cst \neq 0$.
  
  If $\y^0 \in \Char p \cap \TL$ there exists a maximal \gbichar $s \mapsto \gammaG(s)$ with $s \in \R \setminus B$ 
  such that $\gammaG(0) = \y^0$ if $\y^0 \notin \Hb$ and 
$\gammaG(0^\pm) = \y^0$ if $\y^0 \in \Hb^\pm$.
\end{theorem}
Note that there is no uniqueness of such a  maximal 
  \gbichar because of the limited smoothness of $\XG$.
This result is classical in the case of smooth coefficients; see
\cite{MS:78} or \cite[Section 24.3]{Hoermander:V3}. Here, in
the case of the present limited smoothness it can be proven with the
arguments developed in the companion article; see
\cite[Appendix~A]{BDLR2} for a proof.

\subsection{Geometric control conditions}
\medskip
In the present  low regularity framework we state the  geometric control
conditions (GCC)
that coincide with the usual definitions found in the
literature. First, we state the interior geometric control condition.
%%%%%%%%%%%%%%%%%%%%%%%%
% definition           %
%%%%%%%%%%%%%%%%%%%%%%%%
\begin{definition}[interior geometric control]
\label{control-geo-interior}
Let $\omega$ be an open subset of $\M$. One
says that $\omega$ controls geometrically the manifold $\M$
if there exists $T>0$ such that any \gbichar reaches a point
above  $]0,T[\times \omega$. One says that $(\omega, T)$ fulfills GCC.
In such case, one sets
\begin{align*}
  T_{\text{GCC}}(\omega) = \inf \{ T>0; \ (\omega, T) \ \text{fulfills GCC}\}.
\end{align*}
\end{definition}
To state the boundary geometric control we introduce the notion of
boundary escape point.
%%%%%%%%%%%%%%%%%%%%%%%%
% definition           %
%%%%%%%%%%%%%%%%%%%%%%%%
\begin{definition}[boundary escape point]
  \label{def: escape point}
  \begin{enumerate}
    \item a point $\y \in \dTL$ is said to be a boundary escape point in the future if  locally in time all \bichars initiated at $\y$ immediately leave $\TL$ in
the future. One denotes by $\ESC^F$ the set
of all such points
\item a point $\y \in \dTL$ is said to be a boundary escape point in
  the past  if locally in time all \bichars initiated at $\y$ immediately leave $\TL$ in
the past. One denotes by $\ESC^P$ the set
of all such points
\end{enumerate}
Moreover $\ESC=\ESC^F \cup 
\ESC^P$ is called the boundary escape set and points in $\ESC$ are the
boundary escape points.
\end{definition}
The reader should note that the definition of escape points relies on
\bichars, that is, integral curves of $\Hp$ in $\Char p \subset T^*\hL$, and not on
the notion of \gbichars. The latter curves do remain in $\TL$. 
%%%%%%%%%%%%%%%%%%%%%%%%
% lemma                %
%%%%%%%%%%%%%%%%%%%%%%%%
\begin{lemma}
  \label{lemma: glancing non-escape point}
  The following properties hold.
  \begin{enumerate}
  \item $\Hb^- \subset \ESC^F\setminus \ESC^P$ and $\Hb^+ \subset
    \ESC^P\setminus \ESC^F$.
  \item $\sgGb \subset \ESC^F \cap \ESC^P$.
  \item $\sdGb \cap \ESC = \emptyset$.
  \item $\Gb \setminus \ESC \subset \sdGb\cup \glGb$. 
  \end{enumerate}
\end{lemma}

%%%%%%%%%%%%%%%%%%%%%%%%
% definition           %
%%%%%%%%%%%%%%%%%%%%%%%%
\begin{definition}[boundary geometric control]
\label{control-geo-boundary}
Let $\Gamma$ be an open subset of $\d \M$. One
says that $\Gamma$ controls geometrically the manifold $\M$
if there exists $T>0$ such that any \gbichar encounters a boundary  escape point above  $]0,T[\times
\Gamma$.   One says that $(\Gamma, T)$ fulfills GCC.
In such case, one sets
\begin{align*}
  T_{\text{GCC}}(\Gamma) = \inf \{ T>0; \ (\Gamma, T) \ \text{fulfills GCC}\}.
\end{align*}
\end{definition}

%%%%%%%%%%%%%%%%%
% subsection
\subsection{Invariant measure supports}
For a manifold $\M\in \X^1$ we will consider an extension $\tM$ as in
the begining of Section~\ref{geometry}.
The following result is proven in the companion article~\cite{BDLR2}.
\begin{theorem}
  \label{theo-propagation}
Let $(\M, \k, g) \in \X^1$ and let $\mu$ and $\nu$ be two nonnegative
measure densities on $T^* \hL$ and $T^*\d\L \simeq \pdTL$  respectively  that fulfill the
following properties:  
\begin{enumerate}
\item  $\supp \mu \subset \Char p \cap \TL \setminus 0$. 

\item  One has, in the sense of distributions,  
  \begin{align}
    \label{eq: Gerard-Leichtnam equation}
    \Hp \mu  = - \transp{\Hp} \mu = \int_{\y \in \pHb \cup \pGb} 
    \frac{\delta_{\y^+} - \delta_{\y^-}}
    {\dup{\xi^+- \xi^-}{\nx}_{T_x^*\M, T_x\M} 
    } \ d \nu (\y), 
  \end{align}
  where $\y^\pm$ and $\xi^\pm$ are as given in 
   \eqref{eq: relevement H-G-intro}. Here, $\nx$ stands for  the unitary inward
    pointing normal vector in the sense of the metric. 
\end{enumerate}
Then, the support of the measure  $ \mu$ is a union of maximal \gbichars.
\end{theorem}

With the notation of Definitions \ref{def: generalized bichar-intro} and
\ref{def: generalized bichar 2-intro}, the result of Theorem~\ref{theo-propagation} 
means that if $\y \in \supp \mu$, there exists a maximal 
\gbichar $s\to \gammaG(s), \, s \in \R \setminus B$ such that
$\y \in \GammaG \subset \supp \mu$.

The identification $T^*\d\L \simeq \pdTL$ is explained in Section~\ref{sec: A partition of the cotangent bundle at the boundary-intro}.
\begin{remark}\label{integrand}
  \label{remark: integrand GL equation on G-intro}
  If $\y \in \pGb$ then $\y^-$ and $\y^+$ coincide with $\y$ and $\xi^+= \xi^-$.  The value of
  the integrand in \eqref{eq: Gerard-Leichtnam equation} thus requires
  some explanation in this case.  In fact, first consider
  $\y^0 = (\y^{0\prime}, 0) \in \pHb$ with $\y^{0\prime}=(t^0,x^{0\prime},z=0, \tau^0, \xi^{0\prime})$. Then $\y^{0,\pm} \neq \y^0$ and \eqref{eq: relevement
    H-G-intro} give
  $\xi^{0,+} - \xi^{0,-} = 2 \zeta dz$, yielding
  $\dup{\xi^{0,+} - \xi^{0,-}}{\n_{x^0}}_{T_x^*\M, T_x\M} = 2 \zeta$ since
  $\nx= \d_z$ in the coordinates we consider here. Considering a
  $\Con^1$-test function $q(\y)$ one has 
  \begin{align*}
    \dup{\delta_{\y^{0,+} } - \delta_{\y^{0,-} }}{q}
    &= q \big(\y^{0\prime},\zeta\big)  
      -  q \big(\y^{0\prime},-\zeta\big) .
  \end{align*}
  The integrand is thus 
\begin{align*}
    \frac{q \big(\y^{0\prime},\zeta\big)  
      -  q \big(\y^{0\prime},-\zeta\big) }{2 \zeta}.
\end{align*}
  If now a sequence $(\y^{(n)} )_n\subset \pHb$ converges to $\y \in \pGb$
then 
\begin{align*}
  %\label{eq: understanding GL-glancing}
  \frac{\dup{\delta_{\y^{(n),+}} - \delta_{\y^{(n),-}}}{q}}
  {\dup{\xi^{(n),+}-\xi^{(n),-}}{\nx}_{T_x^*\M, T_x\M}}
  \to \d_{\zeta}  q (\y).
\end{align*}
The integrand in \eqref{eq: Gerard-Leichtnam equation} for $\y\in
 \pGb$ is thus to be understood as the derivative with respect to the
 variable $\zeta$ at $\zeta=0$.
Note that this interpretation is very coordinate dependent. We
give a more geometrical interpretation using more intrinsic
coordinates in the companion article~\cite[Section 5.7]{BDLR2}.
\end{remark}

This result was proven in \cite[Th\'eor\`eme 3]{Burq:1997} in the case
$(\M, \k, g) \in \X^2$. The proof of the result of
Theorem~\ref{theo-propagation} in \cite{BDLR2} is more intricate due
to the lower regularity of the metric $g$ and the function $\k$.

%%%%%%%%%%%%%%%%%
% section
%%%%%%%%%%%%%%%%%
\section{A priori estimates and exact controllability}\label{sec.4}

In this section we consider $(\M, \k, g) \in \Y^1$, that is,  $\M$ is
$W^{2, \infty}$ and both $\k$ and the metric $g$ are  Lipschitz. 

First we recall a classical {\em a priori} estimation for the normal
derivative of a solution to the wave equation. Second, we recall the
equivalence between observablity and exact controllability.

%%%%%%%%%%%%%%%%%
% subsection
\subsection{Normal derivative estimation}
\label{sec: Normal derivative estimation}
Denote by $\n$ the unitary normal inward pointing vector field to $\d\M$ in the sense
of the metric. It has the regularity of the metric, that is Lipschitz
here. For a function $w$ and $x \in \d\M$ then $\d_\n w (x) = \n (w)(x) = d
w (x)(\nx)$. In the quasi-normal geodesic coordinates of
Propositon~\ref{prop: quasi-normal coordinates} that can also be
obtained in the $\Y^1$-regularity setting \cite[Appendix~B]{BDLR2} one has $n=
\d_{x_d}$ and thus $\d_{\n} w = \d_d w$. 
%%%%%%%%%%%%%%%%%%%%%%%%
% proposition          %
%%%%%%%%%%%%%%%%%%%%%%%%
\begin{proposition}
  \label{prop: admissibility}
  Assume that $(\M, \k, g) \in \Y^1$.  For any $T>0$ there exists $C>0$
  depending only on $T$, $\M$, $\Norm{\k}{W^{1,\infty}(\M)}$,
  $\Norm{g}{W^{1,\infty}\mathcal T^0_2(\M)}$ such that for any $(\udl{u}^0,
  \udl{u}^1) \in H^1_0(\M)\times L^2 (\M)$ and $f \in L^2_{\loc}(\L)$, if $u$ is the solution to the
  wave equation \eqref{eq: wave equation-intro}
 then 
 \begin{align*}
   %\label{eq: admissibility boundary obs}
   &\Norm{\d_\n u}{L^2(]0,T[\times \d\M)}^2 \\
   &\qquad \leq C \big(
      \int_{-1}^{T+1} \E_{\k, g} (u)(t) \, d t
      + \Norm{\nablag u}{{L^2(-1, T+1; L^2 V(\M))}}
      \Norm{f}{L^2(]-1, T+1[\times\M)}\big).
 \end{align*}
 \end{proposition}
 Below we will use the Neumann trace as an observation operator for
 the wave equation. In this context, with $f=0$, Proposition~\ref{prop:
   admissibility} provides a
 so-called admissibility result; see for instance \cite{TW:09}.

 Note that a more usual and natural form of the estimation is simply
 \begin{equation*}
   %\label{eq: admissibility boundary obs}
   \Norm{\d_\n u}{L^2(]0,T[\times \d\M)}^2 \lesssim
      \int_{-1}^{T+1} \E_{\k, g} (u)(t) \, d t
      +  \Norm{f}{L^2(]-1, T+1[\times\M)}^2.
 \end{equation*}
 This form is however not sufficient in one argument we use in what
 follows; we refer to the use of Proposition~\ref{prop: admissibility} made below
 \eqref{eq: mass2-proof 0}.

 Note that since $u$ vanishes on $\d\M$ one has $\Norm{\nablag
   u_{|\d\M} }{} = |\d_\n u |$. The result of the previous proposition thus can be transferred to $\Norm{\nablag u_{|\R\times \d\M}}{L^2(]0,T[\times \d\M)}$.
%%%%%%%%%%%%%%%%%%%%%%%%
% corollary            %
%%%%%%%%%%%%%%%%%%%%%%%%
%\begin{corollary} 
%  \label{coro-direct}
% Assume that $(\M, \k, g) \in \Y^1$.  For any $T>0$ there exists $C$
%  depending only on $T$, $\M$, $\Norm{\k}{W^{1,\infty}(\M)}$,
%  $\Norm{g}{W^{1,\infty}\mathcal T^0_2(\M)}$ such that for any $(\udl{u}^0,
%  \udl{u}^1) \in H^1_0(\M)\times L^2 (\M)$, if $u$ is the solution to the
%  wave equation \eqref{eq: wave equation-observation-intro}
% then 
% \begin{equation}
%   \label{eq: admissibility boundary obs-bis}
%   \Norm{\nablag u_{|\R\times \d\M}}{L^2(]0,T[\times \d\M)} 
%   \leq C  \E_{\k,g}(u).
% \end{equation}
% \end{corollary}

 The proof of Proposition~\ref{prop: admissibility} follows from an
 examination of the standard proof and a carefull handling of the low
 regularity metric. We will also need to approach the weak solution
 to the wave equation~\ref{eq: wave equation-intro} by a
 sequence of strong solutions.
%%%%%%%%%%%%%%%%%%%%%%%%
% proposition          %
%%%%%%%%%%%%%%%%%%%%%%%
\begin{proposition}
  \label{prop:  strong solution wave equation}
 Suppose that $(\udl{u}^0, \udl{u}^1) \in \big(H^2(\M) \cap H_0^1(\M)\big)
 \times H^1_0(\M)$ and that $f \in L^1_{\loc}(H_0^1(\R; \M))$.  There
 exists a unique
  \begin{align*}
    u \in \Con^0 \big(\R; H^2(\M) \cap H_0^1(\M) \big)
    \cap \Con^1 \big(\R; H^1_0(\M)\big) \cap \Con^2 \big(\R; L^2(\M)\big)
  \end{align*}
  that is a strong solution of \eqref{eq: wave equation-intro}
  meaning  that $(u, \d_t u)\bt = (\udl{u}^0,\udl{u}^1)$ and $P_{\k, g} u =f$
  holds in $L^1_{\loc}\big(\R; L^2(\M)\big)$. 
\end{proposition}
Note that a strong solution is also a weak solution.
Then, if $(\udl{u}^0, \udl{u}^1) \in H_0^1(\M) \times L^2(\M)$, $f \in L^1_{\loc}(\R; L^2(\M))$ and $u$ is
  the weak solution to \eqref{eq: wave equation-intro} given by
  Propostion~\ref{prop: well-posedness nonhomogeneous wave equation}
   and  if $(u_n^0, u_n^1)_n \subset  \big(H^2(\M) \cap H_0^1(\M)\big)
  \times H^1_0(\M)$, $(f_n)_n\subset L^1_{\loc}(H_0^1(\R; \M))$, with $(u_n)_n$  the sequence of
  associated strong solutions,  are such that $(u_n^0, u_n^1) \to (\udl{u}^0, \udl{u}^1)$ in
  $H_0^1(\M) \oplus L^2(\M)$, and $f_n \to f$ in $L^1_{\loc}(\R; L^2(\M))$ then $u_n \to u$ in $\Con^0 \big(\R; H_0^1(\M)  \big)
    \cap \Con^1 \big(\R; L^2(\M)\big)$ from the
continuity of the map \eqref{eq: continuity wave initial conditions}.

\begin{proof}[\bfseries Proof of Proposition~\ref{prop: admissibility}]

  First we consider the case of a strong solution
  \begin{align*}
    u \in \Con^0
  \big(\R; H^2(\M) \cap H_0^1(\M) \big) \cap \Con^1 \big(\R;
  H^1_0(\M)\big) \cap \Con^2 \big(\R; L^2(\M)\big),
  \end{align*}
  with $f \in
  L^2_{\loc}(\L)$.

    Consider a {\em Lipschitz} vector field $X$ that coincides with
    $\n$ on the boundary. We view $X$ as a first-order
    differential operator.  For $\chi \in \Cinfc (-1, T+1)$, nonnegative and  equal to
    $1$ on $]0,T[$, one finds that $[P,  \chi(t) X] u
    \in \Con^0 \big(\R; L^2(\M)\big)$. Set
    $I = \inp{[P,  \chi(t) X] u}{u}_{L^2(\L)}$.
    With the Green formula, that is, two integrations by parts, one finds
    \begin{align}
      \label{eq: estimation dn u-1}
      \Norm{\d_\n u}{L^2(]0,T[\times \d\M)}^2
      &\leq \inp{\chi(t) X u_{|\d\M}}{\d_\n u}_{L^2(\d\L)} \\
      &\leq I - \inp{\chi(t) X u }{f}_{L^2(\L)}
      + \inp{f}{ \chi(t)  X^* u}_{L^2(\L)} \notag\\
      &\lesssim I
      + \Norm{\nablag u}{L^2(-1, T+1; L^2 V(\M))}
      \Norm{f}{L^2(]-1, T+1[\times \M)}.\notag
\end{align}
 %     &\lesssim I
  %    + \E_{\k, g} (u) (0)^{1/2}
  %    \Norm{f}{L^1(-1, T+1; L^2(\M))}
  %    + \Norm{f}{L^1(-1, T+1; L^2(\M))}^2, \notag
  %  \end{align}
  %  using \eqref{eq: evol energy}.
 Writing  $[P,  \chi(t) X]=
    [\d_t^2 ,\chi(t)] X - \chi(t) [A, X]$. One has $I = J -K$ with
    \begin{align*}
      &J = \int_\R\dup{ [\d_t^2 ,\chi(t)] X  u}
      {\bar{u}}_{H^{-1}(\M), H^1_0(\M)} d t,\\
      &K = \int_\R\dup{\chi(t) [A, X] u}{\bar{u}}_{H^{-1}(\M), H^1_0(\M)} d t .
    \end{align*}
    Since $[ \d_t ^2 , \chi]$ is a first-order operator and compactly
    supported (in time)  we can integrate by parts in the time variable
    and obtain the bound, using the Poincaré inequality,
    \begin{align}
      \label{eq: estimation dn u-2}
      |J| &\lesssim  \int_{-1}^{T+1} \Bigl(\Norm{\d_t u}{L^2(\M)}+ \Norm {u} {L^2(\M)}\Bigr)
      \Norm{\nablag u}{L^2 V(\M)}\,  d t 
      \lesssim
      \int_{-1}^{T+1} \E_{\k,g}(u)(t) \, d t.
    \end{align}
    To estimate $|K|$ we use
    a partition of unity subordinated to an atlas on $\M$ and we
    consider the commutator $[A,  X] $ in local coordinates. Recall
    that $A$ takes the form
    \begin{align*}
      A  = \tk^{-1} \d_{x_i} \circ \tk g^{i j} \d_{x_j} 
    \end{align*}
    with $\tk = \k  (\det g)^{1/2}$ where we use Eistein's summation
    convention.
    One thus finds
    \begin{align}
      \label{eq: estimation dn u-3}
      [A,  X] &= \tk^{-1} \d_{x_i} \circ \tk g^{i j} [\d_{x_j},X]
      +\tk^{-1} \d_{x_i} \circ [\tk g^{i j},X]  \circ \d_{x_j}\\
      &\quad + \tk^{-1} \circ [\d_{x_i},X] \circ \tk g^{i j} \d_{x_j}
      + [\tk^{-1},X] \circ \d_{x_i} \circ \tk g^{i j} \d_{x_j} .\notag
    \end{align}
    Write $K = K_1 + \cdots + K_4$ in association with the four terms 
    in \eqref{eq: estimation dn u-3}.
    Since $X$ has Lipschitz coefficients then $ [\d_{x_k},X]$ is a
    vector field with bounded coefficients and $[\tk g^{i
      j},X]$ is a bounded function  in the local
    coordinates. Thus, with an integration by parts in space the
    contribution $|K_1|$, $|K_2|$, and  $|K_3|$ can by estimated by
    \begin{align}
      \label{eq: estimation dn u-4}
      |K_1| +|K_2|+|K_3|\lesssim\int_{\R}\chi(t) \Norm{\nablag u}{L^2 V(\M)}^2\, d t.
    \end{align}
    For the term $K_4$ since  $[\tk^{-1},X] $ is only bounded, an
    integration by parts in space is not possible. 
    Instead, exploiting that $u$ is a solution to the homogenous wave equation one writes
    \begin{align*}
      [\tk^{-1},X] \circ \d_{x_i} \circ \tk g^{i j} \d_{x_j} u
      =  [\tk^{-1},X]\tk A u =  [\tk^{-1},X]\tk \d_t^2 u.
    \end{align*}
    This now allows one to perform an integration by parts with
    respect to the time variable  yielding
    \begin{align}
      \label{eq: estimation dn u-5}
      |K_4| &\lesssim \int_{-1}^{T+1} \Norm{\d_t u}{L^2(\M)}^2 d t
      +
      \int_{-1}^{T+1}\Norm{u}{L^2(\M)}
      \Norm{\d_t u}{L^2(\M)}  d t \\
      &\quad
      + \Norm{u}{L^2(]-1,T+1[\times \M)}
      \Norm{f}{L^2(]-1,T+1[\times \M)}\notag\\
      &\lesssim
      \int_{-1}^{T+1} \E_{\k, g} (u)(t) \, d t
      + \Norm{\nablag u}{L^2(-1, T+1; L^2 V(\M))}
      \Norm{f}{L^2(]-1,T+1[\times \M)}.
      \notag
    \end{align}
    Combining \eqref{eq: estimation dn u-1}, \eqref{eq: estimation dn
      u-2}, \eqref{eq: estimation dn u-4} and \eqref{eq: estimation dn
      u-5} gives
    \begin{align}
      \label{eq: estimation dn u-6}
      \Norm{\d_\n u}{L^2(]0,T[\times \d\M)}^2 &\lesssim
      \int_{-1}^{T+1} \E_{\k, g} (u)(t) \, d t\\
      &\quad + \Norm{\nablag u}{L^2(-1, T+1; L^2 V(\M))}
      \Norm{f}{L^2(]-1,T+1[\times \M)}.\notag
    \end{align}

    If $u$ is now a weak solution, if one approches $u$ by a sequence
    of strong solutions as described below Proposition~\ref{prop:  strong solution wave equation} one
    finds that the normal trace $\d_\n u_{|\d\M}$ makes sense in $L^2(]0,T[\times \d\M)$
    and \eqref{eq: estimation dn u-6} remains true for the weak solution.
\end{proof}

%%%%%%%%%%%%%%%%%
% subsection
\subsection{Exact controllability notions}

We make the classical connexion between the observability properties
of the {\em homogeneous} wave equation given in \eqref{eq: wave equation-observation} and the exact controllability of the
wave equations. We will consider two different wave equations here,
one with an interior source term and a homogeneous Dirichlet boundary
condition and one with a boundary source term through the Dirichlet
boundary condition. In each case we describe what is meant by exact
controllability.
%%%%%%%%%%%%%%%%%
% subsection
\subsubsection{Exact interior controllability}
%Consider the nonhomogeneous wave equation 
%\begin{align}
%  \label{eq: wave equation-interior}
%  \begin{cases}
%     P_{\k, g}  \, y =f 
%     & \dans\ \R\times\M,\\
%     y=0 &  \text{on}\ \R \times\d\M,\\
%      y_{|t=0}  = \udl{y}^0, \ \d_t y_{|t=0}  = \udl{y}^1 
%      & \dans\ \M,
%  \end{cases}
%\end{align}
%with solutions as given by Proposition~\ref.
%Standard results show that it is well-posed.
%%%%%%%%%%%%%%%%%%%%%%%%%
% proposition          %
%%%%%%%%%%%%%%%%%%%%%%%
%\begin{proposition}
%  \label{prop: well-posedness wave equation-interior}
%  Consider $\k$ and $g$ both Lipschitz. 
%  Let $(\udl{y}^0, \udl{y}^1) \in H_0^1(\M) \times L^2(\M)$ and let $f \in
%  L_{\loc}^1\big(\R; L^2(\M)\big)$.   There exists a unique  
%  \begin{align*}
%    y \in \Con^0 \big(\R; H_0^1(\M)  \big)
%    \cap \Con^1 \big(\R; L^2(\M)\big)
%  \end{align*}
%  that is a weak solution of \eqref{eq: wave equation-interior}.
%\end{proposition}

Suppose $\omega$ is an open subset of $\M$.  The notion of exact
interior controllability for the wave equation on $\M$ from $\omega$
in time $T$ is stated as follows.
%%%%%%%%%%%%%%%%%%%%%%%%
% definition           %
%%%%%%%%%%%%%%%%%%%%%%%%
\begin{definition}[exact interior controllability in $H_0^1(\M) \oplus L^2(\M)$]
  \label{def: exact controllability}
  One says that the wave equation is exactly controllable from
  $\omega$ in time $T>0$ if
  for any
  $(\udl{y}^0, \udl{y}^1) \in H_0^1(\M) \times L^2(\M)$, there exists  
  $f\in L^2(]0,T[\times \M)$ such that the weak solution $y$ to
  \begin{equation*}
    %\label{eq.controlint}
    P_{\k,g} y =  \bld{1}_{]0,T[ \times \omega} \, f,
    \quad y_{|\R \times \d\M} =0, \quad (y, \d_t y)\bt = (\udl{y}^0, \udl{y}^1),
  \end{equation*}
  as given by Proposition~\ref{prop: well-posedness nonhomogeneous
    wave equation} satisfies $(y, \d_t y)_{|t= T} = (0,0)$.  The
  function $f$ is called the control function or simply the control.
\end{definition}

%%%%%%%%%%%%%%%%%
% subsection
\subsubsection{Exact boundary controllability}
Consider the nonhomogeneous wave equation with source term given by a
Dirichlet boundary condition. 
\begin{align}
  \label{eq: wave equation-boundary-intro}
  \begin{cases}
     P_{\k, g}  \, y =0
     & \dans\ \R\times\M,\\
     y=f_\d &  \text{on}\ \R \times\d\M,\\
      y_{|t=0}  = \udl{y}^0, \ \d_t y_{|t=0}  = \udl{y}^1 
      & \dans\ \M,
  \end{cases}
\end{align}
Standard results show that it is well-posed.
%%%%%%%%%%%%%%%%%%%%%%%%
% proposition          %
%%%%%%%%%%%%%%%%%%%%%%%
\begin{proposition}
  \label{prop: well-posedness wave equation-boundary}
  Consider $\k$ and $g$ both Lipschitz. 
  Let $(\udl{y}^0, \udl{y}^1) \in L^2(\M) \times H^{-1}(\M)$ and $f_\d \in L_{\loc}^2(\R\times \d\M)$.   There exists a unique  
  \begin{align*}
    y \in \Con^0 \big(\R; L^2(\M)  \big)
    \cap \Con^1 \big(\R; H^{-1}(\M)\big).
  \end{align*}
  that is a weak solution of \eqref{eq: wave equation-boundary-intro}.
\end{proposition}

Let $\Gamma$ be a nonempty open subset of $\d\M$ and $T>0$.  The
notion of exact boundary 
controllability for the wave equation from $\Gamma$ in time
$T$ is stated as follows.
%%%%%%%%%%%%%%%%%%%%%%%%
% definition           %
%%%%%%%%%%%%%%%%%%%%%%%%
\begin{definition}[exact boundary controllability in $L^2(\M) \oplus H^{-1}(\M)$]
  \label{def: exact controllability-boundary}
  One says that the wave equation is exactly controllable from
  $\Gamma$ in time $T>0$ if
  for any
  $(\udl{y}^0, \udl{y}^1) \in L^2(\M) \times H^{-1}(\M)$, there exists  
  $f_\d \in L^2(]0,T[\times \d\M)$ such that the weak solution $y$ to
  \begin{equation*}
    %\label{eq.controlboundary}
    P_{\k,g} y = 0,
    \quad y_{|\R \times \d\M} =\bld{1}_{]0,T[ \times \Gamma} f_\d,
    \quad (y, \d_t y)\bt = (\udl{y}^0, \udl{y}^1),
  \end{equation*}
  as given by Proposition~\ref{prop: well-posedness wave equation-boundary}
  satisfies $(y, \d_t y)_{|t= T} = (0,0)$.
 The function $f_\d$ is called the control function or simply the
 control. 
\end{definition}
%%%%%%%%%%%%%%%%%
% subsection
\subsection{Exact controllability equivalent to observability and corollaries}

The following proposition is standard and states that in the two
cases we consider exact controllability is equivalent to an
obserbability inequality.
%%%%%%%%%%%%%%%%%%%%%%%%
% proposition          %
%%%%%%%%%%%%%%%%%%%%%%%%
\begin{proposition}
  \label{prop: equiv controllability observability}
  \begin{enumerate}
  \item Let $\omega$ be an open subset of $\M$ and $T>0$.
  The wave equation is exactly controllable from
  $\omega$ in time $T$ if and only if the homogeneous wave equation is observable from
  $\omega$ in time $T$.
  \item    Let $\Gamma$ be a nonempty open subset of $\d\M$ and $T>0$. 
     The wave equation is exactly controllable from $\Gamma$ in time $T$  if and only if the homogeneous wave equation is observable from
  $\Gamma$ in time $T$.
    \end{enumerate}
\end{proposition}
With the previous proposition and Theorem~\ref{theo-perturbation} one
deduces the following corollary.
%%%%%%%%%%%%%%%%%%%%%%%%
% corollary            %
%%%%%%%%%%%%%%%%%%%%%%%%
\begin{corollary}[Exact controllability result]
  \label{cor: theo-perturbation}
  Let $(\M, \k, g) \in \X^1$.

   {\em \bfseries Interior exact controllability.}
  Let $\omega$ be an open subset of $\M$ that satisfies the interior geometric control condition associated with the infimum time $T_{GCC}(\omega)$.
  Let $T>T_{GCC}(\omega)$. Then, there exists $\eps>0$ such that if  $(\tilde{\M}, \tilde{\k}, \tilde{ g},
  \tilde{\omega})$ is $\eps$-close to $(\M, \k, g,\omega)$ in the $\Y^1$-topology in the
  sense of Defintion~\ref{def: eps-close} for $(\tilde{\M}, \tilde{\k}, \tilde{ g}) \in \Y^1$
  and $\tilde{\omega}$ an open subset of $\tilde{\M}$, then the wave equation
  associated with $P_{\tilde{\k}, \tilde{ g}}$ on $\tilde{\M}$  is
  exactly controllable 
  from $\tilde{\omega}$ in time $T$.

  {\em \bfseries  Boundary  exact controllability.}
  Let $\Gamma$ be an open subset of $\d \M$ such that
   $\Gamma$ satisfies the boundary geometric control condition
  associated with the infimum time $T_{GCC}(\Gamma)$.
Let $T>T_{GCC}(\Gamma)$. 
   Then, there exists $\eps>0$ such that if  $(\tilde{\M}, \tilde{\k}, \tilde{ g},
  \tilde{\Gamma})$ is $\eps$-close to $(\M, \k, g,\Gamma)$ in the $\Y^1$-topology in the
  sense of Defintion~\ref{def: eps-close} for $(\tilde{\M}, \tilde{\k}, \tilde{ g}) \in \Y^1$
  and $\tilde{\Gamma}$ an open subset of $\d \tilde{\M}$, then the wave equation
  associated with $P_{\tilde{\k}, \tilde{ g}}$ on $\tilde{\M}$  is exactly controllable 
  from $\tilde{\Gamma}$ in time $T$.
\end{corollary}

%%%%%%%%%%%%%%%%%
% section
%%%%%%%%%%%%%%%%%
\section{Semi-classical reduction}
\label{Semi-classical reduction}
In \cite{BDLR0}, on a compact manifold without boundary, we proved an
interior observability estimate for the Klein-Gordon equation by means
of microlocal defect measures. The more general case we consider here,  in the presence of a boundary is technically more involved and requires a semi-classical approach. 
We recall in this section how observability  estimates as in
Definitions~\ref{def: observability interior-intro} and \ref{def: observability boundary-intro}
can be obtained  from counterpart semi-classical observability estimates.

%%%%%%%%%%%%%%%%%
% subsection
\subsection{Dyadic decomposition}
\label{dyadic}

Consider $(\k,g) \in \X^1(\M)$ and the associated operator  $A = A_{\k,g}$ with Dirichlet
boundary conditions. Denote by $\lambda_\nu$ the  nondecreasing sequence
of positive eigenvalues of $-A$ that goes to $+\infty$ and consider
$(e_\nu)_\nu$ a Hilbert basis of $L^2(\M) = L^2(\M, \k \mug)$
of associated {\em real}  eigenfunctions. 
  
  Let $ 0<\alpha<1$,  $ \varrho\in ]1, 1/\alpha[$ and set  $\hk =
  \varrho^{-|k|}$ and 
\begin{align*}
  J_k= \{\nu \in \N;\ \alpha \leq \hk {\sqrt {\lambda _\nu}}<
  \alpha^{-1} \}
  = \{\nu \in \N;\ \alpha \varrho^{|k|}\leq {\sqrt {\lambda _\nu}}<
  \varrho^{|k|} /\alpha\},
\end{align*}
 for  $ k \in \ZZ^*$. Introduce 
\begin{align*}
  E_k=\Span \{  e_\nu; \ \nu\in J_k \},
\end{align*}
equipped with the $L^2$-norm  $\Norm{u}{L^2(\M)}^2 = \Norm{u}{L^2(\M, \k \mug)}^2 = \sum_{\nu\in J_k} | u_\nu
|^2$ for $u = \sum_{\nu \in J_k} u_\nu e_\nu \in E_k$.
Observe that if $u \in E_k$ then $A^n u \in E_k$. Hence, $E_k$ is a
subspace of all the iterated domains of $A$.

At this stage it is important to  note that $J_{-k} = J_k$ implying $E_{-k} = E_k$. 
However, we will identify $u \in E_k$ 
with the  following solution of the wave equation
\begin{align*}
 u=\sum _{\nu \in J_k} 
  e^{\sgn(k) i t{\sqrt {\lambda _\nu}}} u_\nu e_\nu.
\end{align*}
The sign of $k$ here becomes important. Yet, note that $u \in E_k$ if
and only if $\bar{u} \in E_{-k}$, through this identification since
the eigenfunctions $e_\nu$ are chosen real.

Following up, we identify $\d_t^\ell u$ with $u = \sum_{\nu \in J_k}
(i \sgn(k))^\ell
\lambda_\nu^{\ell/2}  u_\nu e_\nu \in E_k$, its value at $t=0$.
Similarly, one identifies $A^s u$  with $ \sum_{\nu \in J_k} 
\lambda_\nu^{s}  u_\nu e_\nu \in E_k$.
%%%%%%%%%%%%%%%%%%%%%%%%
% lemma                %
%%%%%%%%%%%%%%%%%%%%%%%%
\begin{lemma}
  \label{lemma: norm dt uk}
  For $u \in E_k$, the norms 
  \begin{align*} 
    \Norm{\hk \nablag u}{L^2 V(\M)} 
    \ \ \et \ \  
    \Norm{\hk \d_t u}{L^2(\M)}  
  \end{align*}
are both equivalent to $\Norm{u}{L^2(\M)}$, uniformly with
  respect to $k\in \ZZ$.
  More generally, for $\ell \in \N$ and $s\in\R$, the norm
  $\hk^{\ell+2s} \Norm{\d_t^\ell  A^s u}{L^2 (\M)}$ is equivalent to
  $\Norm{u}{L^2(\M)}$ for $u \in E_k$,  uniformly with
  respect to $k\in \ZZ^*$.
\end{lemma}
%% proof of lemma
\begin{proof}
  For the first result one writes
  \begin{align*} 
\Norm{\hk \nablag u}{L^2 V(\M)}^2 = 
\inp{\hk^2 A u }{u}_{L^2(\M)}
=  \Norm{\hk \d_t u}{L^2(\M)}^2.
\end{align*}
Then one concludes using $\hk\lambda_\nu^{1/2} \eqsim 1$ for $\nu \in J_{k}$.
\end{proof}
As a consequence the $L^2$-norm and the square root of the
semi-classical energy $\E^h(u)$
\begin{align}
  \label{eq: equiv H1sc scenergu}
  \E^h(u)=  \frac12 \big( \Norm{\hk\nablag u}{L^2 V(\M)}^2
  + \Norm{\hk\d_t  u}{(\M)}^2\big) = \hk^2 \E(u),
  %\eqsim \Norm{u}{H^1_{sc}}^2,
\end{align}
are equivalent on $E_k$. Note that for $u \in E_k$ both terms in the
semi-classical energy coincide; this is not the case in general for a
solution of the wave equation.

\bigskip
We introduce the
following sets of sequences of functions 
\begin{align*} 
%\label{eq: space B}
&B= \big\{ (u^k)_{k \in {\ZZ}^*};\ u^k \in
E_k  \ \et \  \Norm{u^k}{L^2(\M)}  \leq 1 \big\}, \\
&B^\pm= \big\{ (u^k)_{k \in \pm{\N}^*};\  u^k \in
E_k \ \et \  \Norm{u^k}{L^2(\M)}  \leq 1 \big\}.\notag
\end{align*}

%%%%%%%%%%%%%%%%%
% subsection
\subsection{Semi-classical observation}

The result of Proposition~\ref{3.inegalitesemiclassbord} below for
boundary observation is proven in \cite[Section 4]{Burq:1997} following a
strategy of G.~Lebeau~\cite{Lebeau:92}. The result of
Proposition~\ref{3.inegalitesemiclassint} for interior observation can
be proven similarly. In  \cite{Burq:1997} domains and metrics are
smoother, yet lowering the regularity does not affect the proof that
is only based on semi-classical analysis arguments with respect to the time
 variable. In fact, the proofs of both propositions can be carried out
 within an abstract framework that can be found in \cite{BDLR-abs-arg}.
%%%%%%%%%%%%%%%%%%%%%%%%
% proposition          %
%%%%%%%%%%%%%%%%%%%%%%%%
\begin{proposition}[interior semi-classical observation implies
  classical observation]
\label{3.inegalitesemiclassint}
 Let $\omega$ be a nonempty open subset of $\M$. 
Assume that there exists $C>0$, $ k_0>0$, and $\delta >0$, such that for
any $U= (u^k)_{k\in \N}\in B^+$ and  any $k\geq k_0$ one has
\begin{equation}
\label{3.ineg11int}
\Norm{u^k\bt}{L^2(\M)}
\leq C \Norm{\bld{1}_{I \times \omega}   \hk \d_t u^k}{L^2(\R\times\M)}, \qquad \avec \ I = ]\delta, T-\delta[.
\end{equation} 
Then,  the homogeneous wave equation is observable from
  $\omega$ in time $T>0$ in the sense of Definition~\ref{def: observability interior-intro}.
\end{proposition}
%%%%%%%%%%%%%%%%%%%%%%%%
% proposition          %
%%%%%%%%%%%%%%%%%%%%%%%%
\begin{proposition}[boundary semi-classical observation implies classical observation]
\label{3.inegalitesemiclassbord}
 Let $\Gamma$ be a nonempty open subset of $\d\M$. 
Assume that there exists $C>0$, $ k_0>0$, and $\delta >0$ such that
for any $U= (u^k)_{k\in \N}\in B^+$ and  any $k\geq k_0$, one has 
\begin{equation}
  \label{3.ineg11bord}
  \Norm{u^k\bt}{L^2(\M)}^2 
  \leq C \Norm{\bld{1}_{I \times \Gamma} 
  \hk \d_\n  u^k}{L^2(\R\times \d\M)}, \qquad \avec \ I = ]\delta, T-\delta[.
\end{equation} 
Then, the homogeneous wave equation is observable from
  $\Gamma$ in time $T>0$ (in the sense of Definiton~\ref{def: observability boundary-intro}).
\end{proposition}

Propositions~\ref{3.inegalitesemiclassint}
and~\ref{3.inegalitesemiclassbord} state that if an observability
inequality as in  Definition~\ref{def: observability interior-intro} or Definiton~\ref{def: observability boundary-intro}
holds
in $E_k$ uniformly for large $|k|$, then it also holds for any initial
data with possibly a small  loss (here $\delta$ on each side) in the
time interval  required for observation.

\medskip
The proof presented in \cite{BDLR-abs-arg} is based on
several properties of the observation operator $\Obs$, here  $\Obs =
\bld{1}_{I \times \omega}   \hk \d_t$ in the first case and $\Obs = \bld{1}_{I \times \Gamma} 
  \hk \d_\n$ in the second case:
\begin{enumerate}
\item a unique continuation property, $\Obs u=0$ implying $u=0$ for
  eigenfunctions of the operator $A_{\k,g}$; this
  condition holds in the both cases we consider;  see for instance \cite[Theorem 2.4]{Hoermander:83}
and \cite[Theorems 5.11 and
5.13]{LRLR:V1}.
\item an optional admissibility condition, here given by
  Proposition~\ref{prop: admissibility} in the case $\Obs = \bld{1}_{I \times \Gamma} 
  \hk \d_\n$. In the first case, $\Obs =
\bld{1}_{I \times \omega}   \hk \d_t$ the admissibility condition is trivial.
\end{enumerate}

%%%%%%%%%%%%%%%%%%%%%%%%
% remark               %
%%%%%%%%%%%%%%%%%%%%%%%%
\begin{remark}
\label{rem: equiv norm energy}
In \eqref{eq: equiv H1sc scenergu}, we pointed out that the
$L^2$-norm $\Norm{.}{L^2(\M)}$ is equivalent to the  square root of the semi-classical energy $\E^h(u)$, uniformly in
$k$. Here, $\E^h(u)$ is constant \wrt time $t$, since $u_k$ is solution of the
homogeneous wave equation.  Consequently, one can also replace the
\lhs in \eqref{3.ineg11bord} and \eqref{3.ineg11int} by
\begin{equation*}
  \Norm{u^k}{L^\infty(\R; L^2(\M))}^2
  \ \ \ou \ \ 
\Norm{u^k}{L^2(\interval \times\M)}^2, 
\end{equation*}
for any finite interval $\interval \subset \R$.
\end{remark}

%%%%%%%%%%%%%%%%%
% section
%%%%%%%%%%%%%%%%%
\section{Semi-classical operators and measures}
\label{sec: Semi-classical operators and measures}
%%%%%%%%%%%%%%%%%
% subsection
\subsection{The Schur lemma}
Here, we recall a result that is important  in our
analysis of some semi-classical operators on $\R^d$ in what follows. 
%%%%%%%%%%%%%%%%%%%%%%%%
% lemma                %
%%%%%%%%%%%%%%%%%%%%%%%%
\begin{lemma}[Schur's Lemma]
  \label{lemma: Schur lemma}
  Let $K(x,y)$ be a measurable function on $\R^d\times \R^d$ such that $K(x,.)$ and
  $K(.,y)$ are $L^1$-functions for almost all $x$ and $y$ in $\R^d$
  respectively, with moreover 
  \begin{align*}
    \esssup_{x\in \R^d} \Norm{K(x,.)}{L^1(\R^d)} \leq A
    \ \ \et \ \
    \esssup_{y\in \R^d} \Norm{K(.,y)}{L^1(\R^d)} \leq B,
  \end{align*}
  for some $A\geq 0$ and $B\geq 0$. 
  Then, the operator $\K$ with Schwartz kernel $K(.,.)$ extends as a
  continuous operator on  $L^2(\R^d)$ with 
  $\Norm{\K}{\L(L^2(\R^d))} \leq (AB)^{1/2}$.
\end{lemma}
Assume that the kernel of the operator $\K$ is of the form
\begin{align*}
    K(x,y)= h^{-d} \, k \big(x,\frac{x-y}{h}\big),
\end{align*}
for some measurable function $k$ defined on $\R^d \times \R^d$.
Changes of variables give
\begin{align*}
  %\label{eq: L1 norm kernel}
  &\Norm{K(x,.)}{L^1(\R^d)}
  = \Norm{k(x,.)}{L^1(\R^d)},\\
  &\Norm{K(.,y)}{L^1(\R^d)}
  = \Norm{k(y + h\, .\, ,.)}{L^1(\R^d)}.
  \notag
\end{align*}
The Schur Lemma can be translated accordingly.
%%%%%%%%%%%%%%%%%%%%%%%%
% lemma                %
%%%%%%%%%%%%%%%%%%%%%%%%
\begin{lemma}
  \label{lemma: Schur modified}
 Let the operator $\K$ have Schwartz kernel $K(x,y) = h^{-d} \, k
 \big(x,\frac{x-y}{h}\big)$ with the function $k$ satisfying
 \begin{align*}
    \esssup_{x\in \R^d} \Norm{k(x,.)}{L^1(\R^d)} \leq A
    \ \ \et \ \
    \esssup_{y\in \R^d} \Norm{k(y + h\, . \, ,.)}{L^1(\R^d)} \leq B,
  \end{align*}
  for some $A\geq 0$ and $B\geq 0$. 
  Then, the operator $\K$ with Schwartz kernel $K(.,.)$ extends as a
  continuous operator on  $L^2(\R^d)$ with 
  $\Norm{\K}{\L(L^2(\R^d))} \leq (AB)^{1/2}$.
\end{lemma}
%%%%%%%%%%%%%%%%%%%%%%%%
% Corallary              %
%%%%%%%%%%%%%%%%%%%%%%%%
\begin{corollary}
  \label{cor: Schur modified}
 Let the operator $\K$ have Schwartz kernel $K(x,y) = h^{-d} \, k
 \big(x,\frac{x-y}{h}\big)$ with the function $k$ satisfying
 \begin{align*}
   |k(x,v)|\leq L_0 \jp{v}^{-d-\delta},
   \qquad x \in \R^d, \ v \in\R^d,
 \end{align*}
 for some $\delta>0$ and $L_0>0$. Then, $\K$ extends as a continuous operator on  $L^2(\R^d)$ with
 $\Norm{\K}{\L(L^2(\R^d))} \leq C_{d,\delta} L_0$
 for some $C_{d,\delta}>0$. 
\end{corollary}

%%%%%%%%%%%%%%%%%
% subsection
\subsection{Semi-classical operators on $\R^d$}
\label{sec: Semi-classical operators on Rd}
We recall and develop here some aspects of semi-classical pseudo-differential
operators associated with symbols with fairly low regularity.

Let $h_0>0$. In the semi-classical setting we denote by  $ h\in (0,
h_0]$ a small parameter. 
%%%%%%%%%%%%%%%%%%%%%%%%
% definition           %
%%%%%%%%%%%%%%%%%%%%%%%%
\begin{definition}[symbols]
  \label{def: symbols}
  Let $m,n \in \N \cup \{+\infty\}$, with $n \geq d+1$, and $N \in \R^+$.  Denote by $\Symbol{m,n}{-N}{2d}$ the space of
  all functions $a(x,\xi)$, $x\in \R^d$, $\xi\in \R^d$, such that
  $\d_x^\alpha \d_\xi^\beta a \in L^1_{\loc}(\R^{2d})$ for
  $\alpha, \beta \in \N^d$ with $|\alpha|\leq m$, $|\beta| \leq n$,
  and
  \begin{equation}
  \label{eq: est symbols}
  M_{m,n}^{-N}(a):=\max_{{|\alpha|\leq m}\atop{|\beta| \leq n}}  \esssup_{(x,\xi)}
  \bignorm{\d_x^\alpha \d_\xi ^\beta a(x,\xi)}{} \jp{\xi}^{N} < \infty.
\end{equation}
  % \begin{align}
  %   \bignorm{\d_x^\alpha \d_\xi ^\beta a(x,\xi)}{} 
  %   \leq  C_{\alpha, \beta} \jp{\xi}^{-N}, \qquad x\in \R^d, \ \xi \in \R^d.
  % \end{align}
In addition, one sets  $\Symbolo{m,n}{-N}{2d}$, $n \geq d+1$,  as the set of
all symbols $a \in \Symbol{m,n}{-N}{2d}$ with moreover
$\d_x^\alpha \d_\xi^\beta a \in \Con_0(\R^{2d})$ for $\alpha, \beta \in \N^d$
with
$|\alpha|\leq m$ and $|\beta| \leq n -1-d$.
\end{definition}
Recall that $\Con_0 (\R^{2d})$ is the space of continuous functions on
$\R^{2d}$ that converge to $0$ at infinity.
Both spaces  $\Symbol{m,n}{-N}{2d}$ and $\Symbolo{m,n}{-N}{2d}$ are complete if 
equipped with the norm $M_{m,n}^{-N}(.)$.
The space $\Cinfc(\R^{2d})$ is dense in $\Symbolo{m,n}{-N}{2d}$ for
$N>0$. 

At first, we will be interested in the case $N= d+1$. 
Since
\begin{align*}
  \Symbol{m',n'}{-N'}{2d} \subset \Symbol{m,n}{-N}{2d}
\end{align*}
if $m'
\geq m$, $n'\geq n$, and $N' \geq N$, set
\begin{align*}
  \symb(\R^{2d})=\Symbol{0,d+1}{-(d+1)}{2d}.
\end{align*}
Set
also
\begin{align*}
  \symbo(\R^{2d})=\Symbolo{0,d+1}{-(d+1)}{2d}.
\end{align*}

Faster
decay with respect to $\xi$ will be considered, starting in Section~\ref{sec: Analysis of a commutator}. 
For symplicity,  we will use
\begin{align}
  \label{eq: symbols fast decay}
  \Symbolo{\infty,\infty}{-\infty}{2d} = \bigcap_{N\geq 0}
  \Symbolo{\infty,\infty}{-N}{2d}
\end{align}
in those later sections. 

%%%%%%%%%%%%%%%%%%%%%%%%
% definition           %
%%%%%%%%%%%%%%%%%%%%%%%%
\begin{definition}[semi-classical operators]
  \label{def: semi-classical operators}
 For $u\in \S(\R^d)$ and $a \in \symb(\R^{2d})$ one sets
\begin{equation}
  \label{4.pseudo}
  \OpH(a) u(x)=a(x, hD_x) u(x)=  (2\pi)^{-d} \int e^{ix\cdot \xi}a(x,h\xi){\hat u}(\xi) d\xi.
\end{equation} 
\end{definition}
The Schwartz kernel of $\OpH(a)$ is given
by
\begin{align*}
  %\label{eq: pseudo kernel}
  K_a(x,y)&=(2\pi)^{-d} \int e^{i (x-y)\cdot \xi} a(x,h\xi) d\xi
  = (2\pi h)^{-d} \int e^{i \frac{x-y}{h}\cdot \xi} a(x,\xi) d\xi\\
  &= h^{-d} k_a \big(x,\frac{x-y}{h}\big),\notag
\end{align*}
with
\begin{equation}
  \label{4.k}
  k_a(x,v)=(2\pi)^{-d} \int e^{i v\cdot \xi} a(x,\xi)d\xi.
\end{equation}
Note that \eqref{4.k} is well defined in the sense of classical integrals by the decay
property in the variable $\xi$ of the symbol $a$.
Observe that $L \exp(i v\cdot \xi) = \exp(i v\cdot \xi)$ with $L = (1 -i v
  \cdot \nabla_\xi)/\jp{v}^2$ leading to, with integrations by
parts, 
\begin{align*}
  k_a(x,v) = (2\pi)^{-d} \int e^{i v\cdot \xi}  (\transp L)^N a(x,\xi)d\xi,
\end{align*}
for $N \leq d+1$, with $\transp L= (1 +i v
  \cdot \nabla_\xi)/\jp{v}^2$.  One then obtains
\begin{align}
  \label{eq: estimation k_a}
   | k_a(x,v) | \lesssim M_{0,d+1}^{-(d+1)}(a) \jp{v}^{-(d+1)}, 
  \quad v \in \R^d, \ \ x \in \R^d \ \text{\pp}.
\end{align}
With Corollary~\ref{cor: Schur modified} one deduces the boundedness
of $\OpH(a)$ on $L^2(\R^d)$ with $a$ as above.
%%%%%%%%%%%%%%%%%%%%%%%%
% lemma                %
%%%%%%%%%%%%%%%%%%%%%%%%
\begin{lemma}
\label{4.continuite} Let $a \in \symb(\R^{2d})$. Then
$\OpH(a)$ extends as a uniformly bounded operator on $L^2(\R^d)$ and 
\begin{equation*}
  \Norm{\OpH(a)}{\L(L^2(\R^d))}
  \leq C_d M_{0,d+1}^{-(d+1)} (a).
\end{equation*} 
\end{lemma}

The following remark will be used in what follows.
%%%%%%%%%%%%%%%%%%%%%%%%
% remark               %
%%%%%%%%%%%%%%%%%%%%%%%%
\begin{remark}
  \label{remark: kernel d xi a}
  If $a \in \Symbol{0,d+2}{-(d+1)}{2d}$, note that one has
  \begin{align}
    \label{eq: remark kernel d xi a}
    k_{\d_{\xi_j} a}(x,v) = - i v_j k_a (x,v),\qquad j=1, \dots, d.
  \end{align}
  In fact, with an integration by parts one has
  \begin{align*}
    k_{\d_{\xi_j} a}(x,v)
    &=  (2\pi)^{-d} \int e^{i v\cdot \xi} \d_{\xi_j}  a(x,\xi)d\xi\\
    &= - (2\pi)^{-d}\int \d_{\xi_j} \big(e^{i v\cdot \xi}\big) a(x,\xi)
    d\xi
    = 
    - i v_j k_a (x,v).
  \end{align*}
  
\end{remark}

\bigskip
%%%%%%%%%%%%%%%%%%%%%%%%
% lemma                %
%%%%%%%%%%%%%%%%%%%%%%%%
\begin{lemma}
  \label{lemma: ka(x,v) ->0 as x->infty}
  Let $a \in \symbo(\R^{2d})$. Then, 
  $k_a(x,v) \to 0$ as $|x|\to \infty$ uniformly with respect to $v\in \R^d$.
\end{lemma}
\begin{proof}
  One writes 
  $|k_a(x,v)| \leq g(x) = (2\pi)^{-d} \int |a(x,\xi)| d\xi$.
  Since $|a(x,\xi)| \lesssim \jp{\xi}^{-d-1}$, one finds 
that $g(x) \to 0$ as $|x|\to \infty$ by the \Ldct.
\end{proof}
The following lemma will be of great use in what follows.
%%%%%%%%%%%%%%%%%%%%%%%%
% lemma                %
%%%%%%%%%%%%%%%%%%%%%%%%
\begin{lemma}
  \label{lemma: preparation lemma kernels}
  Let $a \in \symbo(\R^{2d})$.  Let 
 $\rho_h(x,v) \in \Con^0(\R^d\times \R^d)$ be such that 
 \begin{enumerate}
   \item $\Norm{\rho_h}{L^\infty}  \leq C_0$ uniformly in $h$,
   \item $\rho_h(x,v) \to 0$
 as $h\to 0$ uniformly for $(x,v)$ in any compact set. 
\end{enumerate}

Then, one has 
\begin{align*}
   \esssup_{x\in \R^d} \Norm{\rho_h k_a(x,.)}{L^1(\R^d)} \to 0 
   \ \ \et \ \ 
   \esssup_{x\in \R^d} \Norm{\rho_h k_a(x+h.,.)}{L^1(\R^d)} \to 0,
 \end{align*}
 as $h \to 0$, with $k_a(x,v)$ as in \eqref{4.k}.
\end{lemma}
%%%% proof of lemma
\begin{proof}
  Set $m_h = \rho_h k_a $. 
  For $R>0$, by \eqref{eq: estimation k_a} one writes 
\begin{align*}
  \int_{|v| > R} |m_h(x,v)|\, d v \lesssim C_0
  \int_{|v| > R} \jp{v}^{-(d+1)} d v,
\end{align*}
for almost all $x\in \R^d$. 
Let $\eps>0$.  For $R>0$ chosen \suff large one has 
\begin{align}
  \label{eq: lemma 4.10 R large}
  \int_{\R^d} |m_h(x,v)| d v \leq \eps + \int_{|v| \leq R} |m_h(x,v)|\,
  d v, \qquad x\in \R^d \ \text{\pp}.
\end{align}
Next, one writes 
\begin{align*}
  \int_{|v| \leq R} |m_h(x,v)|\, d v 
   \leq C_0 \int_{|v| \leq R} | k_a(x,v) |\, d v.
\end{align*}
Thus, by Lemma~\ref{lemma: ka(x,v) ->0 as x->infty}, for some $R'>0$, one has $\int_{|v| \leq R}
|m_h(x,v)|\, d v \leq \eps/C_0$ for $|x| \geq R'$. One thus has
\begin{align*}
  \esssup_{|x| \geq R'}\int_{\R^d} |m_h(x,v)|\, d v  \leq 2 \eps.
\end{align*}
Consider now the case $|x|\leq R'$.  By hypothesis $|\rho_h(x,v)|\to 0$ as $h\to 0$
uniformly with respect to $x$ and $v$ if $|x|\leq R'$
and $|v| \leq R$. With \eqref{eq: estimation k_a} one finds
\begin{align*}
  \int_{|v| \leq R} |m_h(x,v)|\, d v 
   \lesssim \int_{|v| \leq R} |\rho_h(x,v)|\, d v,
\end{align*}
for almost all $x$ such that  $|x|\leq R'$.
One thus  finds that $\int_{|v| \leq R} |m_h(x,v)|\, d v \leq \eps$
for such $x$ and  for $h>0$ chosen \suff small. With~\eqref{eq: lemma
  4.10 R large} one thus concludes that 
\begin{align*}
  \esssup_{|x| \leq R'}\int_{\R^d} |m_h(x,v)|\, d v  \leq 2 \eps,
\end{align*}
if  $h>0$ is chosen \suff small and thus 
$\esssup_{x \in \R^d}\int_{\R^d} |m_h(x,v)|\, d v \leq 2 \eps$.

\medskip
One obtains {\em mutatis mutandis} that 
$\esssup_{x \in \R^d}\int_{\R^d} |m_h(x+hv,v)|\, d v  \leq 2 \eps$,
for $h$ chosen \suff small. 
\end{proof}

%%%%%%%%%%%%%%%%%%%%%%%%
% proposition          %
%%%%%%%%%%%%%%%%%%%%%%%%
 \begin{proposition}\label{propcommut}
 Let $a \in \symbo(\R^{2d})$. 
 \begin{enumerate}
  \item Consider $\theta \in \Con^0(\R^d) \cap
 L^\infty(\R^d)$. one has
\begin{equation}
  \label{WCM}
  \lim_{h \rightarrow 0 } \Norm{[ a(x, hD_x) , \theta]}{\mathcal{L} (L^2(\R^d))} =0.
\end{equation}
\item 
More generally, if  $(\theta_k)_{k\in \N}  \subset L^\infty(\R^d)$, $\theta \in \Con^0(\R^{d+1})$ is 
such that  $\| \theta_k - \theta\|_{L^\infty} \to 0$ as
$k\to + \infty$, then 
\begin{equation}\label{WCMbis}
  \Norm{ [ a(x, hD_x) , \theta_k] }{\mathcal{L} (L^2(\R^d))}
  =o(1)_{h\to 0\ \et\  k \to \infty}.
\end{equation}
\item Assume in addition that $a \in \Symbolo{0,2+d}{-(1+d)}{2d}$ and
  consider $\theta \in W^{1,\infty}(\R^d)$ then one  has
\begin{align*}%\label{eq: commutator Lipschitz}
  \BigNorm{ [ \OpH(a) , \theta] }{\mathcal{L}(L^2(\R^d))} 
  =O(h).
\end{align*}
\item 
For $a \in \Symbolo{0,2+d}{-(1+d)}{2d}$, if moreover
$\theta \in \Con^1(\R^d) \cap W^{1,\infty}(\R^d)$ then one  has the following properties
\begin{align}\label{WCM2}
\BigNorm{ [ \OpH(a) , \theta] 
  + ih\sum_{j=1}^d 
  \frac{ \d \theta}{ \d x_j} 
  \OpH\bigl(\frac{ \d a} { \d\xi_j}\bigr)
  }{\mathcal{L}(L^2(\R^d))} 
  = o(h),
 \intertext{and}
  \label{WCM3}
 \BigNorm{   [ \OpH(a) , \theta] 
  + i h\sum_{j=1}^d  \OpH\bigl(\frac{ \d a} { \d \xi_j}\bigr)
  \frac{ \d \theta}{ \d x_j}}{\mathcal{L} (L^2(\R^d))} = o(h).
\end{align}
\item 
More generally, if  $(\theta_k)_{k\in \N}  \subset Lip(\R^d)$ is 
such that  $\| \theta_k - \theta\|_{Lip} \to 0$ as
$k\to + \infty$, then 
\begin{align}\label{WCM2bis}
 \BigNorm{   [ \OpH(a) , \theta_k] 
   + i h\sum_{j=1}^d  \frac{ \d \theta_k}{ \d x_j} \OpH\bigl(\frac{ \d a} { \d \xi_j}   \bigr)}{\mathcal{L} (L^2(\R^d))} 
  =h o(1)_{h\to 0 \text{ and } k \to \infty}, 
   \intertext{and}
   %\label{WCM2ter}
 \BigNorm{   [ \OpH(a) , \theta_k] + i h\sum_{j=1}^d  \OpH\bigl(\frac{ \d a} { \d \xi_j}   \bigr)\frac{ \d \theta_k}{ \d x_j}}{\mathcal{L} (L^2(\R^d))} 
  =h o(1)_{h\to 0 \text{ and } k \to \infty}.\notag
\end{align}
\item 
Finally, assume that $a\in \Symbolo{0,N+1+d}{-(d+1)}{2d}$. Let $\phi \in \Cinfc(\R^d)$ be
such that  $\phi a =0$. One has 
\begin{equation}
  \label{WCM4}
  \Norm{ \OpH (a) \circ \phi}{\L (L^2(\R^d))} =o (h^N).
\end{equation}
\end{enumerate} 
\end{proposition}
%%%% proof of proposition
 \begin{proof}
 The kernel of the operator $[ a(x, hD_x) , \theta]$ is given by
 \begin{align}
   \label{eq: kernel commutator}
   K(x,y) = K_a(x,y)  \big( \theta(y) - \theta(x)\big)
   = h^{-d}m_h \big(x,\frac{x-y}{h}\big),
 \end{align}
with $m_h(x,v) = k_a(x,v) \big(\theta(x-hv)-\theta(x)\big)$.
Since $\theta$ is continuous
it is uniformly continuous on any compact set. Thus, one finds $|\theta(x-hv)-\theta(x)|\to 0$ as $h\to 0$
uniformly with respect to $x$ and $v$ if $|x|\leq R'$
and $|v| \leq R$. With Lemma~\ref{lemma: preparation lemma kernels} one obtains  that 
\begin{align*}
  \esssup_{x \in \R^d}\Norm{m_h(x,.)}{L^1(\R^d)} \to 0 
  \ \ \et 
  \esssup_{y \in \R^d}\Norm{m_h(x+h.,.)}{L^1(\R^d)} \to 0,
\end{align*}
as $h\to 0$. 
With Lemma~\ref{lemma: Schur modified} one concludes that the limit in \eqref{WCM}
holds.

\bigskip
To obtain~\eqref{WCMbis} one writes
\begin{align*}
  [ a(x, hD_x) , \theta_k] =   a(x, hD_x) (\theta_k- \theta) -
  (\theta_k- \theta) a(x, hD_x)+  [ a(x, hD_x) , \theta].
\end{align*}
 
\bigskip
Assume now that $a$ and $\theta$ fullfil the requirements of point
(3). The kernel of the operator $\OpH(a) , \theta]$ is given by
\eqref{eq: kernel commutator}. With the first-order Taylor formula one writes
\begin{align*}
  \theta(x-hv)-\theta(x) =  -h \int_0^1 d_x \theta (x-s hv) (v)  d s
  = - h \sum_j v_j \Theta_j(x,hv), 
\end{align*}
with $\Theta_j(x,hv) = \int_0^1 \d_{j} \theta (x-s hv)  d s$. 
With the additional regularity of
$a(x,\xi)$ and Remark~\ref{remark: kernel d xi a} one finds 
\begin{align*}
  m_h(x,v) = - i h \sum_j k_{\d_{\xi_j}a}(x,v) \Theta_j(x,hv),
\end{align*}
yielding
\begin{align*}
  |m_h(x,v)|  
  \lesssim  h \sum_j M_{0,d+1}^{-(d+1)}(\d_{\xi_j}a)\jp{v}^{-(d+1)}
  \lesssim h M_{0,d+2}^{-(d+1)}(a)\jp{v}^{-(d+1)},
\end{align*}
as $|\Theta_j(x,hv)|\lesssim 1$ uniformly in $x$, $v$ and $h$. With Corollary~\ref{cor: Schur
  modified}  one deduces the result of point~(3).
 
\bigskip
Assume now that $a$ and $\theta$ fullfil the requirements of point (4)
of the proposition and denote by ${\mathsf K}(x,y)$ the kernel of $[
  \OpH(a) , \theta] +ih\sum_{j=1}^d \d_{x_j} \theta
\OpH\big(\d_{\xi_j} a\big)$. One has ${\mathsf K}(x,y) =
h^{-d}r_h(x,(x-y)/h)$ with
\begin{align*}
  r_h(x,v)= k_a(x,v) \big(\theta(x-hv)-\theta(x) + h d_x \theta(x) (v) \big),
\end{align*}
using Remark~\ref{remark: kernel d xi a}.
The first-order Taylor formula gives
\begin{align*}
  \theta(x-hv)-\theta(x) + h d_x \theta(x) (v) 
  = h \int_0^1 (d_x \theta(x)  - d_x \theta (x-s hv)) (v)  d s.
\end{align*}
Setting $A_h^j(x,v) = \int_0^1 \big(\d_j \theta(x)  
  - \d_j \theta (x-s h v) \big) d s$,
one finds 
\begin{align*}
  r_h(x,v) = i h \sum_j A_h^j(x,v) k_{\d_{\xi_j} a}(x,v),
\end{align*}
using again Remark~\ref{remark: kernel d xi a}.
Since $\theta \in \Con^1(\R^d) \cap W^{1,\infty}(\R^d)$ one finds  $A_h^j \in \Con^0(\R^{2d})$
and
$\Norm{A_h^j}{L^\infty(\R^{2d})}\leq C_0$ for some $C_0>0$
uniformly with respect to $h$.
Moreover, if $L$ is a compact of $\R^{2d}$, and $0\leq  h
\leq h_0$, if $x, v \in L$ then $x-s h v$ remains
in a compact set of $\R$ where $\d_j \theta$ is uniformly
continuous. One concludes that $A^j_h(x,v)\to 0$ as $h \to 0$ uniformly
with respect to $(x,v)\in L$. Since $\d_{\xi_j} a \in \symbo(\R^{2d})$,  with Lemma~\ref{lemma: preparation lemma
  kernels}
one concludes that
\begin{align*}
  \esssup_{x \in \R^d}\Norm{r_h(x,.)}{L^1(\R^d)} \to 0 
  \ \ \et 
  \esssup_{y \in \R^d}\Norm{r_h(x+h.,.)}{L^1(\R^d)} \to 0,
\end{align*}
as $h\to 0$. 
With Lemma~\ref{lemma: Schur modified} one concludes that the limit in \eqref{WCM2}
holds.

Following the same strategy, one finds  that the kernel of the operator 
\begin{align*}
  [ \OpH(a) , \theta] 
  -i \sum_j \OpH\big(\d_{\xi_j} a \big) \d_j \theta(x)
\end{align*}
  is given by $h^{-d}\tilde{r}_h(x,(x-y)/h)$, with
\begin{align*}
 \tilde{r}_h(x,v)= k_a(x,v)\big( \theta(x-hv)-\theta(x) +h d \theta (x-hv)(v)\big) 
\end{align*}
and applying the Taylor formula as above and Lemma~\ref{lemma:
  preparation lemma kernels} one obtains the limit in~\eqref{WCM3} by
 Lemma~\ref{lemma: Schur modified}.

\bigskip
 To prove~\eqref{WCM2bis}, as above one writes 
\begin{align*}
  [ a(x, hD_x) , \theta_k] =   [ a(x, hD_x),  \theta_k- \theta]  +  [ a(x, hD_x) , \theta].
\end{align*}
and by point (2) of the proposition one observe that  it suffices to prove 
\begin{align}
  \label{eq: to prove WCM2bis}
  \Norm{[ a(x, hD_x),  \theta_k- \theta]}{\mathcal{L}(L^2)}  =
  h o(1),
\end{align}
as $k \to +\infty$. 
Set $\alpha_k = \theta_k- \theta$. The kernel of $[ a(x, hD_x),
\alpha_k]$ is given by ${\mathsf L}(x,y) = h^{-d} q_h (x,(x-y)/h)$
with 
\begin{align*}
  q_h (x,v)
  &= k_a(x,v)\big(\alpha_k (x- h v) -\alpha_k(x)\big) 
  = - h k_a(x,v) \int_0^1 d_x \alpha_k (x-s h v)(v)\, d s\\
  &= -h \sum_{j=1}^d Q^j_{h,k} k_{\d_{\xi_j} a}(x,v),
\end{align*}
where $Q^j_{h,k}(x,v) =\int_0^1 \d_{j} \alpha_k (x-thv)\, d
t$. Arguing as above one obtains \eqref{eq: to prove WCM2bis}.

\bigskip Finally, we
consider the last statement of the proposition, with
$\phi \in \Cinfc(\R^d)$ such that $a\phi = 0$.  The kernel of
$\OpH (a) \phi$ reads $h^{-d}s_h(x,(x-y)/h)$ with
 $s_h(x,v) = k_a(x,v) \phi (x-hv)$.
 With the Taylor formula one writes 
\begin{align*}
  \phi(x - h v) 
  &= \sum_{j\leq N-1}  \frac{(-h)^j}{j!} d^j \phi(x)( v,\dots, v) \\
  &\quad+ \frac{(-h)^N}{(N-1)!} \int_0^1 d^N \phi (x - s h v)(v,
  \dots, v)(1-s)^{N-1} ds,
\end{align*}
which we write 
\begin{align*}
  \phi(x - h v) 
  &= \sum_{j\leq N}  \frac{(-h)^j}{j!} d^j \phi(x)( v,\dots, v) 
    + R^N_h(x,v),
\end{align*}
with 
\begin{align*}
 R^N_h(x,v) &= \frac{(-h)^N}{(N-1)!} \int_0^1 d^N \phi (x - s h v)(v,
 \dots, v)(1-s)^{N-1} ds \\
 &\quad - \frac{(-h)^N}{N!} d^N \phi(x)( v,\dots, v)\\
  &= \frac{(-h)^N}{(N-1)!} \int_0^1 \big( d^N \phi (x - s h v) - d^N
    \phi(x) \big) (v,\dots, v)(1-s)^{N-1} ds.
\end{align*}
Since $\phi a =0$ the same holds for $\d_{x}^\beta \phi a$ for any
$\beta$ and  $s_h(x,v) = k_a(x,v) R^N_h(x,v)$. One has 
\begin{align*}
 R^N_h(x,v) = h^N \sum_{|\beta| = N} v^{\beta} \psi_{h,\beta}^N(x,v),
\end{align*}
with $\Norm{\psi_{h,\beta}^N}{L^\infty(\R^d \times \R^d)}\leq C_0$ for some $C_0>0$
uniformly with respect to $h$ and $\psi_{h,\beta}^N (x,v)\to 0$ as $h \to 0$ uniformly
with respect to $(x,v)$ in a compact set. 
Iterating \eqref{eq: remark kernel d xi a}, one finds
\begin{align*}
 k_a(x,v)  R^N_h(x,v) 
  =  h^N \sum_{|\beta| = N}  i^{|\beta|} k_{\d_\xi^\beta a}(x,v) \psi_{h,\beta}^N(x,v),
\end{align*}
and thus  Lemma~\ref{lemma: preparation lemma
  kernels} and Lemma~\ref{lemma: Schur modified}.
 imply that \eqref{WCM4} holds. 
 \end{proof}
 
 \medskip
 Let $K$ be a compact set of $\R^d$ and $a \in \symbo(\R^{2d})$ such
 that $\supp a \subset K \times\R^d$.
 For these particular symbols,  if $\phi \in \Con_c^0(\R^d)$ is equal
 to $1$ on the $x$-projection of  $\supp a$ then 
 \begin{align*}
   \Norm{\OpH (a) (1 - \phi)}
   {\L (L^2(\R^d))}
   = o(1), \qquad h \to 0,
 \end{align*}
 by Proposition~\ref{propcommut}.
 In fact, we will be inclined to define semi-classical operators up to
 operators in $\L (L^2(\R^d))$ whose norm
 is $o(1)$ as $h \to 0$. Then we denote $[\OpH](a)$ the
 class of operators defined by $\OpH (a) (\phi u)$ where $\phi$ is as
 above. This is further explained by our intention to use
 semi-classical operators on manifolds, here $\M$ or $\L$, that we now
 present. 

%%%%%%%%%%%%%%%%%
% subsection
\subsection{Tangential symbols and operators}
\label{sec: Tangential symbols and operators}
In what follows we also use tangential operator. They are associated
with symbols of the form $a(y, \eta')$ with
$y \in \R^{d}$ and $\eta' \in \R^{d-1}$. 
%%%%%%%%%%%%%%%%%%%%%%%%
% definition           %
%%%%%%%%%%%%%%%%%%%%%%%%
\begin{definition}[tangential symbols]
  \label{def: tangential symbols}
   Let $m,n \in \N \cup \{+\infty\}$, with $n \geq d$, and $N \in \R^+$. 
  Denote by $\Symbolt{m,n}{-N}{d}{d-1}$ the space of
  all functions $a(y,\eta)$, $y\in \R^{d}$, $\eta' \in \R^{d-1}$, such that
  $\d_y^\alpha \d_{\eta'}^\beta a \in L^1_{\loc}(\R^{2d-1})$ for
  $\alpha \in \N^{d}, \beta \in \N^{d-1}$ with $|\alpha|\leq m$, $|\beta| \leq n$,
  and for some $C_{\alpha, \beta}>0$, 
  \begin{align*}
    \bignorm{\d_y^\alpha \d_{\eta'} ^\beta a(y,\eta')}{} 
    \leq  C_{\alpha, \beta} \jp{\eta'}^{-N}, \qquad y\in \R^{d}, \
    \eta' \in \R^{d-1}.
  \end{align*}
In addition, one sets  $\Symbolto{m,n}{-N}{d}{d-1}$, $n \geq d$,  as the set of
all symbols $a \in \Symbolt{m,n}{-N}{d}{d-1}$ with moreover
$\d_y^\alpha \d_{\eta'}^\beta a \in \Con_0(\R^{2d-1})$ for $\alpha \in \N^{d}, \beta \in \N^{d-1}$
with
$|\alpha|\leq m$ and $|\beta| \leq n -d$.
\end{definition}
 Equipped with the norm
\begin{equation*}
  %\label{eq: est tangential symbol}
  M_{\T, m,n}^{-N}(a)=\max_{{|\alpha|\leq m}\atop{|\beta| \leq n}}
  \esssup_{(y,\eta') \in \R^{d}\times \R^{d-1}}
  \bignorm{\d_y^\alpha \d_{\eta'} ^\beta a(y,\eta')}{} \jp{\eta'}^{N},
\end{equation*}
both spaces  $\Symbolt{m,n}{-N}{d}{d-1}$ and $\Symbolto{m,n}{-N}{d}{d-1}$ are complete.

Note that $\Symbolt{m',n'}{-N'}{d}{d-1} \subset \Symbolt{m,n}{-N}{d}{d-1}$ if $m'
\geq m$, $n'\geq n$, and $N' \geq N$.  The case $N=d$ is of interest similarly to symbols defined in
Section~\ref{sec: Semi-classical operators on Rd}. Set
\begin{align*}
  &\symbt(\R^{d}\times \R^{d-1})=\Symbolt{0,d}{-d}{d}{d-1},\\
  &\symbto(\R^{d}\times \R^{d-1})=\Symbolto{0,d}{-d}{d}{d-1},
\end{align*}
and 
\begin{equation*}
  %\label{eq: est tangential symbol 2}
  N_{n}(a)=M_{\T, 0,n}^{-d}(a) = \max_{|\beta| \leq n} \esssup_{(y,\eta')}
  \bignorm{ \d_{\eta'} ^\beta a(y,\eta')}{} \jp{\eta'}^{d}.
\end{equation*}
With $y= (y', z)$, observe that $N_{n}(a)$ correponds to
$M_{0,n}^{-d}(a)$ in \eqref{eq: est symbols} with $z$ acting as a parameter.

\bigskip
For $a\in \symbt(\R^{d}\times \R^{d-1})$, the associated operator is defined by 
\begin{equation*}
  %\label{tangential pseudo}
  \OpH(a) u(y)=a(y, hD_y') u(z,y')
  =  (2\pi)^{1-d} \int_{\R^{d-1}} e^{iy'\cdot \eta'} a(z,y',h \eta'){\hat u}(z,\eta') d\eta',
\end{equation*}
where the Fourier transformation acts in the $y'$ variables. 
In fact,
the action of $\OpH(a)$ is through the Schwartz kernel
\begin{align*}
  K(y,\ty) =  K_a(y',\ty'; z) \otimes \delta_{z-\tilde{z}}, 
\end{align*}
with the tangential kernel
 \begin{align}
   \label{eq: tangential kernel 1}
  K_a(y',\ty'; z)= (2 \pi)^{1-d} \int_{\R^{d-1}} e^{i (y'- \tilde{y'})\cdot \eta'}
   a(y',z, h \eta')\,  d \eta'.
 \end{align}
Then, one has 
\begin{align}
  \label{eq: tangential kernel 2}
  \OpH(a) u (y',z) = \int_{\R^{d-1}} K_a(y',\ty'; z)\,  u(\ty',z)\,
  d \ty'.
\end{align}
If  $a \in \symbt(\R^{d}\times \R^{d-1})$ one finds $K_a(y', \ty'; z) = h^{1-d} k_a
\big(y,\frac{y'-\ty'; z}{h}\big)$ and
\begin{align*}
  %\label{eq: estimation k_a tangential}
   | k_a(y',v; z) | \leq C N_{d}(a) \jp{v}^{-d}, 
  \quad v \in \R^{d-1}, \ z \in \R, \ y' \in \R^{d-1} \ \text{\pp},
\end{align*}
as in \eqref{eq: estimation k_a}.  With Corollary~\ref{cor: Schur
  modified} one has 
\begin{align*}
  \Norm{\OpH(a) u (.,z)}{L^2(\R^{d-1})}
  \leq C_d  N_{d}(a) \Norm{u (.,z)}{L^2(\R^{d-1})},  
  \qquad z \in \R \ \text{\pp},
\end{align*}
for some $C_d>0$ uniform with respect to $z$, yielding 
\begin{align*}
  \Norm{\OpH(a) u }{L^2(\R^{d})}^2 &= 
  \int_\R \Norm{\OpH(a) u (.,z)}{L^2(\R^{d-1})}^2 d z
  \lesssim N_{d}(a)^2 \int_\R \Norm{u (.,z)}{L^2(\R^{d-1})}^2 d z\\
  &\lesssim N_{d}(a)^2  \Norm{ u }{L^2(\R^{d})}^2,
\end{align*}
that is,  the following continuity result. 
%%%%%%%%%%%%%%%%%%%%%%%%
% lemma                %
%%%%%%%%%%%%%%%%%%%%%%%%
\begin{lemma}
\label{lemma: continuity tangential op}
Let $a(y,\eta') \in \symbt(\R^{d}\times \R^{d-1})$. Then
$\OpH(a)$ extends as a uniformly bounded operator on $L^2(\R^{d})$ and 
\begin{equation*}
  \Norm{\OpH(a)}{\L(L^2(\R^{d}))}
  \leq C_d N_{d}(a) .
\end{equation*} 
\end{lemma}

\bigskip
In what follows, we also use symbols of the form $a(y, \eta) = b(y,
\eta')  f(\zeta)$ with $b(y, \eta') \in \symbt(\R^{d}\times \R^{d-1})$ and $f(\zeta)$ a Fourier
mutliplier; see for instance  Proposition~\ref{prop: Euclidean
  symbol division} 
In fact, one has 
\begin{align}
  \label{eq: decomp b f (zeta)}
  \OpH  \big( b(y, \eta') f(\zeta)\big) = \OpH (b) f (h D_z), 
\end{align}
and thus one can write
\begin{align}
  \label{eq: est b f (zeta)}
  \bigNorm{\OpH  \big( b(y, \eta') f(\zeta)\big)}{\L(L^2(\R^{d}))} 
  \leq \Norm{\OpH(b)}{\L(L^2(\R^{d}))} \Norm{ f (h D_z)}{\L(L^2(\R^{d}))},
\end{align}
that we will use several times in what follows. Since 
\begin{align*}
  f (h D_z) = \mathcal{F}^{-1}_{\zeta \rightarrow z} f(h \zeta) \mathcal{F}_{z\rightarrow \zeta}
\end{align*}
if $f$ is bounded one finds 
\begin{align}
  \label{eq: est  f (zeta)}
  \Norm{ f (h D_z)}{\L(L^2(\R^{d}))} \leq \Norm{f}{L^\infty},
\end{align}
since the Fourier transformation
$\mathcal{F}_{z\to\zeta}$ is a an isometry on $L_z^2(\R;
L^2_{y'}(\R^{d-1}))$.

Similarly $f(h D_z)$ has kernel on $\R^{d}$ given by 
\begin{align*}
  \delta_{y'-\ty'}  \otimes K_f(y; z, \tilde{z}),
\end{align*}
with the part only acting in the $z$ variable given by 
\begin{align}
  \label{eq: kernel multiplier 1}
  K_f(z, \tilde{z}) =  h^{-1} k_f \big((z - \tilde{z})/h\big),
\end{align}
with 
\begin{align}
  \label{eq: kernel multiplier 2}
 k_f(v) = (2 \pi)^{-1} \int_{\R} e^{i v\cdot \zeta}
   f(\zeta)\,  d \zeta = \check{f}(v).
\end{align}

%%%%%%%%%%%%%%%%%
% subsection
\subsection{Semi-classical operators on a manifold}
\label{sec: Semi-classical operators on a manifold}
 Let $\mathcal N$ be a $\Con^1$-manifold of dimension $d$ equipped
 with a density measure $\rho$ that allows one to define $L^2(\mathcal N)$.  We denote
 by $\mathcal P (\mathcal N)$ the algebra of bounded operators $B_h$
on $L^2 (\mathcal N)$, depending on $ h\in (0, h_0 ]$ as a
 parameter, and by $\mathcal{R}(\mathcal N)$ the ideal of
 $\mathcal{P} (\mathcal N)$ of the operators $B_h$ such that
 $\Norm{B_h}{\L(L^2)}=o(1) $. Set
 $\mathcal{Q}(\mathcal N) = \mathcal{P} (\mathcal N)/ \mathcal{R}
 (\mathcal N)$.

\bigskip
The following lemma is key towards the 
notion of semi-classical operators on a manifold.
%%%%%%%%%%%%%%%%%%%%%%%%
% lemma                %
%%%%%%%%%%%%%%%%%%%%%%%%
\begin{lemma}[\protect{\cite[Lemme 1.10]{GL:1993}}]
\label{4.changtvar}
Consider $\psi : V\to U$ a $\Con^1$-diffeomorphism
between two open subsets of $\R^d$. Let $a\in \Con^0_c(U \times \R^{d})$
be such that $\d^\beta_\xi a \in \Con^0_c(U \times \R^{d})$ for
$|\beta|\leq d+1$. Set $ b(y,\eta)=a(\psi(y),\transp d\psi^{-1}_y(\eta))\in
\Con^0_c(V\times \R^d)$. Then, for any compact set $ K\subset U$ one has
\begin{equation*}
\Norm{a(x,hD_x) u \circ \psi- b(y,hD_y) (u\circ \psi)}{L^2(V)} 
=o(1) \Norm{u}{L^2(U)},
\qquad h \to 0,
\end{equation*}
uniformly with respect to $u\in L^2(U)$ with support in $K$.  
\end{lemma}
%%%%%%%%%%%%%%%%%%%%%%%%
% definition           %
%%%%%%%%%%%%%%%%%%%%%%%%
\begin{definition}
  \label{def: symbol manifold}
  Let $\mathcal N$ be a $\Con^1$-manifold of dimension $d$.  
 Denote by $\symbc(T^* \mathcal N)$ the space of functions $a\in
\Con_c^0(T^* \mathcal N)$ such that for $|\beta| \leq d+1$, one has  $\d^\beta_\xi a \in \Con_c^0(T^* \mathcal N)$.
\end{definition}

For $a  \in \symbc(T^* \mathcal N)$ and  a chart $\chart = (\O, \cdiff)$ we
denote by $a^\chart$ the local representative of $a$ in this chart. 
Consider two local charts $\chart_1 = (\O_1, \cdiff_1)$
and $\chart_2 = (\O_2, \cdiff_2)$ with $W = \O_1 \cap \O_2 \neq \emptyset$
and $a \in \symbc(T^* \mathcal N)$ supported in $W$. Then,  the
representatives $a^{\chart_1}$ and $a^{\chart_2}$ 
fulfill the assumption of Lemma~\ref{4.changtvar} with $U =
\cdiff_1(W)$, $V = \cdiff_2(W)$ and $\psi = \cdiff_1 \circ
\cdiff_2^{-1}$.

Consider a chart $\chart = (\O, \cdiff)$ as above, $a \in \symbc(T^* \mathcal N)$ and
$\theta, \chi \in \Con_c(\O)$, $\theta\equiv 1$ in a \nhd of
$\supp \chi$. For $ u\in L^2(\mathcal N)$ one may compute
\begin{equation*}
 \chdiff^{*} \circ (\chi a)^\chart (x,hD_x)  \circ
  \big(\chdiff^{-1}\big)^* (\theta u),
\end{equation*}
yiedling an $L^2$-function on $\mathcal N$.

Consider now a locally finite $\Con^1$-partition of unity
$(\chi_i)_{i\in \I}$ subordinated to a given atlas $\mathcal A =
(\chart_i)_{i \in \I}$, $\chart_i = (\O_i, \cdiff_i)$  and a familly of
localisation functions $(\theta_i)_{i \in \I}$ with $\supp \theta_i
\subset \O_i$ and $\theta_i \equiv 1$ on $\supp \chi_i$. We form 
\begin{equation*}
  A u  = \sum_{i\in I} \chdiff_i^{*} \circ (\chi_i a)^{\chart_i} (x,hD_x)  \circ
  \big(\chdiff_i^{-1}\big)^* (\theta_i u).
\end{equation*}

\medskip
From Lemma~\ref{4.changtvar}, the class of the operator $A$ defined above
in  $\mathcal{Q}(\mathcal N)$ is independent of the choice of the
atlas $\mathcal A$, the partition
of unity $(\chi_i)_{i\in \I}$, and the localisation functions
$(\theta_i)_{i \in \I}$. We denote this class by $[\OpH](a)$.
 
\medskip Let $\varphi\in \Con^0(\mathcal N)$. Let $B_h$ and
$\tilde{B}_h$ be two representatives of a class in
$\mathcal{Q}(\mathcal N)$, that is,
$[B_h]= [\tilde{B}_h]$. Observe that $[\varphi B_h]= [\varphi\tilde{B}_h]$, thus defining a multiplication by the function $\varphi$ on
$\mathcal{Q}(\mathcal N)$, which one writes $[\varphi B_h] = \varphi [B_h]$.
If $ a\in \symbc(T^*\mathcal N)$ one has 
\begin{equation*}
%\label{4.32}
[\OpH] (\varphi a)= \varphi   [\OpH](a).
\end{equation*}

%%%%%%%%%%%%%%%%%
% subsection
\subsection{Semi-classical measures}\label{mesureeteq}
This section is borrowed from~\cite{GL:1993} and \cite{Burq:1997}; it
recalls the basic properties of semi-classical measures.

In what follows, we call a sequence of scales $H= (\hk)_k$ a sequence
of positive real numbers that converges to $0$.  If such a sequence of
scales is used we will write $\OpH$ in place of ${\Op^{\hk}}$ for
concision if no confusion can arise. 

%%%%%%%%%%%%%%%%%%%%%%%%
% definition           %
%%%%%%%%%%%%%%%%%%%%%%%%
\begin{definition}[semi-classical measure]
  \label{def: semi-classical measure}
  Let $H= (\hk)_k$ be a sequence of
  scales and $(u_k)_k$ be a bounded sequence  of
  $L^2(\R^d)$.
  Let $\mu$ be a nonnegative Radon measure on $\R^{2d}$.  One says
  that $(u_k)_k$ admits $\mu$ as its semi-classical measure (\scm) at scale $H=(\hk)_k$
  if one has 
  \begin{equation}\label{limite-a}
    \lim_{k\rightarrow + \infty} \inp{\OpH (a) u_k}{u_k}_{L^2(\R^d)}
    = \int_{\R^{2d}} a(x,\xi) d \mu(x,\xi) = \dup{\mu}{a},
 \end{equation}
  for any $a\in \symbo(\R^{2d})$.
\end{definition}

\begin{definition}[mass leakage at infinity]
  \label{def: mass leak at infinity}
  One say that no mass leaks  at infinity at scale $H$ if one has 
  \begin{equation*}
    %\label{eq: mass leak at infinity}
    \lim_{R\rightarrow +\infty}\, \limsup _{k\rightarrow +\infty}
    \Big( \int _{|x|\geq R} |u_{k}(x) |^2 dx
    +
    \int _{h_{k}|\xi|\geq R} |{\hat{u}_{k}}(\xi) |^2 d\xi\Big) =0.
\end{equation*}
One says that there is some mass leakage  at scale $H$ at infinity otherwise.
\end{definition}
%%%%%%%%%%%%%%%%%%%%%%%%
% lemma                %
%%%%%%%%%%%%%%%%%%%%%%%%
\begin{lemma}
  \label{lemma: h-oscillation}
  Suppose that $\big(\hk^s|D_x|^s u_k\big)_k$ is $L^2$-bounded for
  some $s>0$. Then
  \begin{equation*}
    \lim_{R\rightarrow +\infty}\, \limsup _{k\rightarrow +\infty}
    \int _{h_{k}|\xi|\geq R} |{\hat{u}_{k}}(\xi) |^2 d\xi =0.
\end{equation*}
\end{lemma}
%%%% proof of lemma
\begin{proof}
  Write
  $\int_{h_{k}|\xi|\geq R} |{\hat{u}_{k}}(\xi) |^2 d\xi
    \leq R^{-2s}\int_{\R^2} \hk^{2s}|\xi|^{2s}
    |{\hat{u}_{k}}(\xi)|^2 d\xi
    \lesssim R^{-2s}$.
\end{proof}
The following proposition states that up to a subsequence extraction,
every bounded sequence in $L^2(\R^d)$ admits a \scm at some given
scale. It moreover provides a criterium for mass conservation in the
limiting process.
%%%%%%%%%%%%%%%%%%%%%%%%
% proposition          %
%%%%%%%%%%%%%%%%%%%%%%%%
\begin{proposition}[\protect{\cite[Propositions 1.4 and 1.6]{GL:1993}}]
\label{prop: properties scm}
For any  sequence of scales $H= (\hk)_k$, and any bounded sequence
$(u_k)_{k}\subset  L^2(\R^d)$, there exist  a
subsequence $(k_n)_{n\in \mathbb{N}}$ and a nonnegative measure $\mu$ on~$\R^{2d}$
such that the following properties hold: 
\begin{enumerate}
\item $\mu$ is the \scm for the sequence $(u_{k_n})_n$ at scale $(h_{k_n})_n$

\item If no mass leaks  at infinity at scale $H$ in the sense of
  Definition~\ref{def: mass leak at infinity}, then 
\begin{equation}
  \label{eq: pas-de-perte}
  \lim_{n\rightarrow + \infty}\Norm{u_{k_n}}{L^2(\R^{d})}^2 = \mu(\R^{2d}).
\end{equation}
meaning mass is preserved in the limiting process. 
\end{enumerate}
\end{proposition}
%%%%%%%%%%%%%%%%%%%%%%%%
% lemma                %
%%%%%%%%%%%%%%%%%%%%%%%%
 \begin{lemma}
   \label{lem: continu}
 Assume that $\mu$ is the \scm for the sequence $(u_k)_k$ at scale
 $(\hk)_k$. Let $(a_k)_k\subset \symbo( \R^{2d})$ be converging in
 $\symbo(\R^{2d})$ to some $a$, and $(b_k)_k, (b'_k)_k\subset L^\infty(\R^d) $ that
 converges uniformly to some $b, b'\in \Con^0(\R^d)\cap L^\infty(\R^d)$ respectively.  Then   
 \begin{align*}
   \lim_{k\to +\infty}
   \inp{ b'_k \OpH (a_k) b_k\,  u_k}{u_k}_{L^2(\R^d)} = \dup{\mu}{b b' a}. 
 \end{align*}
\end{lemma}
\begin{proof}
  One writes
  \begin{align}
    \label{eq: lem continu}
    b'_k \OpH (a_k) b_k &=   (b'_k -b')  \OpH (a_k) b_k
    + b' \OpH (a_k - a) b_k + b' \OpH (a) (b_k -b)\\ &\quad + b' \OpH (a) b.\notag
  \end{align}
Convergence in $\symbo( \R^{2d})$ shows that the operator
 norms $\Norm{\OpH (a) - \OpH (a_k)}{\mathcal{L}(L^2)}$  converge
 to~$0$ uniformly with respect to $h>0$ by
 Lemma~\ref{4.continuite} and $\OpH (a_k)$ is uniformly bounded in
 $\L(L^2)$. With the convergences of $(b_k)_k$, $(b'_k)_k$ one sees
 that the first three terms in \eqref{eq: lem continu} contribute with
 a  vanishing limit because of the $L^2$-boundedness of $(u_k)_k$. It thus suffices to study the limit of $\inp{ b'
   \OpH (a) b \, u_k}{u_k}_{L^2(\R^d)}$. 
One writes 
\begin{align*}
  \inp{b' \OpH (a) b \, u_k}{u_k}_{L^2(\R^d)} 
  &= \inp{b' [\OpH (a),b] u_k}{u_k}_{L^2(\R^d)} 
  + \inp{b' b  \OpH (a) u_k}{u_k}_{L^2(\R^d)}.
\end{align*}
Since $b' b a \in \symbo(\R^{2d})$ and $b' b  \OpH (a) = \OpH ( b' b  a)$ the
result follows from \eqref{WCM}, the $L^2$-boundedness of $(u_k)_k$,
and \eqref{limite-a}.
\end{proof}
 A consequence is the following result.
%%%%%%%%%%%%%%%%%%%%%%%%
% corollary               %
%%%%%%%%%%%%%%%%%%%%%%%%
\begin{corollary}
  \label{cor: localization measure}
  Assume that $\mu$ is the \scm for a sequence $(u_k)_k$ at scale
 $(\hk)_k$ and let
  $\theta \in \Con^0(\R^{2d}) \cap L^\infty(\R^{2d})$.  Then
  $|\theta|^2 \mu$ is the \scm for the sequence $(\theta u_k)_k$ at scale
  $(\hk)_k$.
\end{corollary}

The convergence in \eqref{limite-a} can be extended to more general symbols
and semi-classical operators.
%%%%%%%%%%%%%%%%%%%%%%%%
% proposition          %
%%%%%%%%%%%%%%%%%%%%%%%%
\begin{proposition}
  \label{prop: extension measure}
  Let $(u_k)_k$ be bounded in $L^2(\R^d)$. Suppose $\varphi \in
  \Cinf(\R^d)$ is such that $(\varphi(\hk D) u_k)_k$ is bounded in $L^2(\R^d)$.
  Suppose $\mu$ is the \scm for the sequence $(u_k)_k$ at scale
  $(\hk)_k$ and there is no mass leakage  at infinity at scale $H=(\hk)_k$ in the sense of
  Definition~\ref{def: mass leak at infinity} for the sequences $(u_k)_k$ and
  $(\varphi(\hk D) u_k)_k$. Suppose  $a(x,\xi)$ (or $a(x,\xi')$, that is,
  a tangential symbol), continuous in $x$ and $(d+1)$-times
  differentiable in $\xi$, is such that $\OpH(a)$ is
  bounded on $L^2(\R^d)$.
  Then, one has
  \begin{align*}
    \lim_{k\rightarrow + \infty} \inp{\OpH (a) \varphi(\hk D) u_k}{u_k}_{L^2(\R^d)}
    = \dup{\mu}{a(x,\xi) \varphi(\xi)}.
  \end{align*}
\end{proposition}
A Typical example is $a\in \symbt(\R^{d}\times \R^{d-1})$ by
Lemma~\ref{lemma: continuity tangential op}. Other examples are $a\in
  S^0(\R^d\times\R^d)$ or $a \in
  S^0(\R^d\times\R^{d-1})$,
with  $S^0$ denoting the usual class of symbols of order
$0$; see \cite[Definition 18.1.1]{Hoermander:V3}. The result also
applies to any $a \in S(1,g)$ for any slowly-varying temperate metric $g$ in the
sense of the Weyl-H\"ormander calculus \cite[Section
18.4-18.5]{Hoermander:V3}; such generality is not needed here.
\begin{remark}
  \label{rem: extension measure}
  An inspection of the proof shows that a sharper assumption is
  \begin{equation*}
    %\label{eq: mass leak at infinity}
    \lim_{R\rightarrow +\infty}\, \limsup _{n\rightarrow +\infty} \Big( \int _{|x|\geq R} |u_{k_n}(x) |^2 dx
    +
    \int _{h_{k_n}|\xi|\geq R} |{\widehat{ \varphi(\hk D) u}_{k_n}}(\xi) |^2 d\xi\Big) =0.
\end{equation*}
Note also that Lemma~\ref{lem: continu} also holds in the tangential
case, for instance for $a\in \symbt(\R^{d}\times \R^{d-1})$.
\end{remark}
%%%% proof of proposition
\begin{proof}
  The proof is the same in both cases and is along the line of Proposition 1.6-(iii) in
  \cite{GL:1993}, yet simpler.
 With the no mass-leakage  assumption and since $\OpH(a)$ is bounded on
  $L^2(\R^d)$ one finds
  \begin{align}
    \label{eq: extension measure1}
    \lim_{R\rightarrow +\infty}\, \limsup _{k\rightarrow +\infty}
    \big| \inp{\OpH (a) \varphi(\hk D) u_k}{u_k}_{L^2(\R^d)}
    - \inp{\OpH (a_R) u_k}{u_k}_{L^2(\R^d)}\big|=0,
  \end{align}
  with $a_R (x,\xi)= \chi(x/R) \chi(\xi/R) a(x,\xi) \varphi(\xi)$. Since $a_R\in
  \symbo(\R^{2d})$ with \eqref{limite-a} one has
  \begin{equation*}
    \lim_{k\rightarrow + \infty} \inp{\OpH (a_R) u_k}{u_k}_{L^2(\R^d)}
    = \dup{\mu}{a_R}.
  \end{equation*}
  With \eqref{eq: extension measure1} one  concludes by  means of the
  \Ldct, since $\mu$ has finite mass by \eqref{eq: pas-de-perte}.
\end{proof}

\bigskip We now extend the notion of semi-classical measures to the
case of manifolds. As above $\mathcal N$ is a $\Con^1$-manifold of
dimension $d$ equipped with a density measure $\rho$ that allows one
to define $L^2(\mathcal N)$. For some basic details on density
measures on manifold we refer for instance to \cite[Section 16.2]{LRLR:V2}.

Set $\lambda=\ell^\infty/c_0$ as the space of
bounded sequences modulo the space of sequences converging to $0$. Let
$U = (u_k)_k$ be a bounded sequence in $L^2(\mathcal N)$ and
$H=(\hk)_k$ be a sequence of scales.  For
$ a\in \symbc(T^*\mathcal N)$, denote by
\begin{align*}
  \big[ \big([\OpH](a) u_k,u_k\big)_{L^2(\mathcal N ,\rho)} \big]_\lambda
\end{align*}
the class in $\lambda$ of
the sequence $\big([\OpH](a) u_k,u_k\big)_{L^2(\mathcal N, \rho)}$.

If now $U=(u_k)_k$ is  bounded in $L^2_{\loc}(\mathcal
N)$ it is sensible to compute
\begin{align*}
  \biginp{[\Op^{\hn}](a)
    \psi u_{k_n}}{u_{k_n}}_{L^2(\mathcal N,\rho)}
\end{align*}
for $a\in \symbc(T^* \mathcal N)$ and $\psi \in \Cinfc(\mathcal N)$ with $\psi =1$ on $\supp a$.
%%%%%%%%%%%%%%%%%%%%%%%%
% definition           %
%%%%%%%%%%%%%%%%%%%%%%%%
\begin{definition}
\label{def: mdm on manifolds}
Let $U = (u_k)_k$ be a bounded sequence in $L^2_{\loc}(\mathcal N)$ and
$H=(\hk)_k$ a sequence of scales.  Denote by $\mathcal{M}^+(U)$ the set
of measures $\mu$ on $T^* \mathcal N$ such that there exists a
subsequence $k_n$ such that 
\begin{align*}
  \lim_{n \to +\infty} \big[ \biginp{[\Op^{\hn}](a)
  \psi u_{k_n}}{u_{k_n}}_{L^2(\mathcal N,\rho)} \big]_\lambda
  =\dup{\mu}{a},
\end{align*}
for any $a\in \symbc(T^* \mathcal N)$ and $\psi\in \Cinfc(\mathcal N)$ with $\psi =1$ on $\supp a$.
\end{definition}
What follows explains that this definition is sensible in the sense that it is independent of the choice made for the function $\psi$. In particular, this coincides with the definition of a \scm in the case of a $L^2$-bounded sequence.

If $\mu$ is the \scm associated with $U = (u_k)_k$, then in any local
chart $\chart=(\O,\cdiff)$, denote by $\mu^\chart$ the local representative of $\mu$, that
is, $(\chdiff^{-1})^* \mu$.  Denote also by $u_k^\chart$ the local
representative of $u_k$, that is, $u_k^\chart = (\chdiff^{-1})^* u_k =
u_k \circ \chdiff^{-1}$. Then, if $K \subset \cdiff(\O)$ is compact, $a \in \symbo(\R^{2d})$ with $\supp a
\subset K\times\R^d$, and $\psi \in \Cinfc\big( \cdiff(\O)\big)$ equal
to $1$ in a \nhd of the $x$-projection of $\supp a$ one has 
\begin{align}
  \label{eq: measure local chart}
  \lim_{k\to +\infty} \inp{\OpH (a) \psi
  u^\chart_k}{u^\chart_k}_{L^2(\R^d, \rho)}
    = \dup{\mu^\chart}{a}.
\end{align}
In what follows $\rho$ will be given by $\k \mug$,  that is, in local coordinates $\rho^\chart = \k^\chart \det
\big(g^\chart\big)^{1/2} d x$. One then has   
\begin{align*}
  \mu^\chart = \k^\chart \det \big(g^\chart\big)^{1/2} m, 
\end{align*}
if $m$ is the
\scm associated with $(u_k^\chart)_k$ yet using the $L^2$-inner product given by the Lebesgue measure
as in Definition~\ref{def: semi-classical measure}, that is, 
\begin{align*}
  \lim_{k\to +\infty} \inp{\OpH (a) \psi
  u^\chart_k}{u^\chart_k}_{L^2(\R^d, dx)}
    = \dup{m}{a}.
\end{align*}
Some of the properties of \scm on $\R^d$ can then be extended to the
case on a manifold. For instance one has the following result. 
%%%%%%%%%%%%%%%%%%%%%%%%
% lemma                %
%%%%%%%%%%%%%%%%%%%%%%%%
 \begin{lemma}
   \label{lem: continu - manifold}
 Assume that $\mu$ is the \scm of a sequence $U=(u_k)_k$ at scale
 $(\hk)_k$ on $\mathcal N$. Let $(a_k)_k\subset \symbc( T^* \mathcal N)$ be converging in
 $\symbc(T^* \mathcal N)$ to some $a$, and $(b_k)_k,
(b'_k)_k\subset L^\infty(\mathcal N) $ that
 converges uniformly to some $b, b'\in \Con^0(\mathcal N)\cap L^\infty(\mathcal N)$ respectively.  Then   
 \begin{align*}
   \lim_{k\to +\infty}
   \big[ \inp{ b'_k [\OpH] (a_k) b_k\,  u_k}{u_k}_{L^2(\mathcal N)}\big]_\lambda
   = \dup{\mu}{b b' a}. 
 \end{align*}
\end{lemma}
The local chart version is
\begin{align}
  \label{eq: lem: continu - manifold - local chart}
   \lim_{k\to +\infty}
   \inp{ b'_k \OpH (a_k) b_k\,  \psi u_k^\chart
  }{u_k^\chart}_{L^2(\R^d, \rho)}
   = \dup{\mu^\chart}{b b' a}, 
 \end{align}
 for $a$ and $\psi$ as given for \eqref{eq: measure local chart} and
 $a_k, b_k, b'_k$ also defined locally accordingly.

\medskip
The following result that yields the existence of \scm. 
%%%%%%%%%%%%%%%%%%%%%%%%
% proposition          %
%%%%%%%%%%%%%%%%%%%%%%%%
\begin{proposition}
  \label{prop: existence measure L2 loc}
  Suppose $H=(\hk)_k$ is a sequence of scales and $U=(u_k)_k$ a
    sequence of functions on $\mathcal N$.
  \begin{enumerate}
    \item 
  If $U=(u_k)_k$ is bounded in $L^2(\mathcal N)$ the set $
  \mathcal{M}^+(U)$ is nonempty.
  \item 
  Suppose $\mathcal N$ is countable at infinity. If $U=(u_k)_k$ is bounded in $L^2_{\loc}(\mathcal N)$ the set $ \mathcal{M}^+(U)$ is nonempty.
  \end{enumerate}
\end{proposition}
%%%% proof of proposition
\begin{proof}
  The result of the first part, that is, if $U=(u_k)_k$ is bounded in
  $L^2(\mathcal N)$, holds by \cite[Section 1]{GL:1993}.
  
  For the second part, as $\mathcal N$ is countable at infinity, there
  exists a sequence of open sets $(\O_n)_n$ with $\O_n \Subset
  \O_{n+1} \Subset \mathcal N$ and $\cup_n \O_n = \mathcal N$.  The
  sequence $(u_k)_k$ is $L^2$-bounded on $\O_1$.  Suppose $a \in
  \symbc(T^* \mathcal N)$ supported in $\O_1$ and $\psi \in\in
  \Cinfc(\mathcal N)$ with $\psi =1$ on $\O_1$. One has
  \begin{align*}
    \big[ \biginp{[\Op^{\hn}](a)
        \psi u_{k_n}}{u_{k_n}}_{L^2(\mathcal N,\rho)} \big]_\lambda
    =\big[ \biginp{[\Op^{\hn}](a)
        \psi u_{k_n}}{\psi u_{k_n}}_{L^2(\mathcal N,\rho)} \big]_\lambda.
\end{align*}
The sequence $(\psi u_k)_k$ is $L^2$-bounded.  By the first part,
there exists an inscreasing function $\varphi_1: \N \to \N$ and a
measure $\mu_1$ on $T^*(\O_1)$ that is the \scm for the subsequence
$(\psi u_{\varphi_1(k)})_k = (u_{\varphi_1(k)})_k$ on $\O_1$.  With
the same reasoning there exists an inscreasing function $\psi_2 : \N
\to \N$ and measure $\mu_2$ on $T^*(\O_2)$ that is the \scm for the
subsequence $(u_{\varphi_2(k)})_k$ on $\O_2$, with $\varphi_2 = \psi_2
\circ \varphi_1$. One has $\mu_2 = \mu_1$ on $T^*(\O_1)$.  One
proceeds by induction yielding two sequences of inscreasing functions
$\varphi_n: \N \to \N$ and $\psi_n: \N \to \N$, with $\varphi_{n+1} =
\psi_{n+1}\circ \varphi_n$, and a sequence of measures $\mu_n$ on
$T^*(\O_n)$, with $\mu_n$ the \scm of $(u_{\varphi_n(k)})_k$ on
$\O_n$. Moreover,for $\ell \in \N$, one has $\mu_{n + \ell} = \mu_n$
on $T^*(\O_n)$.

  There exists a unique measure $\mu$ on $T^*(\mathcal N)$ such that $\mu= \mu_n$ on $T^*(\O_n)$.
  A diagonal extraction yields the subsequence $(u_{\varphi_k(k)})_k$ implying that $\mu_n$ is its \scm on $\O_n$ for any $n \in \N$. Hence, $\mu$ is its   \scm on $\mathcal N$. 
\end{proof}

The notion of \scm can be extended to vector valued sequences.
If $N\in \N^*$, denote by ${\mathcal M}(T^*\mathcal N; \MM_N(\C))$ the space of $N\times N$-matrix
valued Radon measures on $T^* \mathcal N$, and by ${\mathcal M}^+(T^*\mathcal N; \MM_N(\C))$ the subspace fromed by nonnegative Hermitian such measures.
%%%%%%%%%%%%%%%%%%%%%%%%
% definition           %
%%%%%%%%%%%%%%%%%%%%%%%%
\begin{definition}[Hermitian measures]
\label{def: hermitian mdm on manifolds}
Suppose $N \in \N^*$ and $U=(u_k)_k$ is a bounded sequence in $L^2_{\loc}(\mathcal N; \C^N)$ and
$H=(\hk)_k$ a sequence of scales.  Denote by $\mathcal{M}^+(U)$ the set
of measures $\mu \in {\mathcal M}^+(T^*\mathcal N; \MM_N(\C))$ such that there exists a
subsequence $k_n$ such that 
\begin{align*}
  \lim_{n \to +\infty} \big[ \biginp{[\Op^{\hn}](a)
  \psi u_{k_n}}{u_{k_n}}_{L^2(\mathcal N,\rho)} \big]_\lambda
  =  \dup{\trace (a \mu)}{1} = \int_{T^* \mathcal N} \trace (a(x,\xi) d \mu(x,\xi)),
\end{align*}
for any $N\times N$ matrix $a$ with entries in $\symbc(T^* \mathcal
N)$, and $[\Op^{\hn}](a)$ the associated class of matrix valued
operators, and $\psi\in \Cinfc(\mathcal N)$ with $\psi =1$ on $\supp
a$.
\end{definition}
We refer the reader to \cite{Gerard:91,Burq-Bourbaki}.
Each element of the matrix valued measure can also be understood as follows:
\begin{align*}
  \lim_{n \to +\infty} \big[ \biginp{[\Op^{\hn}](a)
  \psi u_{i, k_n}}{u_{j, k_n}}_{L^2(\mathcal N,\rho)} \big]_\lambda
  =  \dup{\mu_{ij}}{a},
  \qquad a \in \symbc(T^* \mathcal N).
\end{align*}
Each diagonal term is nonnegative.  One finds that
\begin{align*}
  %\label{eq: bound off-diagonal terms measure}
  \mu_{ij} \leq \mu_{ii}^{1/2} \mu_{jj}^{1/2}, 
\end{align*}
in the sense that $|\dup{\mu_{ij}}{ a b}|^2 \leq \dup{\mu_{ii}}{|a|^2}\dup{\mu_{ii}}{|b|^2}$.

\medskip
The counterpart to Proposition~\ref{prop: existence measure L2 loc} is the following result.
%%%%%%%%%%%%%%%%%%%%%%%%
% proposition          %
%%%%%%%%%%%%%%%%%%%%%%%%
\begin{proposition}
  \label{prop: existence measure L2 loc vector}
  Suppose $N \in \N^*$ and $H=(\hk)_k$ is a sequence of scales and $U=(u_k)_k$ a
    sequence of function on $\mathcal N$ valued in $\C^N$.
  \begin{enumerate}
    \item 
  If $U=(u_k)_k$ is bounded in $L^2(\mathcal N; \C^N)$ the set $
  \mathcal{M}^+(U)$ is nonempty.
  \item 
  Suppose $\mathcal N$ is countable at infinity. If $U=(u_k)_k$ is bounded in $L^2_{\loc}(\mathcal N; \C^N)$ the set $ \mathcal{M}^+(U)$ is nonempty.
  \end{enumerate}
\end{proposition}

%%%%%%%%%%%%%%%%%
% section
%%%%%%%%%%%%%%%%%

%%%%%%%%%%%%%%%%%
% subsection
\section{The measure propagation equation and proof of observability}

We first state a result that is at the heart of the proof of
Theorems~\ref{theoprinc} and \ref{theo-perturbation}. It expresses how
a \scm $\mu$ associated with solutions to wave equations varies in the direction of the hamiltonian vector field
$\Hp$, in particular at the boundary $\d\L= \R \times \d\M$ where this
variation is connected to a \scm $\nu$ associated with the Neumann
trace. 

\subsection{The measure propagation equation}
\label{sec: The measure propagation equation}
Suppose $(\M, \k,
g) \in\X^1$, $H =(\hk)$ is a scale, and 
$(\kk, \gk)_k$ such that 
$(\M, \kk, \gk)$  converges to  $(\M, \k,
g)$ in the $\Y^1$  topology.

 Suppose $(u_k)_k$ are weak-solutions to
   \begin{align*}
     \d_t^2 u_k  - A_{\kk, \gk}u_k =f_k,
   \end{align*}
   with homogeneous Dirichlet boundary condition, as given in Proposition~\ref{prop: well-posedness nonhomogeneous wave equation}.  
   Extend $u_k$ and $f_k$  by zero to $\hL$. 

   Suppose $(u_k)_k$ is bounded in $L^2_{\loc}(\hL)$, $(\hk
   \d_{\n_k}{u_k}_{|\d\L})_k$ is bounded in $L^2_{\loc}(\d\L)$, and
   $(\hk f_k)_k$ is bounded in $L^2_{\loc}(\hL)$.
   With
   Proposition~\ref{prop: existence measure L2 loc vector}, a Hermitian $2 \times 2$ \scm $M$ on
   $T^*(\hL)$ is associated with a subsequence at scale $H$ of $(u_k, \hk f_k)$.
   Write 
   \begin{align*}
     M = \begin{pmatrix}
       M_{0,0} & M_{0,1}\\
       M_{1,0} & M_{1,1}
     \end{pmatrix}.
   \end{align*}
   Set $\mu = M_{0,0}$.
   Similarly, with Proposition~\ref{prop: existence measure L2 loc},
   there exists a measure $\nu$ on $T^*\d\L$ such that the \scm
   measure associated with (a subsequence of) $\hk \psi(t)\d_{\n_k}{u_k}_{|\d\L}$ is
   $|\psi(t)|^2 \nu$ at scale $H$.
%%%%%%%%%%%%%%%%%%%%%%%%
% theorem              %
%%%%%%%%%%%%%%%%%%%%%%%%
\begin{theorem}
  \label{thm: equationmesure}
  Suppose that 
   \begin{align}
     \label{eq: hyp support mesure}
     \supp \mu \subset \Char p \cap \TL \setminus 0
     \  \ \et \ \
     \supp \nu \subset T^*\d\L  \setminus 0.
     \end{align}
 Then, the two measures $\mu$ and $\nu$ fulfill, in sense of density distributions, 
\begin{align}
    \label{eq: GL equation-propagation theorem}
    \Hp \mu  = - \transp{\Hp} \mu = 2 \Im M_{0,1} + \int_{\y \in \pHb \cup \pGb} 
    \frac{\delta_{\y^+} - \delta_{\y^-}}
    {\dup{\xi^+- \xi^-}{\nx}_{T_x^*\M, T_x\M}} 
  \ d \nu (\y).
\end{align}
The hyperbolic set $\pHb$ and the glancing set $\pGb$ are introduced
in Definition~\ref{def: E', H', G'-intro} and 
$\y^\pm$ and $\xi^\pm$ are as given in 
   \eqref{eq: relevement H-G-intro}. 
The vector field $\nx$   is the unitary inward
    pointing normal vector in the sense of the metric $g$.
\end{theorem}

Here, $p = p_{\k,g}$ and thus $\Hp = \Hamiltonian_{p_{\k,g}}$, the
sets $\pHb$ and $\pGb$ are constructed with respect to the metric $g$
as in Section~\ref{sec: A partition of the cotangent bundle at the boundary-intro}. 
Recall that we identify $T^* \d\L$ and $\pTL$.
Hence, the measure $\nu$ defined on  $T^* \d\L$ is also a measure
on  $\pTL$. The integral performed on $\pHb \cup \pGb$ thus makes
sense. Also, the meaning of the right hand side is explained in
Remark~\ref{integrand}.

Sections~\ref{sec: Proof of the propagation equation1} to \ref{sec:
  Proof of the propagation equation2} are dedicated to the proof of
Theorem~\ref{thm: equationmesure}.  The result of Theorem~\ref{thm:
  equationmesure} is key in the proof the main observability result as
presented in the next section. There, one only considers the case
$f_k=0$ implying $M_{0,1}=0$. The addition of the source term $f_k$
does not provide any complication for the proof Theorem~\ref{thm:
  equationmesure}, hence this slight generalization that can be of use
elsewhere, in particular for the study of stabilization issues.
 
\subsection{Proof of the observability results}
\label{sec: Proof of the observability results}

Here, we provide  the proof of
Theorem~\ref{theo-perturbation} based on the measure equation of
Theoerm~\ref{thm: equationmesure}.
Suppose $(\M, \k, g) \in \X^1$ and $\omega$ is an open subset of
$\M$ (\resp $\Gamma$ an open subset of $\d\M$) 
such that the interior geometric control condition of
Definition~\ref{control-geo-interior} (\resp the boundary geometric
control condition of Definition~\ref{control-geo-boundary}) is
fulfilled and  we consider  some time $T>T_{GCC}(\omega)$ (\resp
$T>T_{GCC}(\Gamma)$). We also consider $\delta>0$ such that $T-2\delta
>T_{GCC}(\omega)$ (\resp
$T-2\delta>T_{GCC}(\Gamma)$).

According to
Propositions~\ref{3.inegalitesemiclassint} (\resp
Proposition~\ref{3.inegalitesemiclassbord}),
to achieve the observability inequalities of
Theorem~\ref{theo-perturbation} for the time interval $]0,T[$
it suffices to prove the
semi-classical observability inequality~\eqref{3.ineg11int}
(\resp~\eqref{3.ineg11bord}) for the time interval $I=]\delta,T-\delta[$ for any $(\tilde{\M}, \tilde{\k}, \tilde{g}, \tilde{\omega})$ (\resp
$(\tilde{\M}, \tilde{\k}, \tilde{g}, \tilde{\Gamma})$) that is $\eps$-close to
$(\M, \k, g, \omega)$ (\resp $(\M, \k, g, \Gamma)$) in
the $\Y^1$-topology in the sense of
Definition~\ref{def: eps-close} and $\eps>0$ chosen \suff small. We preform a contradiction argument based on propagation properties
of semi-classical defect measures.

Below we consider a sequence
$(\M_n, \kn, \gn)$  that converges to  $(\M, \k,
g)$ in the $\Y^1$  topology. For $k \in \ZZ^*$, one denotes by $E_{n, k}$,
the space of functions defined in Section~\ref{dyadic},  here
built on the elliptic operator $A_{\kn, \gn}$ on $\M_n$.

\subsubsection{Initiation of the contradiction argument}
\label{sec: Initiation of the contradiction argument}
In the case of an interior observation, we assume that \eqref{3.ineg11int} does not hold for some 
$(\tilde{\M}, \tilde{\k}, \tilde{g}, \tilde{\omega})$ arbitrary close
to $(\M, \k, g, \omega)$ in the sense recalled above.
Thus, there exists a
sequence $(\M_n, \kn, \gn, \omega_n)_n$ that converges to  $(\M, \k,
g, \omega)$ in the $\Y^1$  topology,
and for each $n\in \N$ and each  $k\in \N^*$ there exists $\ell(n,k) \in
\N$, with $\ell(n,k)\geq k$ and $u_n^{\ell(n,k)} \in E_{n, \ell(n,k)}$,  such that 
\begin{align}
\label{absurde-int}
1
= \Norm{{u_n^{\ell(n,k)}}\bt }{L^2(\M_n, \kn \mugn)} 
\geq  k \, \Norm{\bld{1}_{I \times \omega_n}   h_{\ell(n,k)} \d_t
  u_n^{\ell(n,k)}}{L^2(\R\times\M_n, \kn \mugn dt )},
\end{align}
with   $I = ]\delta, T-\delta[$. 
Note that we have  normalized the \lhs of~\eqref{absurde-int} to be equal to $1$.
The notation $u_n^{\ell(n,k)}$ may seem very cumbersome at this
stage; it will be greatly simplified by a diagonal extraction in what
follows shortly.

\medskip
Similarly, 
in the case of a boundary
observation we assume that there exists a
sequence $(\M_n, \kn, \gn, \Gamma_n)$ that converges  to $(\M, \k, g,
\Gamma)$  in the $\Y^1$  topology,
and for each $n\in \N$ a sequence $ (u_n^{\ell(n,k)})_{k\in \N}$, with
$u_n^{\ell(n,k)} \in E_{n, \ell(n,k)}$ and $\ell(n,k)\geq k$, such that 
\begin{equation}
  \label{absurde-boundary} 
  1= \Norm{{u_n^{\ell(n,k)}}\bt }{L^2(\M_n, \kn \mugn)} 
  \geq  k \, \Norm{\bld{1}_{I \times \Gamma_n} 
  h_{\ell(n,k)} \d_{\n_n}  u_n^{\ell(n,k)}}{L^2(\R\times \d\M_n, \kn \mugbn dt)}, 
\end{equation} 
where  $\n_n$ is the normal to the boundary $\d\M_n$ in the sense  of
the metric $\gn$.

\bigskip
We now proceed with a diagonal extraction along with a zero-extension of
the solutions outside $\L$. 
Set $u_k= u_k^{\ell(k,k)} 1_{\L}$, that is,  the extension by $0$ of the function
$u_k^{\ell(k,k)}$ to $\hL = \R\times \hM$ (see  Section~\ref{geometry}) and $v_k = \hk
\d_{\n_k} {u_k}_{|\d\L}$ its normal partial derivative (in the sense
of $\gk$).
In what follows we will denote $h_{\ell(k,k)}$ and $J_{\ell(k,k)}$ by
$\hk$ and $J_k$ for simpicity. Yet, there will be no possible
confusion.

\medskip
First, with the $W^{2,\infty}$-diffeomorphism of  Definition~\ref{def:
  eps-close} the analysis can be pulled back from $\M_k$ to $\M$ for
each $k \in \N$. By abuse of notation we use the same letters for the
pullbacked functions and metric. Hence, without loss of generality we
may assume that $\M_k = \M$. 

Second, observe that since $ \Norm{{u_k}\bt}{L^2(\M, \kk
  \mugk)} =1$ one has 
\begin{equation*}
  %\label{mass}
  \Norm{{u_k}\bt}{L^2(\M)}
  = 1 + o(1)_{k \to \infty}.
\end{equation*}
If no precision is given, the $L^2$-norm on $\M$ is given by the
density measure $\k \mug$ in what follows. 

 From Lemma~\ref{lemma: norm dt uk} and Remark~\ref{rem: equiv norm energy} one obtains that 
\begin{align}
  \label{eq: mass2}
  1 &\eqsim \Norm{u_k(t,.)}{L^2(\M, \kk \mugk)}  \eqsim  
   \Norm{\hk \d_t u_k(t,.)}{L^2(\M)}  
  \eqsim   \Norm{\hk \nabla_{\!\! \gk} u_k(t,.)}{L^2(\M)} \\
  &\eqsim  \Norm{\hk^2 A_{\kk , \gk} u_k(t,.)}{L^2(\M)}, \notag
\end{align}
for any $t \in \R$ and $k$ large. 
From ellipticity up to the boundary  one deduces~\cite[Theorem 8.12]{GT:77}\footnote{Notice that in~\cite{GT:77} the boundary is assumed $\Con^2$. Still, 
 $W^{1,\infty}$-regularity suffices to reach the conclusion since it is enough to
 make the boundary straight in local coordinates and apply the
 argumentation therein.} 
 \begin{equation}
   \label{eq: mass3}
  \Norm{\hk^2 u_k(t,.)}{H^2(\M)} \eqsim  1, 
 \end{equation}
for any $t \in \R$.

\subsubsection{Measures for the wave equations}
 From Proposition~\ref{prop: admissibility} in the case $f=0$,
 \eqref{eq: mass2}, and \eqref{absurde-int} and \eqref{absurde-boundary}
 (and the fact that $\kn$ converges to $\k$ and $\gn$ to $g$
 in the sense given in Definition~\ref{def: eps-close})
 one obtains the following proposition.
%%%%%%%%%%%%%%%%%%%%%%%%
% proposition          %
%%%%%%%%%%%%%%%%%%%%%%%%
\begin{proposition}
\label{prop: sequence waves}
The sequences $u_k  \in L^\infty(\R; L^2(\hM))$ and  $v_k \in L^2_{loc}(
\d\L)$ satisfy 
\begin{enumerate}
\item For any bounded  interval $\interval \subset \R$ there exists $C>0$ such that 
\begin{equation*}%\label{massbis}
  \Norm{u_k}{L^2(\interval \times \hM)} + \Norm{v_k}{L^2(\interval \times \d\M)} \leq C.
  \end{equation*}
\item With $I = ]\delta, T-\delta[$,  one has
  \begin{itemize}
  \item $\lim_{k\to + \infty } \Norm{u_k}{L^2(I \times \omega)}=0$, if
    the case \eqref{absurde-int} holds, that is, for interior observability,
  \item $\lim_{k\to + \infty } \Norm{v_k}{L^2(I \times \Gamma)} =0$,
    if the case \eqref{absurde-boundary}  holds, that is,  for boundary observability.
  \end{itemize}
\end{enumerate}
\end{proposition}

Recall that we consider on $\d\M$ the density measure $\k \mugb$.
From Propositions~\ref{prop: sequence waves} and \ref{prop: existence
  measure L2 loc} we deduce the following result.
%%%%%%%%%%%%%%%%%%%%%%%%
% proposition          %
%%%%%%%%%%%%%%%%%%%%%%%%
\begin{proposition}
  \label{prop: existence measures wave}
  \begin{enumerate}
  \item 
  There exists a semi-classical measure $\mu$ on $\ThL$ associated with a
  subsequence of $(u_k)_k$. 
  \item There exists a semi-classical measure $\nu$ on $T^* \d\L$
    associated with a subsequence of $(v_k)_k$.
\end{enumerate}
\end{proposition}

By abuse of notation we will use the notation $(u_k)_k$ and
$(v_k)_k$ for both subsequences. Then one has 
\begin{align*}
  &\dup{\mu}{a} 
    = \lim_{k\to +\infty} \big[\inp{[\OpH](a) u_k}{u_k}_{L^2(\hL)} \big]_\lambda, 
    \qquad a \in \symbc(T^*\hL),\\
  &\dup{\nu}{b} 
    = \lim_{k\to +\infty} \big[\inp{[\OpH](b) v_k }{v_k}_{L^2(\d\L)}\big]_\lambda, 
    \qquad b\in \symbc(T^* \d\L),
\end{align*}
where both limits are understood in the sense given in
Definition~\ref{def: mdm on manifolds}. Recall that the spaces of
symbols $\symbc(T^*\hL)$ and $\symbc(T^* \d\L)$ are introduced in Definition~\ref{def: symbol manifold}.
%%%%%%%%%%%%%%%%%%%%%%%%
% proposition          %
%%%%%%%%%%%%%%%%%%%%%%%%
\begin{proposition}
\label{prop: first property measures}
The three following properties hold. 
\begin{enumerate}
\item \label{prop: first property measures 1}  
  If $\interval \subset \R$ is a bounded nonempty open interval,  one
  has $\mu \big(  T^* (\interval   \times \hM) \big) >0$. 
\item One has 
  \begin{align}
    \label{eq: prop: first property measures 2}
    &\supp \mu \subset \Char p \cap \TL \cap \{ \alpha\leq
    \tau\leq  \alpha^{-1}\},\\
     \label{eq: prop: first property measures 3}
  &\supp \nu \subset T^*\d\L \cap \{ \alpha\leq
    \tau\leq  \alpha^{-1}\}.
  \end{align}
\item
    With $I = ]\delta,
T-\delta[$ as above one has
  \begin{itemize}
    \item 
the measure $\mu$ vanishes on
$T^* (I \times \omega)$, in the case of an interior observation,
\item the
measure $\nu$ vanishes on $T^* (I\times \Gamma)$, in the case of a
boundary observation.
\end{itemize}
\end{enumerate}
\end{proposition}
\begin{proof}
  Consider a finite $\Con^2$-partition of unity
$(\chi_i)_{i\in \I}$ subordinated to a given atlas of $\M$; see Section~\ref{geometry}. 
 Let $\varphi \in \Cinfc ( \R)$ be nonvanishing. From \eqref{eq: mass2} 
 one has 
 \begin{align*}
   \Norm{\varphi(t) \chi_i u_k}{L^2(\hL)} \gtrsim 1,
   \quad \text{for some}\ i \in \I.
 \end{align*}
 The semi-classical measure associated with $(\varphi(t) \chi_i
 u_k)_k$ is $|\varphi(t) \chi_i|^2 \mu$ by Lemma~\ref{lem: continu -
   manifold} for the $L^2$-inner product associated with the density
 measure $\k \mug dt$.  In a local chart $\chart=(\O,\cdiffL)$ of
 $\hL$, where $\O = \R \times\hO$ and $\cdiffL(t,x) = (t, \chdiff(x))$
 (see Section~\ref{sec: Local coordinates}), with $\supp \chi_i\subset
 \hO$, from \eqref{eq: mass2} and Lemma~\ref{lemma: h-oscillation} one 
 has  
\begin{align*}
  \bigNorm{\bld{1}_{\hk|(\tau, \xi)|\geq R}\,
  \widehat{ \varphi(t) \chi_i u_k}}{L^2(\R^{d+1})} \lesssim R^{-1},
\end{align*}
There is thus no mass leakage  at infinity at scale $H$ in the
sense of Definition~\ref{def: mass leak at infinity}.

Denote by $\mu^\chart$ the  local representative of  $\mu$. Recall
that $\tk = \k (\det g)^{1/2}$. 
Using that the \scm
of the local representative of  $\varphi(t) \chi_i u_k$
with a $L^2$-inner product associated with the Lebesgue measure $dx$ is
$m = \tk^{-1}|\varphi(t) \chi_i|^2 \mu^\chart$, with Proposition~\ref{prop: properties scm} one finds
\begin{align*}
  m (\R^{2 d+2}) = \lim_{n \to +\infty} \Norm{\varphi(t) \chi_i u_k}{L^2(\R^d)}^2 \gtrsim \lim_{n \to +\infty} \Norm{\varphi(t) \chi_i u_k}{L^2(\hL)}^2 \gtrsim 1,
\end{align*}
hence the first result. 

\medskip
We place ourselves in a local chart $(\O,\cdiffL)$ of $\hL$. Here, $\d\L$
is given by $\{z=0\}$. 
Let $b\in \Cinfc(\R^{2d+2})$ with $\supp b \subset \cdiffL(\O) \times
\R^{d+1}$ and $\psi, \tpsi \in \Cinfc \big(\cdiffL(\O)\big)$ with
$\psi$ equal to $1$
in a \nhd of the $(t,x)$-projection of $\supp b$ and $\tpsi$ equal to 1
in a \nhd  of $\supp \psi$. Here, $\OpH(b) = b(t,x,\hk
D_t, \hk D_x)$.
One has, for any $s\geq 0$,
\begin{equation}
  \label{eq: bornitude  H-s L2}
\Norm{\OpH (b)}{\L(H^{-s}(\R^{d+1}), L^2(\R^{d+1}))}\leq C_s h^{-s} .
\end{equation}
In fact, one uses that $\OpH(b) \OpH(\jp{\xi}^s) = \OpH(b \jp{\xi}^s) $ is bounded on $L^2$ and
$h^s \jp{\xi}^s \leq \jp{h \xi}^s$
yielding
\begin{align*}
  h^s \Norm{\OpH (b) u }{L^2(\R^{d+1})}
  &\lesssim h^s \Norm{\OpH (\jp{\xi}^{-s}) u }{L^2(\R^{d+1})}
  =  h^s \Norm{\jp{h \xi}^{-s} \hat{u} }{L^2(\R^{d+1})}\\
  &\lesssim \Norm{\jp{\xi}^{-s} \hat{u} }{L^2(\R^{d+1})}
  = \Norm{u }{H^{-s}(\R^{d+1})}.
\end{align*}

One has
\begin{align}
  \label{eq: local wave eq}
 \hk^2 P_{\kk, \gk} u_k= - \hk v_k \otimes\delta_{z=0}.
\end{align}
Note that $v_k \otimes
\delta_{z=0}$ is bounded in $H^{-s}(\R_z; L^2_{\loc}(\R^{d}_{y'})) \subset  H_{\loc}^{-s}(\R^{1+d})$ if $s >
1/2$. Thus, with \eqref{eq: bornitude  H-s L2} one finds
$\Norm{\OpH(b) \psi \hk v_k \otimes \delta_{z=0}}{L^2} \leq \hk^{1-s}$.
With $1/2< s < 1$, with \eqref{eq: local wave eq} and since $(u_k)_k$ is bounded in $L^2$, one concludes that
\begin{align*}
  &\lim_{k\to +\infty}
  \hk^2 \biginp{\OpH(b)  \psi P_{\kk, \gk} u_k }
  { u_k}_{L^2(\R^{d+1}, \k \mug dt) }\\
  & \qquad=- \lim_{k\to +\infty} \biginp{\OpH(b)  \psi \hk v_k \otimes
  \delta_{z=0}}{ u_k}_{L^2(\R^{d+1}, \k \mug dt) } =0.
\end{align*}
Now, one has
\begin{align*}
  %\label{develop-1}
  \OpH(b) \psi \hk^2 \d_t^2
  &=- \OpH( \tau^2 b)  \psi
  + \hk^2 \OpH(b)[\psi,  \d^2_t]\\
  &= - \OpH( \tau^2 b)  \psi + \hk^2 O(1)_{\L(H^1,L^2)}.\notag
\end{align*}
With the form of $A_{\kk,\gk}$ in local coordinates one writes, with
 $\tkk=\kk  (\det \gk)^{1/2}$,
\begin{align*}
  %\label{develop-2}
  &\OpH(b) \psi \hk^2 A_{\kk,\gk}\\
  &\qquad = \hk^2 \OpH(b)   \sum_{1\leq i, j \leq d}  \Big( 
    \d_{x_i} \psi \gk^{i j}(x) \d_{x_j} 
  +   [\psi \tkk^{-1}, \d_{x_i} ] 
    \tkk \gk^{i j}(x)\d_{x_j} \Big) 
  \notag\\
  &\qquad  =\sum_{1\leq i, j \leq d}  i \OpH(b\xi_i) \psi \gk^{i j}(x)
    \hk \d_{x_j}  \tpsi
    +  \hk^2\,  O(1)_{\L(H^1,L^2)} \notag\\
  &\qquad  =\sum_{1\leq i, j \leq d}
    \Big( - \psi \gk^{i j}(x)  \OpH(b\xi_i \xi_j ) 
    +i [\OpH(b\xi_i), \psi \gk^{i j}(x)] \hk \d_{x_j}  \Big)  \tpsi\notag\\
  &\qquad \quad   +  \hk^2\,  O(1)_{\L(H^1,L^2)} \notag\\
  &\qquad  =- \sum_{1\leq i, j \leq d}
    \psi \gk^{i j}(x)  \OpH(b\xi_i \xi_j ) \tpsi
    +  \hk^2\,  O(1)_{\L(H^1,L^2)},\notag
\end{align*}
 where in the last line we used Proposition~\ref{propcommut}.
From  Lemma~\ref{lem: continu - manifold} and \eqref{eq: lem: continu - manifold - local chart} one finds
 \begin{equation*}
  0= \lim_{k\to +\infty}
  \hk^2 \biginp{\OpH(b)  \psi P_{\kk, \gk} u_k }{ u_k}_{L^2(\hL, \k \mug dt) }
  =  \dup{\mu^\chart}{ b\,  p_{\k,g}}.
\end{equation*}
  Since  $b$ is arbitrary in $\Cinfc(\R^{2d+2})$ this implies that
  $\supp \mu^\chart \subset \Char \big( p_{\k,g}\big)$, which is the
  first inclusion in \eqref{eq: prop: first property measures 2}. 

  \medskip
To prove
  the second property, that is $\tau \in [ \alpha, \alpha^{-1}]$,
  consider $\varphi \in \Cinfc(\R)$ such that $ \varphi \equiv
0$ in a \nhd of $ [\alpha,\alpha^{-1}]$,  say
  $[(1-\eps)\alpha, (1+\eps) \alpha^{-1} ]$ for some $\eps>0$ to be
  kept fixed, and now $\psi \in
\Cinfc(\R)$. We write 
\begin{equation*}
\varphi(\hk D_t) \psi (t) u_k = \sum _{\nu \in J_k} \varphi(\hk D_t) 
\big( \psi (t) e^{it\sqrt{\lambda_\nu}} \big)  u_\nu e_\nu(x).
\end{equation*}
The Fourier transform of $\varphi(\hk D_t)  \big( \psi (t)
e^{it\sqrt{\lambda_\nu}} \big)$ is 
$\varphi(\hk \tau) \hat{\psi}\big(\tau - \sqrt{\lambda_\nu}\big)$.
As $\hk\sqrt{\lambda_\nu}  \in [\alpha, \alpha^{-1}]$ if $\nu \in
  J_k$ and $\hk \tau \notin [(1-\eps)\alpha, (1+\eps) \alpha^{-1}]$ if
  in the support of $\varphi$, then 
\begin{align*}
  %\label{eq: est diff tau sqrt lambda}
  |\tau - \sqrt{\lambda_\nu}| \gtrsim |\tau| + \hk^{-1},
\end{align*}
in the support of the above Fourier transform. 
With the rapid decay of $\hat{\psi}$
one finds, for any $N \in \N$, 
\begin{align*}
  %\label{eq: decay Fourier localization1}
  \big| \varphi(\hk \tau) \hat{\psi}\big(\tau -
  \sqrt{\lambda_\nu}\big)\big| 
  \leq C_N (|\tau|^{-1} + \hk)^{N},
\end{align*}
leading to 
\begin{align}
  \label{eq: decay Fourier localization2}
 \Norm{ \varphi(\hk D_t)  \big( \psi (t)
e^{it\sqrt{\lambda_\nu}} \big)}{H^\ell(\R)} \leq C_N' \hk^{N}
\end{align}
for any $\ell\geq 0$. With $\ell=0$, one deduces
\begin{align*}
&\Norm{\varphi(\hk D_t) \psi (t) u_k}{L^2(\L, \k \mug  dt)}^2\\
&\qquad \eqsim \Norm{\varphi(\hk D_t) \psi (t) u_k}{L^2(\L, \kk
  \mugk  dt)}^2 
=  \sum_{\nu \in J_k } |u_\nu|^2 \Norm{ \varphi(\hk D_t)  \big( \psi (t)
e^{it\sqrt{\lambda_\nu}} \big)}{L^2(\R)}^2 \\
 &\qquad \leq C_N \hk^N \sum_{\nu \in J_k } |u_\nu|^2 
  = C_N
  \hk^N\Norm{{u_k}\bt}{L^2(\M, \kk \mugk)}^2
   \lesssim C_N \hk^N,
\end{align*}
for any $N \in \N$,
using \eqref{absurde-int} or \eqref{absurde-boundary}, 
implying
$\dup{\mu}{|\varphi(\tau)\psi(t)|^2} = 0$,
which gives the last inclusion in~\eqref{eq: prop: first property
  measures 2} since $\mu$ is nonnegative.

\medskip We now prove the inclusion in \eqref{eq: prop: first
  property measures 3}. One has 
\begin{align*}
   \Norm{\varphi(\hk D_t) \psi (t) u_k}{L^2(\R; H^2(\M, \kk \mugk)) }
      &\lesssim \Norm{\varphi(\hk D_t) \psi (t) u_k}{L^2(\L, \k \mug dt)}\\
      &\quad + \Norm{\mathsf H_g \varphi(\hk D_t) \psi (t) u_k}{L^2(\L, \k \mug dt)},
\end{align*}
where $\mathsf H_g$ denotes the Hessian operator. Using the elliptic
      regularity and  that $u_k$ is solution to the wave equation one obtains
\begin{align*}
  \Norm{\mathsf H_g \varphi(\hk D_t) \psi (t) u_k}
  {L^2(\L, \k \mug dt)}
   &\lesssim \Norm{\varphi(\hk D_t) \psi (t) A_{\kk,\gk} u_k}
  {L^2(\L, \k \mug dt)}\\
  &\lesssim \Norm{\varphi(\hk D_t) \psi (t) \d_t^2  u_k}
    {L^2(\L, \k \mug dt)}\\
  &\lesssim \Norm{\varphi(\hk D_t) \psi (t) \d_t^2  u_k}
    {L^2(\L, \kk \mugk dt)}.
\end{align*}
Then, one writes
\begin{align*}
\Norm{\varphi(\hk D_t) \psi (t) \d_t^2  u_k}
    {L^2(\L, \kk \mugk dt)}^2
  &= \sum_{\nu \in J_k } |u_\nu|^2 \lambda_\nu^2 
  \Norm{ \varphi(\hk D_t)  \big( \psi (t) e^{it\sqrt{\lambda_\nu}}
  \big)}{L^2(\R)}^2 \\
  &\lesssim C_N
  \hk^N\Norm{{u_k}\bt}{L^2(\M, \kk \mugk)}^2
   \lesssim C_N \hk^N,
\end{align*}
for any $N \in \N$, using \eqref{eq: decay Fourier localization2} and that 
$\hk^2 \lambda_\nu \lesssim 1$ for $\nu \in J_k$.
Hence, one has
\begin{align*}
&\Norm{\varphi(\hk D_t) \psi (t) u_k}{L^2(\R; H^2(\M, \kk \mugk)) }^2
  \lesssim C_N \hk^N,
\end{align*}
implying, by the trace formula,
\begin{align*}
&\Norm{\varphi(\hk D_t) \psi (t)  \hk \d_{\n_k}
  {u_k}_{|\d\L}}{L^2(\d\L, \kk \mugbk dt)) }^2
  \lesssim C_N \hk^N.
\end{align*}
This yields $\dup{\nu}{|\varphi(\tau)\psi(t)|^2} = 0$,
which gives \eqref{eq: prop: first property
  measures 3} since $\nu$ is nonnegative.
\end{proof}

%%%%%%%%%%%%%%%%%
% subsection
\subsubsection{Final step of the proof of observability}
\label{sec: Proofs of observability inequalities}

With Theorem~\ref {thm: equationmesure}  at hand we can
complete  the proof of the observability  results of
Theorem~\ref{theo-perturbation}. The two measures $\mu$ and $\nu$ given by
 Proposition~\ref{prop: existence measures wave}, with scale $\hk
 = \varrho^{-|k|}$ of Section~\ref{dyadic}, fulfill the assumptions of 
 Theorem~\ref{thm: equationmesure} by Proposition~\ref{prop: first
   property measures}. Hence, one has
 \begin{align*}
    \Hp \mu  =  \int_{\y \in \pHb \cup \pGb} 
    \frac{\delta_{\y^+} - \delta_{\y^-}}
    {\dup{\xi^+- \xi^-}{\nx}_{T_x^*\M, T_x\M}} 
  \ d \nu (\y),
 \end{align*}
 since the considered wave equations are homogeneous here. 
Theorem~\ref{theo-propagation}
recalled from  the companion article~\cite{BDLR2}  implies that the support of the measure $\mu$ is a union of maximal  \gbichars. 

Recall that $I = ]\delta, T-\delta[$. Let
$\y^0= (t^0, x^0, \tau^0, \xi^0) \in \supp \mu$, with $t^0 \in I$.  According the first
point of Proposition~\ref{prop: first property measures} such a point
exists. Then, there exists a \gbichar
$\gammaG$ with $\gammaG(0) = \y^0$ such that
$\GammaG \subset \supp \mu$. 

\medskip
{\bf Case of an interior observation.}

With the interior geometric control
condition fulfilled by $(\omega, T-2\delta)$ (Definition~\ref{control-geo-interior}) the \bichar $\gammaG$ reaches a point above  $I \times \omega$. Yet, from the last point of  Proposition~\ref{prop: first property
  measures}, the measure $\mu$ vanishes above $I \times \omega$, 
which gives a contradiction and concludes the proof in this case.

\medskip
{\bf Case of a boundary observation.}

With the boundary geometric control condition fulfilled by
$(\Gamma, T-2\delta)$ (Definition~\ref {control-geo-boundary}), 
there exists  $s \in \R$ such that $t(s) \in I$ and 
\begin{enumerate}
\item eiher $\y^1 = \lim_{s\to s^-} \gammaG(s) \in \ESC^F$;
\item or $\y^1 = \lim_{s\to s^+} \gammaG(s) \in \ESC^P$.
\end{enumerate}
The sets $\ESC^F$ and $\ESC^P$ are introduced in Definition~\ref{def:
  escape point}.  With the measure $\nu$ vanishing above $I \times
\Gamma$ by the last point of Proposition~\ref{prop: first property
  measures}, the measure propagation equation is locally $\transp \Hp
\mu=0$. Thus locally, the support of $\mu$ is a union of
$\bichars$. In both possibilities, all such \bichars exit $\TL$
reaching a region where $\mu$ vanishes, which gives a contradiction
and concludes the proof in this case.

%%%%%%%%%%%%%%%%%
% section
%%%%%%%%%%%%%%%%%
\section{Proof of the propagation equation I}
\label{sec: Proof of the propagation equation1}

\subsection{Preliminary remarks and observations}
\label{sec: Preliminary remarks and observations}
Recall that $u_k$ is the zero-extension to $\hL$ of a weak solution to
the wave equation in $\L$. With the homogeneous Dirichlet condition
this extension is $H^1$.

Consider $\chi \in \Cinfc(\R)$ with $0 \notin \supp \chi$.  Since the
coefficient of the wave operator are independent of time $t$ then
$\tilde{u}_k=\chi(h D_t) u_k$ also solves the wave equation in $\L$
with $\tilde{f}_k=\chi(h D_t) f_k$ as source term and the Hermitian
\scm of $\transp (\tilde{u}_k,\hk \tilde{f}_k)$ is $|\chi(\tau)|^2
M$. Similarly, the \scm of $(\hk \d_{\n_k}{\chi(h D_t)
  u_k}_{|\d\L})_k$ is $|\chi(\tau)|^2 \nu$.  If we prove that
\eqref{eq: GL equation-propagation theorem} holds for $M$ and $\nu$
replaced by $|\chi(\tau)|^2 M$ and $|\chi(\tau)|^2 \nu$ then, using
\eqref{eq: hyp support mesure}, one finds that \eqref{eq: GL
  equation-propagation theorem} holds also for $M$ and $\nu$ by the
\Ldct.  Without any loss of generality we may thus replace $u_k$ by
$\tilde{u}_k$ and $f_k$ by $\tilde{f}_k$.  Then, there exists $0<
C_{\mu,0}< 1< C_{\mu,1}< \infty$ such that
\begin{align}
     \label{eq: hyp support mesure-new-mu}
     &\supp \mu \subset \Char p \cap \TL 
     \cap \{ C_{\mu,0}\leq |\tau|\leq  C_{\mu,1}\},\\
     &\supp M_{0,1} \subset \TL 
     \cap \{ C_{\mu,0}\leq |\tau|\leq  C_{\mu,1}\},\notag
   \end{align}
and 
\begin{align}
     \label{eq: hyp support mesure-new-nu}
     \supp \nu \subset T^*\d\L  \cap \{ C_{\mu,0}\leq
    |\tau|\leq  C_{\mu,1}\}.
     \end{align}
If no precision is given, the $L^2$-norm on $\M$ is given by the
density measure $\k \mug$ in what follows.

Suppose $I$ is a time interval. With
the $\tau$-microlocalisation performed above, one has
\begin{align}
   \label{eq: mass2-proof 0}
   \Norm{u_k}{L^2(I \times \M)}
   \eqsim \Norm{\hk^2 \d_t^2 u_k}{L^2(I \times \M)}
   \eqsim \Norm{\hk^2 A_{\kk , \gk} u_k + \hk^2 f_k}{L^2(I \times \M)}.
\end{align}
using the wave equation.
Assume that a subsequence of $u_k$
converges to $0$ in $L^2(I\times \M)$. This gives $\mu=0$ on $T^*(I\times
\M)$. With \eqref{eq: mass2-proof 0}, one finds that $\Norm{\hk^2
  \d_t^2 u_k}{L^2(I \times \M)}\to 0$ and $\Norm{\hk^2 A_{\kk , \gk}
  u_k }{L^2(I \times \M)}\to 0$ also, using that $\hk f_k$ is
$L^2_{\loc}$-bounded.
Ellipticity up to the boundary gives
$\Norm{\hk^2 u_k}{H^2(I\times \M)} \to 0$
and 
interpolation gives 
\begin{align*}
  \Norm{\hk \d_t u_k}{L^2(I\times \M)}  \to 0
  \ \ \et \ \
     \Norm{\hk \nabla_{\!\! \gk} u_k}{L^2(I \times \M)} \to 0.
\end{align*}
With Proposition~\ref{prop: admissibility} one conludes that
$\Norm{\hk \d_\n u_k}{L^2(I \times \d\M)} \to 0$. Hence, all terms in the measure equation vanish, in this case. 
One may thus assume that  $\Norm{u_k}{L^2(I\times \M)}
\gtrsim 1$, for any interval $I$.

With the arguments given just above, one has 
\begin{align}
   \label{eq: mass2-proof}
  1 &\eqsim \Norm{u_k}{L^2(I\times \M)}  \eqsim  
  \Norm{\hk^2 \d_t^2 u_k}{L^2(I\times \M)}
  \eqsim  \Norm{\hk^2 A_{\kk , \gk} u_k}{L^2(I\times \M)}\\
  &\eqsim \Norm{\hk^2 u_k}{H^2(I \times \M)}
  \eqsim \Norm{\hk \d_t u_k}{L^2(I\times \M)}  
      \eqsim  \Norm{\hk \nabla_{\!\! \gk} u_k}{L^2(I\times \M)},\notag
\end{align}
and more generally one has 
\begin{equation}
   \label{eq: mass2-proof4}
  \hk^{\ell+2s} \Norm{\d_t^\ell  A_{\kk , \gk}^s u_k}{L^2 (\M)} \eqsim  1,
\end{equation}
using that $\Norm{\hk^\ell\d_t^\ell f_k}{L^2(I\times \M)}   \eqsim \Norm{f_k}{L^2(I\times \M)}$.

\medskip
Note that equation~\eqref{eq: GL equation-propagation theorem} is
local. Consequently, the proof can be carried out in local charts. It 
greatly simplifies away from the boundary: there \eqref{eq:
  GL equation-propagation theorem} is $\Hp \mu  =2 \Im M_{0,1}$, which follows
from Proposition~\ref{prop: commutator} below.    We will thus
only consider the case of a
local chart $\chart=(\O,\cdiffL)$ near a boundary point, where the
boundary is given by $\{z=x_d=0\}$ and $\M=\{ z>0\}$; see
Section~\ref{sec: Local coordinates}. 
By abuse of notation we will use the notation $A$, $\k$, $g$, $u_k$ and $\mu$
for their representative in the local chart. 

For concision we will write $y=(t,x)$, $y'=(t,x')$,
$\eta=(\tau,\xi)$, and $\eta'= (\tau,\xi')$. 

\medskip
A consequence of \eqref{eq: mass2-proof}--\eqref{eq: mass2-proof4} that we will use is as
follows in $\chart$, for $\psi \in \Cinfc(\R^{d+1})$,
\begin{multline}
  \label{eq: local H2 tangential estimate}
  \sum_{1\leq  i, j \leq d-1} 
  \Norm{\psi \hk^2 D_i D_j \udl{u}^k}{L^2(\R^{d+1})}
  + \sum_{1\leq  j \leq d-1} 
  \Norm{\psi \hk^2 D_j D_d\udl{u}^k}{L^2(\R^{d+1})}\\
  + \Norm{\psi \hk^2 D_d^2 u^k}{L^2(\R^{d+1}_+)}
  \lesssim 1.
\end{multline}
Note that the last term $\hk^2 D_d^2 u^k$ does not lie in
$L^2(\R^{d+1})$ in general as $D_d u_k$ is not continuous across
$z=0$. This explains the computation of its norm only on $\R^{d+1}_+$.

\medskip
As mentioned in Section~\ref{sec: Normal derivative estimation}, in
the quasi-normal geodesic coordinates of Propositon~\ref{prop:
  quasi-normal coordinates} one has $\d_{\n} = \d_d$, if $\n$ is the
unitary normal inward pointing vector field to $\d\M$. Here, we will
not use a different coordinate system if the metric $g_k$ varies. In
the chosen local chart $\chart$, quasi-normal geodesic coordinates adapted to the ``reference''
metric $g$ will be kept fixed. In such coordinates one has $\n_k = \sum_{j} g_k^{dj}
\d_j$. For a function like $u_k$ that vanishes at $z=0$ one has
$\d_{\n_k} {u_k}_{|\d\L} = g_k^{dd} \d_d {u_k}_{|z=0}$.  Hence
\begin{align}
  \label{eq: vk local coordinates}
  v_k = \hk  g_k^{dd} \d_z {u_k}_{|z=0^+},
\end{align}
in local coordinates. Note that ${g_k^{dd}}_{|z=0} = 1 + o(1)$ as $k
\to \infty$, since
${g^{dd}}_{|z=0} = 1$ in the chosen quasi-normal geodesic
coordinates; see Section~\ref{geometry}.

\medskip
From the ``jump formula" the sequences $u_k$ and $v_k$ satisfy 
\begin{equation}\label{eq-approx}
  \hk^2 (\d_t^2 - A_{\kk, \gk})u_k = \hk^2 f_k - \hk v_k  \otimes \delta_{\d \L} .
\end{equation}

\medskip
The proof follows the lines of the proof of~\cite[Theorem 2.3]{GL:1993} (or also~\cite[Théorème 4]{Burq:1997}).
The main differences are as follows.
\begin{itemize}
\item The sequence $(u_k)_k$ is solution to a {\em family of} wave equations associated with Lipschitz metrics.
\item On the \lhs of~\eqref{eq-approx}, the wave operator is $k$
  dependent. In fact, in the application we make in Section~\ref{sec:
    Proof of the observability results}, the sequence $(u_k)_k$ is not
  spectrally localized with respect to a fixed operator, but rather
  with respect to the family of operators $A_k = A_{\kk,\gk}$,
\begin{align*}
  u_k = 1_{- \hk ^2 A_k \in [\alpha^2, \alpha^{-2} [ }u_k.
\end{align*}
  
\item With respect to \cite[Théorème 4]{Burq:1997} the result here
  assumes less smoothness ($W^{1, \infty}$ as opposed to  $\Con^2$),
  while the difference with respect to~\cite[Theorem 2.3]{GL:1993} is
  more subtle: in~\cite{GL:1993}, the authors study only rough ($W^{2, \infty}$)  domains
  of $\R^d$ with the usual flat metric, which ensures the existence of
  local coordinate systems that are {\em regular} with respect to the
  variable tranverse to the boundary, which simplifies greatly the analysis. Note that 
  neither the result of \cite{GL:1993} nor its  proof are preserved by
  change of coordinates. To the opposite the result and proofs in the present
  article are coordinate invariant and thus geometrical.
\end{itemize}
The lack of regularity and the geometrical framework we consider, if compared with
\cite{Burq:1997,GL:1993}, generate technical difficulties.
  In Sections~\ref{sec: More on semi-classical symbols and
  operators}--\ref{sec: Proof of the propagation equation2}, we will
use a particular decomposition of symbols based on the Weierstrass
preparation theorem. This allows one to express symbols as a
first-order polynomial in $\zeta$, the dual variable of $z=x_d$, with
coefficients that are tangential symbols, and a remainder term. An
issue is then the handling of the different terms that lack decay in
$\zeta$, even if the initial symbol is smooth with fast decay. This is a main reason for the introduction of {\em ad hoc}
symbol spaces, taking into account both the low regularity in the $x$
variable (that originates from the  regularity of the  metric we
consider) and this low decay in the $\zeta$ variable. This makes some of
the statements quite technical even though we made an effort to minimize this aspect.

\subsection{Commutator analysis}
\label{sec: Analysis of a commutator}
\bigskip
To establish the propagation equation of Theorem~\ref{thm:
  equationmesure} we carefully compute a commutator.
In fact, assume for a second that $p$ and $b$ are smooth symbols
$-i \hk\Hp b = -i \hk \{p, b\}$ is the principal symbol of the
commutator $[\OpH(p), \OpH(b)]$. Hence, to find the value of
$\dup{\mu}{\Hp b}$ it is very natural to analyse the limit of 
\begin{align*}
  \biginp{\hk [P_{\kk,\gk}, \OpH(b)] u_k}{u_k}_{L^2}.
\end{align*}
Technicalities arise because of the limited smoothness of the
coefficients of $P_{\kk,\gk}$ that prevent one from using standard
semi-classical calculus results. However, there is no restriction on
the smoothness and decay properties of the test symbol $b$.
In fact, in the course of the proof of the
Proposition~\ref{prop: commutator},  differentiations of the symbol
$b$ with respect to $x$ are
needed as well as some decay in $\xi$.   For simplicity, symbols in
$\Symbolo{\infty,\infty}{-\infty}{2d}$
 are thus considered; 
see Definition~\ref{def: symbols} and \eqref{eq: symbols fast decay}.
A classical result is the following one. 
%%%%%%%%%%%%%%%%%%%%%%%%
% lemma                %
%%%%%%%%%%%%%%%%%%%%%%%%
\begin{lemma}
  \label{lemma: smoothing semiclassical operators}
  Let $b \in
\Symbolo{\infty,\infty}{-\infty}{2d}$ and $s, s' \in \R$. There
exists $C_{s,s'}>0$ such that 
\begin{align*}
  %\label{eq: smoothing semiclassical operators}
\Norm{\OpH(\jp{\xi}^{s'})   \OpH(b) u}{L^2(\R^d)}
  \leq C_{s,s'} \Norm{\OpH(\jp{\xi}^{s}) u}{L^2(\R^d)}, \qquad u \in \mathscr S(\R^d).
\end{align*}
\end{lemma}
This means that $\OpH(b)$ sends any {\em semi-classical} Sobolev space in the
intersection of all semi-classical Sobolev spaces. This is due to both the
smoothness in $x$ and $\xi$ and the fast decay in $\xi$ of $b$.

Here, and in what follows one writes $H_z^{\delta}L_{y'}^2$ in place
of $H^\delta(\R_z; L^2(\R_{y'}^{d}))$ in norm indices or duality bracket
indices for the sake of
concision for $\delta \in [-1,1]$. 
For a density measure $\rho$ on $\R^N$ we will denote by
\begin{align*}
  \inp{.}{.}_{H_z^{-\delta}L_{y'}^2, H_z^{\delta}L_{y'}^2}^\rho,
\end{align*}
the complex duality bracket understood with $L^2(\R^{N}, \rho)$
 as a pivot space.
%%%%%%%%%%%%%%%%%%%%%%%%
% proposition          %
%%%%%%%%%%%%%%%%%%%%%%%%
\begin{proposition}
  \label{prop: commutator}
  Suppose $b \in \Symbolo{\infty,\infty}{-\infty}{2d+2}$ with $\supp
  b\subset K \times \R^{d+1}$, for $K$ a compact of $\cdiffL(\O)$, and
  suppose $\psi \in \Cinfc(\cdiffL(\O))$ be equal to $1$ in a \nhd of
  the $y$-projection of $\supp b$.  Set
 \begin{align}
   \label{eq: Lk}
   L_k (b,\psi) &=    i \biginp{ \OpH(b) \psi u_k}
             {v_k \otimes \delta_{z=0}}_{ H_z^{1}L_{y'}^2, H_z^{-1}L_{y'}^2}^{\kk\mugk dt}\\
              &\quad    
             -i \biginp{ v_k \otimes \delta_{z=0}}
             {\psi   \OpH(b)^\star u_k}_{H_z^{-1}L_{y'}^2,
             H_z^{1}L_{y'}^2}^{\kk\mugk dt}.\notag
 \end{align}
 One has $
  \lim_{k\to +\infty }L_k (b,\psi)= - \dup{\mu}{\Hamiltonian_{p_{\k,g}} b} - 2 \dup{\Im M_{0,1}}{b}$.
\end{proposition}
Recall that $ v_k \otimes \delta_{z=0}$ lies in $H_{\loc}^{-s}(\R_z; L^2(\R_{y'}^{d}))$ for
any $s >1/2$, hence one finds the duality brackets appearing in the definition of
$L_k$. Here,  $\OpH(b)^\star$ stands for the adjoint in the sense of the 
$L^2(\R^{d+1}, \kk\mugk dt)$-inner product. 
%% proof
\begin{proof}
 Using $L^2(\R^{d+1}, \kk \mugbk d t)$ as a pivot space and the
 symmetry of $P_{\kk, \gk}$ for the associated inner product and using
 \eqref{eq-approx} one has
   \begin{align*}
     &\biginp{ [  P_{\kk, \gk},\OpH(b) \psi ] u_k}
     { u_k}_{L^2(\R^{d+1};\kk\mugk dt)}\\
     &\quad = \biginp{  \OpH(b) \psi  u_k}
       {  P_{\kk, \gk}  u_k}_{H_z^{1}L_{y'}^2, H_z^{-1}L_{y'}^2}^{\kk\mugk dt}
       - \biginp{  \OpH(b) \psi P_{\kk, \gk} u_k}
       { u_k}_{L^2(\R^{d+1};\kk\mugk dt)}\\
       &\quad = \biginp{  \OpH(b) \psi  u_k}
       {f_k}_{L^2(\R^{d+1};\kk\mugk dt)}
       - \biginp{  \OpH(b) \psi f_k}
       { u_k}_{L^2(\R^{d+1};\kk\mugk dt)}\\
     &\qquad  -\hk^{-1} \biginp{  \OpH(b) \psi  u_k}
      { v_k \otimes\delta_{z=0}}_{H_z^{1}L_{y'}^2,  H_z^{-1}L_{y'}^2}^{\kk\mugk dt}
       \\
     &\qquad+ \hk^{-1}\biginp{ \OpH(b) \psi v_k \otimes\delta_{z=0}}
       {u_k}_{L^2(\R^{d+1};\kk\mugk dt)}\\
     &\quad = \biginp{  \OpH(b) \psi  u_k}
       {f_k}_{L^2(\R^{d+1};\kk\mugk dt)}
       - \biginp{  \OpH(b) \psi f_k}
       { u_k}_{L^2(\R^{d+1};\kk\mugk dt)}\\
     &\qquad   + i \hk^{-1} L_k (b,\psi) .
   \end{align*}
   Since one has
   \begin{align*}
     &\lim_{k\to +\infty}
     \frac{\hk}{i} \big(
     \biginp{\OpH(b) \psi  u_k}{f_k}_{L^2(\R^{d+1};\kk\mugk dt)}
     - \biginp{\OpH(b) \psi f_k}{u_k}_{L^2(\R^{d+1};\kk\mugk dt)}
     \big) \\
     &\qquad = \frac{1}{i} \big( \dup{M_{0,1}}{b} - \dup{M_{1,0}}{b}\big)
     = 2 \dup{\Im M_{0,1}}{b},\notag
   \end{align*}
   the result follows if one proves
   \begin{align*}
     \lim_{k\to +\infty} I_k =  - \dup{\mu}{\Hamiltonian_{p_{\k,g}} b},
     \end{align*}
    with
    \begin{align*}
    I_k &= 
    \frac{\hk }{i}\biginp{ [   P_{\kk, \gk},\OpH(b) \psi ] u_k}
    { u_k}_{L^2(\R^{d+1};\kk\mugk dt)}\\
    &=  
    \frac{1}{i \hk}\biginp{ \tk^{-1} \tkk [  \hk^2  P_{\kk, \gk},\OpH(b) \psi] u_k}
    { u_k}_{L^2(\R^{d+1};\k\mug dt)}.\notag
 \end{align*}
Recall that $\tk = \k \det(g)^{1/2}$ and $\tkk = \kk
\det(\gk)^{1/2}$.
First one writes
\begin{align*}
  [  \hk^2  P_{\kk, \gk}, \OpH(b) \psi]
  = [  \hk^2  P_{\kk, \gk},\OpH(b) ] \psi
  + \OpH(b)  [  \hk^2  P_{\kk, \gk},\psi].
\end{align*}
Since $ [ P_{\kk, \gk},\psi]$ is a differential operator of order one
with Lipschitz coefficients one finds
\begin{align*}
  [  \hk^2  P_{\kk, \gk}, \OpH(b) \psi]
  = [  \hk^2  P_{\kk, \gk},\OpH(b) ] \psi
  + \hk^2\OpH(b)  O(1)_{\L (H^1, L^2)}, 
\end{align*}
and one obtains
\begin{align}
    \label{limite2}
  \lim_{k\to +\infty}I_k = 
   \lim_{k\to +\infty} 
    \frac{1}{i \hk}\biginp{ \tk^{-1} \tkk [  \hk^2  P_{\kk, \gk},\OpH(b) ]\psi u_k}
    { u_k}_{L^2(\R^{d+1};\k\mug dt)}.
 \end{align}
  According to symbolic calculus one has 
  \begin{equation}
    \label{facile}
    [ \hk^2  \d_t^2 , \OpH(b) ] 
    = i   \hk \OpH (2 \tau  \d_tb) + o(\hk) _{\mathcal{L} (L^2)} .
    %= \hk O(1)_{\mathcal{L} (L^2)} .
\end{equation}
The contribution of \eqref{facile} to the limit in \eqref{limite2} is
then 
\begin{equation}
     \label{eq: contrib limite 1}
    \lim_{k\to +\infty} 
    \frac{1}{i \hk}\biginp{  \tk^{-1} \tkk [ \hk^2  \d_t^2 ,\OpH(b)] \psi u_k}
    {  u_k}_{L^2(\R^{d+1}, \k \mug dt)} 
    = \dup{\mu}{2 \tau \d_t b}
    = \dup{\mu}{ \{\tau^2, b \} },
 \end{equation}
by Lemma~\ref{lem: continu - manifold} and \eqref{eq: lem: continu - manifold - local chart}.

Next, with repeated indices convention, one writes
\begin{align*}
  A_{\kk, \gk} =\rho_k^j \d_j + \gk^{ij} \d_i \d_j
  \quad \avec \ \ \rho_k^j =  \tkk^{-1} [\d_i, \tkk \gk^{ij} ],
\end{align*}
with $(\rho_k^j)_k\subset L^\infty$ that  converges to some $\rho^j
\in \Con^0\cap L^\infty$. 
One computes
\begin{align*}
  %\label{eq: computation commutator}
  [\hk^2 A_{\kk, \gk}, \OpH(b) ] = A_1 + A_2 + A_3 + A_4,
\end{align*}
with
\begin{align*}
  &A_1=  \hk^2\rho_k^j [\d_j , \OpH(b) ],
  \quad 
  A_2 = \hk^2[ \rho_k^j  , \OpH(b)] \d_j,\\
  &A_3 =\hk^2\gk^{ij} [ \d_i \d_j , \OpH(b) ], \quad 
  A_4 =  \hk^2 [ \gk^{ij} , \OpH(b)] \d_i \d_j.
\end{align*}
One writes
\begin{align*}
  A_1=  \hk^2\rho_k^j \OpH(\d_{x_j} b) = O(\hk^2)_{\L(L^2)},
\end{align*}
with Lemma~\ref{lemma: smoothing semiclassical operators} since
$\d_{x_j} b \in \Symbolo{\infty,\infty}{-\infty}{2d+2}$ (one can also
argue that $\d_{x_j} b \in \symb(\R^{2d+2})$ and use Lemma~\ref{4.continuite}). With \eqref{WCMbis} in Proposition~\ref{propcommut} one finds
\begin{align*}
  A_2 = \hk^2\sum_j o(1)_{\L(L^2)} \, \d_j =  o(\hk^2)_{\L(H^1, L^2)}. 
\end{align*}
Thus, with \eqref{eq: mass2-proof} one finds that the contributions of $A_1$ and $A_{2}$ to the limit in \eqref{limite2}
are $0$.

\medskip
Next, one writes
\begin{align*}
  A_3 =\hk^2\gk^{ij} [ \d_i \d_j , \OpH(b) ]
  &=
  \hk i \gk^{ij}  \OpH\big(\xi_i\d_{x_j} b +\xi_j\d_{x_i} b   \big)
  + \hk^2 \gk^{ij} \OpH(\d_{x_i}\d_{x_j} b)\\
  &=
  \hk i \gk^{ij}  \OpH\big(\{ \xi_i \xi_j, b\} \big)
  + \hk^2 O(1)_{\L(\L^2)},
\end{align*}
as $\d_{x_i}\d_{x_j} b\in \Symbolo{\infty,\infty}{-\infty}{2d+2}$.
Since $\{ \xi_i \xi_j, b\} \in \symbo(\R^{2d+2})$, with \eqref{eq:
  mass2-proof} one finds that the
contributions of $A_3$ to the limit in \eqref{limite2} is
\begin{align}
  \label{eq: contrib limite 2}
  - \lim_{k\to +\infty} 
    \frac{1}{i \hk}\biginp{ \tk^{-1} \tkk  A_3 \psi u_k}
  {  u_k}_{L^2(\R^{d+1}, \k \mug dt)}
  = - \dup{\mu}{g^{ij} \{ \xi_i \xi_j, b\} },
\end{align}
by Lemma~\ref{lem: continu - manifold} and \eqref{eq: lem: continu - manifold - local chart}.
Finally, with~\eqref{WCM2bis} in Proposition~\ref{propcommut} one writes
\begin{align*}
  A_4 &=  i \hk^3 \sum_\ell \big(\d_{x_\ell} \gk^{ij}\big) \OpH\big(
  \d_{\xi_\ell} b\big) \d_i \d_j
  + o (\hk^3)_{\L(H^2, L^2)}\\
  &= -  i \hk \sum_\ell \big(\d_{x_\ell} \gk^{ij} \big) \OpH\big(
 \xi_i \xi_j  \d_{\xi_\ell} b\big) 
  + o (\hk^3)_{\L(H^2, L^2)},
\end{align*}
implying with \eqref{eq: mass2-proof} that the contributions of $A_4$ to the limit in
\eqref{limite2} is 
\begin{align}
\label{eq: contrib limite 3}
-\lim_{k\to +\infty} \frac{1}{i \hk}
\biginp{\tk^{-1} \tkk  A_4 \psi u_k}
{u_k}_{L^2(\R^{d+1}, \k \mug dt)}
&=\sum_\ell \dup{\mu}
{\xi_i \xi_j  \d_{x_\ell} g^{ij}  \d_{\xi_\ell} b}\\
&= - \bigdup{\mu}{\xi_i \xi_j \{g^{ij}, b\} }.
\notag
\end{align}
Gathering \eqref{eq: contrib limite 1}, \eqref{eq: contrib limite 2}
and \eqref{eq: contrib limite 3} and writing
\begin{align*}
  \{\tau^2, b\} - g^{ij}\{ \xi_i \xi_j, b\} - \xi_i \xi_j \{ g^{ij}, b\}
  = \{ \tau^2 - g^{ij}\xi_i \xi_j, b\}
  = - \{ p_{\k,g} , b\} = - \Hamiltonian_{p_{\k,g}} b ,
\end{align*}
 one obtains the result of the proposition.
\end{proof}

%%%%%%%%%%%%%%%%%
% subsection
\subsection{Time microlocalization}
Above in Proposition~\ref{prop: commutator},  we defined $L_k (b,\psi) $ for symbols $b$ in
$\Symbolo{\infty,\infty}{-\infty}{2d}$ and we are interested in
the limit of $L_k (b,\psi) $ as $k \to +\infty$. With the support
properties of the measure $\mu$ given in \eqref{eq: hyp support mesure-new-mu} one obtains the
following lemma.
%%%%%%%%%%%%%%%%%%%%%%%%
% lemma                %
%%%%%%%%%%%%%%%%%%%%%%%%
\begin{lemma}
  \label{lem: time microlocalization}
  Suppose $\chi\in \Cinfc ( C_{\mu,0}^2, C_{\mu,1}^2)$ be equal to $1$ on a
\nhd of  $[C_{\mu,0}, C_{\mu,1}]$. Let $b \in \Symbolo{\infty,\infty}{-\infty}{2d}$ and $\psi$ be as
  in Proposition~\ref{prop: commutator}. Then, one has 
  \begin{align*}
    \lim_{k\to +\infty} L_k \big((1 - \chi) (\tau) \,  b,\psi \big) =0.
  \end{align*}
\end{lemma}

%%%%%%%%%%%%%%%%%
% section
%%%%%%%%%%%%%%%%%
\section{More on semi-classical symbols and operators}
\label{sec: More on semi-classical symbols and operators}
%%%%%%%%%%%%%%%%%
% subsection
\subsection{Preparation theorem: Euclidean symbol division}
\label{sec: Preparation theorem: Euclidean symbol division}
For technical reasons, it is convenient to consider symbols  with
finer properties here and in what follows.  
%%%%%%%%%%%%%%%%%%%%%%%%
% definition           %
%%%%%%%%%%%%%%%%%%%%%%%%
\begin{definition}
  \label{def: holomorphic symbol}
  One says that $a \in \symbH(\R^{2d+2})$ if $a \in
  \Symbolo{\infty,\infty}{-\infty}{(2d+2)}$ and satisfies moreover the following properties
  \begin{itemize}
  \item $a(y, \eta)$ is compactly supported in the $y$ variable;
  \item $a(y, \eta)$  has a  compactly supported Fourier transform in
    the $\eta$ variable and, consequently, is holomorphic with respect
    to the $\eta$ variable.
\end{itemize}
\end{definition}
This choice is possible and relevant observing that  $\symbH(\R^{2d+2})$ is dense in
$\Symbolo{\infty,\infty}{-\infty}{2d+2}$, the symbol classes  we
consider in Propositon~\ref{prop: commutator}. 
Recall that $y = (t,x)$ and $\eta= (\tau,\xi)$, $x=(z,y')$ (the boundary is given by $\{ z=0\}$. See the beginning of Section~\ref{sec: Proof of the propagation equation1}.

Below we will need the following quantification of the decay of a
symbol $a \in \symbH(\R^{2d+2})$ with the $\eta$ variable allowed to slightly
depart from the real axis: for any $R>0$, $\alpha,
  \beta \in \N^{d+1}$, and  $N\in \N$ there
exists $C_{\alpha,\beta, N,R} >0$ such that 
\begin{align}
  \label{eq: rapid decay symbol SigmaH}
  |\d_y^\alpha \d_{\eta}^{\beta} a(y, \eta)|
  \leq C_{\alpha,\beta, N,R} \jp{\eta}^{-N}, \quad  y \in
  \R^{d+1}, \eta \in \C^{d+1} \ \avec\   |\Im \eta|\leq R.
\end{align}
This is given by the Paley-Wiener theorem; see for instance
\cite[Theorem 7.3.1]{Hoermander:V1}.

\bigskip
The following proposition gives a decomposition of a symbol 
$b \in \symbH(\R^{2d+2})$. For our purpose, that is, using such
symbols in the identity given by Proposition~\ref{prop: commutator},
with Lemma~\ref{lem: time microlocalization} it suffices to work with a time-frequency
truncated symbol. Recall that $\eta = (\eta', \zeta)$ with $\eta' =
(\tau ,\xi')$ and $\zeta$ is the dual variable to $z= x_d$. 
%%%%%%%%%%%%%%%%%%%%%%%%
% proposition          %
%%%%%%%%%%%%%%%%%%%%%%%%
\begin{proposition}[Euclidean symbol division]
  \label{prop: Euclidean symbol division}
  Let $\chi\in \Cinfc ( C_{\mu,0}^2, C_{\mu,1}^2)$ be equal to $1$ on a
\nhd of  $[C_{\mu,0}, C_{\mu,1}]$ and $b(y, \eta)  \in
\symbH(\R^{2d+2})$. 
For $k \in \N$, there exist $b_{0,k}(y,\eta')$, $b_{1,k}(y, \eta')$ and $q_k(y,\eta)$ such
that 
\begin{align}
  \label{eq: symbol decomposition Weierstrass}
  \chi(\tau) b(y,\eta',\zeta) = b_{0,k}(y,\eta') + b_{1,k}(y,\eta') \zeta 
  + q_k(y,\eta',\zeta) \,  p_k(y,\eta',\zeta),
\end{align}
with the following symbol properties
\begin{multline}
  \label{eq: estimation bjk}
  \big| \d_{y}^\alpha \d_{\eta'}^\beta  b_{j,k} (y,\eta') \big| \leq
  C_{N, \beta}
  \jp{\eta'}^{-N}, \\
  \pour\ N \in \N, \ \alpha \in \N^{d+1},  \ |\alpha| \leq 1, \
  \beta \in \N^{d},  \
  j=0,1, \ y \in \R^{d+1}, \ \eta'\in \R^{d},
\end{multline}
and 
\begin{multline}
  \label{eq: estimation rk}
  \big| \d_{y}^\alpha \d_{\eta'}^\beta  \d_\zeta^\delta q_k(y,\eta',\zeta) \big| 
  \leq C_{N, \beta, \delta} 
  \jp{\eta'}^{-N} \jp{\zeta}^{-1-\delta}, \\
  \pour\ N\in \N, \ \alpha \in \N^{d+1},  \ |\alpha| \leq 1, \
  \beta \in \N^{d},  \ \delta \in \N, 
 \ y \in \R^{d+1}, \ (\eta',\zeta) \in \R^{d+1}, 
\end{multline}
uniformly with respect to $k \in \N$. Moreover, $q_k$ admits a
polyhomogeneous development in the~$\zeta$ variable: there 
 exist
$q_k^j (y, \eta')$, $j\in \N^*$, such that 
\begin{multline}
  \label{eq: estimation rjk}
  \big| \d_{y}^\alpha \d_{\eta'}^\beta  q_k^j (y,\eta') \big| \leq
  C_{N, \beta}
  \jp{\eta'}^{-N}, \\
  \pour\ N \in \N,\  \alpha \in \N^{d+1},  \ |\alpha| \leq 1, \
  \beta \in \N^{d},  \ y \in \R^{d+1}, \ \eta'\in \R^{d},
\end{multline}
uniformly with respect to $k \in \N$, and $q_k \sim  \sum_{j\geq 1} q_k^j \zeta^{-j}$ in the following sense:
for $\phi \in \Cinfc(\R)$ equal to $1$ near $0$ one has 
\begin{multline}\label{polyhomog}
  \Big|\d_{y}^\alpha \d_{\eta'}^\beta  \d_\zeta^\delta  
  \Big(q_k (y,\eta',\zeta)  - (1- \phi(\zeta) ) \sum_{j=1}^{M} q_k^j (y, \eta') \zeta^{-j} \Big)  \Big|
  \leq C_{N, M, \beta, \delta}  \jp{\eta'}^{-N} \jp{\zeta}^{-M -1- \delta},\\
  \pour\ M, N \in \N, \ \alpha \in \N^{d+1},  \ |\alpha| \leq 1, \
  \beta \in \N^{d},  \ \delta \in \N,  y \in \R^{d+1}, \ (\eta',\zeta) \in \R^{d+1}.
\end{multline}
\end{proposition}

The decomposition of symbols given in Proposition~\ref{prop: Euclidean
  symbol division} makes tangential symbols appear; they are
introduced in Section~\ref{sec: Tangential symbols and operators}.
Observe that $q_k$ has limited decay in $\zeta$. Yet, the
polyhomogeneous development will be used in what follows. 
%%%% proof of proposition
\begin{proof}
Recall that $p_k(y,\eta)= - \tau^2 +  \sum_{i,j} \gk^{i,j} (x)
  \xi_i\xi_j$ is the principal symbol of $P_{\kk, \gk}$. 
For $\tau \in \supp \chi$, and $\y' =
(y,\eta')$, with $\eta' =(\tau, \xi')$,  having $p_k (\y',\zeta) =0$ reads
\begin{align*}
  \gk^{ij} \xi_i \xi_j = \tau^2 \in (C_{\mu,0}^4, C_{\mu,1}^4), \qquad
  \xi = (\xi',\zeta),
\end{align*}
 meaning that $|\xi'| + |\zeta|  \lesssim C_{\mu,1}^2$ if $y = (t,x)$ remains in
 a bounded domain. Hence for $\supp b \subset K \times
 \R^{d+1}$ with $K$ compact of $\R^{d+1}$ one sees that there exist a bounded domain $L'$
 of $\R^{d-1}$ and $R>0$ 
 \begin{multline}
   \label{eq: localisation zero - zeta}
   y=(t,x',z) \in K, \ \eta' = (\tau,\xi') \in \R^{d}, \ 
   \zeta \in \C,\ \tau \in \supp \chi \ \ \et \ \
   \ p_k(y,\eta',\zeta)=0 \\
   \ \ \imp \ \
   \xi' \in L' \ \et  \ |\zeta| < R.
 \end{multline}
For $r \geq  R$ we will consider the rectangular curve in the complex
plane postively oriented and made with the following pieces
\begin{align*}
  L_{r,R} &= \{ z \in \C; \ -r\leq \Re z\leq r
  \ \et \ \Im z = \pm  R  \} \\
  &\quad \cup  \{ z \in \C; \  \Re z = \pm r \ \et \ 
  -R\leq \Im z  \leq R \},
\end{align*}
that encloses the open ball centred at $0$ with radius $R$. The
important aspect of this contour is that the distance from the real
axis is bounded by $R$ from above allowing one to use the estimation~\eqref{eq: rapid decay symbol SigmaH}

Consider $\tchi \in \Cinfc(\R^{d-1})$ that is  equal to $1$ in  a \nhd of
$L'$.  We decompose  symbols in $\symbH(\R^{2d+2})$ according to Weierstrass preparation Theorem~\cite[Section
7.5]{Hoermander:V1}. Let $b(y,\eta',\zeta) \in \symbH(\R^{2d+2})$.
One may write for
$|\zeta| < r$, with $r\geq R$,
\begin{align*}
 \chi(\tau) \tchi (\xi ') b(y,\eta',\zeta) 
  = \frac{\chi( \tau) \tchi (\xi ') } {2i\pi} \int_{L_{r,R}} 
  b(y,\eta',\tzeta) 
  \frac{d \tzeta}{\tzeta - \zeta}.
\end{align*}
Following the proof of \cite[Theorem 7.5.2]{Hoermander:V1}, using \eqref{eq: localisation zero - zeta}, one 
further writes
\begin{align*}
 \chi(\tau) \tchi (\xi ') b(y,\eta',\zeta) = \tilde{b}_k (y,\eta',\zeta)
  + \tilde{r}_k(y,\eta)\,  p_k(y,\eta),
\end{align*}
with 
\begin{align}
  \label{eq: symbol Euclidean division1}
  \tilde{b}_k (y,\eta',\zeta)
  = \frac{\chi( \tau) \tchi (\xi ')} {2i\pi} 
 \int_{L_{r,R}} \frac{ b( y, \eta', \tzeta )}{ p_k ( y, \eta', \tzeta)} 
\frac {p_k ( y, \eta', \tzeta) - p_k ( y, \eta', \zeta)} {\tzeta -
  \zeta} \, d \tzeta
\end{align}
and 
\begin{align*}
  \tilde{r}_k( y, \eta', \zeta) 
  = \frac{\chi( \tau) \tchi (\xi ')} {2i\pi} \int_{L_{r,R}} 
  \frac{ b( y, \eta', \tzeta )}{ p_k ( y, \eta', \tzeta)} 
  \frac {d \tzeta}{\tzeta - \zeta}.
\end{align*}
Observing that $(p_k ( y, \eta', \tzeta) - p_k ( y, \eta', \zeta))/(\tzeta -
  \zeta)$ is a first-order polynomial in $\zeta$ one finds that
  $\tilde{b}_k$ has the form
\begin{align}
  \label{eq: symbol Euclidean division2}
  \tilde{b}_k (y,\eta',\zeta) = b_{0,k}(y,\eta') + b_{1,k}(y,\eta') \zeta .
\end{align}
It is important to notice that the values of $\tilde{r}_k$, $b_{0,k}$, and 
$b_{1,k}$ are independent of the value of $r$, provided that
$r > |\zeta|$. 
From \eqref{eq: rapid decay symbol SigmaH} and the explicit formula~\eqref{eq: symbol Euclidean division1} 
one
deduces that \eqref{eq: estimation bjk} holds uniformly with respect to $k \in \N$.

Setting 
\begin{align}
  \label{eq: def rk}
  q_k ( y, \eta', \zeta)  = \tilde{r}_k( y, \eta', \zeta)  + (1-
  \tchi)(\xi') \chi(\tau) \frac{b(y,\eta',\zeta)}{ p_k ( y, \eta',
  \zeta)}, 
\end{align}
where the second term is properly defined by \eqref{eq: localisation
  zero - zeta}, one has 
\begin{align}
  \label{eq: symbol decomposition Weierstrass bis}
  \chi(\tau) b(y,\eta',\zeta) = b_{0,k}(y,\eta') + b_{1,k}(y,\eta') \zeta 
  + q_k(y,\eta',\zeta) \,  p_k(y,\eta',\zeta) .
\end{align}
Using that $q_k$ is smooth in the $\eta',\zeta$ variables and that
$p_k(y,\eta',\zeta)$ is invertible for $|(\xi',\zeta)|$ large and
$\tau \in \supp \chi$, with
\eqref{eq: estimation bjk} and \eqref{eq: symbol decomposition
  Weierstrass bis} by induction one
finds that \eqref{eq: estimation rk} holds uniformly with respect to $k \in \N$.

We now consider 
the  polyhomogeneous development of $q_k$ in the $\zeta$ variable.
Observe that the second term on the \rhs of \eqref{eq: def rk} can be
estimated by the remainder in \eqref{polyhomog}.  Hence, it suffices
to consider the term $\tilde{r}_k$. In the support of this term one
has $|\eta'|=|(\xi',\tau)|$ bounded. Observe that it suffices to have
the polyhomogeneous development for $|\zeta|$ large. 
With \eqref{eq: localisation zero
  - zeta}, if $|\zeta| \geq R$ one has $p_k(y, \eta) \neq 0$ and one
can write
\begin{align}
  \label{eq: tilde rk bis}
  \tilde{r}_k = \frac{\chi(\tau) \tchi (\xi ') b(y,\eta',\zeta)
  -  b_{0,k}  (y,\eta') - b_{1,k}  (y,\eta') \zeta}{p_k(y,\eta)},
\end{align}
and $p_k(y,\eta)$ takes the form
\begin{align*}
 p_k(y,\eta)= \big(\zeta - \rho(y,\eta')\big) \big(\zeta - \rho'(y,\eta')\big),
\end{align*}
with the two roots having the same regularity as the coefficients in
the $x$ variable and homogeneous of degree one in $\eta'$, a classical
result based on the Rouch\'e theorem; see for instance \cite[Section 6.A]{LRLR:V2}.  Observe that the first term on
the \rhs of \eqref{eq: tilde rk bis} can be estimated by the remainder
in \eqref{polyhomog}.
For the other terms one writes
\begin{align*}
  \frac{b_{0,k}  (y,\eta') +b_{1,k}  (y,\eta') \zeta}{p_k(y,\eta)}
  = \frac{b_{0,k}  (y,\eta')/\zeta^2  +b_{1,k}  (y,\eta') /\zeta}
  {\big( 1 - \rho(y,\eta')/ \zeta\big) \big( 1- \rho'(y,\eta') / \zeta\big)}.
\end{align*}
Since here $|\eta'|$ is bounded and $y$ remains in a compact domain,
for $|\zeta|$ \suff large one obtains the sought polyhomogeneous
development with a truncated Neumann series. 
  \end{proof}
 
%%%%%%%%%%%%%%%%%%%%%%%%
% remark               %
%%%%%%%%%%%%%%%%%%%%%%%%
\begin{remark}
  \label{remark: Euclidean division of symbol}
  Recall that the symbol $p_k(y,\eta)$ is in fact smooth in $t$ since
  independent of $t$. Hence, estimates \eqref{eq: symbol decomposition Weierstrass}-\eqref{polyhomog} remain valid
  with an arbitrary number of derivatives in $t$. This is however not
  needed in what follows.
\end{remark}

%%%%%%%%%%%%%%%%%
% subsection
\subsection{Low regularity/low conormal decay symbolic calculus}

The limited smoothness with respect to the $x$ variable of the symbols
obtained in Proposition~\ref{prop: Euclidean
  symbol division}  and their limited decay in the $\zeta$ variable force us to
investigate  the symbolic calculus properties of  operators with low
regularity ($W^{1,\infty}$) that we will need in what follows.

\bigskip
In \eqref{eq: estimation rk} and  \eqref{eq: estimation rjk} we have found symbol estimates with distinct decays in the
variables $\eta'$ and $\zeta$.
For a symbol   $a(y,\eta)$ we thus set 
\begin{align}
  \label{eq: estimate full symbol}
  \tilde{N}_{\ell} (a)&=\max_{|\alpha|\leq \ell} \esssup_{(y, \eta)}
                    \big| \d^\alpha_{\eta} a_2(y, \eta)\big|
\jp{\zeta}^{2} \jp{\eta'}^{d+1} .
\end{align}
Observe the difference with $M_{0,\ell}^{-(d+1)}(a)$ in
\eqref{eq: est symbols}.
%%%%%%%%%%%%%%%%%%%%%%%%
% lemma                %
%%%%%%%%%%%%%%%%%%%%%%%%
\begin{lemma}
  \label{lemma: bounded}
  Let $\chi \in \Cinfc(\R)$ be equal to 1 in a \nhd of $0$. 
  Consider a symbol $a(y, \eta)$ that is compactly supported in the
  $y$ variable and of the form
  \begin{equation*}
    %\label{asympt1}
    a(y,\eta)= a_0(y, \eta')
    +a_1(y, \eta')\frac {(1- \chi( \zeta))} {\zeta}
    + a_2(y,\eta),
 \end{equation*}
with $a_j \in \symbt(\R^{d+1}\times \R^d)$, $j=0,1$, and $\tilde {N}_{d+2}(a_2) <
+\infty$.
Then the operator $a(y, h D_y)$ is bounded on $L^2(\R^{d+1})$ and  
\begin{equation*}
%\label{calcsymb1}
\Norm{a(y, h D_y)}{\mathcal{L}(L^2(\R^{d+1}))}
\leq C (N_{d+1}  (a_0)+N_{d+1}  (a_1) +\tilde{N}_{d+2} (a_2)).
\end{equation*}
\end{lemma}
\setcounter{lemmabiss}{\arabic{theorem}}
\begin{proof}
The estimate of the contribution associated with  $a_0$ is given by Lemma~\ref{lemma: continuity tangential op}.
The contribution associated with  $a_1$ is given by \eqref{eq: est b f
  (zeta)}, Lemma~\ref{lemma: continuity tangential op} and \ref{eq: est  f (zeta)}.

 The contribution of the symbol $a_2$ is dealt with by using the same method as in
 the proof of Lemma~\ref{4.continuite}. In fact, the kernel of $a_2(y,
 h D_y) $ is given by $K(y,\ty) = h^{-d-1} k(y,(y-\ty)/h)$ with 
\begin{align*}
  k(y,v) = ( 2\pi)^{-d-1} \int_{\R^{d+1}} e^{ i v \cdot \eta } a_2(y, \eta) d \eta 
 = ( 2\pi)^{-d-1} \int_{\R^{d+1}} e^{ i v \cdot \eta } (\transp L)^{d+2} a_2(y, \eta) d \eta,
\end{align*}
with $L =(1 -i v  \cdot \nabla_\eta)/\jp{v}^2$ and $\transp L= (1 +i v
  \cdot \nabla_\eta)/\jp{v}^2$ since $L \exp(i v\cdot 
  \eta ) = \exp(i v\cdot \eta ) $. 
Using \eqref{eq: estimate full symbol} and that $\jp{\zeta}^{-2}
\jp{\eta'}^{-(d+1)}$ is integrable one finds
\begin{align*}
  |k(y,v)| \lesssim \tilde{N}_{d+2} (a_2) \jp{v}^{-(d+2)} 
  \int_\R \jp{\eta'}^{-d-1} \jp{\zeta}^{-2} d \eta' d \zeta  
  \lesssim \tilde{N}_{d+2} (a_2) \jp{v}^{-(d+2)}.
\end{align*}
One concludes with Corollary~\ref{cor: Schur modified}.
\end{proof}

An inspection of the part of the proof of Lemma~\ref{lemma: bounded}
dedicated to the term $a_1(y,\eta')$ shows that
multiplying $a_0(y,\eta')$ and $a_1(y,\eta')$ by a uniformly bounded
function of $\zeta$ leaves the result unchanged. 
%%%%%%%%%%%%%%%%%%%%%%%%
% lemma                %
%%%%%%%%%%%%%%%%%%%%%%%%
\begin{lemmabis}
  \label{lemmabis: bounded}
  Let $\chi \in \Cinfc(\R)$ be equal to 1 in a \nhd of $0$.
  Let $m(\zeta, h)$ be a bounded function uniformly with respect to $h>0$.
  Consider a symbol $a(y, \eta)$ that is compactly supported in the
  $y$ variable and of the form
  \begin{equation*}
    %\label{asympt1-bis}
    a(y,\eta, h)= a_0(y, \eta') m(\zeta, h)
    +a_1(y, \eta')\frac {(1- \chi( \zeta)) m(\zeta, h)} {\zeta}
 \end{equation*}
with $a_j \in \symbt(\R^{d+1}\times \R^d)$, $j=0,1$.
Then,
the operator $a(y, h D_y, hx)$ is bounded on $L^2(\R^{d+1})$ and  
\begin{equation*}
%\label{calcsymb1-bis}
\Norm{a(y, h D_y, h)}{\mathcal{L}(L^2(\R^{d+1}))}
\leq C (N_{d+1}  (a_0)+N_{d+1}  (a_1)).
\end{equation*}
\end{lemmabis}

\bigskip
%%%%%%%%%%%%%%%%%%%%%%%%
% lemma                %
%%%%%%%%%%%%%%%%%%%%%%%%
\begin{lemma} 
  \label{lemma: bounded 2}
Let $\chi \in \Cinfc(\R)$ be equal to 1 in a \nhd of $0$. 
  Consider a symbol $a(y, \eta)$ that is compactly supported in the
  $y$ variable and of the form
  \begin{equation*}
    %\label{asympt2}
    a(y,\eta)= a_0(y, \eta')
    +a_1(y, \eta')\frac {(1- \chi( \zeta))} {\zeta}
    + a_2(y,\eta).
  \end{equation*}
  \begin{enumerate}
    \item Assume that $N_{d+2} (a_j) < \infty$, $j=0,1$, that
      is,  $a_j\in \Symbolt{0,d+2}{-d-1}{d+1}{d}$, and
      $\tilde{N}_{d+3}(a_2) < \infty$. Then, 
      for $\theta \in W^{1,\infty} (\R_y^{d+1})$, one has 
\begin{align}\label{calcsymb3}
  \bigNorm{[a(y, h D_{y}), \theta]}{\mathcal{L}(L^2(\R^{d+1}))}
  \leq Ch \big(N_{d+2}(a_0)+N_{d+2}(a_1) + \tilde{N}_{d+3}(a_2) \big)
\Norm{\theta}{W^{1, \infty}}.
\end{align}
\item Assume that  $N_{d+2} (\nabla_{y'} a_0) < \infty$, $N_{d+2} (\nabla_{y}
  a_1) < \infty$,  and
$\tilde{N}_{d+3} (\nabla_{y} a_2) < \infty$. Then, 
\begin{multline}\label{calcsymb2}
  \bigNorm{ \bar{a}(y, h D_y)^* - a(y, h D_y)}{\mathcal{L}(L^2(\R^{d+1}))}\\
  \leq C h \big(
   N_{d+2} (\nabla_{y'} a_0)+N_{d+2} (\nabla_{y} a_1)+\tilde{N}_{d+3} (\nabla_{y} a_2)\big).
 \end{multline}
 The adjoint is understood with respect to the inner product
 $L^2(\R^{d+1},d x dt)$.
 \end{enumerate}
\end{lemma}
\setcounter{lemmabiss}{\arabic{theorem}}
%%%%%%%%
\begin{proof}
First, we consider  the contribution of $a_2$  to the  commutator. The
kernel of the commutator is then given by 
$K_2(y, \ty) = h^{-(d+1)} k_2(y, (y-\ty)/h)$
\begin{align*}
  k_2 (y, v)  &=  (2\pi)^{-(d+1)} \int_{\R^{d+1}} e^{ i v \cdot \eta}
                \big(\theta(y -h v) -\theta(y)\big)   a_2 (y, \eta)
                d \eta\\
  &= (2\pi)^{-(d+1)} \int_{\R^{d+1}} e^{ i v \cdot \eta}
                \big(\theta(y -h v) -\theta(y)\big)  
               (\transp L)^{d+3} a_2 (y, \eta)
                d \eta,
\end{align*}
with $L =(1 -i v \cdot \nabla_\eta)/\jp{v}^2$ and $\transp L= (1 +i v
  \cdot \nabla_\eta)/\jp{v}^2$ since $L \exp(i v\cdot 
  \eta ) = \exp(i v\cdot \eta )$. 
Using that $\ell(y, h v) = \big(\theta(y -h v) - \theta(y)\big)/ \Norm{h v}{}$ 
is bounded one can write 
\begin{align*}
  k_2 (y, v)  &=  h (2\pi)^{-(d+1)} \int_{\R^{d+1}} e^{ i v \cdot \eta}
                 \ell(y, h v)
              \Norm{v}{} (\transp L)^{d+3} a_2 (y, \eta)
                d \eta.
\end{align*}
With the form of $\transp L$ and \eqref{eq: estimate full symbol} one
obtains
\begin{align*}
  |k_2(y,v)| &\lesssim h \tilde{N}_{d+3} (a_2) \Norm{\theta}{W^{1,\infty}}\jp{v}^{-(d+2)} 
  \int_{\R^{d}} \jp{\eta'}^{-d-1} d \eta' \int_{\R}  \jp{\zeta}^{-2} d \zeta  \\
  &\lesssim h \tilde{N}_{d+3} (a_2) \Norm{\theta}{W^{1,\infty}}\jp{v}^{-(d+2)}.
\end{align*}
One concludes with Corollary~\ref{cor: Schur modified} as for the
proof of Lemma~\ref{lemma: continuity tangential op}.

\medskip
Second, we consider  the contribution of $a_0$ to the commutator. Since
$[\OpH(a_0),\theta]$ is tangential one can consider its action 
in the $y'$ variable only. As in  \eqref{eq: tangential kernel
  1}--\eqref{eq: tangential kernel 2} one writes
\begin{align*}
  \OpH(a_0) u(y',z) = \int_{\R^d} K_{a_0}(y',\ty'; z)\,  u
  (\ty', z)\,  d \ty',
\end{align*}
with 
\begin{align*} 
K_{a_0}(y',\ty'; z)= (2 \pi)^{-d} \int_{\R^d} e^{i (y'- \tilde{y'})\cdot \eta'}
   a_0(y',z, h \eta')\,  d \eta',
\end{align*}
with $z$ as a parameter.
The associated tangential kernel for the commutator is 
\begin{align*}
  K_0(y', \ty'; z)  &= K_{a_0}(y', \ty'; z)
  \big( \theta(\ty',z) - \theta(y', z)\big)
  =h^{-d} k_0\big(y',     (y' - \ty')/h; z\big),
\end{align*}
with 
\begin{align*}
  k_0(y', v; z) 
  &= (2\pi)^{-d} \int_{\R^d} e^{ i v \cdot \eta'}
  \big(\theta( y' -h v,z) -\theta(y',z)\big) 
  a_0 (y, \eta') \, d \eta'.
\end{align*}
Note that here $v\in \R^d$. With the same argument as above one finds
\begin{align*}
  k_0(y', v; z) 
  &=h  (2\pi)^{-d} \int_{\R^d} e^{ i v \cdot \eta'}
  \ell'(y', h v)  \Norm{v}{}
  (\transp L_0)^{d+2} a_0 (y, \eta') \, d \eta'.
\end{align*}
with 
$\ell'(y', h v,z) = \big(\theta(y' -h v,z) - \theta(y',z)\big)/
\Norm{h v}{}$
and  $\transp L_0= (1 +i v \cdot \nabla_{\eta'})/\jp{v}^2$ yielding
\begin{align*}
  |k_0(y', v; z)| &\lesssim h N_{d+2}(a_0) \Norm{\theta}{W^{1,\infty}}\jp{v}^{-(d+1)} 
   \int_{\R^d} \jp{\eta'}^{-d-1} d \eta' \\
  &\lesssim h N_{d+2}(a_0) \Norm{\theta}{W^{1,\infty}}\jp{v}^{-(d+1)}.\notag
\end{align*}
One concludes with Corollary~\ref{cor: Schur modified}.

\medskip
Third, we consider  the contribution of $a_1$ to the commutator. Set
$f(\zeta) = \big(1- \chi( \zeta)\big)/  \zeta$. As observed in
\eqref{eq: decomp b f (zeta)} one has 
$\OpH(a_1 f(\zeta)) = \OpH(a_1) f(h D_z)$ allowing one to write 
\begin{align*}
  [\OpH(a_1 f(\zeta)), \theta] = \OpH(a_1) [ f(h D_z), \theta]
  + [\OpH(a_1), \theta] f(h D_z).
\end{align*}
By \eqref{eq: est  f (zeta)}  one has $\Norm{ f (h D_z)}{\L(L^2(\R^{d+1}))}
\lesssim 1$ and $[\OpH(a_1), \theta] \lesssim h N_{d+2}(a_1)
\Norm{\theta}{W^{1,\infty}}$ similarly to the treatment of the term associated
with $a_0$, yielding
\begin{align*}
  \bigNorm{  [\OpH(a_1), \theta] f(h D_z) }{\L(L^2(\R^{d+1}))}
  \lesssim h N_{d+2}(a_1)  \Norm{\theta}{W^{1,\infty}}.
\end{align*}
With \eqref{eq: kernel multiplier 1}--\eqref{eq: kernel multiplier 2}
the commutator $[f(h D_z),\theta]$ has tangential kernel acting only
in the $z$ variable
\begin{align*}
  K(y; z, \tilde{z}) = h^{-1} k \big(y; z,(z - \tilde{z})/h\big)
\end{align*}
with 
\begin{align*}
 k (y; z,v)= 
(2 \pi)^{-1} \int_{\R} e^{i v \zeta}
 \big(\theta(y', z -h v) -\theta(y',z)\big) 
 f(\zeta)\,  d \zeta,
\end{align*}
where $v\in \R$ here. One writes 
\begin{align*}
 k (y; z,v)= 
h (2 \pi)^{-1} \int_{\R} v e^{i v \zeta}
 \ell_d (y',  z, h v)  f(\zeta)\,  d \zeta,
\end{align*}
  with $\ell_d (y',
  z, h v) = \big(\theta(y', z -h v) -\theta(y',z)\big) /
(h v)$ with $| \ell_d (y',
  z, h v) | \leq \Norm{\theta}{W^{1,\infty}}$. 
Since $v e^{i v \zeta} = -i \d_\zeta e^{i v\zeta}$, with an
integration by parts, one
finds
\begin{align*}
 k (y; z,v)= 
i h (2 \pi)^{-1} \int_{\R} e^{i v \zeta}
 \ell_d (y',  z, h v)  \d_\zeta f(\zeta)\,  d \zeta.
\end{align*}
Moreover with $\transp L_f= (1 +i v \d_{\zeta})/\jp{v}^2$
 one writes
\begin{align*}
 k (y; z,v)= 
i h (2 \pi)^{-1} \int_{\R} e^{i v\zeta}
 \ell_d (y',  z, h v)  (\transp L_f)^2 \d_\zeta f(\zeta)\,  d \zeta.
\end{align*}
Since $|(\transp L_f)^2 \d_\zeta  f(\zeta)| \lesssim \jp{v}^{-2}\jp{\zeta}^{-2}$
one finds 
$|k (y; z,v)| \lesssim h \Norm{\theta}{W^{1,\infty}} \jp{v}^{-2}$,
implying 
\begin{align*}
   \bigNorm{[ f(h D_z), \theta]}{\mathcal{L}(L^2(\R^{d+1}))}
  \lesssim h \Norm{\theta}{W^{1,\infty}}.
\end{align*}
With Lemma~\ref{lemma: continuity tangential op} one obtains
\begin{align*}
  \bigNorm{ \OpH(a_1) [ f(h D_z), \theta]}{\L(L^2(\R^{d+1}))}
  \lesssim h N_{d+1}(a_1)  \Norm{\theta}{W^{1,\infty}} 
  \lesssim h N_{d+2}(a_1)  \Norm{\theta}{W^{1,\infty}}.
\end{align*}
This concludes the proof of the estimation of the commutator norm.

\bigskip
We now turn to the proof of the estimate for the adjoint. We will observe
that the proof is in fact along the same lines as that for the commutator. We
start with the contribution of $a_2(y,\eta)$.  The kernel
of the operator $\OpH(\bar{a}_2)^* - \OpH(a_2)$ is given by 
$K(y,\ty) = h^{-d-1} k \big(y, (y-\ty)/h\big)$ with 
\begin{align*}
k(y,v) = (2\pi)^{-d-1} \int_{\R^{d+1}} e^{i v
  \cdot \eta } \big( a_2(y - h v,\eta)  - a_2(y,\eta)\big) d \eta,
\end{align*}
with $v \in \R^{d+1}$, 
and one writes 
\begin{align*}
  k(y,v) = (2\pi)^{-d-1} \int_{\R^{d+1}} e^{i v
  \cdot \eta } (\transp L)^{d+3} \big( a_2(y - h v,\eta)  -  a_2(y,\eta)\big) d \eta,
\end{align*}
with $\transp L= (1 +i v
  \cdot \nabla_\eta)/\jp{v}^2$. Since 
\begin{align*}
 \big|  (\transp L)^{d+3} \big( a_2(y - h v,\eta)  -
  a_2(y,\eta)\big)\big| 
  &\leq h \Norm{v}{} 
  \bigNorm {\nabla_y (\transp L)^{d+3}  a_2(y,\eta)}{L^\infty}\\
  &\lesssim  h \jp{v}^{-d-2} \tilde{N}_{d+3} (\nabla_y a_2) \jp{\eta'}^{-d-1} \jp{\zeta}^{-2}.
\end{align*}
The sought estimate follows.

\medskip
For the contribution of the symbol $a_0(y,\eta)$,  the {\em tangential} kernel
of the operator $\OpH(\bar{a}_0)^* - \OpH(a_0)$ is given by 
$K(y',\ty') = h^{-d} k \big(y', (y'-\ty')/h\big)$ with 
\begin{align*}
k(y',v) = (2\pi)^{-d} \int_{\R^{d}} e^{i v
  \cdot \eta' } \big( a_0(y' - h v,\eta')  - a_0(y',\eta')\big) d \eta',
\end{align*}
with here $v\in \R^d$, 
and one writes 
\begin{align*}
  k(y',v) = (2\pi)^{-d} \int_{\R^{d}} e^{i v
  \cdot \eta' } (\transp L)^{d+2} \big( a_0(y' - h v,\eta')  -  a_0(y',\eta')\big) d \eta',
\end{align*}
with $\transp L= (1 +i v
  \cdot \nabla_{\eta'})/\jp{v}^2$ and one finds similarly
\begin{align*}
 \big|  (\transp L)^{d+2} \big( a_0(y - h v,\eta)  -
  a_0(y,\eta)\big)\big| 
  &\lesssim  h \jp{v}^{-d-1} N_{d+2} (\nabla_{y'} a_0) \jp{\eta'}^{-d-1}.
\end{align*}
The sought estimate follows.

\medskip
For the contribution of $a_1(y,\eta)$ one writes
\begin{align*}
  &\OpH \big(\bar{a}_1 \bar{f}(\zeta) \big)^* - \OpH \big ( a_1
    f(\zeta)\big)\\
  &\qquad = f(h D_z) \bar{a}_1 (y, h D_y)^* - a_1 (y, h D_y) f(h D_z) \\
  &\qquad = f(h D_z) \big( \bar{a}_1 (y, h D_y)^* - a_1 (y, h D_y) \big)
  + [f(h D_z) , a_1 (y, h D_y)],
\end{align*}
using that $\bar{f}(h D_z)^* =  f(h D_z)$.
With \eqref{eq: est  f (zeta)}   and applying the argument made for the term associated
with $a_0$ one finds 
 \begin{align*}
  \bigNorm{ f(h D_z) \big( \bar{a}_1 (y, h D_y)^* - a_1 (y, h D_y) \big) }{\L(L^2(\R^{d+1}))}
  \lesssim h N_{d+2}(\nabla_{y'} a_1).
\end{align*}

We now consider the commutator $[f(h D_z) , a_1 (y, h D_y)]$. 
The kernel of $f(h D_z)$ is given by
\begin{align*}
  K_1(y, \ty) = \delta_{y' - \ty'} \otimes K_d
  (z,\tilde{z}), 
\end{align*}
with
\begin{align*}
  K_d (z,\tilde{z}) 
  = (2\pi)^{-1} \int_\R  e^{i (z-\tilde{z}) \zeta} f(h\zeta) d \zeta.
\end{align*}
The kernel of $a_1 (y, h D_y)$ is given by
\begin{align*}
  K_2(y, \ty) =K' (y',\ty'; z) \otimes  \delta_{z - \tilde{z}} , 
\end{align*}
with
\begin{align*}
  K' (y',\ty'; z)
  = (2\pi)^{-d} \int_{\R^{d}} e^{i (y'-\ty') \cdot\eta'} a_1(y',z, h\eta') d \eta'.
\end{align*}
This leads to the following kernel for the commutator $[f(h D_z) , a_1 (y, h D_y)]$
\begin{align*}
  K(y, \ty) =  \int_{\R^{d+1}} \big( K_1(y,\hy) K_2(\hy,\ty) - K_2(y,\hy) K_1(\hy,\ty) \big) d\hy,
\end{align*}
where the integration is understood in the sense of distribution
action. Products here make sense; see \cite[Theorem 8.2.14]{Hoermander:V1}.
This gives
\begin{multline*}
  K(y, \ty) 
  = (2\pi)^{-d-1} \int_\R  e^{i (z-\tilde{z}) \zeta} f(h\zeta) d \zeta\\
  \times \int_{\R^{d}}  e^{i (y'-\ty') \cdot\eta'} \big( a_1(y',\tz, h\eta') -  a_1(y',z, h\eta') \big) d \eta',
\end{multline*}
that we write $K(y, \ty) = h^{-d-1} k \big(y, (y - \ty)/h\big)$ with 
\begin{align*}
  k(y, v) 
  = (2\pi)^{-d-1} \int_\R  e^{i w \zeta} f(\zeta) d \zeta
  \int_{\R^{d}}  e^{i v' \cdot\eta'} \big( a_1(y',z - h w, \eta') -  a_1(y',z, \eta') \big) d \eta',
\end{align*}
for $v = (v',w)$ with $v' \in \R^{d}$ and $w \in \R$. 
Considering the bounded function
\begin{align*}
  \ell (y, z, h w, \eta') = \big( a_1(y',z - h w, \eta') -  a_1(y',z,
  \eta') \big) / (h w),
\end{align*}
one writes
\begin{align*}
  k(y, v) 
  = h (2\pi)^{-d-1} \int_\R  e^{i w \zeta} w f(\zeta) d \zeta
  \int_{\R^{d}}  e^{i v' \cdot\eta'} \ell (y, z, h w, \eta')  d \eta'.
\end{align*}
Set $L_\zeta = (1 - i w \d_\zeta ) /\jp{w}^2$ and $L_{\eta'} = (1 - i
v' \cdot \d_{\eta'} ) /\jp{v'}^2$. 
With $w e^{i w \zeta}  = - i \d_\zeta e^{i w \zeta}$,  $L_\zeta e^{i w
  \zeta}  = e^{i w \zeta}$, and $L_{\eta'} e^{i v' \cdot\eta'} = e^{i
  v' \cdot\eta'}$ one writes 
\begin{align*}
  k(y, v) 
  =i  h (2\pi)^{-d-1} \int_\R  e^{i w \zeta}  (\transp L_\zeta)^2 \d_\zeta f(\zeta) d \zeta
  \int_{\R^{d}}  e^{i v' \cdot\eta'} (\transp L_{\eta'})^{d+1} \ell (y, z, h w,
  \eta')  d \eta',
\end{align*}
and one finds
\begin{align*}
  |k(y, v) | 
 \lesssim   h  \jp{w}^{-2} \jp{v'}^{-d-1} N_{d+1} (\d_z a_1) \int_\R \jp{\zeta}^{-2} d \zeta
  \int_{\R^{d}} \jp{\eta'}^{-d-1} d \eta',
\end{align*}
since $(\transp L_{\eta'})^{d+1} \ell (y, z, h w,
  \eta') \lesssim N_{d+1} (\d_z a_1) \jp{\eta'}^{-d-1}$. This leads to
  the conclusion of the proof.
\end{proof}

\medskip
Multiplying $a_0(y,\eta')$ and $a_1(y,\eta')$ by a \suff rapidely decaying 
function of $\zeta$ leaves the result {\em nearly} unchanged.
%%%%%%%%%%%%%%%%%%%%%%%%
% lemma                %
%%%%%%%%%%%%%%%%%%%%%%%%
\begin{lemmabis} 
  \label{lemmabis: bounded 2}
  Let $\chi \in \Cinfc(\R)$ be equal to 1 in a \nhd of $0$. 
  Let $m(\zeta, h)$ be a bounded function of $\zeta$ and $h>0$ with moreover 
  \begin{align*}
    \d_\zeta^ j m(\zeta, h) \in L^1_\zeta,  \quad
      \ 1\leq j \leq 3 ,
  \end{align*}
 uniformly with respect to $h>0$.
  Consider a symbol $a(y, \eta)$ that is compactly supported in the
  $y$ variable and of the form
  \begin{equation*}
    %\label{asympt2-bis}
    a(y,\eta)= a_0(y, \eta') m(\zeta, h)
    +a_1(y, \eta')\frac {(1- \chi( \zeta)) m(\zeta, h)} {\zeta}
  \end{equation*}
  \begin{enumerate}
    \item Assume that $N_{d+2} (a_j) < \infty$, $j=0,1$, that
      is,  $a_j\in \Symbolt{0,d+2}{-d-1}{d+1}{d}$.  Then, 
      for $\theta \in W^{1,\infty} (\R_y^{d+1})$, one has 
\begin{multline*}%\label{calcsymb3-bis}
  \bigNorm{[a(y, h D_{y}), \theta]}{\mathcal{L}(L^2(\R^{d+1}))}\\
  \leq Ch \big(N_{d+2}(a_0)+N_{d+2}(a_1) \big)
\Norm{\theta}{W^{1, \infty}} \sup_{1\leq j \leq 3} \| \d_\zeta^j m\|_{L^1}.
\end{multline*}
\item Assume that  $N_{d+2} (\nabla_{y} a_0)$, $N_{d+2} (\nabla_{y}
  a_1)$ are finite. Then one has 
\begin{multline*}%\label{calcsymb2-bis}
  \bigNorm{ \bar{a}(y, h D_y)^* - a(y, h D_y)}{\mathcal{L}(L^2(\R^{d+1}))}\\
  \leq C h \big(
   N_{d+2} (\nabla_{y} a_0)+N_{d+2} (\nabla_{y} a_1)\big)\sup_{1\leq j \leq 3} \| \d_\zeta^j m\|_{L^1}.
 \end{multline*}
 The adjoint is understood with respect to the inner product
 $L^2(\R^{d+1},d x dt)$.
 \end{enumerate}
\end{lemmabis}
Note that 
$N_{d+2} (\nabla_{y'} a_0)$ is replaced by $N_{d+2} (\nabla_{y} a_0)$
in the 
second estimate if compated with Lemma~\ref{lemma: bounded 2}.

%%%%%%%%%%%%%%%%%%%%
\begin{proof}
Considering the properties of $m(\zeta, h)$ in both results, only the
contribution associated with the tangential symbol $a_0(y,\eta')$
needs to be analyzed. 
For the commutator, as for the treatment of the term $a_1$ in the proof of
Lemma~\ref{lemma: bounded 2} one estimates the operator norm of
$[m(h D_z, h),\theta]$. Its tangential kernel is 
\begin{align*}
  K(y; z, \tilde{z}) = h^{-1} k \big(y; z,(z - \tilde{z})/h\big)
\end{align*}
with 
\begin{align*}
 k (y; z,v)= 
(2 \pi)^{-1} \int_{\R} e^{i v \zeta}
 \big(\theta(y', z -h v) -\theta(y',z)\big) 
 m (\zeta, h)\,  d \zeta, \qquad v \in \R. 
\end{align*}
Following the proof of Lemma~\ref{lemma: bounded 2} one obtains 
\begin{align*}
 k (y; z,v)= 
i h (2 \pi)^{-1} \int_{\R} e^{i v \zeta}
 \ell_d (y',  z, h v)  (\transp L_f)^2 \d_\zeta  m (\zeta, h)\,  d \zeta.
\end{align*}
where $\ell_d (y',
  z, h v) = \big(\theta(y', z -h v) -\theta(y',z)\big) /
(h v)$ and $\transp L_f= (1 +i v \d_{\zeta})/\jp{v}^2$.
With the properties of the function $m(\zeta,h)$ one obtains 
\begin{align*}
|k (y; z,v)| \lesssim h \Norm{\theta}{W^{1,\infty}}
  \jp{v}^{-2} \sup_{1\leq j \leq 3} \| \d_\zeta^j m\|_{L^1} 
  \lesssim h \Norm{\theta}{W^{1,\infty}}\sup_{1\leq j \leq 3} \| \d_\zeta^j m\|_{L^1} 
  \jp{v}^{-2},
\end{align*}
implying 
\begin{align*}
   \bigNorm{[ f(h D_z), \theta]}{\mathcal{L}(L^2(\R^{d+1}))}
  \lesssim  h \Norm{\theta}{W^{1,\infty}}\sup_{1\leq j \leq 3} \| \d_\zeta^j m\|_{L^1} .
\end{align*}

\medskip
Considering the argument for the adjoint given in the proof of
Lemma~\ref{lemma: bounded 2} one needs to estimate the operator norm
of $[m(h D_z, h) , a_0 (y, h D_y)]$.  Its kernel reads $K(y, \ty) = h^{-d-1} k \big(y, (y - \ty)/h\big)$ with 
\begin{align*}
  k(y, v) 
  = (2\pi)^{-d-1} \int_\R  e^{i w \zeta} m(\zeta,h) d \zeta
  \int_{\R^{d}}  e^{i v' \cdot\eta'} \big( a_0(y',z - h w, \eta') -  a_0(y',z, \eta') \big) d \eta',
\end{align*}
for $v = (v',w)$ with $v' \in \R^{d}$ and $w \in \R$. One obtains 
\begin{align*}
  k(y, v) 
  =i  h (2\pi)^{-d-1} \int_\R  e^{i w \zeta}  (\transp L_\zeta)^2 \d_\zeta m(\zeta,h) d \zeta
  \int_{\R^{d}}  e^{i v' \cdot\eta'} (\transp L_{\eta'})^{d+1} \ell (y, z, h w,
  \eta')  d \eta',
\end{align*}
with $\ell (y, z, h w, \eta') = \big( a_0(y',z - h w, \eta') -  a_0(y',z,
  \eta') \big) / (h w)$. This leads to 
\begin{align*}
  |k(y, v) | 
 \lesssim   h \jp{w}^{-2} \jp{v'}^{-d-1} N_{d+1} (\d_z a_0)\sup_{1\leq j \leq 3} \| d_\zeta^j m\|_{L^1} 
  \int_{\R^{d}} \jp{\eta'}^{-d-1} d \eta',
\end{align*}
and 
  the conclusion of the proof.
\end{proof}
%%%%%%%%%%%%%%%%%%%%%%%%
% remark               %
%%%%%%%%%%%%%%%%%%%%%%%%
\begin{remark}
  \label{remark: additional multiplier cutoff}
A particular choice of Fourier multiplier $m(\zeta, h)$ appearing in
Lemmata~\ref{lemmabis: bounded} and \ref{lemmabis: bounded 2} is 
$m(\zeta, h) = \varphi( h^\beta \zeta)$ for $\varphi \in \mathscr
S(\R; \R)$ and some $\beta>0$.  
Indeed, the following properties hold:
\begin{enumerate}
  \item One has $m(\zeta, h) \lesssim 1$ uniformly with respect to $h>0$ as
    required by Lemmata~\ref{lemmabis: bounded} and \ref{lemmabis: bounded 2}.
  \item  One has obviously
  \begin{align*}
    \sup_h \sup_{1\leq j \leq 3} \| \d_\zeta^j m\|_{L^1} <+\infty.
  \end{align*}
    \end{enumerate}
In what follows, we will use $m(\zeta,h) = \varphi(h^3 \zeta)$ with
$\varphi \in \Cinfc(\R)$. 
\end{remark}

%%%%%%%%%%%%%%%%%
% subsection
\subsection{A class of error terms}
%%%%%%%%%%%%%%%%%%%%%%%%
% definition           %
%%%%%%%%%%%%%%%%%%%%%%%%
\begin{definition}
  \label{def: error operators}
  Let $\mathcal{R}$ be the class of sequences of operators $(R_k)_{k}$
  bounded on $L^2(\R^{d+1})$ by $C \hk^\delta$, $\delta\geq 0$, and from
  $L^2(\R^{d+1})$ to $H^1(\R_z; L^2(\R^d_{y'}))$ by $C\hk^{\rho}$,
  $\rho\geq 0$, with moreover $\delta+\rho >0$. Denote by
  $\mathcal{R}_0$ the class obtained in the case  $(\delta, \rho) =
  (1,0)$. 
\end{definition}
%%%%%%%%%%%%%%%%%%%%%%%%
% lemma                %
%%%%%%%%%%%%%%%%%%%%%%%%
\begin{lemma}
  \label{lem610}
Let $(f_k)_k$ be a bounded sequence of $L^2(\R^{d+1})$ and $(g_k)_k$
be a bounded sequence of $L^2_{y'}(\R^{d})$. Then, if $(R_k)_k \in
\mathcal{R}$ one has 
\begin{align*}
  \lim_{k\to+\infty}  \bignorm{ \biginp{g_k \otimes \delta_{z=0}}
  {R_k f_k}_{H_z^{-1}L_{y'}^2,H_z^{1}L_{y'}^2}}{} =0.
\end{align*}
\end{lemma}
\begin{proof}
As the sequence $g_k \otimes \delta_{z=0}$ is bounded in $H^{-\sigma}(
\R_z; L^2(\R^{d}_{y'}))$, for any $\sigma> 1/2$, the result follows
from the bound
\begin{align*}
  \Norm{R_k}{\L (L^2(\R^{d+1}), H_z^{\sigma}L_{y'}^2)} 
  \leq C \hk ^{\frac{\delta + \rho} 2 + (\sigma- \frac 1 2) ( \rho-
  \delta)},
\end{align*}
obtained by interpolation 
and choosing $\sigma >\frac12$ \suff close to $\frac12$.
\end{proof}
%%%%%%%%%%%%%%%%%%%%%%%%
% corollary            %
%%%%%%%%%%%%%%%%%%%%%%%%
\begin{corollary}
  \label{cor: Lkprime}
Let $b \in
   \Symbolo{\infty,\infty}{-\infty}{2d+2}$ with $\supp b\subset
 K \times \R^{d+1}$, for $K$ a compact of $\cdiffL(\O)$, and let $\psi \in
 \Cinfc(\cdiffL(\O))$ be equal to $1$ in a \nhd of the $y$-projection of
 $\supp b$. Let $L_k (b, \psi)$ be as defined in \eqref{eq: Lk}. One
 has $L_k (b,\psi)  =  L_k' (b,\psi) + o(1)_{k \to +\infty}$ with 
\begin{align}
  \label{eq: Lkprime}
   L_k' (b,\psi) &=    i \biginp{ \OpH(b) \psi u_k}
             { v_k \otimes \delta_{z=0}}_{H_z^{1}L_{y'}^2, H_z^{-1}L_{y'}^2}^{\kk\mugk dt}
            \\
           &\quad
             - i \biginp{ v_k \otimes \delta_{z=0}}
             {\OpH(\bb) \psi u_k}_{H_z^{-1}L_{y'}^2,
            H_z^{1}L_{y'}^2}^{\kk\mugk dt}.\notag
 \end{align}
\end{corollary}
\begin{proof}
  One has to prove that 
  \begin{align*}
  I_k &= \biginp{ v_k \otimes \delta_{z=0}}   
    {\psi\,   \OpH(b)^\star u_k}
    _{H_z^{-1}L_{y'}^2,H_z^{1}L_{y'}^2}^{\kk\mugk dt}\\
    &= \biginp{ v_k \otimes \delta_{z=0}}   
    {\OpH(\bb) \psi u_k}
    _{H_z^{-1}L_{y'}^2,H_z^{1}L_{y'}^2}^{\kk\mugk dt} + o(1)_{k \to +\infty}.
 \end{align*}
 Let $\tpsi \in \Cinfc(\cdiffL(\O))$ be equal to 1
in a \nhd  of $\supp \psi$. One has 
\begin{align*}
  I_k = \biginp{ \tpsi v_k \otimes \delta_{z=0}}   
    {\psi\,   \OpH(b)^\star \tpsi u_k}
    _{H_z^{-1}L_{y'}^2, H_z^{1}L_{y'}^2}^{\kk\mugk dt}.
 \end{align*}
 The adjoint operator with the $\star$-notation is here understood in the sense of the 
inner product $L^2(\R^{d+1}, \kk\mugk dt)$; see Proposition~\ref{prop:
  commutator}. Thus, $\OpH(b)^\star =  \tkk^{-1} \OpH(b)^* \tkk$
where the adjoint with the usual $*$-notation is understood for the inner
product $L^2(\R^{d+1}, dx dt)$.

Since $(\tpsi u_k)_k$ and  $(\tpsi v_k)_k$ are bounded in
 $L^2(\R^{d+1})$ and in $L^2_{y'}(\R^{d})$ respectively, it suffices
 to prove that $\psi\,   \tkk^{-1} \OpH(b)^* \tkk- \OpH(\bb) \psi
 \in \mathcal{R}_0$ by Lemma~\ref{lem610}. One has 
\begin{align*}%\label{prems}
 \psi\,   \tkk^{-1} \OpH(b)^* \tkk- \OpH(\bb) \psi 
 =  \psi\,   \tkk^{-1} \big(\OpH(\tkk b)^*  - \OpH(\tkk  \bb)\big) + [ \psi, \OpH(\bb)].
\end{align*}
From~\eqref{calcsymb3} and~\eqref{calcsymb2} one deduces that 
\begin{align*}
 %\label{eq: cor: Lkprime1}
  \Norm{\psi\,   \tkk^{-1} \OpH(b)^* \tkk- \OpH(\bb)
  \psi}{\L(L^2(\R^{d+1}))}
  \lesssim \hk. 
\end{align*}
To estimate the operator norm from $L^2(\R^{d+1})$ to $H^1(\R_z;
L^2(\R^d_{y'}))$ we compose with $D_z = \hk^{-1} \OpH(\zeta)$ and get
\begin{align*}
  & D_z \big( \psi\,   \tkk^{-1} \OpH(\bb)^* \tkk- \OpH(b) \psi
  \big)\\
  &\quad 
    =  [ D_z, \psi\,   \tkk^{-1}] \OpH(b)^* \tkk
    + \psi\,   \tkk^{-1}  D_z \OpH(b)^* \tkk
    -  [ D_z, \OpH(\bb)] \psi 
    - \OpH(\bb) D_z \psi\notag\\
   &\quad  
     =  D_z \big( \psi\,   \tkk^{-1} \big) \OpH(b)^*\tkk 
     + 
     \hk^{-1} \psi\,   \tkk^{-1} \OpH(b\zeta)^* \tkk
     - \OpH(D_z \bb ) \psi 
     - \hk^{-1} \OpH(\zeta \bb ) \psi\notag.
\end{align*}
One has 
\begin{align*}
  \Norm{
  D_z \big( \psi\,   \tkk^{-1} \big) \OpH(b)^*\tkk
  }{\L(L^2(\R^{d+1}))}
  + \Norm{ \OpH(D_z \bb) \psi }{\L(L^2(\R^{d+1}))}
  \lesssim 1. 
\end{align*}
It thus remains to prove that 
\begin{align*}
  \Norm{
  \psi\,   \tkk^{-1} \OpH(b\zeta)^* \tkk 
  - \OpH(\zeta \bb) \psi
  }{\L(L^2(\R^{d+1}))}
  \lesssim \hk. 
\end{align*}
The argument is the same as for estimate \eqref{cor: Lkprime} with $b$
replaced by $\zeta b$. We conclude using Lemma~\ref{lem610} with $\delta =1$, $\rho =0$.
\end{proof}

%%%%%%%%%%%%%%%%%%%%%%%%
% corollary            %
%%%%%%%%%%%%%%%%%%%%%%%%
\begin{corollary}
  \label{cor: Lsecond}
Let $b \in
   \Symbolo{\infty,\infty}{-\infty}{2d+2}$ with $\supp b\subset
 K \times \R^{d+1}$, for $K$ a compact of $\cdiffL(\O)$, and let $\psi \in
 \Cinfc(\cdiffL(\O))$ be equal to $1$ in a \nhd of the $y$-projection of
 $\supp b$. 
 Let also $\varphi \in \Cinfc(]-2,2[; \R)$ and equal to $1$ on $(-1,1)$.
 Let $L_k (b, \psi)$ be as defined in \eqref{eq: Lk}. One
 has $L_k (b,\psi)  =  L_k''(b,\psi) + o(1)_{k \to +\infty}$ with 
\begin{align*}
  %\label{eq: Lksecond}
   L_k'' (b,\psi) &=    i \biginp{ \OpH(b) \varphi(\hk^3 D_z)\psi u_k}
             { v_k \otimes \delta_{z=0}}_{H_z^{1}L_{y'}^2, H_z^{-1}L_{y'}^2}^{\kk\mugk dt}
             \\
           &\quad
             -i \biginp{ v_k \otimes \delta_{z=0}}
             {\OpH(\bb) \varphi(\hk^3 D_z) \psi u_k}_{H_z^{-1}L_{y'}^2,
             H_z^{1}L_{y'}^2}^{\kk\mugk dt}.\notag
 \end{align*}
\end{corollary}
\begin{proof}
  Arguing as for Corollary~\ref{cor: Lkprime}, starting from the form
  of $L_k' (b,\psi)$ given in \eqref{eq: Lkprime}  it suffices to prove that
  \begin{align*}
    \OpH(b) \big( 1- \varphi(\hk^3 D_z)\big) = \OpH( \gamma_{\hk}) \in
    \mathcal R_0, 
    \qquad \gamma_{\hk}(y,\eta) = b(y,\eta) \big(1 - \varphi(\hk^2
  \zeta)\big)
\end{align*} 
In the support of $1 - \varphi(\hk^2\zeta)$ one has $\hk^2
|\zeta|\gtrsim 1$, which combined with the fast decay of $b$ in $\eta$
yields
\begin{align*}
  |\d_y^\alpha\d_\eta^\beta \gamma_{\hk}(y, \eta) | \lesssim \hk^N
  \jp{\eta}^{-N}, 
\end{align*}
for any $N$. The result follows from Lemma~\ref{lemma: smoothing semiclassical operators}.
\end{proof}

%%%%%%%%%%%%%%%%%
% section
%%%%%%%%%%%%%%%%%
\section{Proof of the propagation equation II: symbol quantization}
\label{sec: Proof of the propagation equation2}

From the support property of the semi-classical measure $\mu$ given in \eqref{eq: hyp support mesure-new-mu}
if  considering the
action of $\mu$ on a symbol in $\symbH(\R^{2d+2})$ it
suffices to  work with a time-frequency
truncated version. That is, for  $\chi\in \Cinfc ( C_{\mu,0}^2, C_{\mu,1}^2)$  equal to $1$ on a
\nhd of  $[C_{\mu,0}, C_{\mu,1}]$ and $b(y, \eta)  \in
\symbH(\R^{2d+2})$ one has $\dup{\mu}{(1-\chi)b} =0$, meaning that  
\begin{align*}
  \lim_{k\to +\infty} L_k'' \big((1-\chi (\tau))  b,\psi \big)
  = \lim_{k\to +\infty} L_k \big((1-\chi (\tau))  b,\psi \big)
  =0.
\end{align*}
With Proposition~\ref{prop: commutator}, we will thus only consider the action of $\mu$ on a symbol of the
form $\chi(\tau) b(y, \eta)$ through the limit of $L_k \big(
\chi (\tau)  b,\psi \big)$ 
and we will now quantize the Euclidean division of Proposition~\ref{prop: Euclidean symbol division}. Even though the symbol $b$ on the \lhs of~\eqref{eq:
  symbol decomposition Weierstrass} exhibits rapid decay in the
variable $\zeta$, it is
not the case for the symbols $b_{0,k}, b_{1,k}$, and $q_k$ on the \rhs
of~\eqref{eq: symbol decomposition Weierstrass}. Following
\cite{GL:1993},  adding a cutoff in the $\zeta$
variable in the form of $\varphi(\hk^3 D_z)$, made possible by
Corollary~\ref{cor: Lsecond},  acts as a
remedy. 

Since $L_k(., \psi)$ and $L_k''(., \psi)$ have the same limit as
$k \to \infty$ by Corollary~\ref{cor: Lsecond}, 
in what follows, we will study sequentially the limits of $L_k''(a,
\psi)$ as $k \to +\infty$ with $a  (y, \eta)= q_k
p_k (y, \eta)$, $a(y,\eta') =b_{0,k}  (y,\eta')$,  and $a  (y, \eta)=
b_{1,k} (y,\eta') \zeta$. 

%%%%%%%%%%%%%%%%%
% subsection
\subsection{Contribution of $q_kp_k$}
We prove that the symbol $q_kp_k (y, \eta)$ yields a
vanishing contribution to the limit of $L_k''(\chi(\tau) b, \psi)$.
%%%%%%%%%%%%%%%%%%%%%%%%
% proposition          %
%%%%%%%%%%%%%%%%%%%%%%%%
\begin{proposition}
  \label{prop: contrib qkpk}
  One has $L''_k (q_kp_k, \psi) = o(1)_{k\to+\infty}$. 
\end{proposition}
Proving this result requires some preliminary results.

\medskip
Set $\varphi_k = \varphi(\hk^3 D_z) $. Naturally,
$\varphi_k$ is  uniformly  bounded on $L^2(\R)$ as a uniformly  bounded Fourier
multiplier. One can view $\varphi_k$ in various manners: one has 
\begin{align*}
  \varphi_k = \Op^{\hk}(\hk^2 \zeta) = \Op^{\hk^3}(\zeta).
\end{align*}
With the second formula, by simply replacing $h$ by $\hk^3$ in the
analysis of Section~\ref{sec: Semi-classical operators and measures},
with point~(3) of Proposition~\ref{propcommut} one has the following
result.
%%%%%%%%%%%%%%%%%%%%%%%%
% lemma                %
%%%%%%%%%%%%%%%%%%%%%%%%
\begin{lemma}
  \label{lem: commutation varphi Lipschitz}
  Let $\theta\in W^{1,\infty}(\R)$. Then,
$\Norm{[\theta, \varphi_k]}{\mathcal{L}(L^2(\R))} \leq C \hk^3$.
\end{lemma}

\medskip
Set
\begin{align*}
  &Q_k =\OpH (q_k), \ \ \bQ_k=\OpH (\bar{q}_k), \ \ 
  P_k =\hk^2 P_{\kk, \gk}, \\
  &G_k = \OpH(q_k p_k), \ \ \et \ \ 
  \bG_k = \OpH(\bar{q}_k p_k).
\end{align*}
Note that 
\begin{align*}
  %\label{eq: Pk}
  P_k  = \OpH(p_k) +   \hk^2 P_k^1,
\end{align*}
where $P_k^1$ is a differential operator of order one with bounded
coefficients.

\medskip
Because of the form of $p_k$ and $P_k$, one writes
  $P_k = P_k^d + P_k^{\mathsf T}$ with
\begin{align*}
  %\label{eq: decomposition Pk}
  &P_k^d = \tkk^{-1} \hk D_z \tkk \gk^{dd}(x) \hk D_z,\\
  &P_k^{\mathsf T}= \hk^2 \d_t^2 + \sum_{{1\leq i, j \leq d} \atop {(i,j) \neq(d,d)}}
  \tkk^{-1} \hk D_i \tkk  \gk^{ij}(x) \hk D_j,\notag
\end{align*}
and $G_k = G_k^d + G_k^{\mathsf T}$ with  
\begin{align*}
  %\label{eq: decomposition Gk}
  &G_k^d = \gk^{dd}(x) \OpH(q_k) \hk^2 D_z^2 \\
  &G_k^{\mathsf T}= \OpH(q_k)\hk^2 \d_t^2 +\sum_{{1\leq i, j \leq d} \atop {(i,j) \neq(d,d)}} 
  \gk^{ij}(x) \OpH(q_k) \hk^2 D_i D_j,\notag
\end{align*}
and 
$\bG_k = \bG_k^d + \bG_k^{\mathsf T}$ with  
\begin{align*}
  %\label{eq: decomposition bar Gk}
  &\bG_k^d = \gk^{dd}(x) \OpH(\bq_k) \hk^2 D_z^2 \\
  &\bG_k^{\mathsf T}= \OpH(\bq_k)\hk^2 \d_t^2 +\sum_{{1\leq i, j \leq d} \atop {(i,j) \neq(d,d)}} 
  \gk^{ij}(x) \OpH(\bq_k) \hk^2 D_i D_j.\notag
\end{align*}
With this notation one has 
\begin{align*}
   L_k''(q_kp_k, \psi)  
  =  L_k^d(q_kp_k, \psi)   + L_k^{\mathsf T}(q_kp_k, \psi),
\end{align*}
with
\begin{align*}
  L_k^d(q_kp_k, \psi) &= i \biginp{ G_k^d \varphi_k  \psi u_k}
    {v_k \otimes \delta_{z=0}}
    _{H_z^{1}L_{y'}^2,H_z^{-1}L_{y'}^2}^{\kk\mugk dt}\\
  &\quad - i \biginp{ v_k \otimes \delta_{z=0}}
    {\bG_k^d \varphi_k\psi u_k}
    _{H_z^{-1}L_{y'}^2,H_z^{1}L_{y'}^2}^{\kk\mugk dt}.
\end{align*}
and 
\begin{align*}
   L_k^{\mathsf T}(q_kp_k, \psi)&= i \biginp{ G_k^{\mathsf T} \varphi_k  \psi u_k}
    {v_k \otimes \delta_{z=0}}
    _{H_z^{1}L_{y'}^2,H_z^{-1}L_{y'}^2}^{\kk\mugk dt}\\
  &- i \biginp{ v_k \otimes \delta_{z=0}}
    {\bG_k^{\mathsf T} \varphi_k\psi u_k}
    _{H_z^{-1}L_{y'}^2,H_z^{1}L_{y'}^2}^{\kk\mugk dt}.
\end{align*}

%%%%%%%%%%%%%%%%%%%%%%%%
% lemma                %
%%%%%%%%%%%%%%%%%%%%%%%%
\begin{lemma}
  \label{lem: qkpk 2}
  \begin{enumerate}
    \item One has $\hk Q_k \in \mathcal  R_0$.
    \item Let $(f_k)_k$ be a bounded sequence in
      $W^{1,\infty}(\R^{d+1)})$. One has $[Q_k, f_k] \in \mathcal
      R_0$.
    
    \item
  The operator $G_k^{\mathsf T} - Q_k P_k^{\mathsf T}$ is a finite sum of operators that lie in 
  \begin{equation}
    \label{eq: qkpk 1}
    \sum_{|\alpha'|=1} \mathcal{R}_0\,  \hk^2 D_{y'}^{\alpha'} D_z
    +\sum_{|\alpha|=2} \mathcal{R}_0\,  \hk^2 D_{y'}^{\alpha}
    + \sum_{|\beta|=1}\mathcal{R}_0\,  \hk  D_y^{\beta},
\end{equation}
with $\mathcal R_0$ as given in
Definition~\ref{def: error operators}. The same holds for $\bG_k^{\mathsf T}  - \bQ_k P_k^{\mathsf T}$
\item The operators $Q_k [P_k^{\mathsf T}, \varphi_k]$ and $\bQ_k [P_k^{\mathsf T},
  \varphi_k]$ are also finite sums of
  operators that lie in the space given by \eqref{eq: qkpk 1}.
\end{enumerate}
\end{lemma}
A proof is given below. 

Observe that
$\mathcal{R} \, \varphi_k \subset \mathcal{R}$ as
$\varphi_k$ is uniformly bounded on $L^2(\R^{d+1})$. Then, 
with Lemma~\ref{lem610}, exploiting \eqref{eq: mass2-proof}--\eqref{eq: mass2-proof4} and the
local estimate \eqref{eq: local H2 tangential estimate}, with the
third item in Lemma~\ref{lem: qkpk 2} 
one obtains 
\begin{align*}
   L_k^{\mathsf T}(q_kp_k, \psi)  
  &= i \biginp{ Q_k P^{\mathsf T}_k \varphi_k  \psi u_k}
    {v_k \otimes \delta_{z=0}}
    _{H_z^{1}L_{y'}^2,H_z^{-1}L_{y'}^2}^{\kk\mugk dt}
  \\
  &\quad
    - i \biginp{ v_k \otimes \delta_{z=0}}
    {\bQ_k P^{\mathsf T}_k\varphi_k\psi u_k}
    _{H_z^{-1}L_{y'}^2,H_z^{1}L_{y'}^2}^{\kk\mugk dt}
    + o(1)_{k\to + \infty}.\notag
\end{align*}
With the fourth item in Lemma~\ref{lem: qkpk 2}, with the same
argumentation one finds
\begin{align*}
    L_k^{\mathsf T}(q_kp_k, \psi)  
  &=  i \biginp{ Q_k \varphi_k P^{\mathsf T}_k \psi u_k}
    { v_k \otimes \delta_{z=0}}_{H_z^{1}L_{y'}^2,H_z^{-1}L_{y'}^2}^{\kk\mugk dt}   
  \\
  &\quad
    -i \biginp{ v_k \otimes \delta_{z=0}}
    {\bQ_k \varphi_k P^{\mathsf T}_k \psi u_k}
    _{ H_z^{-1}L_{y'}^2,H_z^{1}L_{y'}^2}^{\kk\mugk dt}
    + o(1)_{k\to + \infty}.\notag
\end{align*}
 %%%%%%%%%%%%%%%%%%%%%%%%
% lemma                %
%%%%%%%%%%%%%%%%%%%%%%%%
\begin{lemma}
  \label{lemma: commutator Pk psi}
  Let $M_k(y,D_y)$ be a differential operator with
 coefficients that are uniformly bounded with respect to $k$.
 Let $N \in \N$. For some $C_N>0$ one has
 \begin{align*}
   \Norm{Q_k \varphi_k  [M(y,D_y), \psi]
     u_k}{H_z^{1}L_{y'}^2}\leq C \hk^N.
 \end{align*} The same hold for $\bQ_k$ in
  place of $Q_k$.
\end{lemma}
A proof is given below.
Applying Lemma~\ref{lemma: commutator Pk psi}
gives 
\begin{align*}
   L_k^{\mathsf T}(q_kp_k, \psi)  
  &=  i \biginp{ Q_k \varphi_k\psi P^{\mathsf T}_k u_k}
    { v_k \otimes \delta_{z=0}}_{H_z^{1}L_{y'}^2,H_z^{-1}L_{y'}^2}^{\kk\mugk dt}   
  \\
  &\quad
    -i \biginp{ v_k \otimes \delta_{z=0}}
    {\bQ_k \varphi_k\psi P^{\mathsf T}_k u_k}
    _{ H_z^{-1}L_{y'}^2,H_z^{1}L_{y'}^2}^{\kk\mugk dt}
    + o(1)_{k\to + \infty}.\notag
\end{align*}

\bigskip
Our goal is now to handle the terms associated with the operators
$G_k^d$ and $\bG_k^d$ in $L_k^d(q_kp_k, \psi)$.
With the forms of $G_k^d$ and $\bG_k^d$ and Lemma~\ref{lemma: commutator Pk psi},
using that $\gk^{dd}$ is uniformly Lipschitz as $k\to +\infty$, one has 
\begin{align*}
L_k^d (q_kp_k, \psi)   &= i \biginp{ \gk^{dd}  Q_k \varphi_k \psi \hk^2 D_z^2  u_k}
 { v_k \otimes
   \delta_{z=0}}_{H_z^{1}L_{y'}^2,H_z^{-1}L_{y'}^2}^{\kk\mugk dt}   \\
 &\quad   - i \biginp{ v_k \otimes \delta_{z=0}}
    {\gk^{dd}  \bQ_k \varphi_k \psi \hk^2 D_z^2  u_k}
    _{ H_z^{-1}L_{y'}^2,H_z^{1}L_{y'}^2}^{\kk\mugk dt}
   + o(1)_{k\to +\infty}.
\end{align*}
With the "jump formula" one has 
\begin{align*}
  \hk^2 D_z^2  u_k 
  &=  - \hk^2 \d_z  {u_k}_{z=0^+}\otimes \delta_{z=0} + \hk^2 \udl{D_z^2 u_k}\\
   &=- \hk  (\gk^{dd})^{-1} v_k \otimes \delta_{z=0} 
     + \hk^2  \udl{D_z^2 u_k},
  \end{align*}
recalling \eqref{eq: vk local coordinates}, where $\udl{f}$ denotes
the zero-extension of $f_{| z>0}$. 
One writes 
\begin{align*}
 L_k^d (q_kp_k, \psi)  
  =L_k^\delta (q_kp_k, \psi)  + L_k^{\{ z>0\}}(q_kp_k, \psi) + o(1)_{k\to +\infty},
\end{align*}
with 
\begin{align}
\label{eq: L k delta}
 L_k^\delta (q_kp_k, \psi)  
  &=  - i \hk   \biginp{ \gk^{dd}  Q_k \varphi_k \psi  (\gk^{dd})^{-1} v_k \otimes \delta_{z=0}}
 { v_k \otimes \delta_{z=0}}_{H_z^{1}L_{y'}^2,H_z^{-1}L_{y'}^2}^{\kk\mugk dt}   \\
 &\quad   + i \hk \biginp{ v_k \otimes \delta_{z=0}}
    {\gk^{dd}  \bQ_k \varphi_k \psi (\gk^{dd})^{-1} v_k \otimes \delta_{z=0}}
    _{ H_z^{-1}L_{y'}^2,H_z^{1}L_{y'}^2}^{\kk\mugk dt},\notag
\end{align}
and 
\begin{align*}
 L_k^{\{ z>0\}} (q_kp_k, \psi)  
  &=  i \biginp{ \gk^{dd}  Q_k \varphi_k \psi \hk^2  \udl{D_z^2 u_k}}
    { v_k \otimes\delta_{z=0}}_{H_z^{1}L_{y'}^2,H_z^{-1}L_{y'}^2}^{\kk\mugk dt}   \\
 &\quad   - i \biginp{ v_k \otimes \delta_{z=0}}
    {\gk^{dd}  \bQ_k \varphi_k \psi \hk^2  \udl{D_z^2 u_k}}
    _{ H_z^{-1}L_{y'}^2,H_z^{1}L_{y'}^2}^{\kk\mugk dt}.
\end{align*}
With Lemmata~\ref{lem: commutation varphi Lipschitz} and \ref{lem: qkpk 2}  one has
\begin{align*}
  \gk^{dd} Q_k \varphi_k \psi = Q_k \varphi_k \psi \gk^{dd} \mod
  \mathcal R_0 
  \ \ \et \ \ 
  \gk^{dd} \bQ_k \varphi_k \psi = \bQ_k \varphi_k \psi \gk^{dd} \mod \mathcal R_0.
\end{align*}
With  Lemma~\ref{lem610}, as $\hk^2
\udl{D_z^2 u_k}$ is bounded in $L^2$ by \eqref{eq: local H2 tangential
  estimate},  one obtains
\begin{align*}
 L_k^{\{ z>0\}} (q_kp_k, \psi)  
  &=  i \biginp{Q_k \varphi_k \psi\,  \udl{\gk^{dd} \hk^2   D_z^2 u_k}}
    { v_k \otimes\delta_{z=0}}_{H_z^{1}L_{y'}^2,H_z^{-1}L_{y'}^2}^{\kk\mugk dt}   \\
 &\quad   - i \biginp{ v_k \otimes \delta_{z=0}}
    {\bQ_k \varphi_k \psi \, \udl{\gk^{dd}  \hk^2  D_z^2 u_k}}
    _{ H_z^{-1}L_{y'}^2,H_z^{1}L_{y'}^2}^{\kk\mugk dt}.
\end{align*}
 One writes    
\begin{align*}
  \gk^{dd}  \hk^2  D_z^2 = \kk^{-1} \hk  D_z \gk^{dd} \kk \hk  D_z
  + \hk \kk^{-1} \big( D_z (\gk^{dd} \kk)\big)  \hk  D_z. 
\end{align*}
Since $\hk Q_k \in \mathcal R_0$ by Lemma~\ref{lem: qkpk 2} one obtains
\begin{align*}
 L_k^{\{ z>0\}} (q_kp_k, \psi)  
  &=  i \biginp{Q_k \varphi_k \psi\,  \udl{\kk^{-1} \hk  D_z \gk^{dd} \kk \hk  D_z u_k}}
    { v_k \otimes\delta_{z=0}}_{H_z^{1}L_{y'}^2,H_z^{-1}L_{y'}^2}^{\kk\mugk dt}   \\
 &\quad   - i \biginp{ v_k \otimes \delta_{z=0}}
    {\bQ_k \varphi_k \psi \, \udl{\kk^{-1} \hk  D_z \gk^{dd} \kk \hk  D_z u_k}}
    _{ H_z^{-1}L_{y'}^2,H_z^{1}L_{y'}^2}^{\kk\mugk dt}\\
    &\quad + o(1)_{k\to + \infty}.
\end{align*}
One thus obtains 
\begin{align*}
 L_k^{\{ z>0\}} (q_kp_k, \psi)   + L_k^{\mathsf T}(q_kp_k, \psi)  =o(1)_{k\to + \infty}.
\end{align*}
since $P_k = \kk^{-1} \hk  D_z \gk^{dd} \kk \hk  D_z +  P^{\mathsf
  T}_k$ and thus 
\begin{align*}
   L_k''(q_kp_k, \psi)  
  = L_k^\delta (q_kp_k, \psi)  
 +o(1)_{k\to + \infty},
\end{align*}
with $L_k^\delta (q_kp_k, \psi)$ given in \eqref{eq: L k delta}. One
then writes 
\begin{align*}
 L_k^\delta (q_kp_k, \psi)  
  &=  i \hk   
    \bigdup{ N_k  v_k \otimes \delta_{z=0}}
    { v_k \otimes \delta_{z=0}}_{H_z^{\alpha}L_{y'}^2, H_z^{-\alpha}L_{y'}^2}^{\kk\mugk dt},
\end{align*}
for any $\alpha> 1/2$, with 
\begin{align*}
  N_k = (\gk^{dd})^{-1}  \psi  \varphi_k^\star  \bQ_k^\star\gk^{dd} 
    - \gk^{dd}  Q_k \varphi_k \psi  (\gk^{dd})^{-1}, 
\end{align*}
where the adjoints with  the $\star$-notation are understood  in the sense of the 
inner product $L^2(\R^{d+1}, \kk\mugk dt)$, that is, $\varphi_k^\star
\bQ_k^\star   =  \tkk^{-1} \varphi_k \bQ_k^*\tkk$ where the
adjoint with the $*$-notation is understood for the inner
product $L^2(\R^{d+1}, dx dt)$. One thus has 
\begin{align*}
  N_k &= (\gk^{dd})^{-1}  \psi  \tkk^{-1} \varphi_k \bQ_k^*\tkk \gk^{dd} 
    - \gk^{dd}  Q_k \varphi_k \psi  (\gk^{dd})^{-1}\\
  &=\gk^{dd}  \Big(  \psi \big( \tkk (\gk^{dd})^2\big)^{-1}  \varphi_k
    \bQ_k^* \, \tkk (\gk^{dd})^2
    - Q_k \varphi_k \psi \Big)  (\gk^{dd})^{-1} , 
\end{align*}
%%%%%%%%%%%%%%%%%%%%%%%%
% lemma                %
%%%%%%%%%%%%%%%%%%%%%%%%
\begin{lemma}
  \label{lemma: commutator final}
  Let $(f_k)_k$ be a sequence of functions such that 
  \begin{align*}
    \Norm{f_k}{W^{1,\infty}} + \Norm{1/f_k}{W^{1,\infty}} \leq C,
  \end{align*}uniformly with respect to $k$. 
  Let $\eps>0$. 
  For $\alpha >\frac 1 2$ chosen \suff close to $\frac 1 2$. Then,
  \begin{align*}
    \Norm{Q_k\varphi_k\psi - \psi f_k^{-1} \varphi_k \bQ_k^*f_k}
    {\L(H_z^{-\alpha}L_{y'}^2, H_z^{\alpha}L_{y'}^2)} 
    = o(\hk^{-\eps})_{k\to + \infty}. 
\end{align*}
\end{lemma}
A proof is given below. 

As $\tkk$ and $\gk^{dd}$ and their inverses are Lipschitz uniformly
with respect to $k$, 
with Lemma~\ref{lemma: commutator final} one finds that $ L_k^\delta (q_kp_k, \psi)   =
o(\hk^{1-\eps})_{k\to + \infty}$ for any $0 < \eps< 1$,  which concludes the proof of Proposition~\ref{prop: contrib qkpk}.
\hfill \qedsymbol \endproof

\bigskip
%%%% proof of lemma
\begin{proof}[\bfseries Proof of Lemma~\ref{lem: qkpk 2}]

  Let $f_k$ be as in the statement.  Both operators $[Q_k, f_k]$ and
  $\hk Q_k$ are bounded in $\L\big(L^2(\R^{d+1})\big)$ by $C \hk$ by
  Lemma~\ref{lemma: bounded 2} for the first one and Lemma~\ref{lemma:
    bounded} for the second one, recalling the properties of $q_k$
  given in \eqref{eq: estimation rk}--\eqref{eq: estimation rjk}.

  To estimate their operator norm from $L^2(\R^{d+1})$ to $H^1(\R_z;L^2(\R^d_{y'}))$ we compose with $D_z = \hk^{-1} \OpH(\zeta)$. On the one
hand one gets 
\begin{align*}
D_z [Q_k, f_k]
= \OpH(D_z q_k ) f_k
-  \OpH\big( D_z (f_k  q_k)\big) 
+\hk^{-1} [ \OpH (\zeta q_k) , f_k].
\end{align*}
The first two  operators are bounded in  $\L\big(L^2(\R^{d+1})\big)$
uniformly with respect to $k$
by Lemma~\ref{lemma: bounded}. For the third operator, using that $\zeta q_k$ is of
the form given in Lemma~\ref{lemma: bounded 2} by the polyhomogeneous
expansion of $q_k$ given in \eqref{eq: estimation rjk}, one also finds
a  $k$-uniform bound in $\L\big(L^2(\R^{d+1})\big)$. On the other hand, one
has 
\begin{align*}
D_z \hk Q_k= \hk \OpH (D_zq_k) + \OpH (\zeta q_k),
\end{align*}
that also has a $k$-uniform bound in   $\L\big(L^2(\R^{d+1})\big)$ by
Lemma~\ref{lemma: bounded}. The first two points of the lemma are proven.

For the third point, we provide the proof for $G_k'$ and $Q_k$. The proof of $\bG_k'$ and $\bQ_k$ is identical.
One sees that it
suffices to prove that $\OpH(q_k a_k') - Q_k \A_k'$ is a finite sum of
operators that lie in the space given by \eqref{eq: qkpk 1}, with
\begin{align*}
  a_k' = \sum_{{1\leq i, j \leq d} \atop {(i,j) \neq(d,d)}}
  \gk^{ij}\xi_i\xi_j
  \ \ \et \ \  
  \A_k'  = \sum_{{1\leq i, j \leq d} \atop {(i,j) \neq(d,d)}} \tkk^{-1} \hk D_i \tkk  \gk^{ij}(x) \hk D_j.
\end{align*}
One writes 
\begin{align*}
 \A_k' 
  = \sum_{{1\leq i, j \leq d} \atop {(i,j) \neq(d,d)}}  \big( \gk^{ij} \hk^2 D_i D_j 
  + \hk  \tkk^{-1} \big(D_i (\gk^{ij} \tkk) \big)   \hk D_j \big),
\end{align*}
and 
\begin{align*}
  Q_k \gk^{ij} \hk^2 D_i D_j
  &= \gk^{ij} Q_k  \hk^2 D_i D_j  + [Q_k, \gk^{ij}] \hk^2 D_i D_j \\
  &= \OpH(q_k \gk^{ij}\xi_i \xi_j) +  [Q_k, \gk^{ij}] \hk^2 D_i D_j, 
\end{align*}
yielding 
\begin{align*}
  Q_k \A_k' = \OpH(q_k a_k') + \sum_{{1\leq i, j \leq d} \atop {(i,j) \neq(d,d)}}  [Q_k, \gk^{ij}] \hk^2 D_i D_j
  + \hk  Q_k \tkk^{-1} \big(D_i (\gk^{ij} \tkk) \big) \hk D_j. 
\end{align*}
The result thus amounts to having $[Q_k, \gk^{ij}] \in \mathcal{R}_0$
and $\hk Q_k \in \mathcal{R}_0$, which holds by the first two points
of the lemma proven above. This concludes the proof of the third point
of Lemma~\ref{lem: qkpk 2}.

\medskip We now turn to the proof of the fourth point. Since
$[\d_t^2, \varphi_k]=0$ it suffices to consider $Q_k [\A_k', \varphi_k]$. One writes
\begin{align*}
  \A'_k=\sum_{{1\leq i, j \leq d} \atop {(i,j) \neq(d,d)}}\big( \gk^{ij} \hk^2 D_i  D_j 
  + \hk \tkk^{-1} \big(D_i (\gk^{ij} \tkk) \big)  \hk D_j\big),
\end{align*}
yielding
\begin{align*}
  [\A_k', \varphi_k]
  &=\sum_{{1\leq i, j \leq d} \atop {(i,j) \neq(d,d)}}
  \Big( [\gk^{ij}, \varphi_k] \hk^2 D_i D_j
  + \hk \tkk^{-1} \big( D_i (\gk^{ij} \tkk)\big) \varphi_k \hk D_j\\
  &\qquad \qquad \quad 
  -\hk \varphi_k\tkk^{-1} \big(D_i (\gk^{ij} \tkk)\big) \hk D_j
  \big)
  \Big).
\end{align*}
Since $Q_k$ is bounded on $L^2(\R^{d+1})$ and also bounded by $C\hk^{-1}$
from $L^2(\R^{d+1})$ to $H^1(\R_z;L^2(\R^d_y))$ uniformly in $\hk$
it suffices to prove that $[\gk^{ij}, \varphi_k]$
is bounded by
$C\hk$ on $L^2(\R^{d+1})$. This is a consequence of Lemma~\ref{lem: commutation varphi Lipschitz}.
\end{proof}

%%%% proof of lemma
\begin{proof}[\bfseries Proof of Lemma~\ref{lemma: commutator Pk psi}]
  First, we prove 
  \begin{align}
    \label{eq: commutator Pk psi L2}
    \Norm{Q_k \varphi_k  [M_k(y,D_y), \psi] u_k}{L^2(\R^{d+1}))}\lesssim \hk^N.
\end{align}
Second, we prove 
\begin{align}
    \label{eq: commutator Dz Pk psi L2}
    \Norm{D_z Q_k \varphi_k  [M_k(y,D_y), \psi] u_k}{L^2(\R^{d+1}))}\lesssim \hk^N.
\end{align}
Together \eqref{eq: commutator Pk psi L2} and \eqref{eq: commutator Dz
  Pk psi L2} give the sought result. 

\medskip
\paragraph{\em Proof of \eqref{eq: commutator Pk psi L2}.}
  Note that $ [M_k(y,D_y), \psi] u_k$ is bounded in $L^2(\R^{d+1})$ by
  $C \hk^{1-m}$, where $m$ is the order of  $M_k(y,D_y)$, 
  by \eqref{eq: mass2-proof4}.

With the polyhomogeneous development in the $\zeta$ variable of
  $q_k$ given in \eqref{eq: estimation rjk}--\eqref{polyhomog} one writes
  \begin{align*}
    q_k  (y,\eta',\zeta) = \frac{1- \phi(\zeta)}{\zeta}  q_k^1 (y,\eta') +q_k^a (y,\eta',\zeta) ,
  \end{align*}
  with 
  \begin{multline}
    \label{polyhomog:  remainder}
  \Big|\d_{y}^\alpha \d_{\eta'}^\beta  \d_\zeta^\delta  q_k^a (y,\eta',\zeta) \Big|
  \leq C_{N, \beta, \delta}  \jp{\eta'}^{-N} \jp{\zeta}^{-2 - \delta},\\
  \pour\ N \in \N, \ \alpha \in \N^{d+1},  \ |\alpha| \leq 1, \
  \beta \in \N^{d},  \ \delta \in \N,  y \in \R^{d+1}, \ (\eta',\zeta) \in \R^{d+1}.
\end{multline}
One writes $Q_k = \OpH(q_k^1) \OpH\big( \frac{1-
  \phi(\zeta)}{\zeta}\big) + \OpH(q_k^a)$.

Recall that we work in the local chart
  $(\O,\cdiffL)$ at the boundary. Since $\psi = 1$ in a \nhd of the $y$-projection of
  $\supp b$, note that $\supp([M_k(y,D_y), \psi]u_k)$ does not meet the $y$-projection of
  $\supp q_k$ since $\supp q_k \subset \supp b$.
Let $\tpsi, \hpsi \in \Cinfc \big(\cdiffL(\O)\big)$ with $\hpsi$ equal to 1 in a \nhd of
  $\supp \tpsi$ and with $\tpsi$ equal to 1 in a \nhd of the  $y$-projection of
  $\supp b$ and moreover $\psi =1$ in a \nhd of $\supp \hpsi$.
  One has
  $[M_k(y,D_y), \psi] = (1-\hpsi) [M_k(y,D_y), \psi]$ and 
\begin{multline}
  \label{eq: est standard calculus support}
    \Norm{\tpsi \varphi_k (1-\hpsi) }{\L (L^2(\R^{d+1}))} 
    + \Norm{\tpsi 
    \OpH\big( \zeta^{-1}(1- \phi(\zeta)) \big) \varphi_k (1-\hpsi)}
  {\L(L^2(\R^{d+1}))}  \\
  \leq C_N
    \hk^N,\quad N \in \N,
  \end{multline}
  by standard calculus. Since $\OpH(q_k^1)$ and  $\OpH(q_k^a)$ are
  bounded on $L^2(\R^{d+1})$, to obtain \eqref{eq: commutator Pk psi L2} it
  thus suffices to study the $L^2$-boundedness of the operators $\OpH(q_k^1)
  (1-\tpsi)$ and $\OpH(q_k^a)(1-\tpsi)$.

  The tangential kernel of $\OpH(q_k^1) (1-\tpsi)$ is given by 
  \begin{align*}
    K(y', \ty') = (2\pi)^{-d} \int e^{i (y' - \ty')\cdot \eta'} 
    (1-\tpsi(\ty',z)) q_k^1(y',z, \hk \eta') \, d \eta'.
  \end{align*}
  With the joint support properties of $q_k^1$ and $\tpsi$ one finds
  that $\Norm{y' - \ty'}{} \geq C>0$ in the support of the integrand.  Since $L
  \exp(i (y' - \ty')\cdot \eta') = \exp(i (y' - \ty')\cdot \eta')$ with 
  $L = - i \Norm{y' - \ty'}{}^{-2} (y' - \ty')\cdot \d_{\eta'} $ one
  can write 
  \begin{align*}
    K(y', \ty') = (2\pi)^{-d} \int e^{i (y' - \ty')\cdot \eta'} 
    (1-\tpsi(\ty',z))  (\transp L)^N q_k^1(y',z, \hk \eta') \, d \eta'.
  \end{align*}
  With the estimation~\eqref{eq: estimation rjk} for $q_k^1$ one finds 
  \begin{align*}
    |K(y', \ty')| \lesssim  \frac{\hk^N}{\jp{y' - \ty'}^N},
  \end{align*}
  which, by the Schur lemma (Lemma~\ref{lemma: Schur lemma}), gives
  \begin{align}
    \label{eq: eq: commutator Pk psi L2-1}
    \Norm{\OpH(q_k^1) (1-\tpsi)}{\L(L^2(\R^{d+1}))}\lesssim \hk^N.
  \end{align}

  The kernel of $\OpH(q_k^a)(1-\tpsi)$  is given by 
  \begin{align*}
    K(y, \ty) = (2\pi)^{-d-1} \int e^{i (y - \ty)\cdot \eta} 
    (1-\tpsi(\ty)) q_k^a(y, \hk \eta) \, d \eta.
  \end{align*}
  Here, $\Norm{y - \ty}{} \geq C>0$ in the support of the integrand,
  yielding 
  \begin{align*}
    K(y, \ty) = (2\pi)^{-d-1} \int e^{i (y - \ty)\cdot \eta} 
    (1-\tpsi(\ty))  (\transp L)^N  q_k^a(y, \hk \eta) \, d \eta,
  \end{align*}
  with   $L = - i \Norm{y - \ty}{}^{-2} (y - \ty)\cdot \d_{\eta}$,
  implying with \eqref{polyhomog:  remainder}
  \begin{align*}
    |K(y, \ty)| \lesssim  \frac{\hk^N}{\jp{y - \ty}^N},
  \end{align*}
  which, by the Schur lemma, gives
  \begin{align}
    \label{eq: eq: commutator Pk psi L2-2}
    \Norm{\OpH(q_k^a) (1-\tpsi)}{\L(L^2(\R^{d+1}))}\lesssim \hk^N.
  \end{align}
  Together \eqref{eq: eq: commutator Pk psi L2-1} and \eqref{eq: eq:
    commutator Pk psi L2-2} give estimate \eqref{eq: commutator Pk psi
    L2}.

\medskip
\paragraph{\em Proof of \eqref{eq: commutator Dz Pk psi L2}.}
Here, we write 
\begin{align*}
    q_k  (y,\eta',\zeta) 
  = \frac{1- \phi(\zeta)}{\zeta}  q_k^1(y,\eta') 
  + \frac{1- \phi(\zeta)}{\zeta^2}  q_k^2(y,\eta') 
  + q_k^b (y,\eta',\zeta),
  \end{align*}
  with 
  \begin{multline*}
  \Big|\d_{y}^\alpha \d_{\eta'}^\beta  \d_\zeta^\delta  q_k^b (y,\eta',\zeta) \Big|
  \leq C_{N, \beta, \delta}  \jp{\eta'}^{-N} \jp{\zeta}^{-3 - \delta},\\
  \pour\ N \in \N, \ \alpha \in \N^{d+1},  \ |\alpha| \leq 1, \
  \beta \in \N^{d},  \ \delta \in \N,  y \in \R^{d+1}, \ (\eta',\zeta) \in \R^{d+1}.
\end{multline*}
One writes 
\begin{align*}
  Q_k = \OpH(q_k^1) \OpH\big( \frac{1- \phi(\zeta)}{\zeta}\big) 
  + \OpH(q_k^2) \OpH\big( \frac{1- \phi(\zeta)}{\zeta^2}\big) 
  + \OpH(q_k^b).
\end{align*}  
One has 
\begin{align*}
  D_z Q_k &= \OpH(D_z q_k^1) \OpH\big( \frac{1-
  \phi(\zeta)}{\zeta}\big) 
  + \OpH(D_z q_k^2) \OpH\big( \frac{1- \phi(\zeta)}{\zeta^2}\big) 
  + \OpH(D_z q_k^b)\\
  &\quad +\hk^{-1} \OpH(q_k^1) \OpH\big(1- \phi(\zeta)\big) 
  + \hk^{-1} \OpH(q_k^2) \OpH\big( \frac{1- \phi(\zeta)}{\zeta}\big) \\
  &\quad +\hk^{-1}  \OpH(\zeta q_k^b).
\end{align*}  
Similarly to \eqref{eq: est standard calculus support} one has 
\begin{align*}
  \Norm{\tpsi 
    \OpH\big( \zeta^{-2}(1- \phi(\zeta)) \big) \varphi_k (1-\hpsi)}
  {\L(L^2(\R^{d+1}))} 
  \leq C_N
    \hk^N,\quad N \in \N.
  \end{align*}
Observe that $\zeta q_k^b$ has the same symbol properties as
$q_k^a$. The symbol properties of $D_z q_k^1$, $D_z q_k^2$, $D_z
q_k^b$ also allow one to carry out the same kernel estimations as
above yielding \eqref{eq: commutator Dz Pk psi L2}.
\end{proof}

\medskip
\begin{proof}[\bfseries Proof of Lemma~\ref{lemma: commutator final}]
  We claim that
  \begin{align}
     \label{eq lemma: commutator final-claim1}
    \Norm{Q_k \varphi_k\psi - \psi f_k^{-1} \varphi_k \bQ_k^*f_k}
    {\L(H^{-1/2}_zL^2_{y'}, H_z^{1/2} L^2_{y'})}= O(1),
  \end{align}
   and
  \begin{align}
     \label{eq lemma: commutator final-claim2}
    \Norm{Q_k \varphi_k\psi - \psi f_k^{-1} \varphi_k
    \bQ_k^*f_k }{\L(H^{-1}_z(
    L^2_{y'})), (H_z^{1}( L^2_{y'}))}= O(\hk^{-4}).
  \end{align}
  Interpolation of the two estimations then gives the result of the lemma. 
  
\bigskip
We now prove the claimed estimates.
\medskip
\paragraph{\em Proof of  estimate \eqref{eq lemma: commutator final-claim1}.}
First, one has
\begin{align*}
  &Q_k\varphi_k\psi - \psi f_k^{-1} \varphi_k \bQ_k^*f_k\\
  &\quad = Q_k\varphi_k \psi f_k^{-1}f_k - \psi f_k^{-1} \varphi_k \bQ_k^*f_k
  \\
  &\quad = [Q_k \varphi_k, \psi f_k^{-1}] f_k 
  + f_k^{-1} \psi 
    \big( \OpH\big(q_k\varphi(\hk^2 \zeta) \big)
    - \OpH\big(\ovl{q}_k \varphi(\hk^2 \zeta) \big)^* \big)f_k.
\end{align*}
With Lemmata~\ref{lemma: bounded 2} and \ref{lemmabis: bounded 2} one has 
\begin{multline*}
  \Norm{[Q_k \varphi_k, \psi f_k^{-1}]}{\L(L^2(\R^{d+1}))}\\
  + \Norm{\OpH\big(q_k\varphi(\hk^2 \zeta) \big)
    - \OpH\big(\ovl{q}_k \varphi(\hk^2 \zeta) \big)^*}{\L(L^2(\R^{d+1}))}
  \lesssim \hk,
\end{multline*}
yielding 
 \begin{equation}
   \label{eq lemma: commutator final-claim1-a}
   \Norm{Q_k\varphi_k\psi - \psi f_k^{-1} \varphi_k \bQ_k^*f_k}
    {\L(L^2(\R^{d+1}))}\lesssim \hk.
 \end{equation}
Second one writes
\begin{align*}
  &D_z \big( Q_k\varphi_k\psi - \psi f_k^{-1} \varphi_k \bQ_k^*f_k\big)\\
  &\quad = \OpH(D_z q_k) \varphi_k\psi   
    - \big( D_z(\psi f_k^{-1} )\big) \varphi_k \bQ_k^*f_k\\
    &\qquad + \hk^{-1}
  \big( \OpH(\zeta q_k) \varphi_k\psi 
    -  \psi f_k^{-1} \varphi_k \OpH(\zeta \ovl{q}_k)^*f_k\big).
\end{align*}
With Lemma~\ref{lemma: bounded} one finds
\begin{align*}
  \Norm{\OpH(D_z q_k) \varphi_k\psi  }{\L(L^2(\R^{d+1}))}
  + \Norm{\big( D_z(\psi f_k^{-1} )\big) \varphi_k \OpH(\ovl{q}_k)^*f_k}{\L(L^2(\R^{d+1}))}
  \lesssim 1.
\end{align*}
Arguing as for \eqref{eq lemma: commutator final-claim1-a} one finds 
\begin{equation*}
  \Norm{
  \OpH(\zeta q_k) \varphi_k\psi 
    -  \psi f_k^{-1} \varphi_k \OpH(\zeta \ovl{q}_k)^*f_k
  }
    {\L(L^2(\R^{d+1}))}\lesssim \hk.
 \end{equation*}
This gives 
\begin{equation}
   \label{eq lemma: commutator final-claim1-b}
   \Norm{D_z \big( Q_k\varphi_k\psi - \psi f_k^{-1} \varphi_k \bQ_k^*f_k\big)}
    {\L(L^2(\R^{d+1}))}\lesssim 1
 \end{equation}
Together,  \eqref{eq lemma: commutator final-claim1-a} and \eqref{eq
  lemma: commutator final-claim1-b} give
\begin{equation}
   \label{eq lemma: commutator final-claim1-c}
   \Norm{Q_k\varphi_k\psi - \psi f_k^{-1} \varphi_k \bQ_k^*f_k}
    {\L(L^2(\R^{d+1}), H^1_zL^2_{y'})}\lesssim 1.
 \end{equation}
By duality this implies 
\begin{equation}
   \label{eq lemma: commutator final-claim1-d}
   \Norm{Q_k\varphi_k\psi - \psi f_k^{-1} \varphi_k \bQ_k^*f_k}
    {\L(H^{-1}_zL^2_{y'}, L^2(\R^{d+1}))}\lesssim 1.
 \end{equation}
An interpolation of \eqref{eq lemma: commutator final-claim1-c} and
\eqref{eq lemma: commutator final-claim1-d} yields \eqref{eq lemma: commutator final-claim1}.

\medskip
\paragraph{\em Proof of  estimate \eqref{eq lemma: commutator final-claim2}.}
Above we computed $D_z Q_k = \OpH(D_z q_h)+   \hk^{-1} \OpH(q_h
\zeta)$
yielding
\begin{align*}
 %\label{eq: lemma: commutator final-claim2-a}
\Norm{Q_k}{\L(L^2(\R^{d+1}), H_z^{1} L^2_{y'})}
  \lesssim \hk^{-1}, 
\end{align*}
and thus
\begin{align}
  \label{eq: lemma: commutator final-claim2-b}
  \Norm{Q_k \varphi_k\psi}{\L(L^2(\R^{d+1}), H_z^{1}L^2_{y'})}
  \lesssim \hk^{-1}.
\end{align}
One also computes 
\begin{align*}
  D_z Q_k \varphi_k\psi D_z
  = D_z Q_k \varphi_k D_z\psi -  D_z Q_k \varphi_k (D_z\psi).
\end{align*}
Since $\varphi_k D_z = \hk^{-3} \OpH\big( \varphi(\hk^2 \zeta)\hk^2
\zeta\big)$ one finds that   $\varphi_k D_z $ is bounded on
$L^2(\R^{d+1})$ by $C \hk^{-3}$ thus yiedling
\begin{align}
  \label{eq: lemma: commutator final-claim2-c}
  \Norm{Q_k \varphi_k\psi D_z}{\L(L^2(\R^{d+1}), H_z^{1}L^2_{y'})}
  \lesssim \hk^{-4}.
\end{align}
Together \eqref{eq: lemma: commutator final-claim2-b} and \eqref{eq:
  lemma: commutator final-claim2-c} give
\begin{align}
  \label{eq: lemma: commutator final-claim2-d}
  \Norm{Q_k \varphi_k\psi}{\L(H_z^{-1} L^2_{y'}, H_z^{1} L^2_{y'})}
  \lesssim \hk^{-4}.
\end{align}
The same holds for $\bQ_k$ in place of $Q_k$ and 
by duality one obtains  
\begin{align}
  \label{eq: lemma: commutator final-claim2-e}
  \Norm{\psi f_k^{-1} \varphi_k \bQ_k^*f_k}{\L(H_z^{-1} L^2_{y'}, H_z^{1} L^2_{y'})}
  \lesssim \hk^{-4},
\end{align}
and together \eqref{eq: lemma: commutator final-claim2-d} and
\eqref{eq: lemma: commutator final-claim2-e} yield \eqref{eq
  lemma: commutator final-claim2}.
\end{proof}

%%%%%%%%%%%%%%%%%%%%%%%%
% remark               %
%%%%%%%%%%%%%%%%%%%%%%%%
\begin{remark}
Note that the proof we give of \eqref{eq lemma: commutator
  final-claim2} is far from optimal. However, this has no consequence
on the final result of Lemma~\ref{lemma: commutator final}. 
\end{remark}
%In fact,
% additional work could  improve the estimate to
% \begin{align*}
%      \Norm{Q_k \varphi_k\psi 
%   - \psi f_k^{-1} \varphi_k\bQ_k^*f_k }
%   {\L(H^{-1}_z  L^2_{y'}), (H_z^{1} L^2_{y'})}= O(\hk^{-1} ).
%   \end{align*}
%   Such refinement is however not needed for our purpose.
% 

%%%%%%%%%%%%%%%%%
% subsection
\subsection{Contributions of $b_{0,k}$ and $b_{1,k}$}

First, we prove that the symbol $b_{0,k} (y, \eta')$ yields a vanishing
contribution to the limit of $L_k''(\chi(\tau) b, \psi)$.  Second, we
prove that the symbol $b_{1,k} (y, \eta')$ yields a contribution to
the limit of $L_k''(\chi(\tau) b, \psi)$ as opposed to the other
symbols appearing in the Euclidean division of Proposition~\ref{prop:
  Euclidean symbol division}. This contribution implies the action of
the semi-classical measure $\nu$ at the boundary.

The tangential nature of $\OpH(b_{0,k})$ and  $\OpH(b_{1,k})$  allows one to consider traces
through the action of the Dirac measure $\delta_{z=0}$. A key point of
the proof of this section is the understanding of traces
after the action of the regularizing operator $\varphi_k$.

\medskip 
Consider $w\in L^2(\R^{d+1})$ such that $w^+=w_{|\{z>0\}} \in H^1\big(
\R^+_z; L^2(\R^{d}_{y'})\big)$ and $w^-=w_{|\{z<0\}} \in H^1\big(
\R^-_z; L^2(\R^{d}_{y'})\big)$. One the one hand,  
$w^+ \in \Con^0 \big( [0,+\infty[_z;  L^2(\R^{d}_{y'})\big)$ and
$w _{|z=0^+}  = w^+_{|z=0^+}=\lim_{z \to 0^+} w(z)$ makes sense in
$L^2(\R^{d}_{y'})$ classically. Similarly $w _{|z=0^-}  = w^-_{|z=0^-}=\lim_{z \to 0^-} w(z)$
makes sense. On the other hand, 
the trace of $(\varphi_k w)_{|z=0}$ can be approximated by the mean
of the two traces of $w$ as in the
following lemma. 
%%%%%%%%%%%%%%%%%%%%%%%%
% lemma                %
%%%%%%%%%%%%%%%%%%%%%%%%
\begin{lemma}
  \label{lemma: trace varphik}
  Let $\varphi \in \Cinfc(\R)$ be real valued and equal to $1$ near
  $0$. There exists $C>0$ such that 
  \begin{multline*}
    %\label{eq: trace varphik}
    \bigNorm{\varphi(h^3 D_z)\udl{w}_{|z=0} 
      -  \frac 1 2 \big( w_{|z=0^-}+ w_{|z=0^+}\big)}{L^2(\R^{d})}\\
    \leq C h^{\frac 3 2 } \big( \Norm{\d_z w}
    {L^2( \R^-_z ;L^2(\R^{d}_{y'}))} + \Norm{\d_z w}
    {L^2( \R^+_z ;L^2(\R^{d}_{y'}))} \big),
\end{multline*}
for $h>0$ and $w\in L^2(\R^{d+1})$ such that $w^\pm \in H^1\big( \R^\pm_z; L^2(\R^{d}_{y'})\big)$.
\end{lemma}
%%%% proof of lemma
\begin{proof}
  By linearity and symmetry it suffices to consider $w$ such
  that $w_{\{z<0\}} =0$  and prove\footnote{In what follows, we will actually use
    Lemma~\ref{lemma: trace varphik} in the case of a function
    vanishing in $\{ z<0\}$.} 
  \begin{equation}
    \label{eq: trace varphik proof}
    \bigNorm{\varphi(h^3 D_z)w_{|z=0} 
      -  \frac 1 2 w_{|z=0^+}}{L^2(\R^{d})} 
    \leq C h^{\frac 3 2 } \Norm{\d_z w}
    {L^2( \R^+_z ;L^2(\R^{d}_{y'}))}.
\end{equation}
  Denote by $\hvarphi$ the inverse Fourier transform of
  $\varphi$. The Parseval formula gives
  \begin{align}
    \label{eq: trace varphik-1}
    2\pi \varphi(h^{3} D_z) w_{|z=0} 
    &= \int_\R \varphi(h^3 \zeta)  \hat{w} (\zeta)\,
    d\zeta\\
    &= \int_\R \hvarphi(z) w(h^3 z)\,  dz
    =  \int_{\R^+} \hvarphi(z) w (h^3 z)\,  dz.\notag
 \end{align}
 For $z\geq 0$, with the Cauchy-Schwarz inequality one finds
 \begin{align}  
   \label{eq: trace varphik-2}
   \Norm{w(z) - w_{|z=0^+}}{L^2(\R^{d})} 
   = \BigNorm{\int_0^z \d_z w (s) ds }{L^2(\R^{d})}  
   \leq z^{1/2} \Norm{\d_z w}{L^2( \R^+_z ; L^2(\R^{d}_{y'}))}.
 \end{align}
 Using that $\hvarphi$ is even since $\varphi$ is real valued one has 
 $\int_{\R^+}  \hvarphi = \pi$ since $\int_\R \hvarphi = 2 \pi \varphi(0)
 =2 \pi$. With \eqref{eq: trace varphik-1} and \eqref{eq: trace varphik-2} one thus obtains
 \begin{align*}  
   &2 \pi \bigNorm{\varphi(h^{3} D_z) w_{|z=0}  - \frac 1 2 w_{|z=0^+}}{L^2(\R^{d})} 
   = \BigNorm{\int_{\R^+} \hvarphi (z) \big( w(h^3 z)- w_{|z=0^+}\big)\, dz 
   }{L^2(\R^{d})} \\
  &\qquad \leq h^{3/2} \Norm{\d_z w}{L^2( \R^+_z ; L^2(\R^{d}_{y'}))} 
    \int_{\R^+} |\hvarphi (z)| z^{1/2}\,  dz 
\lesssim h^{3/2} \Norm{\d_z w}{L^2( \R^+_z ; L^2(\R^{d}_{y'}))},
 \end{align*}
which is the sought result \eqref{eq: trace varphik proof}. 
\end{proof}

\bigskip
%%%%%%%%%%%%%%%%%%%%%%%%
% proposition          %
%%%%%%%%%%%%%%%%%%%%%%%%
\begin{proposition}
  \label{prop: contrib b0}
  One has $L''_k (b_{0,k}, \psi) = o(1)_{k\to+\infty}$. 
\end{proposition}
%%%%
\begin{proof}
Using that $\Norm{\psi u_k}{H^1( \R_z^+; L^2(\R _{y'}^d))} = O(\hk^{-1})$ by
\eqref{eq: mass2-proof} one writes 
 \begin{align*}
   L''_k (b_{0,k}, \psi)
   &=  i  \biginp{\OpH \big({b_{0,k}}_{|z=0}\big)
     (\varphi_k\psi u_k)_{|z=0^+}}{v_k}
     _{L^2(\R^d), L^2(\R^d)}^{\kk\mugbk dt}\\
 &\quad -i  \biginp{v_k}{\OpH \big({b_{0,k}}_{|z=0}\big) 
     (\varphi_k\psi u_k)_{|z=0^+}}
     _{L^2(\R^d), L^2(\R^d)}^{\kk\mugbk dt} \\
   &=  \frac i 2  \biginp{\OpH \big({b_{0,k}}_{|z=0}\big)
   (\psi u_k)_{|z=0^+}}{v_k}
     _{L^2(\R^d), L^2(\R^d)}^{\kk\mugbk dt}
   \\
   &\quad - \frac i 2  \biginp{v_k}{\OpH \big({b_{0,k}}_{|z=0}\big) 
     (\psi u_k)_{|z=0^+}}
     _{L^2(\R^d), L^2(\R^d)}^{\kk\mugbk dt}
   + O(\hk^{1/2})\\
   &=O(\hk ^{1/2}),
\end{align*}
using the homogeneous Dirichlet boundary condition,
that is, ${u_k}_{|z=0^+}=0$. 
\end{proof}

\bigskip
%%%%%%%%%%%%%%%%%%%%%%%%
% proposition          %
%%%%%%%%%%%%%%%%%%%%%%%%
\begin{proposition}
  \label{prop: contrib b1}
  One has 
  \begin{align*}
    %\label{eq: contrib b1}
    L''_k (b_{1,k}\zeta , \psi)
    &=  \biginp{\OpH \big({b_{1,k}}_{|z=0}\big) 
      \big((g_k^{dd} )^{-1} \psi\big)_{|z=0}  v_k}{v_k}
     _{L^2(\R^d), L^2(\R^d)}^{\kk\mugbk dt}
     + o(1)_{k\to+\infty}
    \end{align*}
\end{proposition}
\begin{proof}
With \eqref{eq: local H2 tangential
  estimate} one has $\Norm{\hk D_z \psi  u_k}{H^1( \R^+_z;
  (L^2(\R^d_{y'}))}= O(\hk^{-1})$, which gives 
\begin{align*}
   L''_k (b_{1,k} \zeta, \psi)
   &= i  \biginp{\OpH \big({b_{1,k}}_{|z=0}\big)
   (\varphi_k \hk D_z \psi u_k)_{|z=0^+}}{v_k}
   _{L^2(\R^d), L^2(\R^d)}^{\kk\mugbk dt}
    \\
 &\quad -  i  \biginp{v_k}{\OpH \big( \text{$\bb_{1,k}$}_{|z=0}\big) 
     (\varphi_k \hk D_z \psi u_k)_{|z=0^+}}
     _{L^2(\R^d), L^2(\R^d)}^{\kk\mugbk dt}\\
  &=\frac 1 2  \biginp{\OpH \big({b_{1,k}}_{|z=0}\big)
   (\hk \d_z \psi u_k)_{|z=0^+}}{ v_k}
   _{L^2(\R^d), L^2(\R^d)}^{\kk\mugbk dt}\\
    &\quad +\frac 1 2  \biginp{v_k}{\OpH \big(\text{$\bb_{1,k}$}_{|z=0}\big) 
     (\hk \d_z \psi u_k)_{|z=0^+}}
     _{L^2(\R^d), L^2(\R^d)}^{\kk\mugbk dt}
   + O(\hk^{1/2}),
\end{align*}
by Lemma~\ref{lemma: trace varphik}. 
With \eqref{eq: vk local coordinates} one has $v_k = \hk  g_k^{dd} \d_z {u_k}_{|z=0^+}$ yielding
\begin{align*}
  \hk \d_z \psi {u_k}_{|z=0^+} 
  =\hk  (\d_z \psi)_{|z=0} {u_k}_{|z=0^+} + \hk  \psi_{|z=0} (\d_z u_k)_{|z=0^+}
  = (g_k^{dd} )^{-1} \psi_{|z=0}v_k.
\end{align*} 
One then obtains
\begin{align*}
   L''_k (b_{1,k} \zeta, \psi)
   &=\frac 1 2  \biginp{\OpH \big({b_{1,k}}_{|z=0}\big)
    \big((g_k^{dd} )^{-1} \psi\big)_{|z=0} v_k}{v_k}
   _{L^2(\R^d), L^2(\R^d)}^{\kk\mugbk dt}\\
 &\quad 
   + \frac 12  \biginp{v_k}{\OpH \big(\text{$\bb_{1,k}$}_{|z=0}\big) 
      \big( (g_k^{dd} )^{-1} \psi\big)_{|z=0} v_k}
     _{L^2(\R^d), L^2(\R^d)}^{\kk\mugbk dt}
   + O(\hk^{1/2}).
\end{align*}
One writes 
\begin{align*}
   &\biginp{v_k}{\OpH \big(\text{$\bb_{1,k}$}_{|z=0}\big) 
       \big( (g_k^{dd} )^{-1} \psi\big)_{|z=0} v_k}
     _{L^2(\R^d), L^2(\R^d)}^{\kk\mugbk dt}\\
  &= \biginp{  \big( (g_k^{dd} )^{-1} \psi\big)_{|z=0}\OpH \big(\text{$\bb_{1,k}$}_{|z=0}\big)^\star   v_k}
  {v_k}_{L^2(\R^d), L^2(\R^d)}^{\kk\mugbk dt},
\end{align*}
with the adjoint operator with the $\star$-notation understood in the
sense of the inner product $L^2(\R^{d}, (\kk) _{|z=0} \mugbk dt)$,
that is,
\begin{align*}
  \OpH \big(\text{$\bb_{1,k}$}_{|z=0}\big)^\star =
(\tkk^{-1})_{|z=0} \OpH
  \big(\text{$\bar{b}_{1,k}$}_{|z=0}\big)^*(\tkk)_{|z=0},
\end{align*}
where the adjoint with the $*$-notation is understood in the sense of
the inner product $L^2(\R^{d}, dx' dt)$. With the two points of
Lemma~\ref{lemma: bounded 2} one finds
\begin{multline*}
  \big\|\OpH \big({b_{1,k}}_{|z=0}\big) \big((g_k^{dd} )^{-1} \psi\big)_{|z=0}\\
  -
  \big( \psi\tkk^{-1} 
  (g_k^{dd} )^{-1}\big)_{|z=0} \OpH
  \big(\text{$\bar{b}_{1,k}$}_{|z=0}\big)^* (\tkk)_{|z=0}
  \big\|_{\L(L^2(\R^{d}))}
  \lesssim \hk, 
\end{multline*}
yielding the result. 
\end{proof}

%%%%%%%%%%%%%%%%%
% subsection
\subsection{Proof conclusion and further support property of the measure $\nu$}
With Proposition~\ref{prop: commutator}, Corollary~\ref{cor: Lsecond},
and Propositions~\ref{prop: contrib qkpk}, \ref{prop: contrib b0}, and
\ref{prop: contrib b1} one now has
\begin{align}
  \label{eq: Hp mu = nu b1-A}
  - \dup{\mu}{\Hp b} &= 2 \dup{\Im M_{0,1}}{b}\\
  & \quad + \lim_{k \to \infty}
  \biginp{\OpH \big({b_{1,k}}_{|z=0}\big) 
      \big((g_k^{dd} )^{-1} \psi\big)_{|z=0}  v_k}{v_k}
     _{L^2(\R^d), L^2(\R^d)}^{\kk\mugbk dt}.\notag
\end{align}
Let $b \in \symbH(\R^{2d+2})$. 
With Proposition~\ref{prop: Euclidean symbol division} one writes
\begin{align*}
  \chi(\tau) b(y,\eta',\zeta) = b_{0,k}(y,\eta') + b_{1,k}(y,\eta') \zeta 
  + q_k(y,\eta',\zeta) \,  p_k(y,\eta',\zeta),
\end{align*}
and 
\begin{align}
  \label{eq: Euclidean division final}
  \chi(\tau) b(y,\eta',\zeta) = b_{0}(y,\eta') + b_{1}(y,\eta') \zeta 
  + q(y,\eta',\zeta) \,  p(y,\eta',\zeta).
\end{align}
With \eqref {eq: symbol Euclidean division1} and \eqref{eq: symbol
  Euclidean division2} one finds that $N_{d+1} (b_{1} - b_{1,k}) =
o(1)_{k\to \infty}$. Since also $\Norm{\tkk - \tk}{L^\infty}=
o(1)_{k\to \infty}$ one gets
\begin{align*}
  &\biginp{\OpH \big({b_{1,k}}_{|z=0}\big) 
    \big((g_k^{dd} )^{-1} \psi\big)_{|z=0}  v_k}{v_k}
     _{L^2(\R^d), L^2(\R^d)}^{\kk\mugbk dt}\\
  &\qquad = \biginp{\OpH \big({b_{1}}_{|z=0}\big) 
    \big((g_k^{dd} )^{-1} \psi\big)_{|z=0}  v_k}{v_k}
     _{L^2(\R^d), L^2(\R^d)}^{\k\mugb dt}
    + o(1)_{k\to+\infty}.
\end{align*}
Since the \scm of $(v_k)$ is $\nu$, with \eqref{eq: Hp mu = nu b1-A}, Lemma~\ref{lem: continu - manifold}
and \eqref{eq: lem: continu - manifold - local chart}, one obtains
\begin{align}
  \label{eq: Hp mu = nu b1-B}
  - \dup{\mu}{\Hp b} = 2 \dup{\Im M_{0,1}}{b} +\dup{\nu}{{b_{1}}_{|z=0}},
\end{align}
as $\psi _{|z=0}  = 1$ in a \nhd of $\supp({b_{1}}_{|z=0})$ and 
as $\Norm{(g_k^{dd} )^{-1}_{|z=0}  - 1}{L^\infty}\to 0$ as $k \to
+\infty$ since $ g^{dd}_{|z=0} =1$ in the chosen quasi-normal
geodesic coordinates associated with the metric $g$; see Proposition~\ref{prop: quasi-normal coordinates}.

\medskip
With the results obtained above, the description of $\supp \nu$
in \eqref{eq: hyp support mesure-new-nu}
can be refined. 
%%%%%%%%%%%%%%%%%%%%%%%%
% proposition          %
%%%%%%%%%%%%%%%%%%%%%%%%
\begin{proposition}
  \label{prop: supp nu 2}
  One has $\supp \nu \subset \big( \pHb \cup \pGb\big) \cap \{ C_{\mu,0}\leq \tau\leq  C_{\mu,1}\}$.
\end{proposition}
%%%% 
\begin{proof}
  The inclusion $\tau \in [C_{\mu,0}, C_{\mu,1}]$ is given in \eqref{eq: hyp support mesure-new-nu}. 

Consider
  $a(y,\eta') \in \Cinfc(\R^{d+1}\times \R^d)$ supported in a \nhd of
  $\{z=0\}$ and with $a(y',z=0,\eta')$
  supported in the elliptic region $\pEb= (\pHb \cup \pGb)^c$. 
  One has 
  \begin{align*}
    \ppi(\Char p \cap \{ z=0\}) \cap \supp a_{|z=0} = \emptyset.
  \end{align*}
  See
  Section~\ref{sec: A partition of the cotangent bundle at the
    boundary-intro}. If the support of $a$ is chosen \suff small in
  the $z$ variable, then this remains true away from $\{ z=0\}$, in
  the sense that 
  \begin{align*}
    a (y', z, \eta') \neq 0  \ \ \imp \ \ p(y',z, \eta', \zeta) \neq 0 \quad
    \forall \zeta \in \R. 
  \end{align*}
  Because of the homogeneity of $p= p_{\k, g}$ and $p_k= p_{\kk,\gk}$,
  by compactness, the same property holds for $p_k$ in place of $p$ for $k$ chosen \suff large.
  For such choice one can set
  \begin{align*}
    q_k (y,\eta) = - \frac{a(y, \eta') \zeta}{p_k (z, \zeta)}.
  \end{align*}
  It is Lipschitz in $y$, smooth and compactly supported in $\eta'$
  and admits a
polyhomogeneous development in the $\zeta$ variable as in \eqref{eq:
  estimation rjk}--\eqref{polyhomog}. In fact, it reads
\begin{align*}
  0 = a(y, \eta') \zeta + q_k (y,\eta) p_k (z, \zeta), 
\end{align*}
precisely of the form given by the Euclidean symbol division of
Proposition~\ref{prop: Euclidean symbol division}, with $b=0$,
$b_{0,k}=0$, and $b_{1,k}=a$. 

With~\eqref{eq: Hp mu = nu b1-A} and the definition of the
semi-classical measure $\nu$, one finds
\begin{align*}
  \dup{\nu}{a_{|z=0}} = - \dup{\mu}{0}
  =  0.
\end{align*}
This gives the result considering the support property of $a_{|z=0}$. 
\end{proof}

\medskip
Let $b \in \symbH(\R^{2d+2})$ as above, with $b_1 (y,\eta') $ defined by \eqref{eq: Euclidean division final}.
Let $\y' = (y, z=0, \eta') \in \supp \nu$.  With
Proposition~\ref{prop: supp nu 2} one has $\y' \in \pHb\cup \pGb$ and
$\tau \in [C_{\mu,0}, C_{\mu,1}]$.  Let $\zeta^\pm$ be
defined as in \eqref{eq: relevement H-G-intro} and
$\y^\pm = (\y', \zeta^\pm) \in \Hb^\pm \cup \Gb$.  With \eqref{eq:
  Euclidean division final} one finds
\begin{align*}
  b(\y^+)- b(\y^-) = b_{1} (\y') ( \zeta^+ - \zeta^-) \quad
  \si \  \y' \in \supp \nu.
 \end{align*}
Hence, in $\supp \nu$ the function 
$\big(b(\y^+)- b(\y^-)\big) / (\zeta^+ - \zeta^-)$ 
is well defined even for points $\y' \in \pGb$. One has 
\begin{align*}
  b_{1} (\y') = \frac{ b(\y^+)- b(\y^-)}{\zeta^+ - \zeta^-}
  = 
  \frac{\dup{\delta_{\y^+} - \delta_{\y^-}}{b} }
    {\dup{\xi^+- \xi^-}{\nx}_{T_x^*\M, T_x\M}},\qquad \y' \in \supp \nu,
\end{align*}
since $\nx$ is  here the unitary inward
    pointing normal vector in the sense of the metric $g$. 
Hence, from \eqref{eq: Hp mu = nu b1-B} and Proposition~\ref{prop: supp nu 2} one concludes the proof of Theorem~\ref {thm: equationmesure}.
\hfill \qedsymbol \endproof

%%%%%%%%%%%%%%%%%
% section
%%%%%%%%%%%%%%%%%
\section{Measure equation at isochrones and necessary geometric
  control condition}
\label{sec: Measure equation at isochrones and necessary geometric control condition}
In this section and in Section~\ref{sec: Proof the measure equation at
  isochrones}, we suppose $(\M,\k,g) \in \X^1$ chosen fixed, and write $p =
p_{\k,g}$. Thus, $\Hp = \Hamiltonian_{p_{\k,g}}$
   % subsection
%%%%%%%%%%%%%%%%%
\subsection{Necessary geometric control condition}
\label{sec: Necessary geometric control condition}
The geometric conditions we formulate
here state that given any point
of $\y^0 \in T^* \L$ at least one \bichar that goes though $\y^0$ reaches a point above the observation region.
%%%%%%%%%%%%%%%%%%%%%%%%
% definition           %
%%%%%%%%%%%%%%%%%%%%%%%%
\begin{definition}[weak interior geometric control condition]
\label{weak control-geo-interior}
Let $\omega$ be an open subset of $\M$ and $T>0$. One
says that $(\omega,T)$ fulfills the weak interior geometric
control condition if for any $\y^0 \in \Char p \cap \TL$ and for any \nhd $V$  of $[0,T]\times \ovl{\omega}$, at least one
\gbichar that goes through $\y^0$ reaches a point
above $V$. 
\end{definition}
%%%%%%%%%%%%%%%%%%%%%%%%
% definition           %
%%%%%%%%%%%%%%%%%%%%%%%%
\begin{definition}[weak boundary geometric control condition]
\label{weak control-geo-boundary}
Let $\Gamma$ be an open subset of $\d \M$ and $T>0$. One says that
$(\Gamma,T)$ fulfills the weak boundary geometric control
condition if for any $\y^0 \in \Char p \cap \TL$ and any \nhd $V_\d$ of $[0,T]\times \ovl{\Gamma}$, at least one \gbichar
that goes through $\y^0$ encounters a boundary escape point (see
Definition~\ref{def: escape point}) above
$V_\d$.
\end{definition}

The following theorem states the result of Theorem~\ref{theorem:
  necessary control conditions-intro} in the framework of the precise
Definitions~\ref{weak control-geo-interior} and \ref{weak control-geo-boundary}.
%%%%%%%%%%%%%%%%%%%%%%%%
% theorem         %
%%%%%%%%%%%%%%%%%%%%%%%%
\begin{theorem}
  \label{theorem: necessary control conditions}
  \begin{enumerate}
  \item Interior observability (Definition~\ref{def:
      observability interior-intro}) implies the weak interior
    geometric control condition.
  \item Boundary observability (Definition~\ref{def:
      observability boundary-intro}) implies the weak boundary
    geometric control condition. 
  \end{enumerate}
\end{theorem}
This theorem is proven in Section~\ref{sec: Proof of the necessary
  GCC}. Its proof uses a measure equation similar to that of
Theorem~\ref{thm: equationmesure}, yet across isochrones $\{t=\cst\}$,
and the construction of concentrating initial conditions.
%%%%%%%%%%%%%%%%%%%%%%%%
% remark               %
%%%%%%%%%%%%%%%%%%%%%%%%
\begin{remark}
  \label{remark: case of uniqueness}
  \begin{itemize}
  \item 
  In the case of uniqueness of \gbichars, the weak geometric control
  condition stated here coincides with the usual necessary condition
  for observability to hold. 
  \item If one replaces the rough cut-off
    $\unitfunction{[0,T]\times\omega}$ (\resp
    $\unitfunction{[0,T]\times\Gamma}$) by the smoother version
    $1_{[0,T]}\Theta(x)$, then the (properly modified) geometric
    control condition is a necessary and sufficient condition in the
    case of uniqueness of \gbichars (see~\cite{BG:1997}). However, when
    uniqueness does not hold this is no more the case as there is still
    the discrepancy between the two conditions ({\em necessary}: at
    least one \gbichar reaches the set $[0,T] \times \{\Theta
    >0\}$; sufficient: {\em all }\gbichars reach this
    set).
  \end{itemize}
  \end{remark}
 For a
 given point $\y^0= (x^0, t^0, \tau^0, \xi^0) \in \TL$, the proof of
 Theorem~\ref{theorem: necessary control conditions} requires the
 construction of a sequence of initial data spectrally localized such
 that the measure of the associated sequence of solutions is supported
 on \gbichars passing through $\y^0$. This is performed in several
 steps:
\begin{itemize}
\item By an explicit calculation, it is possible to do so if one
  forgets the spectral localization (See Section~\ref{sec:
    Concentration at a point}).
\item We then apply the spectral dyadic projector. Here, the difficulty
  comes from the low regularity assumptions on the coefficients (see
  Section~\ref{sec: Dyadic projection}).
\item We prove a transport equation which allows one to transfer the
  information on the traces of the solutions at $t=t^0$ to $\{t>t^0\}$ (see
  Section~\ref{sec: Proof of the necessary GCC}).
\end{itemize}

% subsection
  %%%%%%%%%%%%%%%%% 
 \subsection{Measure equation at isochrones}
\label{sec: Measure equation at isochrone}

 With $\ubt\in \R$, we consider the isochrone $\I = \{ t=\ubt\}$  In
 $\L$.  We naturally identify $\I$ with $\M$, and $\TM$ with $T^*\I$. For 
 $(x,\xi)\in \TM$, identified with $\y =(\ubt,x,\tau=0,\xi)$,  the
 polynomial  $\tau \mapsto p(\ubt, x,\tau, \xi)$ has exactly two roots
 $\tau^+(\y)>0$ and $\tau^- (\y)= - \tau^+ (\y)<0$. If compared to Section~\ref{sec: A partition of the cotangent bundle at
   the boundary-intro} one only faces hyperbolic points in the
 present situation. Set
    \begin{align*}
      %\label{eq: y plus moins}
      \y^{\oplus} = (\ubt, x,\tau^+(\y), \xi), \qquad 
      \y^{\ominus} = (\ubt, x,\tau^-(\y), \xi).
    \end{align*}
Denote by $a_{\k,g}(x,\xi)$ the principal symbol of $A_{\k,g}$, that
is, $a_{\k,g}(x,\xi) = - g_x^{ij} \xi_i \xi_j$ in local coordinates. 
Suppose $H =(\hk)$ is a scale.
For each $k$, suppose $\udl{u}_k^0 \in H^1_0(\M)$,  $\udl{u}_k^1 \in
L^2(\M)$, $f_k \in L_{\loc}^2(\L)$, and $u_k$ is a weak solution to 
\begin{align*}
  %\label{eq: wave equation-measure-eq}
  \begin{cases}
     P_{\k, g}  \, u_k =f_k
     & \dans\ \R\times\M,\\
     u_k=0 &  \dans\ \R \times\d\M,\\
      {u_k}_{|t=\ubt}  = \udl{u}_k^0, \ \d_t  {u_k}_{|t=\ubt}  = \udl{u}_k^1 
      & \dans\ \M.
  \end{cases}
\end{align*}
One extends the diffent functions  by zero outside $\M$ and $\L$.
Suppose the following holds.
 \begin{assumption}
   \label{assum: isochrone measure equation}
      \begin{enumerate}
      \item 
      The sequences $(\udl{u}_k^0)_k$ and $(\hk \udl{u}_k^1)_k$ are both bounded
      in $L^2(\hM)$ and $\udl{U}_k= \transp(\udl{u}_k^0, \hk
      \udl{u}_k^1)$ admits at scale $H$ the Hermitian \scm on $\ThM$ 
      \begin{align*}
        \nu^0 = \begin{pmatrix}
        \nu^0_{0,0} & \nu^0_{0,1}\\
        \nu^0_{1,0} & \nu^0_{1,1}
      \end{pmatrix}
    \end{align*}
   supported away from $\d\M$.
  \item The sequences $(u_k)_k$ and $(\hk f_k)_k$ are both bounded in
    $L_{\loc}^2(\hL)$, and $\transp(u_k, \hk f_k)_k$ admits at scale
    $H$ the Hermitian \scm on $\ThL$
    \begin{align*}
      M = \begin{pmatrix}
        M_{0,0} & M_{0,1}\\
        M_{1,0} & M_{1,1}
      \end{pmatrix}.
    \end{align*}
Set $\mu = M_{0,0}$. 
  \item No mass leaks at infinity at scale $H$ for $(\psi(t) u_k)_k$ and $(\hk
    \psi(t) f_k)_k$, for any $\psi \in \Cinfc(\R)$, and there exists $C>0$ such that, for any interval $I
    \subset \R$,
    \begin{align}
      \label{eq: non concentration assumption}
      \Norm{u_k}{L^2(I \times \M)} + \Norm{\hk f_k}{L^2(I \times \M)}
      \leq C |I|, \qquad k \in \N. 
    \end{align}
    \item One has 
      \begin{align}
        \label{eq: hyp isochrone support mesure mu t=0}
        \supp \mu \subset \Char p \cap \TL \setminus 0
        \ \ \et \ \
        \supp \nu^0 \subset \TM  \setminus 0.
      \end{align}
    \end{enumerate}  
 \end{assumption}
 The sequence $ \transp( \unitfunction{t>\ubt}u_k,  \unitfunction{t>\ubt} \hk f_k)_k$ admits at scale
    $H$ a Hermitian \scm $M^+$ on $\ThL$, with the following 
 natural connection with $M$.
%%%%%%%%%%%%%%%%%%%%%%%%
% lemma                %
%%%%%%%%%%%%%%%%%%%%%%%%
\begin{lemma}
  \label{lemma: step function measure}
  One has $M^+ = \unitfunction{t>\ubt} \, M$.
\end{lemma}
A proof is given below.
One sets $\mu^+ = M_{0,0}^+= \unitfunction{t>\ubt}\,  \mu$. 

\medskip
At $t=\ubt$ and away from $\d\L$ the measure $\mu$ is solution to $\Hp \mu=0$.
The measure equation we establish concerns $\mu^+$ and involves $M_{0,1}^+$ and 
the Hermitian measure $\nu^0$.
%%%%%%%%%%%%%%%%%%%%%%%%
    % theorem              %
    %%%%%%%%%%%%%%%%%%%%%%%% 
    \begin{theorem}
      \label{theorem: measure equation at t=0}  
      Suppose $\Omega$ is an open subset  of $\M$ with $\ovl{\Omega} \cap
      \d\M = \emptyset$.
      In $T^* (\R\times\Omega)$ one has
  \begin{align}
    \label{eq: GL equation t=0}
    \Hp \mu^+  = - \transp{\Hp} \mu^+ 
    &= 2 \Im M_{0,1}^+
    + \int_{\y \in \TM} 
    \frac{\delta_{\y^{\oplus}} - \delta_{\y^{\ominus}}}
    {\tau^+- \tau^-} \ d (a_{\k,g}\,\nu^0_{0,0} - \nu^0_{1,1})(\y)\\
    &\quad + \int_{\y \in \TM} 
    (\delta_{\y^{\oplus}} + \delta_{\y^{\ominus}})
   \ d  \Im \nu^0_{0,1}(\y),\notag
  \end{align}
   in the sense of distributions.
\end{theorem}
A proof is given in Section~\ref{sec: Proof the measure equation at isochrones}.

   \begin{remark}
     \label{rem: support assumption mu t=0}
     The open subset $\Omega$ is introduced as the measure
       equation~\eqref{eq: GL equation t=0} is only proven to hold away
       from the boundary $\d\L$.
     
     % $\phantom{-}$
   %   \begin{enumerate}
   %   \item The open subset $\Omega$ is introduced as the measure
   %     equation~\eqref{eq: GL equation t=0} is only proven to hold away
   %     from the boundary $\d\L$.
   %   \item
   % Since $|\tau| = |a_{\k,g}(x,\xi)|^{1/2}$ on $\Char p$ and
   % $|a_{\k,g}(x,\xi)|^{1/2}\eqsim |\xi|$ uniformly in $x$, note that the
   % support assumption for $\mu^+$ also reads
   % \begin{align*}
   %      \supp(\mu^+) \subset \Char p \cap \TL 
   %      \cap \{ C_{\mu,0}'\leq |\tau| \leq  C_{\mu,1}'\}.
   %    \end{align*}
   %    However, note that the cotangent variable $\tau$ has no meaning for the measure
   %    $\nu_0$ that is defined on $T^* \M$.
   %    \end{enumerate}
   \end{remark}

   In the simpler context of the wave coefficients with constant
   coefficients, one can find in \cite[Proposition~4.4]{PG:96} a
   result expressing the measure $\mu$ by means of measures associated
   with intial conditions. In the more general context we have here,
   deriving a formula for $\mu$ or $\mu^+$ is not possible. Yet, the
   result of Theorem~\ref{theorem: measure equation at t=0}   provides
   a transport equation solved by $\mu^+$.

\medskip
%%%% proof of lemma
\begin{proof}[\bfseries Proof of Lemma~\ref{lemma: step function measure}]
  For simplicity we consider $\ubt =0$ here without loss of
  generality. 

  We prove that $M_{0,1}^+ = \unitfunction{t>\ubt}\, M_{0,1}$. The
  proof is the same for the other matrix entries. Set $v_k = \hk f_k$.
  Suppose that $\beta \in \Cinfc(\R)$ with $0\leq \beta \leq 1$ and
  $\beta(0)=1$. Then, for $\beta_n (t) = \beta(n t)$, with
  Proposition~\ref{prop: extension measure} (adapted to Hermitian
  measures) and dominated convergence, one obtains
  \begin{align*}
    \lim_{k \to +\infty} \inp{\beta_n u_k}{v_k}_{L^2(\hL)}
    = \dup{M_{0,1}}{\beta_n} \ 
    \mathop{\longrightarrow}_{n \to \infty} \ 
    \dup{M_{0,1}}{\unitfunction{\{t=0\}}},
  \end{align*}
  using that no mass leaks at infinity for both $u_k$ and $v_k$ by
  Assumption~\ref{assum: isochrone measure equation}.  By \eqref{eq:
    non concentration assumption} one has$\big| \inp{\beta_n
    u_k}{v_k}_{L^2(\hL)} \big| \lesssim 1/n$ uniformly in $k$.  Thus
  one finds
  \begin{align}
  \label{eq: no mass concentration t=0-bis}
  \unitfunction{\{t=0\}}M_{0,1}=0.
  \end{align}
Suppose $\chi\in\Cinf(\R)$ is such that $0 \leq \chi \leq 1$ and
$\chi(t)=0$ if $t<0$ and $\chi(t) =1$ if $t>1$. Set $\chi_n (t) =
\chi( n t )$.  Consider $b \in \symbc(T^* \hL)$, $B^h$ a
representative of $[\OpH(b)]$, and $\psi \Cinfc(\hL)$ with $\psi=1$ on
the $(t,x)$-projection of $\supp(b)$. One writes
  \begin{align*}
    &\inp{B^h \psi \unitfunction{t>0} u_k}
    {\unitfunction{t>0} v_k}_{L^2(\hL)}\\
    &\quad = \inp{B^h \psi \unitfunction{t>0}  u_k}
    {(\unitfunction{t>0} -\chi_n) v_k}_{L^2(\hL)}
    + \inp{B^h \psi (\unitfunction{t>0} -\chi_n)  u_k}
    {\chi_n  v_k}_{L^2(\hL)}\\
    &\qquad + \inp{B^h \psi \chi_n  u_k}
    {\chi_n  v_k}_{L^2(\hL)}
  \end{align*}
Let $\eps>0$.  Since $B^h$ is bounded on $L^2(\hL)$, $(u_k)_k$ and
$(v_k)_k$ are bounded in $L^2_{\loc}(\hL)$, by \eqref{eq: non
  concentration assumption} there exists $n_0\in \N$ such that,
\begin{align*}
    \big| \inp{B^h \psi\unitfunction{t>0} u_k}
    {(\unitfunction{t>0} -\chi_{n}) v_k}_{L^2(\hL)}
    + \inp{B^h \psi (\unitfunction{t>0} -\chi_{n}) u_k}
    {\chi_n  v_k}_{L^2(\hL)}\big| \leq \eps,
  \end{align*}
uniformly in $k$, for $n \geq n_0$.
There exists also $n_1 \geq n_0$ such that 
\begin{align*}
  |\dup{M_{0,1}}{b  (\unitfunction{t>0} -\chi_{n}^2) }| \leq \eps
\end{align*}
for $n \geq n_1$ by dominated convergence using \eqref{eq: no mass
  concentration t=0-bis}. One thus concludes that 
\begin{multline*}
  \big| \inp{B^h \psi \unitfunction{t>0} u_k}
    {\unitfunction{t>0} v_k}_{L^2(\hL)}
  - \dup{M_{0,1}}{\unitfunction{t>0} b}\big|\\
  \leq 2 \eps 
  + \big| \inp{B^h \psi \chi_n  u_k}
    {\chi_n  v_k}_{L^2(\hL)}
  - \dup{M_{0,1}}{ \chi_{n}^2  b}\big|,
\end{multline*}
for $n \geq n_1$ and $k \in \N$. Set $n=n_1$. Then, there exists $k_0\in \N$ such
that 
\begin{align*}
  \big| \inp{B^h \psi \unitfunction{t>0} u_k}
    {\unitfunction{t>0} v_k}_{L^2(\hL)}
  - \dup{M_{0,1}}{\unitfunction{t>0} b}\big|
  \leq 3 \eps,
\end{align*}
for $k \geq k_0$ implying the result. 
\end{proof}

%%%%%%%%%%%%%%%%%%%%%%%%
% remark               %
%%%%%%%%%%%%%%%%%%%%%%%%
\begin{remark}
  \label{remark: step function measure}
  Similarly one proves that the off-diagonal entries of the measures of $\transp(\unitfunction{t>0} u_k, \hk f_k)$ and $\transp( u_k, \unitfunction{t>0} \hk f_k)$ are also given by
  $M^+_{0,1} = \unitfunction{t>0} M_{0,1}$ and $M^+_{1,0}=\unitfunction{t>0} M_{1,0}$.
\end{remark}

% subsection
%%%%%%%%%%%%%%%%%
\subsection{Concentration at a point}
\label{sec: Concentration at a point}
Pick $\psi \in \mathcal \S(\R^d)$, $x^0 \in \R^d$, and $\xi^0 \in \R^d
\setminus 0$. Set
\begin{align*}
  w_h (x) = h^{-d/4} e^{i\scp{x}{\xi^0}/h}
  \psi\big(h^{-1/2}(x-x^0)\big).
\end{align*}
One has
$\Norm{w_h}{L^2} = \Norm{\psi}{L^2}$ and $(w_h)_h$ admits the measure
$\Norm{\psi}{L^2}^2 \, \delta_{(x^0,\xi^0)}$ as its semi-classical
measure (at scale $h$); this follows from  computations
based on oscillatory-integral arguments.
%%%%%%%%%%%%%%%%%%%%%%%%
% lemma                %
%%%%%%%%%%%%%%%%%%%%%%%%
\begin{lemma}
  \label{lemma: Hs norm wh}
  For all $s \in \R$ one has $\Norm{(-\Delta)^{s/2} w_h}{L^2}
  \sim h^{-s} |\xi^0|^{s}\Norm{\psi}{L^2}$.
\end{lemma}
%%%% proof of lemma
\begin{proof}
  Write
  $\hw_h (\xi) = h^{d/4}\hpsi \big(h^{1/2}\xi-h^{-1/2}\xi^0\big)$
  assuming $x^0=0$ without any loss of generality.  One then computes
  \begin{align*}
    \Norm{|\xi|^s \hw_h}{L^2}^2 = h^{-2s} \int
  |h^{1/2}\xi+\xi^0|^{2s}\, |\hpsi(\xi)|^2 d \xi \sim h^{-2s} |\xi^0|^{2s}
  \Norm{\hpsi}{L^2}^2,
  \end{align*}
  by the \Ldct.
\end{proof}

If $b(x,\xi) \in \Symbol{0,0}{-N}{2d}$ for some $N$, note that $\OpH(b)$ as in
\eqref{4.pseudo} makes sense if acting on $\S(\R^d)$.
%%%%%%%%%%%%%%%%%%%%%%%%
% Proposition              %
%%%%%%%%%%%%%%%%%%%%%%%%
\begin{proposition}
  \label{prop: prep Helffer-Sjostrand}
  Suppose $b(x,\xi) \in \Symbol{0,d+2}{-N}{2d}$ for some $N\in \R$, and 
  \begin{align*}
    b_j(x,\xi) = \jp{\xi}^{-m_j} \int_0^1 \d_{\xi_j} b(x,s\xi + \xi^0)\,
    d s \in \symb(\R^{2d}) = \Symbol{0,d+1}{d+1}{2d},
  \end{align*}
  $j=1, \dots, d$, for some $m_j\geq 0$. 
  One has 
  \begin{align*}
   \OpH(b) w_h
   = b(x,\xi^0) w_h+\sum_{1\leq j\leq d}M_{0,d+1}^{-(d+1)} (b_j ) \, O(h^{1/2}) 
   \ \dans\ L^2(\R^d) \quad \text{as}\ \ h \to 0.
  \end{align*} 
\end{proposition}
\begin{proof}
  Set $\psi_h(x) = \psi\big(x-h^{-1/2} x^0)\big)$ and 
  \begin{align*}
    v_h =  h^{-d/4} \psi\big(h^{-1/2}(x-x^0)\big)
    = h^{-d/4} \psi_h (h^{-1/2}x).
    \end{align*} Then, $w_h = b(x,h D + \xi_0) v_h$ and
  $\OpH(b) w_h = e^{i\scp{x}{\xi^0}/h} b(x,h D + \xi_0) v_h$ by standard
  computations, yielding with $q(x,\xi) = b(x,\xi + \xi^0) -b(x, \xi^0)$,
  \begin{align*}
    \tilde{w}_h = \big( \OpH(b) - b(x,\xi^0)\big) w_h
    &=  e^{i\scp{x}{\xi^0}/h} \big( b(x,h D + \xi_0) - b(x,\xi^0)
      \big)v_h\\
      &=\frac{h^{-d/4}e^{i \scp{x}{ \xi^0}/h} }{(2\pi)^d} \int e^{i \scp{h^{- 1/2} x}{\xi}}
    q(x, h^{1/2}\xi) \hpsi_h (\xi)\, d \xi.
 \end{align*}
  Note that $\Norm{\tilde{w}_h}{L^2} = \Norm{u_h}{L^2}$ with $u_h$ given by 
  \begin{align*}
    u_h(x)  =\frac{1}{(2\pi)^d} \int e^{i \scp{x}{\xi}} q(h^{1/2} x,h^{1/2} \xi)
    \hpsi_h (\xi)\, d \xi. 
  \end{align*}
  Write 
  $q(x,\xi) = \sum_j  \xi_j  \int_0^1 \d_{\xi_j} b(x,s\xi + \xi^0)\,  d s$.
  This gives $u_h(x) = h^{1/2}\sum_j u_{j,h}(x)$ with 
  \begin{align*}
    u_{j,h}(x) & =\frac{1}{(2\pi)^d} \int e^{i \scp{x }{\xi}}
                 b_j(h^{1/2} x, h^{1/2}\xi)
    \jp{h^{1/2} \xi}^{m_j} \xi_j \hpsi_h (\xi)\, d \xi,
  \end{align*}
  with $b_j$ as in the proposition statement. Set $\psi_{h,j} = \jp{h^{1/2} D}^{m_j} D_{x_j} \psi_h$.
  Then $u_{j,h}(x) = b_j (h^{1/2}x,h^{1/2}D) \psi_{h,j} (x)$, that is, a
  semi-classical operator acting on $\psi_{h,j}$, yet with $h$ replaced by
  $h^{1/2}$. First, observe that $\Norm{\psi_{h,j}}{L^2}$ is bounded uniformly in
  $h$.  Second, with Lemma~\ref{4.continuite} one finds
  \begin{align*}
    \Norm{u_{j,h}}{L^2}\lesssim M_{0,d+1}^{-(d+1)} \big(b_j (h^{1/2}x,\xi)\big),
  \end{align*}
  which concludes the proof since $M_{0,d+1}^{-(d+1)} \big(b_j
  (h^{1/2}x,\xi)\big) = M_{0,d+1}^{-(d+1)} (b_j)$. 
\end{proof}
Set $b_z(x,\xi) = (z+a_{\k,g}(x,\xi))^{-1}$
with $z\in \C \setminus \R$. Then, $b_z \in \Symbol{0,d+2}{-2}{2d}$. Define
$b_j$ as in the statement of Proposition~\ref{prop: prep
  Helffer-Sjostrand} with $b_z$ in place of $b$ and $m_j  =m= 2d+3$.
%%%%%%%%%%%%%%%%%%%%%%%%
% lemma                %
%%%%%%%%%%%%%%%%%%%%%%%%
\begin{lemma}
  \label{lemma: estimate qj}
  There exists $C>0$ such that $M_{0,d+1}^{-(d+1)} (b_j)\leq C |\Im z|^{-3-d}$. 
\end{lemma}
%%%% proof of lemma
\begin{proof}
First, one has
$\d_{\xi_j} b_z(x,\xi) = - (z+a_{\k,g}(x,\xi))^{-2} \d_{\xi_j} a_{\k,g}(x,\xi)$, 
implying 
\begin{align*}
  \jp{\xi}^{d+1} \big| b_j (x,\xi)\big| 
  \lesssim  \jp{\xi}^{d+1-m} |\Im z|^{-2} (|\xi^0| + |\xi|) 
  \lesssim |\Im z|^{-2},
\end{align*}
as $m= 2 d+3$.
Second, note that $\d_{\xi}^{\beta} \big( \jp{\xi}^{-m}\d_{\xi_j}
b_z(x,s\xi + \xi^0)\big)$ is equal to a linear combination of terms 
\begin{align*}
\d_{\xi}^{\beta_1} \jp{\xi}^{-m} \,  s^{|\beta_2|}
  \big(\d_{\xi}^{\beta_2} \d_{\xi_j} b_z\big) (x,s\xi + \xi^0),
  \ \ \avec \ \ \beta_1+ \beta_2 = \beta.
\end{align*}
As $\big|\big(\d_{\xi}^{\beta_2} \d_{\xi_j} b_z\big) (x,s\xi + \xi^0)\big|
\lesssim |\Im z|^{-|\beta_2|-2} (|\xi^0| + |\xi|)^{|\beta_2|+1}$ one
obtains
\begin{align*}
  \jp{\xi}^{d+1} \big|\d_{\xi}^{\beta} b_j (h^{1/2}x,\xi)\big| 
  \lesssim \jp{\xi}^{d+2+ |\beta| -m}  |\Im z|^{-|\beta|-2}
  \lesssim |\Im z|^{-d-3},
\end{align*}
for $|\beta| \leq d+1$ as $m = 2d+3$. 
\end{proof}
%%%%%%%%%%%%%%%%%%%%%%%%
% Corollary            %
%%%%%%%%%%%%%%%%%%%%%%%%
\begin{corollary}
  \label{cor: prep Helffer-Sjostrand}
 Set $b_z(x,\xi) = (z+a_{\k,g}(x,\xi))^{-1}$ with $z\in \C \setminus \R$. One has 
 \begin{align*}
   \OpH(b_z) w_h
   = b_z(x,\xi^0) w_h+h^{1/2} |\Im z|^{-d-3} O(1)
   \ \dans\ L^2(\R^d) \quad \text{as}\ \ h \to 0.
 \end{align*}
\end{corollary}

%%%%%%%%%%%%%%%%%
% subsection
\subsection{Dyadic projection}\label{sec: Dyadic projection}

Consider $(x^0, \xi^0)\in T^*\M$ with $x^0 \notin \d\M$ and
$\hchart=(\hO, \chdiff)$ a local chart with $x^0 \in \hO$. In this
local chart, introduce $(w_h)_h$ as above:
\begin{align}
  \label{eq: wh before dyadic proj}
  w_h (x) = h^{-d/4} e^{i\scp{x}{\xi^0}/h}
  \psi\big(h^{-1/2}(x-x^0)\big),
\end{align}
with $\psi \in \Cinfc(\R^d)$ with $\psi =1$ in a \nhd of $0$ here.
Consider a scale $H=(h_{k})_k$ and choose $k$ \suff large so that
$\supp v_k\subset \hO$, with $v_k = \chdiff^* w_k$ for $w_k= w_{\hk}$. 

As above, denote by $a_{\k,g}(x,\xi)$ the symbol of the operator
$A_{\k,g}$, that is, in local coordinates, $a_{\k,g}(x,\xi) =
-g^{ij}(x) \xi_i \xi_j$.  Suppose $0< \alpha<1$, and $\chi \in
\Cinfc(]\alpha^2, \alpha^{-2}[)$, with $\chi =1$ on $[\alpha,
      \alpha^{-1}]$.  Here, we prove that $\chi(-a_{\k,g}(x,\xi^0))
    v_{k}$ and $\chi(-\hk^2 A_{\k,g}) v_{k}$ coincide in $L^2(\M)$ up
    to a $o(1)$ remainder. One can view $\chi(-\hk A_{\k,g}) v_{k}$ as
    some ``projection'' of $v_k$ onto the dyadic subspace $E_{k}$
    introduced in Section~\ref{dyadic}.
%%%%%%%%%%%%%%%%%%%%%%%%
% lemma                %
%%%%%%%%%%%%%%%%%%%%%%%%
\begin{lemma}
  \label{lemma: approx resolvent}
  Suppose $\theta, \ttheta
  \in \Cinfc(\chdiff(\hO))$, with $\ttheta=1$ on a \nhd of
  $\supp \theta$.
  Set $\htheta = \chdiff^* \theta \in \Conc^2(\hO)$.
For $z \in \C \setminus\R$ set $b_z(x,\xi) = \big(a_{\k,g}(x,\xi)+z\big)^{-1}$. One has
\begin{align*}
  (h^{2}A_{\k,g}+z)^{-1} \htheta
  = \chdiff^{*}\theta \OpH(b_z)  \ttheta \big(\chdiff^{-1}\big)^{*}
  + R_{h},
\end{align*}
with $\Norm{R_{h}}{\L(L^{2})} = O(h) |\Im(z)|^{-1}\big(1 + |z|^{1/2} \big)$.
\end{lemma}
A proof is given below. Note that if $z \in \C \setminus\R$,
the operator $\OpH(b_z)$ is  well defined and bounded on $L^2(\R^d)$ by Lemma~\ref{lemma: refined L2-continuity} given below. 

\bigskip With $\chi$  as above, consider 
$\tchi \in \Cinfc(\C)$ an almost analytic extension of
$\chi$. The Helffer-Sj\"ostrand formula~\cite{DS:99} gives
\begin{align}
  \label{eq: Helffer-Sjostrand}
  \chi (-h^{2}A_{\k,g}) = \frac{1}{2i\pi}  \lim_{\epsilon \to 0^+}
  \int_{|\Im z|\geq \epsilon} \bar{\d}\tchi(z) (h^{2}A_{\k,g}+z)^{-1}
  dz\wedge  d\bar{z}.
\end{align}
The function $\tchi$ has the following properties: $\tchi_{|\R} =
\chi$, there exists $C>0$ such that $\supp \tchi \subset \supp \chi
+ i [-C, C]$, and for any $n \in \N$ there exists $C_n>0$ such that 
\begin{align*}
  %\label{eq: almost analytic extension}
  |\bar{\d} \tchi (z)| \leq C_n |\Im z|^n.
\end{align*}

Choose $\theta, \ttheta$ as in Lemma~\ref{lemma: approx resolvent} and
$\htheta = \chdiff^* \theta$.  One obtains
\begin{align*}
 \chi (-h^{2}A_{\k,g}) \htheta
  &= \frac{1}{2i\pi}  \lim_{\epsilon \to 0^+}
  \int_{|\Im z|\geq \epsilon} \bar{\d}\tchi(z) (h^{2}A_{\k,g}+z)^{-1}\htheta
  dz\wedge  d\bar{z}\\
  &=\frac{1 }{2i\pi}  \chdiff^{*} \theta
    \lim_{\epsilon \to 0^+}
    \int_{|\Im z|\geq \epsilon} \bar{\d}\tchi(z) \OpH(b_z)
    dz\wedge  d\bar{z} \ttheta \big(\chdiff^{-1}\big)^{*} + O(h)_{\L(L^2)}\\
  &=\frac{1 }{2i\pi}  \chdiff^{*} \theta \OpH\Big(
    \lim_{\epsilon \to 0^+}
    \int_{|\Im z|\geq \epsilon} \bar{\d}\tchi(z) b_z
    dz\wedge d\bar{z} \Big) \ttheta \big(\chdiff^{-1}\big)^{*} + O(h) _{\L(L^2)}\\
  &=\chdiff^{*} \theta \OpH \big( \chi( -a_{\k,g})\big)
    \ttheta \big(\chdiff^{-1}\big)^{*} + O(h) _{\L(L^2)},
\end{align*}
meaning that
\begin{align}
  \label{eq: chi A-Helffer-Sjostrand}
  \chi (-h^{2}A_{\k,g}) =  [\OpH] \big( \chi( -a_{\k,g})\big),
\end{align}
with the notation introduced in Section~\ref{sec: Semi-classical operators on a manifold}.

Consider now $\theta$ such that 
$\htheta =1$ on $\supp v_k$, for $k$ \suff large.
One has 
\begin{align*}
  \chi (-h^{2}A_{\k,g}) v_k = \chdiff^{*} \theta \OpH \big( \chi( -a_{\k,g})\big)
    w_k
  + O(h)_{L^2},
\end{align*}
yielding, with Proposition~\ref{prop: prep Helffer-Sjostrand},
\begin{align}
  \label{eq: dyadic projection approx}
  \chi (-h^{2}A_{\k,g}) v_k
  = \chi\big( -a_{\k,g}(x,\xi^0) \big) v_k
  + O(h^{1/2})_{L^2(\M)}.
\end{align}
Since $(v_k)_k$ has $\Norm{\psi}{L^2}^2 \, \delta_{(x^{0},\xi^{0})}$
for \scm, one has the following result.
%%%%%%%%%%%%%%%%%%%%%%%%
% lemma                %
%%%%%%%%%%%%%%%%%%%%%%%%
\begin{lemma}
  \label{lemma: scm after projection}
  The two sequences $\chi (-h^{2}A_{\k,g}) v_k$ and $\chi\big(-
  a_{\k,g}(x,\xi^0)\big) v_k$ have the same \scm, that is,
$\big|\chi\big(-a_{\k,g}(x^{0},\xi^{0})\big)\big|^{2} \Norm{\psi}{L^2}^2 \, \delta_{(x^{0},\xi^{0})}$.
\end{lemma}

\bigskip
%%%%%%%%%%%%%%%%%%%%%%%
%%%%%%%%%%%%%%%%%%%%%%%
\begin{proof}[\bfseries Proof of Lemma~\ref{lemma: approx resolvent}]
For $z \in  \C \setminus\R$ and $w \in \S(\R^{d})$ one has 
\begin{equation*}
  \theta \OpH(b_z)  (\ttheta w)(x)
  = (2\pi)^{-d} \int e^{i \scp{x}{\xi}} \theta(x) b_z(x,h \xi)\widehat
  {\ttheta w}(\xi)\,  d\xi .
\end{equation*}
With the form of the differential operator $A_{\k,g}$ in local coordinates
compute 
\begin{align*}
  &h^2 A_{\k,g}
    \Big( e^{i \scp{x}{\xi}} \theta(x) b_z(x,h \xi) \Big)\\
  &\quad
    = i\tk^{-1}\sum_{i,j} h \d_{x_i}
    \Big(e^{i \scp{x}{\xi}} \tk  g^{ij} \theta  b_z(x,h\xi) h \xi_j \Big)\\
    &\qquad +  h\tk^{-1}\sum_{i,j} h \d_{x_i} \Big(e^{i \scp{x}{\xi}}\tk g^{ij}
    \d_{x_j} \big( \theta b_z\big) (x,h\xi) \Big)\\
  &\quad
    =e^{i \scp{x}{\xi}} (\theta a \, b_z)(x,h\xi)
    +i  h \tk^{-1}e^{i \scp{x}{\xi}} m(x,h \xi) 
  + h^2\tk^{-1} e^{i \scp{x}{\xi}}
    \sum_{i} \d_{x_i} \ell_i (x,h \xi),
\end{align*}
where
\begin{align*}
  m(x,\xi) = \sum_{i,j}\Big(
    \d_{x_{i}}(\tk  g^{ij}\theta b_z) (x,\xi) \xi_j
    + \tk  g^{ij} \d_{x_{j}}(\theta b_z) (x,\xi) \xi_i
    \Big),
\end{align*}
with $\ell_i = \sum_{j} \tk g^{ij} \d_{x_j} \big( \theta b_z\big)$. 
We deduce that
\begin{equation}
  \label{eq: resolvent}
(h^{2}A_{\k,g}+z) \theta \OpH(b_z) \ttheta =   \theta+ i h \tk^{-1}\OpH(m) \ttheta+
h^{2}\sum_{i}\tk^{-1} \d_{x_i}\OpH(\ell_i) \ttheta.
\end{equation}
One checks that $m$ fulfills the assumptions of
Lemma~\ref{lemma: refined L2-continuity} with $\delta =1$ and so do
the symbols $\ell_{i}$, $1\leq i \leq d$,  with $\delta =2$,
implying that $\OpH(m)$, and $\OpH(\ell_i)$ are bounded on $L^{2}(\R^{d})$. 

The following bounds hold for the resolvent
\begin{align*}
  &\Norm{(h^{2}A_{\k,g}+z)^{-1}}{\L(L^{2}(\M))}\leq |\Im z|^{-1},\\
  &\Norm{(h^{2}A_{\k,g}+z)^{-1}}{\L(L^{2}(\M), H_0^1 (\M))}\leq |h \Im z|^{-1}
  (|\Re z| + |\Im z|)^{1/2}.
\end{align*}
From the second estimate one deduces also that
\begin{align*}
  \Norm{(h^{2}A_{\k,g}+z)^{-1}}{\L(H^{-1}(\M), L^{2}(\M))}\leq |h \Im z|^{-1}
  (|\Re z| + |\Im z|)^{1/2}.
\end{align*}
One thus obtains
\begin{align*}
  &\Norm{(h^{2}A_{\k,g}+z)^{-1} \chdiff^{*}
    \OpH(m) \ttheta \big(\chdiff^{-1}\big)^{*}}{\L(L^{2}(\M)}
    \lesssim |\Im(z)| ^{-1} ,
\end{align*} 
and
\begin{align*}
  &\Norm{(h^{2}A_{\k,g}+z)^{-1} \tk^{-1} \chdiff^{*}
    \d_{x_i}\OpH(\ell_i)
    \ttheta\big(\chdiff^{-1}\big)^{*}}{\L(L^{2}(\M))}\\
  &\qquad \lesssim |h \Im(z)|^{-1} 
    |z|^{1/2} \Norm{\tk}{W^{1,\infty}}, \quad i=1, \dots, d.
\end{align*} 
If one applies the resolvent $(h^{2}A_{\k,g}+z)^{-1}$ to the left of
identity~\eqref{eq: resolvent}, one then obtains
\begin{align*}
  &\Norm{(h^{2}A_{\k,g}+z)^{-1} \htheta
    - \chdiff^{*} \theta \OpH(b_z) \ttheta
    \big(\chdiff^{-1}\big)^{*}}{\L(L^{2}(\M))}\\
  &\qquad \lesssim h | \Im(z)| ^{-1} \big(1 + |z|^{1/2} \Norm{\tk}{W^{1,\infty}}\big),
\end{align*}
which gives the result.
\end{proof}

\bigskip
For $\delta \geq 0$ set
\begin{align*}
  L^{-\delta}_{0,d+1}(a) = \max_{|\beta|\leq d+1}\sup_{(x,\xi)}
  \big| \d_{\xi}^{\beta}a(x,\xi)| \jp{\xi}^{|\beta| + \delta} .
\end{align*}
Compare $ L^{-\delta}_{0,d+1}$ and
$M^{-(d+1)}_{0,d+1}$. Here, less decay is expected on $a(x,\xi)$; yet
decay improves with differentiations with respect to $\xi$. 
%%%%%%%%%%%%%%%%%%%%%%%%
% lemma                %
%%%%%%%%%%%%%%%%%%%%%%%%
\begin{lemma}
  \label{lemma: refined L2-continuity}
Suppose $a(x,\xi) \in L^{\infty}(\R^{2d})$ is smooth in $\xi$ and 
$L^{-\delta}_{0,d+1}(a) < \infty$
for some  $\delta > 0$.
Then, $\OpH(a)$ is bounded on $L^{2}(\R^{d})$ and
\begin{align*}
  \Norm{\OpH(a)}{\mathcal{L}(L^{2}(\R^d))} \leq C_{\delta,d} L^{-\delta}_{0,d+1}(a).
  \end{align*}
\end{lemma}
Compare with Lemma~\ref{4.continuite}.
%%%%%%%%%%%%%%%%%%%%%%%
%%%%%%%%%%%%%%%%%%%%%%%
\begin{proof}%[Proof of Lemma~\ref{lemma: refined L2-continuity}]
  Consider $\theta \in \Cinfc(\R^d)$such that $0\leq \theta\leq 1$, 
  $\theta(\xi) =1$ if $|\xi| \leq 1/2$, and  $\theta(\xi) =0$ if
  $|\xi| \geq 1$.  Set $\psi_0 = \theta$ and 
  \begin{equation*}
	\psi( \xi) = \theta(\xi) - \theta(2\xi) \ \  \et \ \
        \psi_j(\xi) = \psi(2^{-j} \xi ) \ \pour\  j \in \N^*, 
\end{equation*}
yielding a dyadic partition of unity
$1= \sum_{j \in \N}\psi_j$. Set  $a_{j}(x,\xi)=\psi_{j}(\xi)a(x,\xi)$.
With Lemma~\ref{4.continuite} one finds
\begin{align*}
  \Norm{\OpH(a_0)}{\mathcal{L}(L^{2}(\R^d))} \lesssim
  M^{-(d+1)}_{0,d+1}(a_0) \lesssim L^{0}_{0,d+1}(a) \lesssim L^{-\delta}_{0,d+1}(a),
\end{align*}
since $\psi_0$ has compact support. Consider now $j \geq 1$. With
$\tilde{h}_{j}=2^{-j} \hk$ one writes 
\begin{align*}
  \OpH(a_{j}) v(x)
  &= (2\pi )^{-d}\int e^{i \scp{x}{\xi}}
    \psi(2^{-j} \hk \xi) a(x,\hk \xi) \hat{v}(\xi)\, d\xi\\
   &= (2\pi )^{-d}\int e^{i \scp{x}{\xi}}
     \psi(\tilde{h}_{j}\xi) a(x, \tilde{h}_{j} 2^j \xi) \hat{v}(\xi)\, d\xi\\
     &=\Op^{\tilde{h}_{j}} (b_j) v(x).
\end{align*}
with $b_j (x,\xi) = \psi(\xi) a(x, 2^j \xi)$. The symbol $b_j$ is
compactly supported in $\xi$ and for $\beta \in \N^d$, with $|\beta|
\leq 1+d$, one finds
\begin{align*}
  \jp{\xi}^{d+1}\big| \d_{\xi}^{\beta} b_{j}(x,\xi)\big|
  &\lesssim  \jp{\xi}^{d+1}\sum_{\beta'+\beta''=\beta}
  \big|\d_{\xi}^{\beta'}\psi(\xi) \big|
  \big|\d_{\xi}^{\beta''}a(x,2^{j}\xi) \big|\\
  &\lesssim L^{-\delta}_{0,d+1}(a)  \jp{\xi}^{d+1} \sum_{\beta'+\beta''=\beta}
     2^{|\beta''| j} \jp{2^j \xi}^{-|\beta''| - \delta} \big|\d_{\xi}^{\beta'}\psi(\xi) \big|.
\end{align*}
Since $|\xi| \gtrsim 1$ in the compact  $\supp \psi$, one obtains
\begin{align*}
  \jp{\xi}^{d+1}\big| \d_{\xi}^{\beta} b_{j}(x,\xi)\big|
  &\lesssim 2^{- \delta j }  L^{-\delta}_{0,d+1}(a)  \sum_{\beta'+\beta''=\beta}
   \jp{\xi}^{d+1}  \big|\d_{\xi}^{\beta'}\psi(\xi) \big|
    \lesssim 2^{-\delta j}  L^{-\delta}_{0,d+1}(a).
\end{align*}
Lemma~\ref{4.continuite} implies
$\Norm{\Op^{\tilde{h}_{j}} (b_j)}{\mathcal{L}(L^{2}(\R^d))} \lesssim 2^{-\delta j }  L^{-\delta}_{0,d+1}(a)$.
Convergence of $\sum_{j} 2^{-\delta j}$ gives the conclusion.
\end{proof}

In what follows we will also need the following results.
%%%%%%%%%%%%%%%%%%%%%%%%
% lemma                %
%%%%%%%%%%%%%%%%%%%%%%%%
\begin{lemma}
  \label{lemma: boundedness constructed solution}
  There exists $C>0$ such that
  $\Norm{\hk \nablag \chi (-\hk^{2}A_{\k,g}) v_k}{L^2(\M)} \leq C$.
  If $(x^0,\xi^0)$ in the definitions of $w_k$ in \eqref{eq: wh
    before dyadic proj} and $v_k= \chdiff^* w_k$ is chosen such that
  $\chi\big(- a_{\k,g}(x^0,\xi^0) \big) \neq 0$, then there exists $C'>0$
  such that 
  \begin{align*}
    &1 /C \leq  \Norm{\chi (-\hk^{2}A_{\k,g}) v_k}{L^2(\M)} \leq C
     \ \ \et  \ \
    1/C \leq  \Norm{\hk \nablag \chi (-\hk^{2}A_{\k,g}) v_k}{L^2(\M)} \leq C,
  \end{align*}
  for $k$ \suff large. 
\end{lemma}
%%%%%%%%%%%%%%%%%%%%%%%%
% lemma                %
%%%%%%%%%%%%%%%%%%%%%%%%
\begin{lemma}
  \label{lemma: H1sc conv concentration}
  One has
  \begin{align*}
    \hk \nablag \, \chi (-\hk^{2}A_{\k,g}) v_k
    = \hk \nablag \, \big(\chi\big(-
    a_{\k,g}(x,\xi^0)\big) v_k\big) + O(\hk^{1/4})_{L^2(\M)}.
    \end{align*}
\end{lemma}

%%%% proof of lemma
\begin{proof}[\bf Proof of Lemma~\ref{lemma: boundedness constructed solution}]
  Set $\tilde{v}_k = \chi (-\hk^{2}A_{\k,g}) v_k$.  One writes
  \begin{align}
    \label{eq: proof boundedness constructed solution}
    \Norm{\hk \nablag \tilde{v}_k }{L^2(\M)}^2 = \biginp{\tchi
      (-\hk^{2}A_{\k,g}) v_k}{\tilde{v}_k }_{L^2(\M)},
  \end{align}
  with $\tchi(\lambda) = \lambda\, \tchi(\lambda)$. The same analysis
 used for $\tilde{v}_k$ applies to $\tchi (-\hk^{2}A_{\k,g}) v_k$. In
  particular $\tchi (-\hk^{2}A_{\k,g}) v_k = \tchi\big(-
  a_{\k,g}(x,\xi^0)\big) v_k + O(\hk^{1/2})_{L^2}$. One thus obtains the
  first result
  
  There exists a \nhd $V$ of $x^0$ such that $\big|\chi\big(
  -a_{\k,g}(x,\xi^0)\big)\big| \gtrsim 1$ for $x \in V$. For $k$ \suff
  large $\supp v_k \subset V$ implying
  \begin{align*}
    \Norm{\chi\big(- a_{\k,g}(x,\xi^0)\big)  v_k }{L^2(\M)} \gtrsim \Norm{v_k
    }{L^2(\M)} \gtrsim 1.
  \end{align*}
  With \eqref{eq: dyadic projection approx} one concludes that $\Norm{\tilde{v}_k}{L^2(\M)} \gtrsim 1$. 

  Arguing the same with
  \eqref{eq: proof boundedness constructed solution} and using
  that $\big| \tchi \, \chi\big(- a_{\k,g}(x,\xi^0)\big) \big| \gtrsim 1$ in a \nhd of
  $x^0$, one obtains that $\Norm{\hk \nablag \tilde{v}_k}{L^2(\M)} \gtrsim 1$.
\end{proof}
\begin{proof}[\bf Proof of Lemma~\ref{lemma: H1sc conv concentration}]
  Set $z_k = \chi (-\hk^{2}A_{\k,g}) v_k - \chi\big(-
  a_{\k,g}(x,\xi^0)\big) v_k$.  With \eqref{eq: dyadic projection
    approx} one has $\Norm{z_k}{L^2(\M)} = O(\hk^{1/2})$.
  Lemma~\ref{lemma: boundedness constructed solution} gives a
  $L^2$-bound for the sequence $\hk \nablag \chi (-\hk^{2}A_{\k,g}) v_k$ and a
  simple computation gives $\hk \nablag \chi\big(-
  a_{\k,g}(x,\xi^0)\big) v_k$ also $L^2$-bounded. Hence, a
  preliminary estimate is $\Norm{\hk \nablag z_k}{L^2(\M)} = O(1)$.

  Compute $\Norm{\hk \nablag z_k}{L^2(\M)}^2 = N_1 + N_2$ with 
  \begin{align*}
    &N_1 = \inp{- \hk^2A_{\k,g} \chi (-\hk^{2}A_{\k,g}) v_k}{z_k}_{L^2(\M)},
    \\
    &N_2 =
    \biginp{ \hk \nablag \chi\big(- a_{\k,g}(x,\xi^0)\big) v_k}
        {\hk \nablag z_k}_{L^2(\M)}.
    \end{align*}
  Note that $- \hk^2 A_{\k,g}\chi (-\hk^{2}A_{\k,g}) v_k = \tchi
  (-\hk^{2}A_{\k,g}) v_k$, with $\tchi(\lambda) = \lambda\,
  \tchi(\lambda)$, is $L^2$-bounded since the same analysis used for
  $\chi (-\hk^{2}A_{\k,g}) v_k$ applies. Hence, $N_1 = O(\hk^{1/2})$.
  Writing
  \begin{align*}
    \hk \nablag \chi\big(- a_{\k,g}(x,\xi^0)\big)) v_k
    = \hk \big[\nablag, \chi\big(- a_{\k,g}(x,\xi^0)\big)\big]  v_k
    + \chi\big(- a_{\k,g}(x,\xi^0)\big) \hk \nablag v_k,
  \end{align*}
  one finds
  \begin{align*}
    N_2 =
    \biginp{\chi\big(- a_{\k,g}(x,\xi^0)\big)) \hk \nablag  v_k}{\hk \nablag z_k}_{L^2(\M)}
    + O(\hk).
  \end{align*}
  With a similar commutator computation one further obtains
  \begin{align*}
    N_2 =
    - \biginp{\chi\big(- a_{\k,g}(x,\xi^0)\big) \hk^{2}A_{\k,g} v_k}{z_k}_{L^2(\M)}
    + O(\hk).
  \end{align*}
  With Lemma~\ref{lemma: Hs norm wh} one concludes that $N_2 = O(\hk^{1/2})$. 
\end{proof}

 %%%%%%%%%%%%%%%%%%%%%%%
 %%%%%%%%%%%%%%%%%%%%%%%
 
% \todo{Nicolas et Belhassen: Continuité de $Op^{h}a(x,\xi)$ sans
% décroissance en $\xi$} \begin{lemma}\label{L3-continuity} Suppose
% $a(x,\xi)$ is $\Con^{1}$ in $x$ and $\Cinf$ in $\xi$
% and \begin{align} M_{0}(a) = \max_{\vert\beta \vert \leq 1,
% \vert\alpha\vert \leq d+1}\sup_{(x,\xi)}\vert
% \d_{x}^{\beta}\d_{\xi}^{\alpha}a(x,\xi)\vert(1+\vert\xi\vert)^{\vert\alpha\vert}
% < \infty .\end{align} Then, $Op^{h}(a)$ is bounded on
% $L^{2}(\R^{d})$ with $\Vert Op^{h}(a)\Vert_{\mathcal{L}(L^{2})} \leq
% C_{d} M_{0}(a)$.
%\end{lemma}

%%%%%%%%%%%%%%%%%
% subsection
\subsection{Proof of the necessary geometric control condition}
\label{sec: Proof of the necessary GCC}
Here, we prove Theorem~\ref{theorem: necessary control conditions}.
Assume that observability holds and yet the condition of
Definition~\ref{weak control-geo-interior} (\resp Definition~\ref{weak
  control-geo-boundary}) does not hold. This section aims to reach a
contradiction.

If the weak interior geometric control condition does not hold, there
exist $\y^0= (t^0, x^0, \tau^0, \xi^0) \in \Char p \cap \TL$ and $V$
an open 
\nhd of $[0,T]\times \ovl{\omega}$ such that no \gbichar going through
$\y^0$ reaches a point above $V$.  If the weak boundary geometric
control condition does not hold, there exist $\y^0$ and $V_\d$ an open
\nhd
of $[0,T]\times \ovl{\Gamma}$ such that no \gbichar going through
$\y^0$ reaches a boundary escape point above $V_\d$.  
\subsubsection{Interior initial point} We first treat the case $\y^0 \in \Char(p)\cap (\TL\setminus
\d\TL)$. The case $\y^0 \in \Char(p)\cap \d\TL$ is treated in a second
round.

One has $\tau^0\neq 0$. With some scaling in the cotangent variables, one may assume $|\tau^0| \in
[\alpha, \alpha^{-1}]$ for some  $0< \alpha< 1$.
One has $(\tau^0)^2 = - a_{\k,g}(x^0,\xi^0) = g^{ij}(x^0) \xi_i^0 \xi_j^0$. 

Suppose $\chi\in \Cinfc(\R)$ is such that $\chi \geq 0$,
$\supp \chi \subset ]\alpha^2, \alpha^{-2}[$, and 
$\chi\big((\tau^0)^2\big) =\Norm{\psi}{L^2}^{-1}$, with $\psi$ used
in \eqref {eq: wh before dyadic proj}. With the sequence $(v_k)_k\subset L^2(\M)$
constructed above set
\begin{align*}
  \udl{u}^0_k= \chi (-h^2 A_{\k,g}) v_k \ \ \et \ \
  \udl{u}^1_k= i \hk^{-1}\tau^0\udl{u}^0_k
\end{align*}
and denote by $u_k$ the solution to the {\em homogeneous} wave equation
\begin{align}
  \label{eq: wave equation-measure-eq-homog}
  \begin{cases}
     P_{\k, g}  \, u_k =0
     & \dans\ \R\times\M,\\
     u_k=0 &  \dans\ \R \times\d\M,\\
      {u_k}_{|t=t^0}  = \udl{u}_k^0, \ \d_t  {u_k}_{|t=t^0}  = \udl{u}_k^1 
      & \dans\ \M.
  \end{cases}
\end{align}
Since $\udl{u}^0_k\subset E_k$,
with $\udl{u}^0_k= \sum _{\nu \in J_k} u^0_{k,\nu} e_\nu$, one finds
\begin{align*}
  u_k = \sum _{\nu \in J_k}  \big( e^{i (t-t^0) \sqrt {\lambda _\nu}}u_{k,\nu} 
   + e^{-i  (t-t^0)  \sqrt {\lambda _\nu}}  u_{k,-\nu}\big)  e_\nu,
\end{align*}
with $u_{k,\pm \nu} = u^0_{k,\nu}\big( 1 \pm
\hk^{-1}\lambda_\nu^{-1/2}\tau^0 \big) /2$.
With Lemma~\ref{lemma: boundedness constructed solution} observe that 
\begin{align}
  \label{eq: boundedness constructed solution2} 
  \E^h(u_k) = \frac12 \big( \Norm{\hk \nablag \udl{u}^0_k}{L^2}^2 
  + \Norm{\hk \udl{u}^1_k}{L^2}^2\big) 
  =
  \frac12 \big( \Norm{\hk \nablag \udl{u}^0_k}{L^2}^2 
  + (\tau^0)^2 \Norm{\udl{u}^0_k}{L^2}^2\big) 
  \eqsim 1.
\end{align}
The solution $(u_k)_k$ is  bounded in $L_{\loc}^2(\L)$ as in
Section~\ref{sec: Initiation of the contradiction argument}, and can
be associated with a \scm $\mu$ (up to a possible subsequence
extraction). Associated with $\hk \d_\n {u_k}_{|\d\L}$ is a \scm
$\nu$. Arguing as in Proposition~\ref{prop: first property measures} one finds 
\begin{align}
  \label{eq: measure support properties isochrone}
    &\supp \mu \subset \Char p \cap \TL \cap \{ \alpha\leq
  |\tau| \leq  \alpha^{-1}\},
  \\
  &\supp \nu \subset T^*\d\L \cap \{ \alpha\leq
    |\tau| \leq  \alpha^{-1}\}.\notag
\end{align}
Hence, Theorem~\ref{thm: equationmesure} applies and 
both measure $\mu$ and $\nu$ satisfy the measure propagation
equation~\eqref{eq: GL equation-propagation theorem}.
As in the proof of Proposition~\ref{prop: first property measures} one
finds that no mass leaks at infinity at scale $H$, in the sense of Definition~\ref{def: mass leak at infinity}.

As the \scm of $(v_k)_k$ at scale $H=(\hk)_k$ is
$\Norm{\psi}{L^2}^{2} \delta_{(x^0,\xi^0)}$, with Lemma~\ref{lemma:
  scm after projection} one finds the following
results.
%%%%%%%%%%%%%%%%%%%%%%%%
% lemma                %
%%%%%%%%%%%%%%%%%%%%%%%%
\begin{lemma}
  \label{lemma: scm y1k}
  The sequence  $(\udl{u}^0_k, \hk \udl{u}^1_k)_k$ admits the Hermitian \scm 
  \begin{align*}
    \nu^0
    = \begin{pmatrix}
        \nu^0_{0,0} & \nu^0_{0,1}\\
        \nu^0_{1,0} & \nu^0_{1,1}
      \end{pmatrix}
                    = \begin{pmatrix}
        1 & -i \tau^0 \\
         i \tau^0 & (\tau^0)^2
      \end{pmatrix} \delta_{(x^0,\xi^0)}
  \end{align*}
   on $\TM$   at scale $H=(\hk)_k$.
\end{lemma}

As in Section~\ref{sec: Measure equation at isochrone} denote by
$\mu^+$ the \scm associated with  $( \unitfunction{t>t^0} u_k)_k$.
By Lemma~\ref{lemma: step function measure} one has $\mu^+ =
\unitfunction{t>t^0} \, \mu$. Observe that
\begin{align*}
  a\,\nu^0_{0,0} - \nu^0_{1,1}
  = - 2 (\tau^0)^2 \delta_{(x^0,\xi^0)}
  \ \ \et \ \
  \Im \nu^0_{0,1}
  = -\tau^0 \delta_{(x^0,\xi^0)}.
\end{align*}
Theorem~\ref{theorem: measure equation at t=0}
applies and, as $\tau^+- \tau^- = 2 \big| a_{\k,g}(x^0, \xi^0)\big|^{1/2} = 2 |\tau^0|$,
one obtains
\begin{align}
  \label{eq: transport near t=0}
    \Hp \mu^+  = - \transp{\Hp} \mu^+ 
    = -2\tau^0\ \delta_{\y^0},
\end{align}
away from $\dTL$,
using that 
for 
$\y = (t^0,x^0,0, \xi^0)$ one has $\y^0 = \y^{\oplus}$ if $\tau^0>0$ and $\y^0 = \y^{\ominus}$ if $\tau^0<0$.
%%%%%%%%%%%%%%%%%%%%%%%%
% lemma                %
%%%%%%%%%%%%%%%%%%%%%%%%
\begin{lemma}
  \label{lemma: mu vanishes near boundary}
  The measure $\mu$ vanishes in a \nhd of $\{ t=t^0\} \cap \dTL$.
\end{lemma}
A proof is given below. With Lemma~\ref{lemma: mu vanishes near boundary} and \eqref{eq: transport near t=0}, 
one concludes that $\supp \mu^+ \cap \{t = t^0\} = \{ \y^0\}$. As $\mu^+ = \unitfunction{t>t^0}\mu$ one
also has $\supp \mu \cap \{t = t^0\} = \{ \y^0\}$.
With Theorem~\ref{theo-propagation} one obtains the following lemma. 
%%%%%%%%%%%%%%%%%%%%%%%%
% lemma                %
%%%%%%%%%%%%%%%%%%%%%%%%
\begin{lemma}
  \label{lemma: necessary condition supp mu y0}
  The support of $\mu$ is a union of maximal \gbichars that go through $\y^0$. 
\end{lemma}

\medskip
\noindent {\bf Case 1: interior observation.}
If interior observablity holds, then inequality \eqref{eq: interior
  observability-intro} is valid for the sequence $(u_k)_k$. By
\eqref{eq: boundedness constructed solution2} one has
$\Norm{\bld{1}_{]0,T[ \times \omega}\, \hk \d_t u_k}{L^2(\L)}
\gtrsim 1$,
implying 
\begin{align}
  \label{eq: nonzero measure interior observation}
  \supp \mu \cap T^*V\neq \emptyset.
\end{align}
The open set $V$ is introduced in the beginning of the proof.
In fact, consider
$\varphi \in \Cinfc(\L)$ nonnegative such that $\supp \varphi \subset V$ and
$\varphi=1$ in a \nhd of $[0,T] \times \ovl{\omega}$. With
Proposition~\ref{prop: extension measure} 
one finds
$\dup{\mu}{\varphi \tau^2}
  =\lim_{k\rightarrow + \infty} \inp{\varphi \hk \d_t u_k}
  {\hk \d_t u_k}_{L^2(\L)}
  \gtrsim 1$,
yielding \eqref{eq: nonzero measure interior observation}.
With Lemma~\ref{lemma: necessary condition supp mu y0} however, the existence
of a point in $\supp \mu \cap T^*V$ yields a contradiction with the
choice of the point $\y^0$ made at the beginning of the proof.

\bigskip \noindent {\bf Case~2: boundary observation.}  If boundary
observablity holds, then inequality \eqref{eq: boundary
  observability-intro} is valid for $(u_k)_k$. With \eqref{eq:
  boundedness constructed solution2} one has $\Norm{\bld{1}_{]0,T[
    \times \Gamma}\, \d_{\n} u_{|\R\times \d\M}}{L^2(\d\L)} \gtrsim
1$, implying that $\supp \nu \cap T^*V_\d\neq \emptyset$; the open
set $V_\d$ is introduced in the beginning of the proof. Suppose $\py^1
= (t^1, x^1, \tau^1, \xi^1) \in T^*V_\d$.
  \begin{description}
  \item[Case $\bld{\py^1 \in \pHb}$, a hyperpoblic point] Denote
    by $\y^{1,\pm} \in \Hb^\pm$ the points such that $\ppi(\y^{1,\pm}) = \py^1$.
    They are boundary escape points. With Lemma~\ref{lemma: necessary
      condition supp mu y0} the existence of such a point in
    $\supp \mu$ yields a contradiction. Thus $\mu=0$ locally near
    these points.  With Theorem~\ref{thm: equationmesure}, near a
    hyperbolic point one has $\transp \Hp \mu = \tilde{\mu} \otimes
    \delta_{z=0}$, for $\tilde{\mu}$ some
    measure on $\dTL$.  Here, one has $\y^{1,\pm} \notin
    \supp(\tilde{\mu} \otimes \delta_{z=0})$ implying $\py^1 \notin
    \supp \nu$.  One concludes that $\pHb\cap \supp \nu =
    \emptyset$.

  \item[Case $\bld{\py^1= \y^1 \in \Gb \cap \ESC}$, a glancing escape point]
    $\y^1 \in
    \pGb = \Gb$. If $\y^1 \in \supp \mu$ one reaches a contradiction with
    Lemma~\ref{lemma: necessary condition supp mu y0} as $\y^1$ is boundary escape point.
     Thus, locally $\mu=0$.  In local coordinates, in a \nhd $W$ of $\y^1$,
    Theorem~\ref{thm: equationmesure} and Remark~\ref{remark:
      integrand GL equation on G-intro} give
    \begin{align*}
      %\label{eq: action nu escape point}
      \dup{\nu}{\d_\zeta q_{|z=\zeta=0}}=0,
    \end{align*}
    for any $q\in \Cinfc(\R^{2d+2})$ supported in $W$, since there is
    no hyperbolic point in $\supp \nu \cap W$. As
    any compactly supported function $\tilde{q}$ on $\{z=\zeta=0\}$
    can be written in the form $\d_\zeta q_{|z=\zeta=0}$, this implies
    that $\nu$ vanishes in a \nhd of $\y^1$.  One concludes that $\Gb \cap \ESC \cap \supp \nu = \emptyset$.
  \end{description}
  With Proposition~\ref{prop: supp nu 2}, Lemma~\ref{lemma: glancing
    non-escape point}, and the two cases above, one concludes that
  $\supp \nu \subset \Gb \setminus \ESC \subset \sdGb\cup
  \glGb$. Yet, the measure $\nu$ has no mass on this set by
  Proposition~3.5 in the companion article~\cite{BDLR2}, that is, $\dup{\nu}{\unitfunction{\sdGb \cup \glGb}}=0$, implying that
  $\nu$ vanishes; a contradiction.
  
  \bigskip
  \subsubsection{Boundary initial point}
  We now treat the case $\y^0 \in \Char(p)\cap \d \TL$.
  
\noindent   {\bf Case~3: $\bld{\y^0 \in \Char p\cap \d \TL}$ for a interior observation.}
 Suppose  that $V$ is a \nhd of
$[0,T]\times \ovl{\omega}$ in $\d\TL$ such that no \gbichar going through $\y^0$
reaches a  point above $V$. 
Consider $\tilde{V}$ a \nhd of
$[0,T]\times \ovl{\omega}$ in $\d\TL$ such that $\tilde{V} \Subset V$ and $\eps=
\dist (\tilde{V}, V^c)$. For $\udl{\y} = (\udl{t},\udl{x},\udl{\tau},\udl{\xi}) \in \TL$ and $\mathsf T>0$,  set
\begin{align*}
  \Gamma^{\mathsf T}(\udl{y})
  = \{|t-\udl{t}| \leq {\mathsf T} \} \cap \bigcup_{\udl{y} \in \GammaG}
  \GammaG,
\end{align*}
that is, the union of all \gbichar that pass through $\udl{y}$, restricted to the  time interval $[\udl{t}-{\mathsf T}, \udl{t}+{\mathsf T}]$.
%\begin{multline*}
%  \Gamma^S(\y)
%  = \{ \ty =(\tilde{t},\tilde{x},\tilde{\tau},\tilde{\xi}) \in \Char p \cap \TL; \ 
%  |\tilde{t} -t|  \leq S, \ \ty =\gammaG(s), \ \et \ \y =\gammaG(0) \\
%  \pour \ \gammaG \ \text{a maximal \gbichar}\}.
%\end{multline*}

With the continuity result of Proposition~2.11 in the companion
article~\cite{BDLR2}, for $\mathsf T>0$ there exists $\delta>0$ such
for any $\tilde{\y}^0 \in \TL$ one has
\begin{align*}
  \dist(\tilde{\y}^0, \y^0) \leq \delta \ \et \ 
    \y \in \Gamma^{\mathsf T}(\tilde{\y}^0)
    \ \ \imp \ \
    \dist\big(\y, \Gamma^{\mathsf T}(\y^0)\big) \leq \eps/2.   
\end{align*}
Thus, for $\mathsf T$ chosen \suff large there exists
$\tilde{\y}^0 \in \Char p \cap (\TL \setminus \d\TL)$ such that no
\gbichar going through $\tilde{\y}^0$ reaches a point above
$\tilde{V}$, meaning we are back to the configuration
considered above.

\bigskip \noindent {\bf Case 4: $\bld{\y^0 \in \Char p\cap \d \TL}$ for
  a boundary observation.} Note that one cannot argue as in the case of an interior observation
 since $\Gb \setminus \ESC \subset \ovl{\Hb}$. However, the method used here applies to the case of an interior observation; the argument is yet much more involved. Note also that the argument simplifies if \gbichars are uniquely defined, that is, in the presence of a \gbichar flow.

Suppose $V_\d$ is a \nhd of $[0,T]\times
\ovl{\Gamma}$ in $\d\TL$ such that no \gbichar going through $\y^0$
reaches a boundary escape point above $V_\d$.  Write $\y^0 = (t^0,
x^0, \tau^0, \xi^0)$, where $x^0=(x^{\prime,0},z^0)$ with $z^0=0$.
One considers a sequence $(\y^n)_n \subset \Char p \cap \TL
 \setminus \dTL$ such that $\y^n = (t^0, x^{\prime,0},z^n, \tau^n,
 \xi^n) \to \y^0$ as $n \to +\infty$, that is, $z^n \to 0^+$ and
 $(\tau^n, \xi^n) \to (\tau^0, \xi^0)$. With each $\y^n$,
 construct a sequence of solutions $(u_{n,k})_k$ to the wave equation
 as done above, that is, with a \scm $\mu_n$ whose support is a union
 of maximal \gbichars that go through $\y^n$.  One has $\supp \mu_n \subset
 \Char p \cap \{ \alpha \leq |\tau| \leq \alpha^{-1}\}$. With this
 construction, the mass of $\mu_n$ on $T^* \big( (-\mathsf T,\mathsf T)\times
 \hM\big)$ is uniformly bounded for any $\mathsf T>0$.  This implies, that
 $(\mu_n)_n$ is a bounded sequence of measure in the sense of the
 measure topology. Consequently, there exists a measure $\mu$ such
 that $\mu_n \to \mu$ for a subsequence, still denoted by $\mu_n$, in
 the measure topology on $T^* \hL$. One has $\supp \mu \subset
 \Char p \cap \{ \alpha \leq |\tau| \leq \alpha^{-1}\}$.

 Consider $\tM$ a bounded \nhd of $\ovl{\M}$ in $\hM$ and set $\tL = \R \times \tM$.
 Denote by $U$ a \nhd of $\Char p \cap T^*\tL \cap \{ \alpha \leq |\tau| \leq \alpha^{-1}\}$ in $T^* \hL$, such that $U_{\mathsf T} = U \cap \{ |t-t^0| \leq {\mathsf T}\}$ is compact.
 %Note that $\supp \mu_n
 %\cap U_n$ and $\supp \mu\cap U_n$ are compact sets.
 There exists a
 increasing function $\varphi: \N \to \N$ such that 
 \begin{align*}
   %\label{eq: necessary condition conv measure1}
   \big| \dup{\mu-\mu_{\varphi(n)}}{a}\big| \leq \frac{1}{n} \Norm{a}{L^\infty},
   \qquad a \in \Con_c^0(U_n),
 \end{align*}
 recalling that the strong topology is equivalent to the weak topology
 for a converging sequence of measures; see \cite[Section
   34.4]{Treves:67}. Denote by $\mu_n$ this extracted sequence for concision:
 \begin{align}
   \label{eq: necessary condition conv measure1}
   \big| \dup{\mu-\mu_n}{a}\big| \leq \frac{1}{n} \Norm{a}{L^\infty},
   \qquad a \in \Con_c^0(U_n).
 \end{align}

 Consider a sequence $(\psi_n)_n \subset
 \Cinfc(\R)$ such that $\psi_n=1$ on a \nhd of $[-n,n]$ and
 $\supp \psi_n \subset [-n-1, n+1]$. We write $\psi_n$ in place of
 $\psi_n(t)$ for concision. The measure $|\psi_n|^2 \mu_n$ is the limit of
 $\mu_{n,k}= |W \psi_n u_{n,k}|^2$ as $k \to +\infty$ in the measure
 sense, where $W \psi_n u_{n,k}$ is (a variant of) the Wiegner transform
 of $\psi_n(t) u_{n,k}(t,x)$; see \cite{GL:1993}. Thus, there exists a
 increasing function $\tvarphi: \N \to \N$ such that
 \begin{align}
   \label{eq: necessary condition conv measure2}
   \big| \dup{\mu_{n,\tvarphi(n)}-|\psi_n|^2\mu_n}{a}\big|
   =\big| \dup{\mu_{n,\tvarphi(n)}-\mu_n}{a}\big|
   \leq \frac{1}{n} \Norm{a}{L^\infty},
   \qquad a \in \Con_c^0(U_n). 
 \end{align}
 From \eqref{eq: necessary condition conv measure1}--\eqref{eq:
   necessary condition conv measure2}, one finds that
 $\mu_{n,\tvarphi(n)} \to \mu$ on $U$ in the measure sense as $n
 \to \infty$.  It follows that $\mu$ is the \scm of $v_n = u_{n,
   \varphi(n)}$ on $\hL$ at scale $h_{\varphi(n)}$ by
 \cite[Proposition~1.4]{GL:1993}. Denote by $\nu$ the \scm of
 $h_{\varphi(n)} \d_\nu v_n$, by potentially performing yet another
 subsequence extraction. One has $\supp \nu \subset \{ \alpha \leq
 |\tau| \leq \alpha^{-1}\}$. Theorems~\ref{thm: equationmesure} and
 \ref{theo-propagation} apply, implying that $\supp \mu$ is a union
 of maximal \gbichars.

 Suppose ${\mathsf T}>0$ and 
 $a \in \Con_c^0(U_{\mathsf T})$ is such that $\supp a \cap \Gamma^{\mathsf T}(\y^0) =
 \emptyset$. In particular, set $\eps = \dist\big(\supp a,
 \Gamma^{\mathsf T}(\y^0)\big)$. There exists $N\in \N^*$ such that
\begin{align*}
  n \geq N \ \et \ 
    \y \in \Gamma^{\mathsf T}(\y^n)
    \ \ \imp \ \
    \dist\big(\y, \Gamma^{\mathsf T}(\y^0)\big) \leq \eps/2,  
\end{align*}
by Proposition~2.11 in the companion article~\cite{BDLR2}. Because of
the description of $\supp \mu_n$ given above one finds that
$\dup{\mu_n}{a}=0$ if $n \geq N$. With \eqref{eq: necessary condition
  conv measure1} one obtains $\big| \dup{\mu}{a}\big| \leq \frac{1}{n}
\Norm{a}{L^\infty}$ if $n \geq N$ thus giving $\dup{\mu}{a}=0$. Hence,
\begin{align*}
  %\label{eq: supp mu necessary cond boundary case}
  \supp \mu \cap U_{\mathsf T} \subset \Gamma^{\mathsf T}(\y^0).
\end{align*}
One concludes that $\supp \mu$ is a union of maximal \gbichars that all go through $\y^0$. 

One is now in the same position as in the proof of the case of a
boundary observation where $\y^0 \notin \d\TL$. The proof can be
carried out {\em mutatis mutandis}: first, $\supp \nu \cap \pHb=
\emptyset$, second, $\Gb \cap \ESC \cap \supp \nu = \emptyset$
implying that $\supp \nu \subset \sdGb\cup \glGb$ yielding a contradiction.

 This concludes the proof of
Theorem~\ref{theorem: necessary control conditions}.
\hfill \qedsymbol \endproof

\begin{proof}[\bf Proof of Lemma~\ref{lemma: mu vanishes near boundary}]
  Consider $y_k$ solution to the homogeneous wave equation \eqref{eq:
    wave equation-measure-eq-homog} with ${y_k}_{|t=t^0}= \udl{y}^0_k
  = \chi( -a_{\k,g}(x,\xi^0)v_k$ and $\d_t {y_k}_{|t=t^0}=
  \udl{y}^1_k= i \hk^{-1} \tau^0 \udl{y}^0_k$.  Since
  $\supp({y_k}_{|t=t^0})$ and $\supp(\d_t {y_k}_{|t=t^0})$ are away
  from $\d\M$, by finite-speed propagation $y_k$ vanishes in a fixed
  open \nhd $W$ of $\{ t=t^0\} \cap \d\L$.  By \eqref{eq: dyadic
    projection approx} and Lemma~\ref{lemma: H1sc conv concentration} one has 
      \begin{align*}
        \Norm{\hk \nablag (\udl{u}^0_k -
          \udl{y}^0_k)}{L^2(\M)} \to 0
        \ \ \et \ \
        \Norm{\hk \udl{u}^1_k - \hk \udl{y}^1_k}{L^2(\M)} \to 0.
        \end{align*}
      One concludes that the semi-classical energy of $u_k -y_k$ converges to 0.
      Hence, one finds
      \begin{align*}
        \Norm{\hk \nablag u_k}{L^2(W)}^2 + \Norm{\hk \d_t u_k}{L^2(W)}^2 \to 0,
      \end{align*}
      yielding $(|\tau|^2 + |\xi|_x^2) \mu=0$ in $T^*
      W$. In particular, this implies that the support of $\mu$ is
      restricted to the null section in $T^* W$. With \eqref{eq:
        measure support properties isochrone} one obtains that $\mu$
      vanishes in $T^* W$.
\end{proof}
%%%%%%%%%%%%%%%%%
% section
%%%%%%%%%%%%%%%%%
\section{Proof of the measure equation at an isochrone}
\label{sec: Proof the measure equation at isochrones}

Here, we prove Theorem~\ref{theorem: measure equation at t=0}.  We
treat the case $\ubt =0$ without any loss of generality.  At the
hypersurface $t=0$, there is no boundary condition. The two traces
${u_k}_{|t=0} = \udl{u}_k^0$ and $\d_t {u_k}_{|t=0} = \udl{u}_k^1$
have to be taken into account in the analysis. Proceeding as is done
for the measure equations at the boundary $\d\L$ in Sections~\ref{sec:
  Proof of the propagation equation1} and \ref{sec: Proof of the
  propagation equation2} makes a double-layer potential appear, and it
cannot be handeled by the method used therein. We choose to proceed
differently here, letting the measure act on tangential symbols. This
approach can for instance be found in \cite{Robbiano:22} for the
treatment of Zaremba boundary conditions at a boundary.

\subsection{Preliminary filtering}
As in the proof of the measure equation of Theorem~\ref{thm:
  equationmesure} we first apply some filtering to reduce the support
of the measures. The principle is very similar to what is done in the
beginning of Section~\ref{sec: Preliminary remarks and observations},
yet more technical because of the low regularity of the coefficients
of the operator $A_{\k,g}$.

Consider $\chi \in \Cinfc(\R)$ with $0 \notin \supp \chi$. Set
$\tilde{u}_k=\chi(-h^2 A_{\k,g}) u_k$, $\tilde{f}_k=\chi(-h^2
A_{\k,g}) f_k$, $\udl{\tilde{u}}_k^0 = \chi(-h^2 A_{\k,g})
\udl{u}_k^0$, and $ \udl{\tilde{u}}_k^1 = \chi(-h^2 A_{\k,g})
\udl{u}_k^1$. One has
\begin{align*}
  %\label{eq: filtered wave equation-measure-eq}
  \begin{cases}
    P_{\k, g}  \, \tilde{u}_k =\tilde{f}_k
     & \dans\ \R\times\M,\\
     \tilde{u}_k =0 &  \dans\ \R \times\d\M,\\
     {\tilde{u}_k }_{|t=0}  = \udl{\tilde{u}}_k^0, \
     \d_t  {\tilde{u}_k }_{|t=0}  = \udl{\tilde{u}}_k^1 
      & \dans\ \M.
  \end{cases}
\end{align*}
%%%%%%%%%%%%%%%%%%%%%%%%
% proposition          %
%%%%%%%%%%%%%%%%%%%%%%%%
\begin{proposition}
  \label{prop: isochrone filtering}
  The sequence $\transp(\tilde{u}_k, \hk \tilde{f}_k)_k$ admits $|\chi(-a_{\k,g})|^2 M$ as its Hermitian \scm  on $\TL\setminus \dTL$ at scale $H$. The sequence
  $\transp(\udl{\tilde{u}}_k^0, \hk
  \udl{\tilde{u}}_k^1)$  admits $|\chi(-a_{\k,g})|^2 \nu^0$ as its  Hermitian \scm  on $\TM\setminus \dTM$
   at scale $H$.
\end{proposition}
The proof of this intuitive result is given in Section~\ref{sec: proof
  prop: isochrone filtering} below.  Note that the \scm of
$\tilde{u}_k$ also reads $|\chi(\tau^2)|^2 \mu$ from the assumed
support properties.

\medskip
With Lemma~\ref{lemma: step function measure}, the sequence $\transp
(\unitfunction{t>0} \tilde{u}_k, \unitfunction{t>0} \hk
\tilde{f}_k)_k$ has $|\chi(-a_{\k,g})|^2 M^+$ for measure.  If we
prove that the measure equation \eqref{eq: GL equation t=0} holds for
$M^+$ and $\nu^0$ replaced by $|\chi(-a_{\k,g})|^2 M^+$ and
$|\chi(-a_{\k,g})|^2 \nu^0$, then using \eqref{eq: hyp isochrone
  support mesure mu t=0} one finds that \eqref{eq: GL equation t=0}
holds also for $\mu^+$ and $\nu^0$ by the \Ldct.  Without any loss of
generality we may thus replace $u_k$ by $\tilde{u}_k$, $f_k$ by $\tilde{f}_k$, $\udl{u}_k^0$
by $\udl{\tilde{u}}_k^0$, and $\udl{u}_k^1$ by $\udl{\tilde{u}}_k^1$.
Then, there exists $0< C_{\mu,0}< 1< C_{\mu,1}< \infty$ such that
\begin{align}
     \label{eq: hyp support mesure mu t=0}
  \supp \mu^+ \subset \Char p \cap \TL 
  \cap \{ C_{\mu,0}\leq |\xi|\leq  C_{\mu,1}\},
   \end{align}
and 
\begin{align*}
     %\label{eq: hyp support mesure nu t=0}
     \supp \nu^0 \subset \TM  \cap \{ C_{\mu,0}\leq |\xi|\leq  C_{\mu,1}\}.
\end{align*}

Suppose $I$ is a time interval. With the filtering used above, one has
\begin{align}
   \label{eq: mass11-proof 0}
   \Norm{u_k}{L^2(I \times \M)}
   \eqsim \Norm{\hk^2 A_{\k , g} u_k }{L^2(I \times \M)}
   \eqsim \Norm{\hk^2 \d_t^2 u_k - \hk^2 f_k}{L^2(I \times \M)}.
\end{align}
Assume that a subsequence of $u_k$
converges to $0$ in $L^2(I\times \M)$. This gives $\mu=0$ on $T^*(I\times
\M)$. With \eqref{eq: mass11-proof 0}, one finds that $\Norm{\hk^2
  \d_t^2 u_k}{L^2(I \times \M)}\to 0$ and $\Norm{\hk^2 A_{\k , g}
  u_k }{L^2(I \times \M)}\to 0$ also, using that $\hk f_k$ is
$L^2_{\loc}$-bounded.
Then, ellipticity up to the boundary gives
$\Norm{\hk^2 u_k}{H^2(I\times \M)} \to 0$
and 
interpolation gives 
\begin{align*}
  \Norm{\hk \d_t u_k}{L^2(I\times \M)}  \to 0
  \ \ \et \ \
     \Norm{\hk \nabla_{\!\! \gk} u_k}{L^2(I \times \M)} \to 0.
\end{align*}
Since \eqref{eq: non concentration assumption} implies the time continuity
of the semi-classical energy uniformly in $k$ one obtains that
\begin{align*}
  \Norm{\hk \d_t u_k (0,.)}{L^2(\M)}  \to 0
  \ \ \et \ \
  \Norm{\hk \nabla_{\!\! \gk} (0,.)}{L^2(\M)}  \to 0.
\end{align*}
One concludes that $\nu^0=0$.  Hence, all terms in the measure equation vanish, in this case. 
One may thus assume that  $\Norm{u_k}{L^2(I\times \M)}
\gtrsim 1$, for any interval $I$.
Then, one finds that
\begin{align*}
   %\label{eq: mass11-proof 1}
  1 \eqsim  \Norm{u_k}{L^2(I \times \M)}
   \eqsim \Norm{\hk^2 A_{\k , g} u_k }{L^2(I \times \M)}
   \eqsim \Norm{\hk^2 \d_t^2 u_k}{L^2(I \times \M)},
\end{align*}
and one further obtains
\begin{align}
  \label{eq: mass11-proof 2}
  % \label{eq: isochrone - bound after filtering}
  \Norm{\hk^{\ell+ 2\ell'} D_t^\ell A_{\k,g}^{\ell'}
  u_k}{L^2(I \times \M)} \eqsim 1, \qquad \ell, \ell'\in \N,
\end{align}
and
\begin{align*}
  %\label{eq: mass11-proof 3}
  \Norm{\hk \d_t u_k}{L^2(I\times \M)}  \eqsim 
     \Norm{\hk \nabla_{\!\! \gk} u_k}{L^2(I \times \M)} \eqsim 1.
\end{align*}
%Suppose $\varphi \in \Cinfc(\R)$ with $\varphi =1$ in a \nhd of $0$.
By \eqref{eq: mass11-proof 2} one finds that $\unitfunction{t>0}
\psi(t) u_k \in H^s_{sc}(\L)$ for any $s \in [0,1/2[$ and $\psi\in
    \Cinfc(\R)$. From Lemma~\ref{lemma: h-oscillation} one deduces the following result. 
%    With Lemma~4.3 in \cite{Robbiano:22} one obtains, for
%    such $s$,
%\begin{align*}
%  \Norm{\OpH\big(1 - \varphi(\tau/R)\big) \unitfunction{t>0} \psi(t) 
%  u_k}{L^2} \leq C_{s} R^{-s},
%\end{align*}
    %uniformly in $k$. One deduces the following result.
%%%%%%%%%%%%%%%%%%%%%%%%
% lemma                %
%%%%%%%%%%%%%%%%%%%%%%%%
\begin{lemma}
  \label{lemma: no mass leaks 1t u}
   No mass leaks at infinity at scale $H$ for $\unitfunction{t>0} (\psi(t) u_k)_k$.
\end{lemma}

\subsubsection{Proof of Proposition~\ref{prop: isochrone
    filtering}}
\label{sec: proof prop: isochrone filtering}
We prove that the \scm of $\tilde{u}_k$ is $|\chi(-a_{\k,g})|^2 \mu$ on $\TL\setminus \dTL$. The proof for the other sequences and measures are the same. 

Suppose $j \in  \Cinfc(\TL)$ and $\psi \in \Cinfc(\L)$ with $\psi=1$ on the $(t,x)$-projection of $\supp j$ and $\supp(\psi) \cap \d\L = \emptyset$. 
Arguing as in Proposition~\ref{prop: commutator} one proves the
following lemma.
%%%%%%%%%%%%%%%%%%%%%%%%
% lemma                %
%%%%%%%%%%%%%%%%%%%%%%%%
\begin{lemma}
  \label{lemma: commutator A}
  One has $\big[[\OpH](j)\psi, h^2 A_{\k,g}\big] \in h \L(L^2_{\loc}(\L))$.
\end{lemma}
One deduces the following result.
%%%%%%%%%%%%%%%%%%%%%%%%
% lemma                %
%%%%%%%%%%%%%%%%%%%%%%%%
\begin{lemma}
  \label{lemma: commutator A bis}
  One has $\big[[\OpH](j)\psi, (z+h^2 A_{\k,g})^{-1}\big]
  \in h |\Im z|^{-2} \L(L^2_{\loc}(\L))$. 
\end{lemma}
%%%% proof of lemma
\begin{proof}
  Suppose $\Im z\neq 0$. With Lemma~\ref{lemma: commutator A} one
  writes
  \begin{align*}
    (z+h^2 A_{\k,g}) [\OpH](j) \psi- [\OpH](j) \psi  (z+h^2 A_{\k,g})  \in h
    \L(L^2_{\loc}(\L)).
  \end{align*}
  Letting $(z+h^2 A_{\k,g})^{-1}$ act, both from the left and the
  right one obtains the result using that $\bigNorm{(z+h^2
    A_{\k,g})^{-1}}{\L(L^2(\M))} \lesssim |\Im z|^{-1}$.
\end{proof}
With the Helffer-Sj\"ostrand formula \eqref{eq: Helffer-Sjostrand} 
the result extends  to $ \chi (-h^{2}A_{\k,g})$.
%%%%%%%%%%%%%%%%%%%%%%%%
% lemma                %
%%%%%%%%%%%%%%%%%%%%%%%%
\begin{lemma}
  \label{lemma: commutator A ter}
  Suppose $\chi \in \Cinfc(\R)$.
  One has
  \begin{align*}\big[[\OpH](j)\psi, \chi (-h^{2}A_{\k,g})\big]
    \in h \L(L^2_{\loc}(\L)).
  \end{align*}
\end{lemma}
With the previous lemma, one writes
\begin{align*}
  &\biginp{[\OpH](j)\psi \chi (-h^{2}A_{\k,g}) u_k}
  {\chi (-h^{2}A_{\k,g}) u_k}_{L^2(\L)}\\
  &\quad
    =\biginp{\chi (-h^{2}A_{\k,g}) [\OpH](j)\psi \chi (-h^{2}A_{\k,g}) u_k}
    {u_k}_{L^2(\L)}\\
    &\quad= \biginp{\chi (-h^{2}A_{\k,g})^2
      [\OpH](j) \psi u_k}
      {u_k}_{L^2(\L)} + O(\hk)\\
  &\quad= \biginp{\chi^2 (-h^{2}A_{\k,g}) [\OpH](j)  \psi u_k}
    {u_k}_{L^2(\L)} + O(\hk)\\
  &\quad= \biginp{  [\OpH](\chi^2 (-a_{\k,g})) [\OpH](j)\psi  u_k}
    {u_k}_{L^2(\L)} + O(\hk),
\end{align*}
using \eqref{eq: chi A-Helffer-Sjostrand}.

It suffices to prove the result of the proposition with a test
function supported in a local chart. Moreover, it can be chosen of the
form $j(t,x,\tau,\xi) = j_1(t,x) j_2(\tau,\xi)$ with $j_1, j_2
\in\Cinfc(\R^d)$. Then, $\OpH(j) = j_1 \OpH(j_2)$. With a partition of
unity and using Proposition~\ref{propcommut} one has
$\big[[\OpH](\chi^2 (-a_{\k,g})), j_1\big] \in h \L(L^2(\M))$ yielding
\begin{align*}
  &\biginp{[\OpH](j)\psi \chi (-h^{2}A_{\k,g}) u_k}
    {\chi (-h^{2}A_{\k,g}) u_k}_{L^2(\L)}\\
    &\qquad = \biginp{ [\OpH](j_1 j_2 \chi^2 (-a_{\k,g}) ) \psi u_k}
    {u_k}_{L^2(\L)} + O(\hk).
\end{align*}
One thus obtains
\begin{align*}
  &\lim_{k\to +\infty }\biginp{[\OpH](j) \chi (-h^{2}A_{\k,g}) \psi u_k}
    {\chi (-h^{2}A_{\k,g}) u_k}_{L^2(\L)}
    = \dup{\mu}{\chi^2 (-a_{\k,g}) j },
\end{align*}
which is the result of  Proposition~\ref{prop: isochrone
    filtering}.
\hfill \qedsymbol \endproof

\subsection{Symbol decomposition}
\label{sec: isochrone Symbol decomposition}

The measure equation is local. Consequently, the remaining of the proof can be carried
out in local charts. Suppose $\chart = (\O, \cdiffL)$ is such a
local chart,  with $\O$ \nhd of  a point $\y^0 \in \L\setminus \d\L$.

As in Section~\ref{sec: Preparation theorem: Euclidean symbol
  division} we consider $b(\y) \in \symbH(\R^{2d+2})$ and
proceed with a Euclidean symbol division. The symbol $b$ has compact
support in the $y=(t,x)$ variables, $\supp b \subset K \times
\R^{d+1}$, $K$ a compact subset of $\cdiffL(\O)$ and fast decay in the $\eta=
(\tau,\xi)$ variables.

Consider $0< C_0 < C_0'< C_1'<C_1$ and $\chi\in \Cinfc (\R^d)$ such
that $\chi(\xi) =1 $ in a \nhd of $\{C_{0}'\leq |\xi| \leq C_{1}'\}$ and supported in $\{C_{0}\leq |\xi| \leq C_{1}\}$. One has 
%  ]C_0, C_1[)$ equal to $1$ on a \nhd of
%    $[C_{0}', C_{1}'] \subset ]C_0, C_1[$ for some $0< C_0 < C_0'<
%    C_1'<C_1$ one has
\begin{align*}
  \chi(\xi) b(\y)
  = b_0(t,x,\xi) + b_1(t,x,\xi) \tau
  + q p(\y),
\end{align*}
by Proposition~\ref{prop: Euclidean symbol division},  here, with
the role of the variable $\xi_d$ played by $\tau$. The symbols $b_0$
and $b_1$ fulfill (an adapted version of) Property~\eqref{eq:
  estimation bjk}, that is,
\begin{multline}
  \label{eq: isochrone estimation bj}
  \big| \d_{y}^\alpha \d_{\xi}^\beta  b_{j} (y,\xi) \big| \leq
  C_{N, \beta}
  \jp{\xi}^{-N}, \\
  \pour\ N \in \N, \ \alpha \in \N^{d+1},  \ |\alpha| \leq 1, \
  \beta \in \N^{d},  \
  j=0,1, \ y \in \R^{d+1}, \ \xi\in \R^{d},
\end{multline}
and $q$ fulfills (an adapted version of) Property~\eqref{eq: estimation rk}.

Recall that $\y= (t,x,\xi,\tau) \in \supp \mu^+$ implies $|\xi| \in
[C_{\mu,0}, C_{\mu,1}]$; see \eqref{eq: hyp support mesure mu t=0}. Without loss of generality one can assume $C_{\mu,0} < 1
< C_{\mu,1}$. The constants  $C_0, C_0', C_1', C_1$ are chosen such that
\begin{align*}
  &\supp \chi \subset \{ C_{\mu,0}^2< |\xi| < \C_{1,\mu}^2\}\\
 &\et \ \ \chi=1 \ \dans \ \text{a \nhd of} \
  \{C_{\mu,0}\leq |\xi| \leq C_{\mu,1}\}.
\end{align*}
One 
writes
\begin{align*}
  \dup{\transp\Hp \mu}{b} = \dup{\mu}{\Hp b}= \dup{\mu}{\chi \Hp b}
  =\dup{\mu}{\Hp (\chi b)}
\end{align*}
since $\chi = 1$ on $\supp \mu$ and $\supp d_\xi \chi \cap \supp \mu =
\emptyset$.
One has
\begin{align*}
 \Hp (\chi b) = \Hp b_0 + (\Hp b_1) \tau +  (\Hp q) p,
\end{align*}
implying
\begin{align}
  \label{eq: isochrone simplification symbol}
  \dup{\transp\Hp \mu}{b}
  =\dup{\mu}{\Hp b_0 + (\Hp b_1) \tau}.
\end{align}
With the support properties of $\mu$ given in \eqref{eq: hyp support
  mesure mu t=0} this last duality bracket makes sense.

\subsection{Commutator and limits}
As in Section~\ref{sec: Analysis of a commutator} we consider the
commutator given by $[\hk^2 P_{\k,g}, \OpH(b) \psi]$. With \eqref{eq: isochrone
  simplification symbol},  the symbol $b$ can be  chosen of the form
$b(\y) = b_0(y,\xi) + b_1(y,\xi)\tau$, with
$\supp b_j \subset K \times \R^{d+1}$ for $K$ a compact of
$\cdiffL(\O)$ and $b_j$ fulfilling Property~\eqref{eq: isochrone
  estimation bj}, $j=0,1$, 
meaning fast decay in the $\xi$ variable and only one derivative in
the $x$ variable.
Recall that $y = (t,x)$ and $\y= (t,x,\tau,\xi)$. The symbols $b_j$  are
tangential with respect to time $t$.

Suppose $\psi \in
 \Cinfc(\cdiffL(\O))$ is  equal to $1$ in a bounded \nhd of $K$.
Arguing as in the proof of Proposition~\ref{prop: commutator} one
obtains the followins result.
%%%%%%%%%%%%%%%%%%%%%%%%
% lemma                %
%%%%%%%%%%%%%%%%%%%%%%%%
\begin{lemma}
  \label{lemma: isochrone commutator bj P} 
  One has
  \begin{align*}
    [\hk^2 P_{\k,g},\OpH(b)\psi] = - i \hk
    \OpH(\Hp b)\psi+ o(\hk)_{\L(L^2)}.
  \end{align*}
\end{lemma}
One computes
\begin{align*}
  \Hp b (\y) = \{ p,  b\} (\y)
  = \tb_0(y,\xi) + \tb_1(y,\xi)\tau
  + \tb_2(y,\xi)\tau^2.
\end{align*}
From \eqref{eq: isochrone estimation bj} one has the estimations
\begin{align*}
  %\label{eq: isochrone estimation tbj}
  \big| \d_{\xi}^\beta  \tb_{j} (y,\xi) \big| \leq
  C_{N, \beta}
  \jp{\xi}^{-N}, \quad
  \pour\ N \in \N, \ 
  \beta \in \N^{d},  \
  \ y \in \R^{d+1}, \ \xi\in \R^{d},
\end{align*}
that is, $\tb_j \in \Symbolto{0,n}{-N}{d+1}{d}$, for any $n$ and $N$.
Note that $\psi=1$ on a \nhd of the $(t,x)$-projection of $\supp \tb_j$, $j=0,1, 2$.
By Lemma~\ref{lem:
  continu} in the tangential case (see Remark~\ref{rem: extension measure}) and Proposition~\ref{prop: extension measure},
one has
\begin{align*}
  \lim_{k \to +\infty}\inp{\OpH(\tb_0)\psi u_k}{u_k}_{L^2(t>0)} =
  \lim_{k \to +\infty}
  \inp{\OpH(\tb_0)\psi\unitfunction{t>0}u_k}
  {\unitfunction{t>0}u_k}_{L^2(\R^d)}
  = \dup{\mu^+}{\tb_0},
\end{align*}
using that that no mass leaks at infinity at scale $H$ by
Lemma~\ref{lemma: no mass leaks 1t u}.

The following two results also hold.
%%%%%%%%%%%%%%%%%%%%%%%%
% lemma                %
%%%%%%%%%%%%%%%%%%%%%%%%
\begin{lemma}
  \label{lemma: isochrone mesure b1tau b2tau2}
  One has
  \begin{align*}
  \lim_{k \to +\infty}\inp{\OpH(\tb_1)\hk D_t \psi u_k}{u_k}_{L^2(t>0)} 
  = \dup{\mu^+}{\tb_1 \tau}
  \end{align*}
  and
  \begin{align*}
  \lim_{k \to +\infty}\inp{\OpH(\tb_2)\hk^2 D_t^2 \psi u_k}{u_k}_{L^2(t>0)} 
  = \dup{\mu^+}{\tb_2 \tau^2}.
\end{align*}
\end{lemma}
With Lemmata~\ref{lemma: isochrone commutator bj P}  and \ref{lemma: isochrone mesure b1tau b2tau2} one concludes
that
\begin{align}
  \label{eq: isochrone limit commutator}
  \lim_{k \to +\infty} 
  i \inp{[\hk P_{\k,g},\OpH(b)\psi]  u_k}{u_k}_{L^2(t>0)} 
  = \dup{\mu^+}{\Hp b} = \dup{\transp \Hp \mu^+}{b}.
\end{align}

\medskip
\begin{proof}[\bf Proof of Lemma~\ref{lemma: isochrone mesure b1tau
    b2tau2}]
  We treat the second limit. As
  $[h D_t^2, \psi] \in h \L(L^2) h D_t + h \L(L^2)$, with  \eqref{eq: mass11-proof 2} one has
  \begin{align*}
  \inp{\OpH(\tb_2)\hk^2D_t^2 \psi u_k}{u_k}_{L^2(t>0)} 
    &= \inp{\OpH(\tb_2) \psi \hk^2D_t^2  u_k}{u_k}_{L^2(t>0)} + O(\hk) \\
    &= - \inp{\OpH(\tb_2) \psi \hk^2A_{\k,g}u_k}{u_k}_{L^2(t>0)}+
      O(\hk)\\
    &= - \inp{\OpH(\tb_2) \psi \hk^2A_{\k,g} \unitfunction{t>0}u_k}{\unitfunction{t>0}u_k}_{L^2}+
      O(\hk).
\end{align*}
Then, one obtains
\begin{align*}
  \lim_{k \to +\infty}
  \inp{\OpH(\tb_2)\hk^2D_t^2 \psi u_k}{u_k}_{L^2(t>0)}
  = - \dup{\mu^+}{\tb_2 a_{\k,g}},
\end{align*}
arguing as in Proposition~\ref{prop: commutator} and using Lemma~\ref{lem:
  continu} in the tangential case (see Remark~\ref{rem: extension measure}) and 
Proposition~\ref{prop: extension measure}, using that that no mass leaks at infinity at scale $H$ by Lemma~\ref{lemma: no mass leaks 1t u}. Then, using that
$\supp \mu^+ \subset \supp \mu \subset \Char p$ where $\tau^2 =
-a_{\k,g}$, one obtains
\begin{align*}
  \lim_{k \to +\infty}
  \inp{\OpH(\tb_2)\hk^2D_t^2 \psi u_k}{u_k}_{L^2(t>0)}
  = \dup{\mu^+}{\tb_2 \tau^2}.
\end{align*}

Second, suppose $\varphi \in \Cinfc(\R)$ with $\varphi =1$ in a \nhd of $0$.  One writes
\begin{align*}
  \inp{\OpH(\tb_1)\hk D_t \psi u_k}{u_k}_{L^2(t>0)}
  = \inp{\OpH(\tb_1)\unitfunction{t>0} \hk D_t \psi
  u_k}{\unitfunction{t>0}u_k}_{L^2}
  = I_1 + I_2,
  \end{align*}
with 
\begin{align*}
  I_1 = \inp{\OpH(\tb_1) \OpH\big(1 - \varphi(\tau/R)\big) \unitfunction{t>0} \hk D_t \psi
  u_k}{\unitfunction{t>0}u_k}_{L^2},
\end{align*}
and
\begin{align*}
  I_2 = \inp{\OpH(\tb_1) \OpH\big(\varphi(\tau/R)\big) \unitfunction{t>0} \hk D_t \psi
  u_k}{\unitfunction{t>0}u_k}_{L^2},
\end{align*}
Since $\hk^{\ell+ 2\ell'} D_t^\ell A_{\k,g}^{\ell'}\psi u_k$ is bounded in $L^2$ by
\eqref{eq: mass11-proof 2} for
any $ \ell, \ell'\in \N$, one finds that  $\unitfunction{t>0} \hk D_t \psi
  u_k \in H^s_{sc}(\L)$ for any $s \in [0,1/2[$. With Lemma~4.3
  in \cite{Robbiano:22} one obtains, for such $s$,
\begin{align*}
  \Norm{\OpH\big(1 - \varphi(\tau/R)\big) \unitfunction{t>0} \hk D_t \psi
  u_k}{L^2} \leq C_{s} R^{-s},
\end{align*}
uniformly in $k$. 
Thus, for any $s\in [0,1/2[$, one has 
\begin{align}
   \label{eq: isochrone mesure b1tau b2tau2-1}
  |I_1| \leq  C'_s R^{-s}.
\end{align}
Next, to treat the second term one writes
\begin{align*}
  \unitfunction{t>0} \hk D_t \psi u_k = \hk D_t \big(
\psi \unitfunction{t>0} u_k \big) +i \hk \psi _{|t=0} \udl{u}_k^0
\otimes \delta_{t=0},
\end{align*}
yielding $I_2 = I_2' + I_2''$ with 
\begin{align*}
  I_2' = \inp{\OpH(\tb_1) \OpH\big(\varphi(\tau/R)\big) \hk D_t \big(
\psi \unitfunction{t>0} u_k \big) }{\unitfunction{t>0}u_k}_{L^2},
\end{align*}
and 
\begin{align*}
  I_2'' &= i \hk \inp{\OpH(\tb_1) \psi _{|t=0}
  \udl{u}_k^0  \otimes \OpH\big(\varphi(\tau/R)\big)
  \delta_{t=0}}{\unitfunction{t>0}u_k}_{L^2}\\
  &= i \hk \inp{\OpH(\tb_1) \psi _{|t=0}
  \udl{u}_k^0  \otimes \OpH\big(\varphi(\tau/R)\big)
  \delta_{t=0}}{\unitfunction{t>0} \psi u_k}_{L^2}.
\end{align*}
Set $\hat{b}_1 (t,x,\tau,\xi)= \tau \varphi(\tau/R) \tb_1 (x,\xi)$; it
is smooth, with compact support in $(t,x)$ and has fast decay in
$(\tau,\xi)$. One has $I_2' = \inp{\OpH(\hat{b}_1)\big(\psi 
\unitfunction{t>0} u_k \big) }{\unitfunction{t>0}u_k}_{L^2}$ and 
\begin{align*}
  %\label{eq: isochrone mesure b1tau b2tau2-2}
  \lim_{k \to \infty} I_2' = \dup{\mu^+}{\hat{b}_1}
 \  \mathop{\to}_{R \to \infty}  \dup{\mu^+}{\tau \tb_1},
\end{align*}
by dominated convergence as  $\supp \mu^+$ is compact in the
$(\tau,\xi)$ variables by \eqref{eq: hyp support mesure mu t=0}.

One has $\OpH\big(\varphi(\tau/R)\big)
  \delta_{t=0} = h^{-1} R\check{\varphi}(h^{-1} R t)$, with
  $\check{\varphi}$ the inverse Fourier transform of
  $\varphi$. Thus one obtains
  \begin{align}
    \label{eq: isochrone mesure b1tau b2tau2-3}
  |I_2''| \lesssim R 
    \Norm{\udl{u}_k^0}{L^2(\M)}\Norm{\check{\varphi}(h^{-1} R .)}{L^2(\R)}
    \Norm{\psi \unitfunction{t>0} u_k}{L^2(\L)} 
    \lesssim \hk^{1/2}  R^{1/2} \Norm{\check{\varphi}}{L^2(\R)}
    \lesssim \hk^{1/2} R^{1/2}.
\end{align}
With \eqref{eq: isochrone mesure b1tau b2tau2-1}--\eqref{eq: isochrone
  mesure b1tau b2tau2-3} one concludes the proof. 
\end{proof}

\subsection{End of the proof of the measure equation at an isochrone}

With integrations by parts and one computes 
\begin{align*}
  i \inp{[\hk P_{\k,g},\OpH(b)\psi]  u_k}{u_k}_{L^2(t>0)} = I_k + J_k,
\end{align*}
with 
\begin{align*}
  &I_k = i \inp{\OpH(b)\psi  u_k}{\hk f_k}_{L^2(t>0)}
  - i \inp{\OpH(b)\psi \hk f_k}{u_k}_{L^2(t>0)},
  \\
  &J_k =  -i \hk \inp{\d_t \OpH(b)\psi {u_k}_{|t=0}}{{u_k}_{|t=0}}_{L^2(\M)}
  + i \inp{\OpH(b)\psi {u_k}_{|t=0}}{{\hk \d_t u_k}_{|t=0}}_{L^2(\M)}.
\end{align*}
With Remark~\ref{remark: step function measure} one obtains
\begin{align}
  \label{eq: lim Ik}
  \lim_{k \to \infty} I_k
  = i \dup{M_{0,1}^+}{b} - i \dup{M_{1,0}^+}{b}
  = 2 \dup{\Im M_{0,1}^+}{b}.
\end{align}
One writes
$\d_t \OpH(b)\psi
  = \OpH(b) \psi \d_t + \big( \OpH(\d_t b) \psi + \OpH(b)
  (\d_t\psi)\big)$.
  As $\d_t b = \d_t b_0 + (\d_tb_1)\tau$, 
this gives
\begin{align*}
  -i \hk \inp{\d_t \OpH(b)\psi {u_k}_{|t=0}}{{u_k}_{|t=0}}_{L^2(\M)}
  &=-i \inp{\OpH(b)\psi {\hk \d_t  u_k}_{|t=0}}{{u_k}_{|t=0}}_{L^2(\M)}\\
  &\quad + O(\hk),
\end{align*}
using that $\big( {u_k}_{|t=0}\big)_{k}$ and $\big(\hk \d_t
{u_k}_{|t=0}\big)_{k}$  are bounded sequences in $L^2(\M)$.  One thus finds
\begin{align*}
  J_k
  &= -i \inp{\OpH(b)\psi {\hk \d_t  u_k}_{|t=0}}{{u_k}_{|t=0}}_{L^2(\M)}
  + i\inp{\OpH(b)\psi {u_k}_{|t=0}}{{\hk \d_tu_k}_{|t=0}}_{L^2(\M)}\\
  &\quad   + O(\hk).
\end{align*}
One further writes
\begin{align*}
  \OpH(b)\psi &=\OpH(b_0)\psi + \OpH(b_1)D_t\psi\\
  &=\OpH(b_0)\psi
  + \OpH(b_1)\psi D_t
  -i \hk \OpH(b_1)(\d_t \psi),
\end{align*}
implying
\begin{align*}
  &\inp{\OpH(b)\psi {\hk \d_t  u_k}_{|t=0}}{{u_k}_{|t=0}}_{L^2(\M)}\\
  &\quad=\inp{\OpH(b_0)\psi {\hk \d_t  u_k}_{|t=0}}{{u_k}_{|t=0}}_{L^2(\M)}
  -i \inp{\OpH(b_1)\psi {\hk^2 \d_t^2
  u_k}_{|t=0}}{{u_k}_{|t=0}}_{L^2(\M)}\\
   &\qquad  + O(\hk)\\
  &\quad=\inp{\OpH(b_0)\psi {\hk \d_t  u_k}_{|t=0}}{{u_k}_{|t=0}}_{L^2(\M)}
  -i \inp{\OpH(b_1)\psi {\hk^2 A_{\k,g}
  u_k}_{|t=0}}{{u_k}_{|t=0}}_{L^2(\M)}\\
    &\qquad + O(\hk),
\end{align*}
and
\begin{align*}
  &\inp{\OpH(b)\psi {u_k}_{|t=0}}{{\hk \d_t u_k}_{|t=0}}_{L^2(\M)}\\
  &\quad =\inp{\OpH(b_0)\psi {u_k}_{|t=0}}{{\hk \d_t u_k}_{|t=0}}_{L^2(\M)}
 -i \inp{\OpH(b_1)\psi \hk \d_t{u_k}_{|t=0}}{{\hk \d_t u_k}_{|t=0}}_{L^2(\M)}\\
   &\qquad + O(\hk).
\end{align*}
This yields
\begin{align*}
  J_k &= -i \inp{\OpH(b_0)\psi {\hk \d_t  u_k}_{|t=0}}{{u_k}_{|t=0}}_{L^2(\M)}
  - \inp{\OpH(b_1)\psi \hk^2 A_{\k,g}  {u_k}_{|t=0}}{{u_k}_{|t=0}}_{L^2(\M)}\\
  &\quad 
    +i \inp{\OpH(b_0)\psi {u_k}_{|t=0}}{{\hk \d_t u_k}_{|t=0}}_{L^2(\M)}
    + \inp{\OpH(b_1)\psi \hk \d_t{u_k}_{|t=0}}{{\hk \d_t u_k}_{|t=0}}_{L^2(\M)}\\
    &\quad + O(\hk).
\end{align*}
Arguing as in the proof of Proposition~\ref{prop: commutator}
one finds 
\begin{align*}
  \inp{\OpH(b_1)\psi \hk^2 A_{\k,g}
  {u_k}_{|t=0}}{{u_k}_{|t=0}}_{L^2(\M)}
  = \inp{\OpH(b_1 a_{\k,g})\psi 
  {u_k}_{|t=0}}{{u_k}_{|t=0}}_{L^2(\M)} + O(\hk).
\end{align*}
With \eqref{eq: isochrone limit commutator} and \eqref{eq: lim Ik}, this gives
\begin{align}
  \label{eq: isochrone pre-measure equation}
  \dup{\transp \Hp \mu^+}{b}
  &= 2 \dup{\Im M_{0,1}^+}{b}- i \dup{\nu^0_{1,0}}{b_0} - \dup{\nu^0_{0,0}}{b_1 a_{\k,g}}
  +  i \dup{\nu^0_{0,1}}{b_0} + \dup{\nu^0_{1,1}}{b_1}\\
  &= 2 \dup{\Im M_{0,1}^+}{b}- 2 \dup{\Im \nu^0_{0,1}}{b_0}
  +   \dup{\nu^0_{1,1} - a_{\k,g} \nu^0_{0,0}}{b_1}.\notag
\end{align}

With the analysis of 
Section~\ref{sec: isochrone Symbol decomposition}, Equation~\eqref{eq:
  isochrone pre-measure equation} holds for
a symbol $b(\y) \in \symbH(\R^{2d+2})$ that we write
$b(\y) = b_0(t,x,\xi) + b_1(t,x,\xi) \tau + q p_{\k,g}(\y)$.
Suppose $(x,\xi) \in T^* \M$ and set $\y = (t=0,x,\xi,0)$, and 
 \begin{align*}
      \y^{\oplus} = (0, x,\tau^+(\y), \xi), \qquad 
      \y^{\ominus} = (0, x,\tau^-(\y), \xi).
 \end{align*}
 as in Section~\ref{sec: Measure equation at isochrone}. Since
 $\y^{\oplus},  \y^{\ominus}\in \Char p$ and $\tau^+(\y) = -
 \tau^-(\y)$, one has
 \begin{align*}
   b_0 (t=0,x,\xi)
   = \frac12 \big( b (\y^{\oplus}) + b(\y^{\ominus})
   \big),
   \quad
    b_1 (t=0,x,\xi)
   = \frac{ b (\y^{\oplus}) - b(\y^{\ominus}}{\tau^+(\y) - \tau^-(\y)}.
 \end{align*}
Plugged in \eqref{eq: isochrone pre-measure equation}, this gives the
result of Theorem~\ref{theorem: measure equation at t=0}.
\hfill \qedsymbol \endproof

%%%%%%%%%%%%%%%%%%
% references
%%%%%%%%%%%%%%%%%%
\bibliographystyle{plain}

\end{document}